\RequirePackage{fix-cm}

\documentclass[leqno, fleqn, centertags, 12pt]{article}

\usepackage{latexsym}
\usepackage{amsmath}
\usepackage{amscd}
\usepackage{amssymb}
\usepackage{verbatim}

\usepackage[T1]{fontenc}

\usepackage[upright]{fourier}

\usepackage{fullpage}

\usepackage{tikz-cd}

\tikzcdset{column sep/ttiny/.initial=0.2em}
\tikzcdset{row sep/hhuge/.initial=5.2em}

\makeatletter
\renewcommand{\@makefntext}[1]{\vspace*{0.5ex}\parindent=0em
\hspace*{-0.4em}
\hbox to 0.4em{\hss\@makefnmark}\hspace*{0.4em}{#1}
}
\makeatother

\newcounter{mysectionnumber}
\newcommand{\mysection}[2]{
\setcounter{equation}{0}
\refstepcounter{mysectionnumber}
\section*{ \textnormal{{\themysectionnumber.}\oss {#1}}}\label{#2}}

\numberwithin{equation}{section}

\newcommand{\mynonumbersection}[1]{
\vspace{-0.0ex}
\section*{{}\hspace*{0.00em}$\phantom{1.}$\textnormal{{#1}}}}

\newcommand{\myit}[1]{\textbf{\textit{#1}}\hspace{0.0em}}

\newcounter{myparnum}[mysectionnumber]
\renewcommand{\themyparnum}{\themysectionnumber.\arabic{myparnum}}
\newcommand{\mypar}[2]{\refstepcounter{myparnum}{\vspace{\medskipamount}\textbf{{\themyparnum. #1}\label{#2}}\hspace{0.5em}}}

\newcounter{mylemmanum}[myparnum]
\newcommand{\myuppar}[1]{\vspace{\medskipamount}\textbf{#1}\hspace*{0.5em}}

\newcounter{myappendnumber}
\newcounter{myaparnum}[myappendnumber]
\newcommand{\myappend}[2]{\setcounter{footnote}{0}
\setcounter{myaparnum}{0}
\refstepcounter{myappendnumber}
\section*{\textnormal{A\halfff\fff.\fff{\themyappendnumber.}\oss {#1}}}\label{#2}}

\newcommand{\myapar}[2]{\refstepcounter{myaparnum}{\vspace{\medskipamount}\textbf{{\themyaparnum. #1}\label{#2}}\hspace{0.5em}}}

\renewcommand{\themyaparnum}{A\halfff\fff.\fff\themyappendnumber.\arabic{myaparnum}}

\newcommand{\proof}{\vspace{\medskipamount}{\textbf{{\emph{Proof}.}}\hspace*{0.7em}}}
\newcommand{\subproof}{\vspace{\medskipamount}{{\emph{Proof}.}\hspace*{0.7em}}}
\newcommand{\eproof}{ $\blacksquare$}
\newcommand{\esubproof}{ $\square$}

\newcommand{\dis}{\displaystyle}

\def\sss{\hspace{0.05em}\ }
\def\dss{\hspace{0.1em}\ }
\def\trs{\hspace{0.15em}\ }
\def\qss{\hspace{0.2em}\ }
\def\pss{\hspace{0.3em}\ }
\def\oss{\hspace{0.4em}\ }

\def\halfff{\hspace*{0.025em}}
\def\fff{\hspace*{0.05em}}
\def\dff{\hspace*{0.1em}}
\def\trf{\hspace*{0.15em}}
\def\qff{\hspace*{0.2em}}
\def\pff{\hspace*{0.3em}}
\def\off{\hspace*{0.4em}}

\def\ttff{{\hspace*{-0.05em}--\hspace*{0.15em}}}

\newcommand{\nsp}{\hspace*{-0.1em}}
\newcommand{\nnsp}{\hspace*{-0.15em}}

\newcommand{\dnsp}{\hspace*{-0.2em}}

\renewcommand{\leq}{\leqslant}
\renewcommand{\geq}{\geqslant}
\newcommand{\id}{\mathop{\mbox{id}}\nolimits}

\newcommand{\zzz}{\mathbf{Z}}

\newcommand{\rrr}{\mathbf{R}}
\newcommand{\nnn}{\mathbf{N}}

\newcommand{\image}{\operatorname{Im}}
\newcommand{\kernel}{\operatorname{Ker}}
\DeclareMathOperator{\pr}{pr}
\DeclareMathOperator{\const}{const}

\newcommand{\toto}{\longrightarrow}
\newcommand{\ttoo}{\hspace*{0.2em}\longrightarrow\hspace*{0.2em}}

\begin{document}

\setlength{\baselineskip}{12pt plus 0pt minus 0pt}
\setlength{\parskip}{12pt plus 0pt minus 0pt}
\setlength{\abovedisplayskip}{12pt plus 0pt minus 0pt}
\setlength{\belowdisplayskip}{12pt plus 0pt minus 0pt}

\newskip\smallskipamount \smallskipamount=3pt plus 0pt minus 0pt
\newskip\medskipamount   \medskipamount  =6pt plus 0pt minus 0pt
\newskip\bigskipamount   \bigskipamount =12pt plus 0pt minus 0pt

\author{Nikolai\qss V.\qss Ivanov}
\title{Notes\qss on\qss the\qss bounded\qss cohomology\qss theory}
\date{}

\footnotetext{\hspace*{-0.65em}\copyright\oss 
Nikolai\qss V.\qss Ivanov,\oss 2017,\qss 2020\qss (revision).\trs 
The present paper is partially based on 
author's papers\qss \cite{i1},\oss \cite{i2}.\oss
The results of these papers were obtained in\dss 1984\dss at the
Leningrad Branch of\dss Steklov Mathematical Institute.\oss
Neither the work on the new results and proofs included in the present paper\halfff,\qss
nor its preparation were supported by any governmental 
or non-governmental agency,\qss 
foundation,\qss 
or institution.}

\footnotetext{\hspace*{-0.65em}The author is grateful to\dss M.\qss Prokhorova 
for careful reading  
of this paper and correction of many infelicities.}

\maketitle

\vspace*{6ex}

\renewcommand{\baselinestretch}{1}
\selectfont

\myit{\hspace*{0em}\large Contents}\vspace*{1ex} \vspace*{\bigskipamount}\\ 
\hbox to 0.8\textwidth{\myit{Preface}\hspace*{0.5em}  2}\hspace*{0.5em} \vspace*{1ex}\\
\hbox to 0.8\textwidth{\myit{\phantom{1}1.}\hspace*{0.5em} Introduction\hfil  4}\hspace*{0.5em} \vspace*{0.25ex}\\
\hbox to 0.8\textwidth{\myit{\phantom{1}2.}\hspace*{0.5em} Compactly generated spaces\hfil 11}\hspace*{0.5em} \vspace*{0.25ex}\\
\hbox to 0.8\textwidth{\myit{\phantom{1}3.}\hspace*{0.5em} Weakly principal and principal bundles\hfil 16}\hspace*{0.5em} \vspace*{0.25ex}\\
\hbox to 0.8\textwidth{\myit{\phantom{1}4.}\hspace*{0.5em} McCord classifying spaces and principal bundles\hfil 22}\hspace*{0.5em} \vspace*{0.25ex}\\
\hbox to 0.8\textwidth{\myit{\phantom{1}5.}\hspace*{0.5em} Spaces with amenable fundamental group\hfil 25}\hspace*{0.5em} \vspace*{0.25ex}\\
\hbox to 0.8\textwidth{\myit{\phantom{1}6.}\hspace*{0.5em} Weak equivalences\hfil 38}\hspace*{0.5em} \vspace*{0.25ex}\\
\hbox to 0.8\textwidth{\myit{\phantom{1}7.}\hspace*{0.5em} Elements of\dss homological algebra\hfil 42}\hspace*{0.5em} \vspace*{0.25ex}\\
\hbox to 0.8\textwidth{\myit{\phantom{1}8.}\hspace*{0.5em} Bounded cohomology and the fundamental group\hfil 54}\hspace*{0.5em}\vspace*{0.25ex}\\
\hbox to 0.8\textwidth{\myit{\phantom{1}9.}\hspace*{0.5em} The covering theorem\hfil 64}\hspace*{0.5em} \vspace*{0.25ex}\\
\hbox to 0.8\textwidth{\myit{10.}\hspace*{0.5em} An algebraic analogue of the mapping theorem\hfil 74}\hspace*{0.5em} \vspace*{1ex}\\
\myit{Appendices}\hspace*{0.5em}  \hspace*{0.5em} \vspace*{1ex}\\
\hbox to 0.8\textwidth{\myit{\phantom{1}A.1.}\hspace*{0.5em} The differential of the standard resolution\hfil 84}\hspace*{0.5em} \vspace*{0.25ex}\\
\hbox to 0.8\textwidth{\myit{\phantom{1}A.2.}\hspace*{0.5em} The complex\dss $B\dff(\dff G^{\fff \bullet},\pff U\dff)^{\fff G}$\hfil 86}\hspace*{0.5em} \vspace*{0.25ex}\\
\hbox to 0.8\textwidth{\myit{\phantom{1}A.3.}\hspace*{0.5em} The second bounded cohomology group\hfil 87}\hspace*{0.5em} \vspace*{0.25ex}\\
\hbox to 0.8\textwidth{\myit{\phantom{1}A.4.}\hspace*{0.5em} Functoriality with coefficients\hfil 89}\hspace*{0.5em}  \vspace*{0.25ex}\\
\hbox to 0.8\textwidth{\myit{\phantom{1}A.5.}\hspace*{0.5em} Straight and Borel straight cochains\hfil 90}\hspace*{0.5em}  \vspace*{0.25ex}\\
\hbox to 0.8\textwidth{\myit{\phantom{1}A.6.}\hspace*{0.5em} Double complexes\hfil 92}\hspace*{0.5em} \vspace*{1ex}\\
\hbox to 0.8\textwidth{\myit{References}\hspace*{0.5em}\hfil 93}\hspace*{0.5em}  \vspace*{0.25ex}

\vspace*{2ex}

\mynonumbersection{Preface}

\vspace*{6pt} 
\myuppar{The sources of\dss the bounded cohomology theory.}
The bounded cohomology theory was introduced by Gromov in his fundamental paper\qss \cite{gro}.\oss
It was born in a rich and diverse geometric context\halfff.\oss 
This context includes the works of\dss Milnor\halfff,\pss Hirsch--Thurston,\pss
Sullivan,\pss Smillie,\pss and others on flat bundles and their characteristic classes,\oss
especially the Euler class.\oss
It includes also Mostow's rigidity for hyperbolic manifolds and
Thurston's ideas\qss \cite{t}\qss about the geometry of $3$\dnsp-manifolds.\oss
But the main motivation came from riemannian geometry and,\pss in particular\halfff,\pss
from Cheeger's finiteness theorem.\oss
To quote Gromov\qss \cite{gro},\vspace*{-6pt}
\begin{quote}
The main purpose of this paper is to provide new estimates from
below for the minimal volume in term of the\qss \emph{simplicial\dss volume}\qss
defined in section\qss (0.2).
\end{quote}
\vspace*{-6pt} 
Here the\qss \emph{minimal\dss volume}\qss of a riemannian manifold $M$
is the infimum of the total volume of complete riemannian metrics on $M$
subject to the condition that the absolute value of all sectional curvatures is\qss $\leq\qff 1$\nnsp.\oss
The\qss \emph{simplicial\dss volume}\qss of $M$ is a topological invariant of $M$
defined in terms of\dss the singular homology theory.\oss
The bounded cohomology theory emerged as the most efficient tool to deal with the simplicial volume.\oss
For further details of\dss this story we refer the reader to\qss \cite{gro},\pss
Introduction,\pss and\qss \cite{gro-book},\oss Sections\qss 5.34 -- 5.43.

\myuppar{The origins of\dss the present paper\halfff.}
These are much less glamorous.\oss
Gromov's paper\qss \cite{gro}\qss consists of\dss two types of results:\oss
the geometric and topological results motivating and applying the bounded cohomology theory\halfff;\oss
and\dss the results dealing with the bounded cohomology theory proper\halfff.\oss
Back in the early\dss 1980-is\dss I\dss found the results of\dss the first type fairly accessible,\oss
but failed in my attempts to understand the proofs of\dss the general results
about the bounded cohomology,\oss
such as the vanishing of\dss the bounded cohomology
of simply-con\-nected spaces.\oss

Still,\oss I was fascinated by this theorem and wanted\dss to know a proof\halfff.\oss
The only option available was to prove it myself\halfff.\oss
I found a proof based on a modification of\dss Cartan--Serre killing homotopy group process\qss \cite{cs}\qss
and\dss Dold--Thom construction\qss \cite{dt}.\oss
Using the Dold--Thom results apparently required to limit the theory 
by spaces homotopy equivalent to countable\dss CW-complexes,\pss
but this was sufficient for Gromov's applications.\pss
Emboldened by this success,\oss I found proofs of most of\dss 
Gromov's basic theorems.\oss 
These proofs were based on a strengthening of the vanishing theorem,\pss
the language of\dss the homological algebra,\pss 
as suggested\dss by\dss R.\dss Brooks\qss \cite{br},\pss
and the sheaf\dss theory.\oss
I did not attempt to deal with the last part of\qss \cite{gro},\oss
which is devoted to relative bounded cohomology and applications to complete riemannian manifolds.\oss

This work was reported in my\dss 1985\dss paper\qss \cite{i1}.\oss
The present paper grew out of a modest project to correct the fairly inadequate English translation
of\qss \cite{i1}\qss and\dss to typeset it in LaTeX\halfff.\oss
This project was started in early\dss 1990-is\dss and was abandoned two times
due to the lack of convenient tools to code moderately complicated commutative diagrams.\oss
The third attempt quickly got out of control\dss and\dss lead to a major revision of\dss the theory.\oss

A stimulus for the revision 
was provided by a remark of\trs Th.\dss B\"{u}hler in the Introduction
to his monograph\qss \cite{bu}.\oss
Referring to the theorem about vanishing of\dss the bounded cohomology of a
simply connected space $X$ homotopy equivalent to a countable\dss CW-complex\halfff,\oss he wrote\vspace*{-6pt}
\begin{quote}
The proof\dss of\dss this result is quite difficult and not very well understood
as is indicated by the strange hypothesis on $X$\dss
(Gromov does not make this assumption explicit\halfff,\oss
his proof is however rather sketchy to say at least).\oss
The reason for this is the fact that the complete proof given by Ivanov
uses the Dold--Thom construction which necessitates the countability assumption.\oss
\end{quote}
\vspace*{-6pt}
In the present paper all results are proved without 
the countability assumption.\oss
While keeping the core ideas intact\halfff,\pss
this requires reworking the theory from the very 
beginning.\footnote{Some remarks of\qss 
M.\dss Blank\qss \cite{bl1},\pss \cite{bl2}\qss 
apparently imply that\dss later\dss Th.\dss B\"{u}hler 
realized that the countability assumption is not crucial\dss for the methods of\qss \cite{i1}\halfff,\oss
but did not wrote down any details.}
An argument going back to Eilenberg\qss \cite{e}\qss
allows to extend most of\dss the results to arbitrary topological spaces.\oss

In a separate development\halfff,\pss
I realized that the power of\dss the relative homological algebra is rather limited.\oss
The bounded cohomology groups are real vector spaces
equipped with a canonical semi-norm,\pss
which is the\dss \emph{raison\dss d'\^{e}tre}\pss of\dss
Gromov's theory.\oss
While relative homological algebra provides a convenient framework
for discussing the bounded cohomology groups
as topological vector spaces,\oss
it does not allow to efficiently recover the canonical semi-norms.\oss
By this reason the relative homological algebra is
moved in the present paper from the forefront to the background.\oss
In par\-tic\-u\-lar\halfff,\pss
the bounded cohomology of groups are defined in terms of\dss
the standard resolutions,\pss
in contrast with\qss \cite{i1},\pss
where they were defined in terms of
strong relatively injective resolutions,\pss
and a trick
was used to recover the canonical semi-norms.\oss

\myuppar{Spaces vs. groups.}
In an agreement with the spirit of\dss Gromov's paper\qss \cite{gro},\oss
the focus of the present paper is on the bounded cohomology of\dss spaces,\oss
and\dss the bounded cohomology of\dss discrete groups
are treated as a tool\dss to study the bounded cohomology of spaces.\oss

By now the bounded cohomology theory of\dss groups is a well established
subject in its own right\halfff,\pss
dealing with locally compact groups and especially with
the Lie groups and lattices in them.\oss
Although some of\dss the ideas of\qss \cite{i1}\qss 
found another home in this theory\halfff,\oss
the present paper is intended\dss to provide the foundations for\dss the
bounded cohomology theory of spaces but not of groups.\oss
The monographs by\dss N.\dss Monod\qss \cite{mo}\qss and\dss
Th.\dss B\"{u}hler\qss \cite{bu}\qss discuss the foundation of\dss
the bounded cohomology theory of groups from two rather different points of view.

The focus of\dss the cohomology theory of groups 
is on the cohomology with non-trivial coefficients.\oss
By this reason in this paper the bounded cohomology theory of groups 
is discussed in the case of
general coefficients,\pss
in contrast with\qss \cite{i1},\pss 
where the case of non-trivial coefficients was left to the reader\halfff.\oss
Given the results of the paper\halfff,\oss
it is a routine matter to generalize the bounded cohomology theory of spaces to\dss
cohomology with local coefficients.\oss
But in order not to obscure the main ideas,\pss
the exposition is limited by the case of\dss
trivial coefficients.\oss

\mysection{Introduction}{1}

\vspace*{6pt}
\myuppar{Singular cohomology.}
Let us recall the definition of\dss the singular cohomology of\dss topological spaces.\oss
Anticipating the theory of bounded cohomology,\oss
we will restrict our attention by cohomology with coefficients in $\rrr$\nnsp.\oss
Let\dss $X$\dss be a topological space.\oss 
For each integer\qss $n\dff \geq\dff 0$\qss let\dss $\Delta_{\fff n}$\dss 
be the standard\dss $n$\dnsp-di\-men\-sion\-al simplex and\dss let\qss $S_{n}\fff(\dff X \dff)$\qss 
be the set of continuous maps\qss 
$\Delta_{\fff n}\toto X$\nnsp,\oss 
i.e. the set of\dss $n$\dnsp-dimensional singular simplices in\dss $X$\nnsp.\oss 
The real\dss $n$\dnsp-dimensional singular cochains are defined as 
arbitrary functions\qss 
$S_{n}\fff(\dff X \dff)\toto \rrr$\nnsp.\oss 
They form a group and even a vector space over $\rrr$\nnsp,\oss 
which is denoted by\dss $C^{\fff n}\fff(\dff X \dff)$\nnsp.\oss 
The formula\vspace*{3pt}
\[
\quad
\partial\fff f\dff(\sigma)
\off =\off
\sum^{\dff n\qff +\qff 1}_{\dff i\qff =\qff 0}\qff
(\dff -\qff 1)^{\halfff i}\qff f\dff(\dff \partial_{\fff i}\fff \sigma \dff)\dff,
\] 

\vspace*{-9pt}
where\qss $\partial_{\fff i}\fff \sigma$\qss is the $i$\dnsp-th 
face of the singular simplex\dss $\sigma$\dnsp,\oss 
defines a map\vspace*{3pt} 
\[
\quad
\partial
\dff \colon\dff
C^{\fff n}\fff(\dff X \dff)
\qff \ttoo\qff 
C^{\fff n\dff +\dff 1}\fff(\dff X \dff)\dff.
\] 

\vspace*{-9pt}
As is well known,\pss $\partial\dff \circ\dff \partial\off =\off 0$\dnsp.\oss 
The real singular cohomology groups\dss $H^{\fff *}(\dff X \dff)$\qss 
of the space\dss $X$\dss are defined as the cohomology groups of\dss the complex\vspace*{3pt}
\begin{equation*}
\quad
\begin{tikzcd}[column sep=large, row sep=normal]\dis
0 \arrow[r]
& 
C^{\fff 0}\fff(\dff X \dff) \arrow[r, "\dis \partial\off"]
& 
C^{\fff 1}\fff(\dff X \dff) \arrow[r, "\dis \partial\off"]
&   
C^{\fff 2}\fff(\dff X \dff) \arrow[r, "\dis \partial\off"]
&
\off \ldots \off.
\end{tikzcd}
\end{equation*}

\vspace*{-6pt}
\myuppar{Bounded cohomology.}
Gromov\qss \cite{gro}\qss modified the above definition in a minor\halfff,\pss
at the first sight\halfff,\pss way.\oss
He replaced the space\dss $C^{\fff n}\fff(\dff X \dff)$\dss
of all\dss functions\qss
$S_{n}\fff(\dff X \dff)\toto \rrr$\qss
by the space\dss $B^{\fff n}\fff(\dff X\dff)$\dss
of\dss bounded\dss functions\qss
$S_{n}\fff(\dff X \dff)\toto \rrr$\nnsp.\oss
Such functions
are called\qss \emph{bounded\dss $n$\dnsp-cochains}.\oss
Clearly\halfff,\qss\vspace*{3pt} 
\[
\quad
\partial\dff \bigl(\dff B^{\fff n}\fff(\dff X \dff) \dff\bigr)
\off \subset\off 
B^{\fff n\dff +\dff 1}\fff(\dff X \dff)
\] 

\vspace*{-9pt}
and\dss hence\vspace*{3pt}
\begin{equation*}
\quad
\begin{tikzcd}[column sep=large, row sep=normal]\dis
0 \arrow[r]
& 
B^{\fff 0}\fff(\dff X \dff) \arrow[r, "\dis \partial\off"]
& 
B^{\fff 1}\fff(\dff X \dff) \arrow[r, "\dis \partial\off"]
&   
B^{\fff 2}\fff(\dff X \dff) \arrow[r, "\dis \partial\off"]
&
\off \ldots \off,
\end{tikzcd}
\end{equation*}

\vspace*{-9pt}
is a complex of real vector spaces.\oss
The\qss \emph{bounded cohomology spaces}\qss of\dss $X$\dss are defined as 
the cohomology spaces of\dss this complex
and are denoted\dss by\qss
$\widehat{H}^{\fff *}(\dff X \dff)$\dnsp.\oss 
In more details,\vspace*{3pt}
\[
\quad
\widehat{H}^{\fff n}(\dff X \dff) 
\off =\off
\kernel\dff 
\left(\dff 
\partial\qff \left|\qff B^{\fff n}\fff(\dff X \dff) \right.
\dff\right)
\Bigl/\qff
\image\dff
\left(\dff
\partial\qff \left|\qff B^{\fff n\dff -\dff 1}\fff(\dff X \dff) \right.
\dff\right)\dff.
\Bigr.
\]

\vspace*{-9pt}
The bounded cohomology spaces\qss $\widehat{H}^{\fff n}(\dff X \dff)$\qss
are real vector spaces carrying a\qss \emph{canonical\dss semi-norm},\oss
their\qss \emph{raison\dss d'\^{e}tre}.\oss
In order to define this semi-norm,\oss let us recall\dss that\dss 
$B^{n}\fff(\fff X\fff)$\dss is a Banach space with the
norm\qss $\|\qff \bullet \qff\|$\qss defined as\vspace*{3pt}
\[
\|\qff f \qff\|
\off =\off
\sup\nolimits_{\qff \sigma\qff \in\qff S_{n}\fff(\dff X\dff)}\qff 
\left|\qff f\fff(\fff \sigma\fff) \qff\right|\dff.
\]

\vspace*{-9pt}
The cohomology\qss $\widehat{H}^{\fff n}\fff(\dff X \dff)$\qss
inherits a
semi-norm\qss $\|\qff \bullet \qff\|$\qss
from the norm on\qss $B^{n}\fff(\fff X\fff)$\qss
in an obvious manner\halfff.\oss
Namely\halfff,\oss if\qss $c\qff \in\qff \widehat{H}^{\fff n}\fff(\dff X \dff)$\dnsp,\oss
then\vspace*{3pt}
\[
\quad 
\|\qff c \qff\|
\off =\off
\inf\qff \|\qff f \qff\|\dff,
\] 

\vspace*{-9pt}
where the infimum is taken over all cochains\qss
$f\qff \in\qff B^{\fff n}\fff(\fff X\fff)$\qss representing the cohomology class\dss $c$\nnsp.\oss 
It may happen that\qss $c\qff \neq\qff 0$\qss
but\qss
$\|\qff c \qff\|
\off =\off
0$\dnsp.\oss
This is possible if and only if the image of\qss 
$\partial
\dff \colon\dff  B^{\fff n\dff -\dff 1}\fff(\fff X\fff)
\ttoo
B^{\fff n}\fff(\fff X\fff)$\qss 
is not closed.\oss 

The bounded cohomology have the same functorial properties as the usual ones.\oss 
Thus each continuous map\qss $X\toto Y$\qss induces a map\qss 
$\widehat{H}^{\fff *}\fff(\dff Y \dff)
\ttoo 
\widehat{H}^{\fff *}\fff(\dff X \dff)$\dnsp,\oss 
and the maps induced by homotopic maps\qss $X\ttoo Y$\qss are equal.\oss 
In particular\halfff,\pss 
$\widehat{H}^{\fff *}\fff(\dff X \dff)
\off =\off 
0$\qss if\dss $X$\dss is contractible.\oss
The proofs are the same as in the singular cohomology theory,\oss
since the chain maps and the chain homotopies 
used in these proofs
map bounded cochains to bounded cochains.\oss

\myuppar{Bounded cohomology of groups.}
The\qss \emph{bounded cohomology}\qss
$\widehat{H}^{\fff n}\fff(\dff G \dff)$\dss of a discrete group $G$
are defined by Gromov\qss (see\qss \cite{gro},\oss Section\qss 2.3)\qss
as the bounded cohomology of an Eilenberg--MacLane space\dss $K\dff(\dff G\fff,\pff 1 \fff)$\dnsp.\oss
But this definition is hardly used by Gromov directly.\oss 

There is also an algebraic definition
in the same spirit as the definition of\dss the bounded cohomology of spaces.\oss
One starts with the definition of the classical cohomology groups\dss
${H}^{\fff n}\fff(\dff G\fff,\pff \rrr \dff)$\dss
based on the standard resolution of\dss the trivial\dss $G$\dnsp-module $\rrr$ 
and replaces arbitrary real-valued functions by the bounded ones.\oss 
In more details,\oss
let\dss $B\fff(\dff G^{\fff n} \dff)$\dss 
be the vector space of\dss bounded\dss real-valued\dss functions\qss
$G^{\fff n}\ttoo \rrr$\nnsp.\oss
The usual supremum norm turns\dss $B\fff(\dff G^{\fff n} \dff)$\dss
into a Banach space,\oss
and the action of\dss $G$\dss defined\dss by\vspace*{3pt} 
\[
\quad
(\fff h\cdot\nsp f \dff)\dff
(\fff g_{\fff 1}\fff,\pff \ldots\fff,\pff g_{\fff n\dff -\dff 1}\fff,\pff  g_{\fff n}\fff)
\off =\off
f\dff(\fff g_{\fff 1}\fff,\pff \ldots\fff,\pff g_{\fff n\dff -\dff 1}\fff,\pff  g_{\fff n}\halfff h\fff)
\]

\vspace{-9pt} 
turns\dss $B\fff(\dff G^{\fff n} \dff)$\dss into a\dss $G$\dnsp-module.\oss
Let us consider the sequence\vspace*{3pt}
\begin{equation*}
\quad
\begin{tikzcd}[column sep=large, row sep=normal]\dis
0 \arrow[r]
& 
B\dff(\fff G \dff) \arrow[r, "\dis d_{\dff 0}\off"]
&   
B\dff(\fff G^{\dff 2} \dff) \arrow[r, "\dis d_{\dff 1}\off"]
&
B\dff(\fff G^{\dff 3} \dff) \arrow[r, "\dis d_{\dff 2}\off"]
&
\off \ldots \off,
\end{tikzcd}
\end{equation*}

\vspace*{-6pt}
where the\qss \emph{differentials}\dss $d_{\dff n}$\dss are defined\dss by the formula\vspace*{6pt}
\[
\quad
d_{\dff n}\dff(\fff f \fff)
(\fff g_{\fff 0}\fff,\pff 
g_{\fff 1}\fff,\pff 
\ldots\fff,\pff g_{\fff n\dff +\dff 1}\fff)
\off =\off
(\dff -\qff 1\dff)^{n\dff +\dff 1}\fff
f\dff(\fff g_{\fff 1}\fff,\pff 
g_{\fff 2}\fff,\pff 
\ldots\fff,\pff g_{\fff n\dff +\dff 1}\fff)
\]

\vspace*{-30pt}
\[
\quad
\phantom{d_{\dff n}\dff(\fff f \fff)
(\fff g_{\fff 0}\fff,\pff 
g_{\fff 1}\fff,\pff 
\ldots\fff,\pff g_{\fff n\dff +\dff 1}\fff)
\off =\off}
+\off
\sum\nolimits_{\dff i\qff =\qff 0}^{\dff n}\qff (-\qff 1)^{\dff n\dff -\dff i}\qff 
f\dff(\dff g_{\fff 0}\fff,\pff \ldots\fff,\pff
g_{\dff i}\fff g_{\fff i\dff +\dff 1}\dff,\pff \ldots\fff,\pff   g_{\dff n\dff +\dff 1} \dff)\dff.
\]

\vspace{-6pt}
A standard calculation shows that\oss 
$\dis
d_{\dff n\dff +\dff 1}\dff \circ\trf  d_{\dff n}\off =\off 0$\oss 
for all\qss $n\qff \geq\qff 0$\dnsp,\oss
i.e.\qss that\dss $B\dff(\fff G^{\dff \bullet\dff +\dff 1} \dff)$\dss
together with the differentials $d_{\dff \bullet}$ is a complex\halfff.\oss
A motivation for the formula defining $d_{\dff n}$\nnsp,\oss
which also leads to a non-calculational proof\dss of\dss the identity\qss 
$\dis
d_{\dff n\dff +\dff 1}\dff \circ\trf  d_{\dff n}\off =\off 0$\dnsp,\oss
is discussed in Appendix\qss \ref{categories-classifying}.

The differentials $d_{\dff n}$ commute with the action of\dss $G$\nnsp,\oss 
and\dss hence the subspaces of\dss $G$\dnsp-invariant vectors form a sub\-comp\-lex\halfff.\oss
The\qss \emph{bounded cohomology spaces}\qss
$\widehat{H}^{\fff n}\fff(\dff G \dff)$\dss 
of\dss $G$\dss are defined
as the cohomology spaces of\dss this sub\-comp\-lex\halfff.\oss
Such a definition is also contained in\qss \cite{gro},\oss
albeit\dss hidden inside of\dss a proof\halfff.\oss
See\qss \cite{gro},\oss the last page of\dss Section\qss 3.3.\oss

\myuppar{Amenable groups.}
The notion of an amenable group was introduced by von Neumann\qss \cite{n}.\oss
Morally speaking,\oss a group $G$ is amenable if
one can assign to each bounded function\qss
$G\ttoo \rrr$\qss a real number deserving to be called its\qss
\emph{mean}\qss or the\qss \emph{average value}\qss in a way 
invariant under the translations by the elements of\dss $G$\nnsp.\oss
More formally,\pss
$G$\dss is called\qss \emph{amenable}\qss
if\dss there exists a bounded\dss linear functional\qss
$B\fff(\dff G \dff)\ttoo \rrr$\qss invariant under the natural\dss right action of\dss $G$\dss
on the space of such functionals and having the norm\qss $\leq\qss 1$\nnsp.\oss
Amenable groups 
play a central role in the bounded cohomology theory.\oss
The following theorem is the first illustration of\dss this role.

\vspace*{2pt}
\myuppar{Theorem.} 
\emph{Suppose that\qss $X$\dss 
is a path-connected space homotopy equivalent to a\dss CW-complex and\dss having amenable fundamental group.\oss
Then\oss 
$\widehat{H}^{\fff i}\fff(\dff X\dff)
\off =\off 
0$\oss 
for all\qss $i\qff \geq\qff 1$\nnsp.\oss}

\vspace*{6pt}
See Theorem\qss \ref{simply-connected-homology}.\oss
In particular\halfff,\pss
$\widehat{H}^{\fff i}\fff(\dff X\dff)
\off =\off 
0$\oss 
for all\qss $i\qff \geq\qff 1$\qss
if $X$ is simply-connected and\dss
homotopy equivalent to a\dss CW-complex\halfff.\oss
While this theorem is remarkable by itself\halfff,\oss
the following stronger form of\dss it is crucial\dss for the methods of\dss the present paper\halfff.\oss
Recall that\qss \emph{contracting chain homotopy}\pss for a complex is
a chain homotopy between its identity map and the zero map.

\vspace*{3pt}
\myuppar{Theorem.}  
\emph{Suppose that\qss $X$\dss 
is a path-connected space homotopy equivalent to a\dss CW-com\-plex\halfff.\oss 
If\qss $\pi_{\fff 1}\fff(\dff X \dff)$\qss is amenable,\oss
then the complex\qss
$B^{\fff \bullet}\fff(\dff X \dff)$\qss
admits a contracting chain homotopy
consisting of\qss bounded maps with the norm\qss $\leq\qff 1$\nnsp.\oss}

\vspace*{6pt}
See Theorem\qss \ref{simply-connected-homotopy}.\oss
These two theorems are the main results of\dss Section\qss 
\ref{amenable}.\oss
Sections\qss \ref{cg-spaces}\qss --\qss \ref{classifying-spaces}\qss
are devoted\dss to the tools used in their proofs\halfff.\oss
Section\qss \ref{classifying-spaces}\qss is a review of\dss
McCord's theory of classifying spaces\qss \cite{mcc},\oss
which replaced\dss the Dold--Thom construction used in\qss \cite{i1}.\oss
McCord\qss \cite{mcc}\qss works in the category of compactly generated spaces,\pss
which are discussed in Section\qss \ref{cg-spaces}.\oss
Section\qss \ref{w-principal}\qss is devoted\dss to a modification of\dss
the theory of principal\dss bundles used in\dss Section\qss
\ref{amenable}.

\vspace*{3pt}
\myuppar{Theorem.} 
\emph{Suppose that\qss $X\fff,\pff A$\dss are path-connected spaces and\qss
$\varphi\dff \colon\dff
A\ttoo X$\qss is\dss a weak equivalence,\oss 
i.e.\qss a map inducing isomorphisms of\dss all homotopy groups.\oss
Then the induced map}\vspace*{3pt}
\[
\quad
\varphi^{*}
\qff \colon\qff
\widehat{H}^{\fff n}\fff(\dff X \dff)
\qff \ttoo\qff
\widehat{H}^{\fff n}\fff(\dff A \dff)
\]

\vspace*{-9pt}
\emph{of\dss the bounded cohomology groups is an isometric isomorphism for all\dss $n$\nnsp.\oss}

\vspace{6pt}
See Theorem\qss \ref{weak-equivalence-maps}.\oss
By well\dss known theorems of\qss J.H.C.\dss Whitehead,\pss
this theorem allows to extend the two previous theorems to arbitrary topological spaces
and reduce most of\dss the questions about\dss bounded cohomology
to the case of spaces homotopy equivalent to CW-complexes.

\myuppar{Homological algebra.}
Section\qss \ref{algebra}\qss begins with a fragment of homological algebra 
needed to work with the bounded cohomology of\dss discrete groups.\oss
While most of\dss the results of\dss this sort appear to be direct adaptations
of\dss the classical results,\oss we included\dss detailed proofs.\oss
The reason was very well spelled out by\dss N.\dss Monod in the Introduction
to his monograph\qss \cite{mo}.\vspace*{-6pt}
\begin{quote}
Well,\oss first of all,\pss
such ``obvious transliterations'' 
often just fail\dss to hold\dss true.\pss
Furthermore,\pss
it happens also that usual proofs do not yield\dss
the most accurate statement in our setting.\oss
This is illustrated e.g. by the fundamental lemma on comparison of\dss the resolutions. 
\end{quote}
\vspace*{-6pt}
The difficulties with the fundamental\dss lemma\dss encountered\dss by\dss N.\dss Monod\qss
\cite{mo}\qss are exactly the same as in\trs \cite{i1}\trs and\dss the present paper\halfff.\oss
The main result of\dss Section\qss \ref{algebra}\qss is Theorem\qss \ref{comparing-to-standard},\oss
one of\dss the two main tools to deal with these difficulties.\oss
It is too technical to state it in the introduction,\oss but\dss the point is that
it allows sometimes to replace
the existence statement of\dss 
the fundamental\dss lemma by an explicit construction
with an adequate control of\dss the norms.\oss 
Theorem\qss \ref{comparing-to-standard}\qss
extends Theorem\qss 3.6\qss of\qss \cite{i1}\qss
to the case of\dss twisted coefficients.\oss
In his situation,\pss N.\dss Monod deals with this difficulty by proving an analogue
of\dss Theorem\qss 3.6\qss of\qss \cite{i1}.\oss
See\qss \cite{mo},\oss Theorem\qss 7.3.1.

\myuppar{The bounded cohomology of\dss spaces and of\dss groups.}
This is the topic of\dss Section\qss \ref{spaces}.\oss
Suppose that $G$ is a discrete group 
and\qss 
$p\dff \colon\dff \mathcal{X}\ttoo X$\pss
is a locally trivial principal\dss $G$\dnsp-bundle.\oss
Then $G$ acts freely on $\mathcal{X}$\nnsp,\pss
the quotient space $\mathcal{X}/G$ is equal to $X$\nnsp,\oss
and $p$ is a covering space projection.\oss
One can construct a\dss $G$\dnsp-equivariant morphism of complexes\vspace*{3pt}
\[
\quad
r_{\dff \bullet}
\qff \colon\qff
B\dff(\fff G^{\dff \bullet\dff +\dff 1} \dff)
\qff \ttoo\qff
B^{\fff \bullet}\dff(\dff \mathcal{X} \dff)\dff.
\]

\vspace{-9pt}
consisting of maps with the norm\qss $\leq\qff 1$\nnsp.\oss
This construction is the second main tool for overcoming
difficulties with the fundamental\dss lemma.\oss
It is fairly flexible and works also in other situations.\oss
See,\oss for example,\oss Appendix\qss \ref{borel}.\oss 
The morphism $r_{\dff \bullet}$ 
induces a map\oss\vspace*{3pt}
\[
\quad
\widehat{H}^{\dff *}\fff(\dff G \dff)
\ttoo
\widehat{H}^{\dff *}\fff(\dff X \dff)
\]

\vspace{-9pt}
with the norm\qss $\leq\qff 1$\nnsp.\oss
If\qss $\pi_{\fff 1}\fff(\dff \mathcal{X} \dff)$\qss
is amenable,\oss
then\qss
$B^{\fff \bullet}\fff(\dff \mathcal{X} \dff)$\qss
admits a contracting chain homotopy
consisting of\dss bounded maps with the norm\qss $\leq\qff 1$\nnsp.\oss
This allows to apply Theorem\qss \ref{comparing-to-standard}\qss
and conclude that there exists a\dss $G$\dnsp-equivariant morphism of complexes\vspace*{3pt}
\[
\quad
u_{\dff \bullet}
\qff \colon\qff
B^{\fff \bullet}\dff(\dff \mathcal{X} \dff)
\qff \ttoo\qff 
B\dff(\fff G^{\dff \bullet\dff +\dff 1} \dff)
\]

\vspace{-9pt}
consisting of maps with the norm\qss $\leq\qff 1$\nnsp.\oss

The morphism $u_{\dff \bullet}$ 
induces a map\oss 
$\widehat{H}^{\dff *}\fff(\dff X \dff)
\ttoo
\widehat{H}^{\dff *}\fff(\dff G \dff)$\oss
with the norm\qss $\leq\qff 1$\nnsp.\oss
The two maps between $\widehat{H}^{\dff *}\fff(\dff X \dff)$
and $\widehat{H}^{\dff *}\fff(\dff G \dff)$
turn out to be mutually inverse.\oss
This leads to the following fundamental theorem.\oss

\myuppar{Theorem.} 
\emph{If\oss 
$\pi_{\fff 1}\fff(\dff \mathcal{X} \dff)$\qss
is amenable,\oss then the maps}\vspace*{0pt}
\begin{equation*}
\quad
\begin{tikzcd}[column sep=normal, row sep=normal]\dis
u_{\dff *}
\qff \colon\qff
\widehat{H}^{\dff *}\fff(\dff X \dff)
\off 
\arrow[r, shift left=3pt]
&
\off 
\widehat{H}^{\dff *}\fff(\dff G \dff) 
\off\dff \colon\dff
r_{\fff *}\qff,
\arrow[l, shift left=3pt]
\end{tikzcd}
\end{equation*}

\vspace*{-9pt}
\emph{induced by\dss $u_{\dff \bullet}$\dss 
and\dss $r_{\dff \bullet}$\dss respectively,\oss
are mutually inverse isometric isomorphisms.\oss}

\vspace*{6pt}
See Theorem\qss \ref{acyclic-bundles-isomorphism}\qss 
for a formally stronger result\halfff.\oss
This theorem,\oss in particular\halfff,\oss
immediately implies that there exists an isometric isomorphism\oss
$\widehat{H}^{\dff *}\fff(\dff X \dff) 
\qff \ttoo \qff
\widehat{H}^{\dff *}\fff(\dff \pi_{\fff 1}\halfff(\dff X \dff) \dff)$\dnsp.\oss
With these tools at hand,\oss
it is an easy matter to prove the functoriality of this isomorphism
and other expected properties.\oss 
See Section\qss \ref{spaces}.\oss
As the first application of\dss these tools to the bounded cohomology of spaces,\oss
we prove the first of\dss two main Gromov's theorems about bounded cohomology,\oss
called\dss by him the\qss \emph{Mapping\dss theorem}.

\myuppar{The Mapping\dss theorem.} 
\emph{Let\qss $X\fff,\pff Y$\qss 
be two path-connected spaces and\dss let\oss
$\varphi\dff \colon\dff Y\ttoo X$\oss be a continuous map.\oss
If\dss the induced\dss homomorphism of the fundamental groups}\oss\vspace*{3pt} 
\[
\quad
\varphi_*
\dff \colon\dff 
\pi_{\dff 1}\fff(\dff Y \dff)\ttoo \pi_{\dff 1}\fff(\dff X \dff)
\]

\vspace*{-9pt} 
\emph{is surjective with amenable kernel,\oss 
then\oss 
$\varphi^*
\dff \colon\dff
\widehat{H}^{\dff *}\fff(\dff X \dff)\ttoo \widehat{H}^{\dff *}\fff(\dff Y \dff)$\oss
is an isometric isomorphism.}

\vspace*{9pt}
See Theorem\qss \ref{mapping-theorem}.\oss
This is a far reaching generalization of\dss the vanishing theorem for spaces
with ame\-na\-ble fundamental group.\oss
Another generalization deals with coverings of spaces. 
In order to state it\halfff,\oss
let\dss us\sss call\sss a path-connected subset $Y$ of\dss a topological space $X$\dss \emph{a\-me\-na\-ble}\qss
if\trs the image 
of\trs the inclusion homomorphism\dss
$\pi_{\dff 1}\fff(\trf Y \trf)
\ttoo 
\pi_{\dff 1}\fff(\trf X \trf)$\sss
is\dss an amenable subgroup of\dss $\pi_{\dff 1}\fff(\trf X \trf)$\dnsp.\oss
Let\dss us call a covering\sss $\mathcal{U}$\sss of\dss a space $X$ by its subsets\dss \emph{nice}\sss\qss
if\dss $X$\sss is\dss paracompact\sss and\sss $\mathcal{U}$\sss is\dss either open,\oss
or\sss is\dss closed and\dss locally\sss finite,\oss and\dss finite
intersections of\dss elements of\dss $\mathcal{U}$\qss
(including\sss $X$\sss itself\dff)\qss
are path connected and,\oss in\dss the case of\dss closed\sss $\mathcal{U}$\nnsp,\oss 
behave nicely\sss with\sss respect\dss to\sss singular\sss homology.\oss
See the assumption\dss ({\fff}C{\fff})\dss at\dss the beginning of\qss Section\qss \ref{covering}.\oss

\myuppar{The\dss Covering\dss theorem.} 
\emph{Let\dss
$\mathcal{U}$\dss is a covering of a space\qss $X$\qss by amenable subsets.\oss
Suppose that\dss
$\mathcal{U}$\dss 
is\dss a nice covering.\oss 
Let\qss $N$\dss be the nerve of\dss the covering\trs $\mathcal{U}$\dss and\dss $|\fff N \fff|$\dss
be the geometric realization of\qss $N$\nnsp.\oss
Then the canonical homomorphism\qss
$\widehat{H}^{\fff *}\fff(\dff X \dff)\ttoo {H}^{\fff *}\fff(\dff X \dff)$\qss
can be factored through the canonical homomorphism\qss
${H}^{\fff *}\fff(\dff |\fff N \fff | \dff)\ttoo {H}^{\fff *}\fff(\dff X \dff)$\dnsp.\oss}

\vspace*{6pt}
See Theorem\qss \ref{covering-theorem}.\oss
If\dss $\mathcal{U}$\sss is\dss open,\oss then\sss the conclusion\sss holds
for arbitrary\sss $X$\nnsp.\oss
See\qss \cite{i3}.\oss
The covering theorem is a generalization and a much more 
precise version of\dss the second of\dss two 
main Gromov's theorems about\dss bounded cohomology,\oss
called by him the\qss \emph{Vanishing\dss theorem}.\oss

\myuppar{The\qss Vanishing\qss theorem.}
\emph{Suppose that a manifold\pss $X$\dss 
can be covered by open amenable subsets
in such a way that every point of\qss $X$\dss
is contained in no more than $m$ elements of this covering.\oss
Then the canonical homomorphism\oss
$\widehat{H}^{\fff i}\fff(\dff X \dff)
\ttoo
{H}^{\fff i}\fff(\dff X \dff)$\oss
vanishes for\qss $i\qff \geq\qff m$\nnsp.\oss}

\vspace{6pt}
See Theorem\qss \ref{vanishing}.\oss
Here we implicitly assume that manifolds are required to be paracompact\qss
(this is needed to ensure the niceness of\dss the covering\fff).\oss
Since Gromov developed the bounded cohomology theory for the sake of
applications to riemannian manifolds,\oss
and all riemannian manifolds are paracompact\halfff,\oss
this seems to be a reasonable assumption.\oss

The proofs of\dss the last two theorems are independent of
the theory of\dss bounded cohomology of groups,\oss
i.e.\qss of Sections\qss \ref{algebra}\qss and\qss \ref{spaces}.\oss
Using the sheaf\dss theory,\oss
these theorems are deduced from the vanishing
of\dss the bounded cohomology of spaces with amenable fundamental groups.\oss

\myuppar{The mapping theorem for groups.} 
\emph{Suppose that\dss $A$\dss is a normal amenable subgroup of\dss
a group\trs $\Gamma$\dnsp.\oss
Then the quotient map\oss
$\Gamma\ttoo \Gamma/A$\qss induces an isometric isomorphism}\oss
\[
\quad
\alpha^{\fff *}
\qff \colon\qff
\widehat{H}^{\dff *}\fff(\dff \Gamma/A \dff)
\qff \ttoo\qff
\widehat{H}^{\dff *}\fff(\dff \Gamma \dff)\dff.
\]

\vspace*{-3pt}
By using Eilenberg--MacLane spaces
and\dss  isometric isomorphisms\oss
$\widehat{H}^{\dff *}\fff(\dff X \dff) 
\qff \ttoo \qff
\widehat{H}^{\dff *}\fff(\dff \pi_{\fff 1}\halfff(\dff X \dff) \dff)$\oss
one can easily deduce this theorem from the Mapping theorem.\oss
Nevertheless,\oss a proof\dss based only on the results of\dss 
Section\qss \ref{algebra}\qss is given in Section\qss \ref{a-mapping-theorem}.\oss
The main reason is that the algebraic proof
is a drastically simplified version of\dss 
the main part of\dss the proof of\dss the vanishing
theorem of\dss the bounded cohomology of spaces with amenable fundamental groups.\oss
As such,\oss it is quite instructive.\oss
The algebraically-minded readers may prefer to start with\dss
Sections\qss \ref{algebra}\qss and\qss \ref{a-mapping-theorem}.\oss

\myuppar{Twisted coefficients.}
Another reason for presenting an algebraic proof of the mapping theorem for groups
is that this is the best place to go beyond the definition of\dss the bounded cohomology
with twisted coefficients and\dss to prove something non-trivial about them.\oss
The main theorem of\dss Section\qss \ref{a-mapping-theorem},\oss
namely Theorem\qss \ref{normal-amenable},\oss
is a version of\dss the mapping theorem for groups
for the bounded cohomology with twisted coefficients,\oss
although not with arbitrary ones.\oss
In order for the proof\dss to work\halfff,\oss
the module of coefficients $U$ should admit an invariant averaging procedure
for $U$\dnsp-valued functions on amenable groups.\oss
In a related context of\dss the continuous cohomology of\dss Banach algebras,\oss
B.\dss Johnson\qss \cite{j}\qss noticed that this is the case when $U$ 
is the Banach dual of some other Banach module,\oss
and in Section\qss \ref{a-mapping-theorem}\qss we work only with such modules of coefficients.\oss
This class of coefficients plays a central role in the work of\dss N.\dss Monod,\oss
and one may think that this is\dss ``the right''\dss class of coefficients for 
the bounded cohomology theory.\oss

In the context of\dss the bounded cohomology of\dss discrete groups
this class of coefficients was used for the first time by\dss
G.A.\dss Noskov\qss \cite{no},\oss
who proved a weaker version of\dss Theorem\qss \ref{normal-amenable}.\oss
Contrary to his claim,\oss his proof cannot be adjusted to prove that the induced map is an isometry.\oss
Cf.\qss the above discussion of\dss the difficulties with the fundamental\dss lemma.

\myuppar{Means of\dss vector-valued functions.}
They play a central role in the 
bounded cohomology theory with twisted coefficients.\oss
One of\dss the tools used may be of independent interest\halfff.\oss

Let\sss 
$U\off =\off V^{\fff *}$\sss
be a Banach space dual to some Banach space $V$\dnsp.\oss
Let\sss $A$\sss be a set and\dss let\sss
$B\fff(\dff A\fff,\pff U \dff)$\dss
be the Banach space of\dss maps\qss
$f\dff \colon\dff A\ttoo U$\qss
such that the real-valued\dss function\dss $\|\qff f \qff\|$\dss is bound\-ed.\oss
Then every mean\qss 
$M\dff \colon\dff B\fff(\dff A \dff)\ttoo \rrr$\qss
leads to a mean\qss
$m\dff \colon\dff B\fff(\dff A\fff,\pff U \dff)\ttoo U$\dnsp.\oss
The mean $m$ 
commutes with all operators\qss
$L^*\dff \colon\dff U\ttoo U$\qss
ad\-joint to bounded operators\qss
$L\dff \colon\dff V\ttoo V$\dnsp.\oss

If\dss $A$\dss is a group and $M$ is an invariant mean,\oss
then $m$ is also invariant\halfff.\oss
See Lemma\qss \ref{vector-means}.\oss
The construction\qss
$M\qff \longmapsto\qff m$\qss
commutes with push-forwards by surjective homomorphisms\qss
(see Sections\qss \ref{amenable}\qss and\qss \ref{a-mapping-theorem}\qss 
for the 
definition).\oss
This is an obvious,\oss
but crucial observation.\oss

\myuppar{Bounded cohomology of\dss spaces with twisted coefficients.}
The theory of\dss bounded cohomology of spaces can be fairly routinely
extended\dss to the case of\dss twisted\qss ({\fff}better known as\qss \emph{local}\trf)\qss
coefficients in $\pi_{\fff 1}\halfff(\dff X \dff)$\dnsp-modules
dual to Banach $\pi_{\fff 1}\halfff(\dff X \dff)$\dnsp-modules.\oss
The properties of\dss the above construction of means of\dss vector-valued\dss functions
ensure that all arguments in Section\qss \ref{amenable}\qss still work\halfff.
Apparently,\oss there is no averaging procedure for more general coefficients,\oss
except in hardly interesting special cases,\oss such as for spaces with finite fundamental group.\oss

\myuppar{Appendices.}
There are six appendices,\oss
complementing the main part of\dss the text\halfff.\oss

\myuppar{A.\ref{categories-classifying}.}
This appendix is devoted\dss to a conceptual 
motivation of\dss the definition of\dss the cohomology of groups,\oss
both the bounded and the usual ones.\oss
There is no reason to overload\dss the main part of\dss the text by this
categorical approach.\oss

\myuppar{A.\ref{invariants}.}
This appendix goes in the opposite direction and gives an explicit description
of\dss the subcomplex of\dss invariant subspaces of\dss $B\dff(\fff G^{\dff \bullet} \dff)$\dnsp.\oss

\myuppar{A.\ref{second}.}
This appendix is devoted\dss to a proof\dss of\dss a theorem to the effect\qss
that\dss $\widehat{H}^{\dff 2}\fff(\dff G \dff)$\dss 
is always a Hausdorff\qss
(and\dss hence a Banach)\qss space.\oss
See Theorem\qss \ref{second-is-norm}.\oss
This result is due independently to\dss
Sh.\dss Matsumoto and\dss Sh.\qss Morita\qss \cite{mm}\qss
and\dss the author\qss \cite{i2}.\oss

\myuppar{A.\ref{functoriality-coefficients}.}
In this appendix the functorial properties of\dss
the bounded cohomology of groups with coefficients are presened in full generality,\oss
complementing Sections\qss \ref{spaces}\qss
and\qss \ref{a-mapping-theorem}.

\myuppar{A.\ref{borel}.}
In this appendix the methods of\dss 
Section\qss \ref{spaces}\qss 
are applied to the complex of\qss \emph{straight\qss Borel\qss cochains},\oss
playing an important role in Gromov's paper\qss \cite{gro}.\oss

\myuppar{A.\ref{double-complexes}.}
This appendix\sss is\sss devoted\dss to\sss two well\dss known\dss theorems
about double complexes used in Section\qss \ref{covering}.\oss 
An old-fashioned reference\sss
to\sss the proof\dss of\dss one of\trs them\dss is\dss given,\oss
as also a hint\sss for a\sss proof\dss of\trs the other.\oss

\mysection{Compactly\qss generated\qss spaces}{cg-spaces}

\vspace*{6pt}
\myuppar{Test spaces and test maps.}
Following\qss \cite{td}\qss
(see\qss \cite{td},\oss Section\qss 7.9),\oss
a\dss \emph{test space}\qss is defined as
a compact Hausdorff space,\pss 
and a\dss \emph{test map}\qss
is defined as a continuous map\qss 
$f\dff \colon\dff C\toto X$\qss from a test space to a topological space $X$\nnsp.\pss

\myuppar{Compactly generated spaces.}
A subset $A$ of a topological space $X$ is called\qss \emph{compactly closed}\qss
if\dss for every test map\qss $f\dff \colon\dff C\ttoo X$\qss
the preimage\dss $f^{\dff -\dff 1}(\fff A\fff)$\dss is closed.\oss
If $A$ is closed in $X$\dnsp,\oss then $A$ is compactly closed.\oss
For a topological space $X$ the space
$k\fff(\dff X\dff)$ is defined as the topological space 
having the same set of points as $X$ and
compactly closed subsets of $X$ as the closed sets. 
The space $X$ is called\qss \emph{compactly generated}\pss if\oss
$k\fff(\dff X\dff)\off =\off X$\oss as topological spaces,\oss i.e.\qss 
if every compactly closed subset of $X$ is closed.\oss

\mypar{Lemma.}{test-maps}
\emph{Suppose that $C$ is a test space.\oss
A set-theoretic map\qss $f\dff \colon\dff C\ttoo X$\qss
is continuous as a map into\qss $X$\qss if and only if\dss it is continuous as map into\qss $k\fff(\dff X\dff)$\dnsp.}

\proof\qss
If $f$ is continuous as a map into $X$\nnsp,\oss then $f$ is a test map.\oss 
If $A$ is closed in $k\fff(\dff X\dff)$\dnsp,\oss
then $A$ is compactly closed in $X$ and hence\dss $f^{\dff -\dff 1}(\fff A\fff)$\dss is closed.\oss
Hence $f$ is continuous as a map into $k\fff(\dff X\dff)$\dnsp.\oss
If$f$ is continuous as a map into $k\fff(\dff X\dff)$ and 
$A$ is closed in $X$\dnsp,\oss then $A$ is compactly closed and\dss 
hence\dss $f^{\dff -\dff 1}(\fff A\fff)$\dss is closed.\oss
It follows that   
$f$ is continuos as a map into $X$\nnsp.\oss  \eproof

\mypar{Lemma.}{cg}
\emph{Suppose that\qss $X\fff,\pff Y$\qss are topological spaces.}\vspace{-9pt}
\begin{itemize}
\item[(i)] \emph{The identity map\qss $k\fff(\dff X\dff)\toto X$\qss is continuous.}
\item[(ii)] \emph{\dnsp$k\fff(\dff X\dff)$ is a compactly generated space,\oss
i.e.\oss $k\fff(\dff k\fff(\dff X\dff)\dff)\off =\off k\fff(\dff X\dff)$\dnsp.\oss} 
\item[(iii)] \emph{If\qss $F\dff \colon\dff X\ttoo Y$\qss is continuous,\oss
then $F$ is continuous as a map\qss 
$k\fff(\dff X\dff)\ttoo k\fff(\dff Y\dff)$\dnsp.}
\end{itemize}

\vspace{-6pt}
\proof\qss
The part\qss (i)\qss is obvious.\oss 
Lemma\qss \ref{test-maps}\qss implies that the properties of being compactly closed in $X$
and in $k\fff(\dff X\dff)$ are e\-quiv\-a\-lent\halfff.\oss 
The part\qss (ii)\qss follows.\oss
If\qss
$f\dff \colon\dff C\ttoo X$\qss 
is a test map,\oss
then\qss $F\dff \circ f$\qss is a test map\qss $C\ttoo Y$\dnsp.\oss
Therefore,\oss
if $A$ is compactly closed in $Y$\dnsp,\oss
then\vspace*{3pt}
\[
\quad
f^{\dff -\dff 1}(\dff F^{\dff -\dff 1}(\dff A \dff) \dff)
\off =\off
(\dff F\dff \circ f \dff)^{-\dff 1}(\dff A \dff)
\]

\vspace{-9pt}
is closed\dss in $C$\nnsp.\oss
It follows that\dss $F^{-\dff 1}(\fff A\fff)$\dss is compactly closed in $X$\nnsp.\oss
The part\qss (iii)\qss follows.\oss  \eproof

\mypar{Lemma.}{cg-quotient}
\emph{Suppose that\qss $p\dff \colon\dff X\ttoo Y$\qss is a quotient map\qss
({\fff}i.e.\dss $p$\dss is a continuous surjective map and\dss $Y$\dss
has the quotient topology with respect to\dss $p$\nnsp).\oss
If\pss $X$\trs is compactly generated,\oss then\qss $Y$\trs is compactly generated.}

\proof\qss
Let\qss $B\qff \subset\qff Y$\qss be compactly closed subset\halfff.\oss
Since $p$\sss is a quotient map,\pss $B$\dss is closed if\qss
$p^{\fff -\dff 1}(\dff B\dff)$\qss is closed.\oss
If\qss
$f\dff \colon\dff C\toto X$\qss
is a test map,\oss
then\qss $p\dff \circ\dff f\dff \colon\dff C\toto Y$\qss 
is a test map and\dss hence\vspace*{3pt}
\[
\quad
f^{\dff -\dff 1}(\dff p^{\dff -\dff 1}(\dff B \dff) \dff)
\off =\off
(\dff p\dff \circ f \dff)^{-\dff 1}(\dff B \dff)
\]

\vspace{-9pt} 
is closed.\oss
It follows that\dss
$p^{\dff -\dff 1}(\dff B \dff)$\dss is compactly closed,\oss and\dss hence is closed in $X$\nnsp.\oss
It follows that $B$ is closed in $Y$\dnsp.\oss 
The lemma follows.\oss  \eproof

\mypar{Lemma.}{cg-subspace}
\emph{Suppose that\pss $Y\qff \subset\qff X$\qss is a subspace of a compactly generated space $X$\nnsp.\oss
If\pss $Y$\qss is either closed or open,\oss then\pss $Y$\qss is compactly generated.\oss}

\proof\qss
Suppose that $Y$ is closed.\oss
Let $B$ be a compactly closed subset of $Y$\nnsp.\oss
We need to prove that $B$ is closed in $Y$\nnsp,\oss
or\halfff,\oss equivalently,\oss that $B$ is closed in $X$\nnsp.\oss
If\qss
$f\dff \colon\dff C\toto X$\qss
is a test map,\oss
then\qss
$D\qff =\qff f^{\fff -\dff 1}(\dff Y\dff)$\qss is closed in $C$\nnsp.\oss
Since $C$ is compact\halfff,\pss $D$ is also compact and\dss hence the map\qss
$g\dff \colon\dff D\ttoo Y$\qss is a test map.\oss
Since $B$ is compactly closed in $Y$\dnsp,\oss
it follows that\qss $g^{\fff -\dff 1}(\dff B\dff)$\qss is closed in $D$ and\dss hence in $C$\nnsp.\oss
But\qss $B\qff \subset\qff Y$\qss implies that\qss
$f^{\fff -\dff 1}(\dff B\dff)
\qff =\qff
g^{\fff -\dff 1}(\dff B\dff)$\qss
and hence\dss $f^{\fff -\dff 1}(\dff B\dff)$\dss is closed in $C$\nnsp.\oss
It follows that $B$ is compactly closed in $X$ and hence is closed in $X$\nnsp.\oss
Therefore,\pss $B$ is closed in $Y$\dnsp.\oss
Since $B$ was an arbitrary compactly closed subset of $Y$\dnsp,\oss 
it follows that the space $Y$
is compactly generated.\oss

Suppose now that $Y$ is open.\oss
Let $B$ be a compactly closed subset of $Y$\nnsp.\oss
In order to prove that $B$ is closed in $Y$ it is sufficient to prove
that\qss $Y\dff \smallsetminus\dff B$\qss is open in $Y$\nnsp,\oss
or\halfff,\oss equivalently,\oss that\qss $Y\dff \smallsetminus\dff B$\qss is open in $X$\nnsp.\oss
Since $X$ is compactly generated,\oss it is sufficient to prove that\qss
$f^{\fff -\dff 1}(\dff Y\dff \smallsetminus\dff B \dff)$\qss is open for every test map\qss
$f\dff \colon\dff C\toto X$\nnsp.\oss
Suppose that\qss $y\qff \in\qff f^{\fff -\dff 1}(\dff Y\dff \smallsetminus\dff B \dff)$\dnsp.\oss
Then\qss $f^{\fff -\dff 1}(\dff Y \dff)$\qss is an open neighborhood of $y$ in $C$\nnsp.\oss
Since $C$ is compact Hausdorff and hence is regular\halfff,\oss
there exists an open set\qss $U\qff \subset\qff C$\qss such that\qss $y\qff \in\qff U$\qss
and the closure\dss $\overline{U}$\dss is contained in\qss $f^{\fff -\dff 1}(\dff Y \dff)$\dnsp.\oss
Then\dss $\overline{U}$\dss is a compact Hausdorff space and hence $f$ induces a test map\qss
$g\dff \colon\dff \overline{U}\ttoo Y$\dnsp.\oss
Since $B$ is compactly closed,\oss the preimage\qss $g^{\fff -\dff 1}(\dff B \dff)$\qss
is closed in\dss $\overline{U}$\dss and hence\vspace*{3pt}
\[
\quad
f^{\fff -\dff 1}(\dff Y\dff \smallsetminus\dff B \dff)\dff \cap\dff  \overline{U}
\off =\off
\overline{U}\dff \smallsetminus\dff g^{\fff -\dff 1}(\dff B \dff)
\]

\vspace*{-9pt}
is open in\dss  $\overline{U}$\dnsp.\oss
It follows that\qss
$f^{\fff -\dff 1}(\dff Y\dff \smallsetminus\dff B \dff)\dff \cap\dff U$\qss is open in $U$ and\dss
hence is open in $C$\dnsp.\oss 
Therefore\qss $f^{\fff -\dff 1}(\dff Y\dff \smallsetminus\dff B \dff)$\qss
is a neighborhood of $y$ in $C$\nnsp.\oss
It follows that\qss  $f^{\fff -\dff 1}(\dff Y\dff \smallsetminus\dff B \dff)$\qss
is open in $C$\nnsp,\oss
and\dss hence\qss $Y\dff \smallsetminus\dff B$\qss is open in $X$
and $B$ is closed in $Y$\dnsp.\oss  \eproof

\myuppar{Products.}
The usual cartesian product of topological spaces\qss $X\fff,\pff Y$\qss will be denoted by\qss
$X\times_c Y$\dnsp.\oss
Even if\qss $X\fff,\pff Y$\qss
are compactly generated,\oss
then\qss $X\times_c Y$\qss does not need to be compactly generated.\oss
In order to deal with this issue,\oss
a new product\pss 
$X\times_{k} Y
\off =\off
k\fff(\dff X\times_c Y\dff)$\pss is introduced.\oss

\mypar{Lemma.}{k-product}
\emph{If\pss $X\fff,\pff Y$\qss
are compactly generated,\oss
then\qss $X\times_{k} Y$\qss is the product in the category of all compactly generated
spaces and continuous maps.\oss}

\proof\qss
By\dss Lemma\qss \ref{cg}\dff(i)\qss
the map\qss $X\times_{k} Y\ttoo X\times_{c} Y$\qss 
is continuous.\oss
Since the projections\qss $X\times_{c} Y\ttoo X\fff,\pff Y$\qss are continuous,\oss
the same projections considered as maps\qss $X\times_{k} Y\ttoo X\fff,\pff Y$\qss are also continuous.\oss
Suppose that $Z$ is a compactly generated space and that $f$ and $g$
are continuous maps\qss $Z\toto X$\qss and\qss $Z\toto Y$\qss respectively.\oss
Since\qss $X\times_{c} Y$\qss is the product in the category of all topological spaces,\oss
there is a unique map\qss $(\fff f,\qff g \fff)\dff \colon\dff Z\ttoo X\times_c Y$\qss
with the components\qss $f\fff,\pff g$\nnsp.\oss
By Lemma\qss \ref{cg}\dff(iii)\qss
the map\dss $(\fff f,\qff g \fff)$\dss leads to a continuous map\qss
$k\fff(\fff f,\qff g \fff)\dff \colon\dff k\fff(\dff Z \dff)\ttoo k\fff(\dff X\times_c Y \dff)$\qss
with the components\qss $f\fff,\pff g$\nnsp.\oss
Since\qss $Z\qff =\qff k\fff(\dff Z \dff)$\dnsp,\oss
this completes the proof\halfff.\oss \eproof

\mypar{Lemma.}{times-lc}
\emph{If\pss $X$\qss is compactly generated and\qss $Y$\qss is locally compact Hausdorff\halfff,\oss
then\qss $X\dff \times_{c}\dff Y$\qss is compactly generated and hence\pss
$X\dff \times_k\dff Y\off =\off X\dff \times_c\dff Y$\dnsp.\oss}

\proof\qss
The following proof\dss largely follows the proof of Theorem\qss 4.3\qss in\qss \cite{st}.\oss 
Let $A$ be a compactly closed subset of\qss $X\times_c\dff Y$\qss
and\qss $(x\fff,\pff y)$\qss be a point in the complement of $A$\nnsp.\oss
Since $Y$ is locally compact Hausdorff\halfff,\oss there is a compact\qss 
(and Hausdorff\fff)\qss 
neighborhood $U$ of $y$ in $Y$\dnsp.\oss
The inclusion\qss 
$i
\dff \colon\dff 
x\dff \times_c\dff U\ttoo X\times_c\dff Y$\qss is a test map
and\dss hence\oss
$A\dff \cap\dff (\fff x\dff\times_c\dff U\dff)
\off =\off
i^{\fff -\dff 1}\fff(\fff A \fff)$\qss is closed.\oss
Therefore,\oss replacing $U$ by a smaller compact neighborhood of $y$ if necessary,\oss
we can assume that\oss $A\dff \cap\dff (\fff x\dff\times_c\dff U\dff)\off =\off \emptyset$\nnsp.\oss
Let $B$ be the projection of\qss
$A\dff \cap\dff (\dff X\dff\times_c\dff U\dff)$\qss
to\dss $X$\nnsp.\oss
Then\qss $x\qff \not\in\qff B$\qss by the choice of\dss $U$\nnsp.\oss
If\qss
$f\dff \colon\dff C\ttoo X$\qss
is a test map,\oss
then\qss\vspace*{3pt}
\[
\quad
f\dff \times_c\dff \id_{\dff U}
\qff \colon\qff
C\dff \times_c\dff U\ttoo X\dff \times_c\dff Y
\]

\vspace*{-9pt}
is also a test map.\oss 
Since $A$ is compactly closed,\oss 
$(\dff f\dff \times_c\dff \id_{\dff U}\fff)^{\fff -\dff 1}\fff(\fff A\fff)$\qss
is closed in\qss $C\dff \times_c\dff U$\dnsp.\oss
Since $U$ is compact\halfff,\oss the projection of\dss this preimage to $C$ is closed in $C$\nnsp.\oss
But this projection is equal to\qss $f^{\fff -\dff 1}\fff(\fff B \fff)$\dnsp.\oss
It\dss follows that\qss $f^{\fff -\dff 1}\fff(\fff B \fff)$\qss is closed for any test map $f$
and\dss hence $B$ is closed.\oss
Since\qss $x\qff \not\in\qff B$\qss
and\qss $y\qff \in\qff U$\dnsp,\oss
the set\qss $(\dff X\dff \smallsetminus\dff B\dff)\dff \times_c\dff U$\qss
is a neighborhood of\qss $(x\fff,\pff y)$\qss disjoint from $A$\nnsp.\oss
Since\qss $(x\fff,\pff y)$\qss is an arbitrary point in the complement of $A$\nnsp,\oss
it follows that $A$ is closed in\qss $X\dff \times_c\dff Y$\dnsp.\oss  \eproof

\mypar{Lemma.}{times-quot}
\emph{Suppose that\pss $X\fff,\pff Y\fff,\pff Z$\qss are compactly generated spaces and\qss
$p\dff \colon\dff X\ttoo Y$\qss is a quotient map\qss
({\fff}i.e.\dss $p$\dss is a continuous surjective map and\dss $Y$\dss
has the quotient topology with respect to\dss $p$\nnsp).\oss
Then\oss \emph{$p\times_k\dff \id_{\dff Z}\qff \colon\qff X\dff \times_k\dff Z\ttoo Y\dff\times_k\dff Z$}\oss
is a quotient map.\oss}

\proof\qss 
If spaces\qss $X\fff,\pff Y\fff,\pff Z$\qss are Hausdorff\halfff,\oss
there is a direct proof in the same spirit as the proof of Lemma\qss \ref{times-lc}.\oss
See\qss \cite{st},\oss Theorem\qss 4.4.\oss
In the general case it is more convenient to reduce this lemma 
to a theorem about function spaces.\oss
See\qss \cite{td},\oss Theorem\qss 7.9.19.\oss  \eproof

\myuppar{Weakly Hausdorff spaces.}
A topological space $X$ is called\dss \emph{weakly Hausdorff}\pss
if\sss $f\fff(\fff C\fff)$\sss is closed for every test map $f\dff \colon\dff C\toto X$\nnsp.\oss
In the category of compactly generated spaces weakly Hausdorff spaces play the same
role as Hausdorff spaces in the category of all topological spaces.\oss
By Lemma\qss \ref{diagonal}\qss below
weakly Hausdorff compactly generated spaces
can be defined in the same way as separated schemes in\dss the algebraic geometry.\oss
The impatient readers may skip the next three lemmas\qss
(which will not be used directly)\qss and\dss go
to the last subsection of\dss this section.

\mypar{Lemma.}{wh}\vspace{-12pt}
\begin{itemize}
\item[(i)] \emph{A Hausdorff space is weakly Hausdorff\halfff.}
\item[(ii)] \emph{A weakly Hausdorff space is a $T_{\fff 1}$\dnsp-space.}
\item[(iii)] \emph{A subspace of a weakly Hausdorff space is weakly Hausdorff\halfff.}
\item[(iv)] \emph{If\pss $X$\qss is weakly Hausdorff\halfff,\oss then\qss $k\fff(\dff X\dff)$\qss 
is also weakly Hausdorff\halfff.}
\item[(v)] \emph{If\pss $X$\qss is a weakly Hausdorff space,\oss 
then\dss 
$f\fff(\fff C\fff)$\dss is Hausdorff\qss 
for every test map\qss $f\dff \colon\dff C\toto X$\nnsp.}
\item[(vi)] \emph{If\pss $X\fff,\pff Y$\qss are weakly Hausdorff spaces,\oss
then\qss $X\times_c Y$\qss is weakly Hausdorff\halfff.}
\item[(vii)] \emph{If\pss $X\fff,\pff Y$\qss are weakly Hausdorff spaces,\oss
then\qss $X\times_k Y$\qss is weakly Hausdorff\halfff.}
\end{itemize}

\vspace{-6pt}
\proof\qss
The properties\qss (i)--(iii)\qss are trivial.\oss
Let us prove\qss (iv).\oss
Suppose that $X$ is weakly Haus\-dorff\halfff.\oss
If\qss
$f\dff \colon\dff C\ttoo k\fff(\dff X\dff)$\qss 
is a test map,\oss
then by Lemma\qss \ref{test-maps}\qss $f$ is also a test map\qss $C\ttoo X$\nnsp.\oss
Since $X$ is weakly Hausdorff\halfff,\oss
the image\dss $f\fff(\fff C\fff)$\dss is closed in $X$\nnsp.\oss
It follows that\dss $f\fff(\fff C\fff)$\dss is closed in $k\fff(\dff X\dff)$\nnsp.\oss
This proves that $k\fff(\dff X\dff)$ is weakly Hausdorff\halfff.\oss 
This proves\qss (iv).\oss

Let us prove\qss (v).\oss
Suppose that\qss $x\fff,\pff y\qff \in\qff  f\fff(\fff C\fff)$\qss and\qss $x\qff \neq\qff y$\nnsp.\oss
By Lemma\qss \ref{wh}\dff(ii)\qss the sets\qss $\{\dff x\dff\}$\qss and\qss $\{\dff y\dff\}$\qss are closed
and hence the preimages\qss $f^{\dff -\dff 1}(x)$\qss and\qss $f^{\dff -\dff 1}(y)$\qss are closed.\oss
Since these preimages are disjoint and the space $C$\nnsp,\oss being compact and Hausdorff\halfff,\oss
is normal,\oss there exist disjoint open sets\qss $U\fff,\pff V\qff \subset\qff C$\qss
containing\qss $f^{\dff -\dff 1}(x)\fff,\pff f^{\dff -\dff 1}(y)$\qss respectively.\oss
Then\qss $C\dff \smallsetminus\dff U$\qss
and\qss $C\dff \smallsetminus\dff V$\qss are closed in $C$ and hence are compact and Hausdorff\halfff.\oss
Since $X$ is weakly Hausdorff\halfff,\oss the images\qss $f\fff(\fff C\dff \smallsetminus\dff U\fff)$\qss
and\qss $f\fff(\fff C\dff \smallsetminus\dff V\fff)$\qss of these complements are closed\dss
and\dss hence the complements\qss
$f\fff(\fff C\fff)\dff \smallsetminus\dff f\fff(\fff C\dff \smallsetminus\dff U\fff)$\qss
and\qss
$f\fff(\fff C\fff)\dff \smallsetminus\dff f\fff(\fff C\dff \smallsetminus\dff V\fff)$\qss
of these images
in\dss $f\fff(\fff C\fff)$\dss are open in\dss $f\fff(\fff C\fff)$\dnsp.\oss
Since these complements are disjoint and contain $x$ and $y$ respectively,\oss
this proves that there are two disjoint open sets in\dss $f\fff(\fff C\fff)$\dss 
containing the points\qss $x\fff,\pff y$\qss respectively.\oss
This proves\qss (v).\oss  

The claim\qss (vi)\qss follows from\qss (v)\qss and 
the standard properties of\dss the usual product\dss
$\times_c$\nsp,\oss
and\qss (vi)\qss immediately follows from\qss (iv)\qss and\qss (vi).\oss  \eproof

\mypar{Lemma.}{cg-for-wh}
\emph{Suppose that\pss $X$\qss is weakly Hausdorff\halfff.\oss
Then\trs $A$\sss is compactly closed in\qss $X$\qss if and only if\dss
for every compact Hausdorff subspace\trs $D$\sss of\pss $X$\qss the intersection\qss
$A\dff \cap\dff D$\qss is closed\dss in\trs $D$\nnsp.\oss}

\proof\qss
Suppose that $A$ is compactly closed.\oss
If $D$ is a compact Hausdorff subspace of $X$\nnsp,\oss
then the inclusion map\qss $D\toto X$\qss is a test maps.\oss
Since\dss $A\dff \cap\dff D$\dss is the preimage of $A$ under this inclusion map,\pss
$A\dff \cap\dff D$\dss is closed in $D$\nnsp.\oss

Conversely,\oss suppose that $A$ satisfies the condition from the lemma.\oss
If\qss $f\dff \colon\dff C\ttoo X$\qss is a test map,\oss
then\dss $f\fff(\fff C\fff)$\dss is Hausdorff\dss by Lemma\qss \ref{wh}\dff(v)\qss
and hence is a compact Hausdorff subspace of $X$\nnsp.\oss
Therefore\qss $f\fff(\fff C\fff)\dff \cap\dff A$\qss is closed in\dss $f\fff(\fff C\fff)$\dss
and\dss hence\pss
$f^{\dff -\dff 1}(\dff A\dff)
\off =\off
f^{\dff -\dff 1}(\dff f\fff(\fff C\fff)\dff \cap\dff A \dff)$\pss
is closed in $C$\nnsp.\oss
It follows that $A$ is compactly closed.\oss  \eproof

\mypar{Lemma.}{diagonal}
\emph{If\pss $X$\qss is compactly generated,\oss
then\pss $X$\qss is weakly Hausdorff\dss if and only if the diagonal\oss
$D_{\dff X}\off =\off \{\qff (\fff x\fff,\pff x\fff)\qff \mid\qff x\qff \in\qff X \qff\}$\oss
is closed in\pss $X\times_k X$\nnsp.\oss}\qss

\proof\qss
Suppose that $X$ is weakly Hausdorff\halfff.\oss
Let\qss
$f\dff \colon\dff C\ttoo X\times_k X$\qss
be an arbitrary test map.\oss
Let\qss $f_{\dff 1}\fff,\pff f_{\dff 2}$\qss be the components of $f$\dnsp.\oss
By Lemma\qss \ref{wh}\dff(v)\qss the images\qss
$f_{\dff 1}\fff(\fff C\fff)$\qss and\qss $f_{\dff 2}\fff(\fff C\fff)$\qss
are compact Hausdorff subspaces,\oss
and\dss hence\pss
$K\off =\off 
f_{\dff 1}\fff(\fff C\fff)\qff \cup\qff f_{\dff 2}\fff(\fff C\fff)$\pss
is also compact Hausdorff\halfff.\oss
It follows that\dss $D_{\dff K}$\dss is closed in\dss $K\times_c\dff K$\dss
and\dss hence in\dss  $K\times_k\dff K$\nnsp.\oss
It follows that\dss $f^{\fff -\dff 1}\fff(\trf D_{\dff K} \dff)$\dss is closed.\oss
But\oss
$f^{\fff -\dff 1}\fff(\trf D_{\dff X} \dff)
\off =\off
f^{\fff -\dff 1}\fff(\trf D_{\dff K} \dff)$\oss 
and\dss therefore\dss $f^{\fff -\dff 1}\fff(\trf D_{\dff X} \dff)$\dss is closed.\oss
Since\dss $X\times_k X$\dss is compactly generated,\oss 
it follows that\dss $D_{\dff X}$\dss is closed in\dss $X\times_k\dff X$\nnsp.\oss

Suppose now that\dss $D_{\dff X}$\dss is closed in\dss $X\times_k\dff X$\nnsp.\oss
Let\qss
$f\dff \colon\dff C\ttoo X$\qss
be a test map.\oss
We need to show that the image $f\fff(\fff C\fff)$ is closed in $X$\nnsp.\oss
Let\qss
$g\dff \colon\dff D\ttoo X$\qss
be some other test map.\oss
Then\oss\vspace*{3pt}
\[
\quad
f\times g
\dff \colon\dff 
C\times_c D
\off =\off 
C\times_k D\ttoo X\times_k X
\]

\vspace*{-9pt}
is a test map into\dss $X\times_k X$\nnsp.\qff\oss
Therefore\dss
$(\dff f\times g \dff)^{\fff -\dff 1}\fff(\dff D_{\dff X} \dff)$\dss
is closed in\dss $C\times_c D$\nnsp.\oss
Since\vspace*{3pt}
\[
\quad
g^{\fff -\dff 1}\fff(\trf f\fff(\fff C\fff) \dff)
\off =\off
\mbox{pr}_{\dff 2}\fff
\left(\dff
(\dff f\times g \dff)^{\fff -\dff 1}\fff(\dff D_{\dff X} \dff)
\dff\right)
\]

\vspace*{-9pt}
and\dss the spaces\qss $C\fff,\pff D$\qss are compact Hausdorff space,\oss
this implies that\dss
$g^{\fff -\dff 1}\fff(\trf f\fff(\fff C\fff) \dff)$\dss
is closed in\dss $C$\nnsp.\oss
Since $X$ is compactly generated,\oss 
it follows that $f\fff(\fff C\fff)$ is closed in $X$\nnsp.\oss
Since $f$ is an arbitrary test map,\oss this proves that $X$ is weakly Hausdorff\halfff.\oss  \eproof

\myuppar{CW-complexes.}
At the end of the day,\oss all spaces used in this paper turn out to be Hausdorff\halfff,\oss
at least if the bounded cohomology theory is limited to Hausdorff spaces.\oss
But the proofs depend on the theory
of\dss topological groups in the category of compactly generated spaces and their
classifying spaces developed by\dss McCord\qss \cite{mcc}\dss 
in the context of weakly Hausdorff compactly generated spaces.\oss
See Section\qss \ref{classifying-spaces}\qss for an outline of this theory.\oss

The topological groups needed are,\qss in general,\qss topological groups in the category
of compactly generated spaces\qss
but not the topological groups in the usual sense.\oss
The reason is that they are\dss CW-complexes
and\dss the product $\times_c$ of\dss two\dss CW-complexes 
is not necessarily a\dss CW-complex\halfff.\oss

We will need the classifying spaces of these topological groups only when they happen to be CW-com\-plexes.\oss
While usually\dss CW-com\-plexes are assumed to be Hausdorff\dss from the very beginning,\oss
this assumption is superfluous.\oss
If\qss CW-complexes are defined in terms of the consecutive attaching of cells,\oss
then one can prove that they are normal,\pss and\dss in particular\halfff,\pss Hausdorff spaces.\oss
See\qss \cite{h},\oss Proposition\qss A.3.\oss
If\qss CW-complexes are defined in terms of the characteristic maps of cells,\oss
it is sufficient to assume that the space is only weak Hausdorff\halfff.\oss
Lemma\qss \ref{wh}\dff(v)\qss assures that
the proof of the equivalence of this definition with the other one
applies to weakly Hausdorff spaces.\oss
See\qss \cite{h},\oss the proofs of\dss Propositions\qss A.1\qss and\qss A.2.\oss

\newpage
\mysection{Weakly\qss principal\qss and\qss principal\qss bundles}{w-principal}

\vspace*{6pt}
\myuppar{Weakly principal bundles.}
Suppose that $G$ is simultaneously a group and a topological space.\oss
Initially no conditions relating the group structure with the topological structure are imposed.\oss
The product in $G$ is denoted by the dot\halfff:\pss
$(a\fff,\pff b)\qff \longmapsto\qff a\dff \cdot\dff b$\nnsp.\oss
A map\qss $p\dff \colon\dff E\ttoo B$ is said to be a\qss
\emph{weakly principal bundle}\qss with the fiber $G$\sss
if the following conditions hold.\vspace*{-6pt}
\begin{itemize}
\item[(i)] \dnsp$p\dff \colon\dff E\ttoo B$\qss is a Serre fibration,\oss
i.e. the homotopy covering property holds for cubes,\oss
or\halfff,\oss what is the same,\oss
for standard simplices $\Delta_{\fff n}$\nnsp,\oss $n\qff \geq\qff 0$\nnsp.\oss
\item[(ii)] A free action of $G$ on $E$ is given.\oss
No continuity requirements  
are imposed at this stage.\oss
The fibers\dss $p^{\fff -\dff 1}\fff(y)$\dnsp,\oss
$y\qff \in\qff E$\nnsp,\oss are orbits of this action 
and\dss hence\qss $B\qff =\qff E/G$\qss as a\qss \emph{set}\halfff.\oss
\item[(iii)] If\qss $\tau\dff \colon\dff \Delta_{\fff n}\ttoo E$\qss
is a singular simplex in $E$\nnsp,\oss
then every singular simplex\qss $\sigma\dff \colon\dff \Delta_{\fff n}\ttoo E$\qss
such that\oss $p\dff \circ\dff \sigma\off =\off p\dff \circ\dff \tau$\oss
has the form\oss 
\begin{equation}
\label{lifts}
\quad
\sigma
\qff \colon\qff
x\qff \longmapsto\qff 
g\fff(x)\dff\cdot\dff \tau\fff(x)
\end{equation}
for some\qss \emph{continuous}\qss map\qss $g\dff \colon\dff \Delta_{\fff n}\ttoo G$\nnsp.\oss
We will abbreviate\qss (\ref{lifts})\qss as\qss
$\sigma\off =\off g\dff \cdot\dff  \tau$\nnsp.\oss
\end{itemize}

\vspace*{-6pt}
By the condition\qss (ii)\qss 
for any two set-theoretic maps\qss
$\sigma\fff,\pff \tau\qff \colon\dff \Delta_{\fff n}\ttoo E$\qss
such that\oss 
$p\dff \circ\dff \sigma\off =\off p\dff \circ\dff \tau$\oss
there is a unique map\qss $g\dff \colon\dff \Delta_{\fff n}\ttoo G$\qss
such that\qss $\sigma\qff =\qff g\dff \cdot\dff  \tau$\nnsp.\oss
The point of the condition\qss (iii)\qss is in the requiring $g$ to be continuous if
$\tau$ and $\sigma$ are continuous.\oss

For\qss $n\qff =\qff 0$\qss the continuity requirement in\qss (iii)\qss is vacuous,\oss
and this special case of\qss (iii)\qss follows from\qss (ii).\oss
In fact\halfff,\oss this special case means that
every non-empty fiber is an orbit of $G$\nnsp.\oss

The property\qss (ii)\qss implies 
that all fibers are non-empty.\oss
In view of this
the property\qss (i)\qss implies that
for every singular simplex\qss $\rho\dff \colon\dff \Delta_{\fff n}\ttoo B$\qss
there exists a\qss \emph{lift}\pss of $\rho$ to $E$\nnsp,\oss
i.e.\qss a singular simplex\qss $\tau\dff \colon\dff \Delta_{\fff n}\ttoo E$\qss
such that\qss $\rho\qff =\qff p\dff \circ\dff \tau$\nnsp.\oss

\myuppar{What is really needed.}
The definition of weakly principal bundles presents an intermediate ground between
the classical definition of principal $G$\dnsp-bundles for topological groups $G$ 
and the properties needed for the applications we have in mind.\oss
The continuity requirement in\qss (iii)\qss isn't crucial.\oss
What is really needed is the following.\vspace*{-3pt}
\begin{itemize}
\item[(a)]
There are sets\dss $L_{\fff n}$\dss
of maps\qss $\Delta_{\fff n}\ttoo G$\qss with the following property.\oss
If\qss $\tau\dff \colon\dff \Delta_{\fff n}\ttoo E$\qss is a singular simplex\halfff,\oss
then every lift of $p\dff \circ\dff \tau$ has the form\qss
$g\dff \cdot\dff  \tau$\qss
with\dss $g\qff \in\qff L_{\fff n}$\nnsp.\oss
\item[(b)] Each set\dss $L_{\fff n}$\dss is invariant under the action of\dss the group\qss
$\Sigma_{\fff n\dff +\qff 1}$\qss of\dss symmetries of\dss $\Delta_{\fff n}$\nnsp.
\item[(c)] The sets\dss $L_{\fff n}$\dss form a\qss 
\emph{coherent system}\pss in the following sense.\oss
For every simplicial embedding\qss 
$i\dff \colon\dff \Delta_{\fff n\dff -\dff 1}\ttoo \Delta_{\fff n}$\qss
the set $L_{\fff n\dff -\dff 1}$ is equal to the set\qss
$\{\qff g\dff \circ\dff i \qff \mid\qff g\qff \in\qff L_{\fff n} \qff\}$\nnsp.\oss
\end{itemize}

\vspace{-3pt}
\myuppar{Locally trivial principal bundles.}
Suppose that $G$ is a
topological group in the sense of the category of compactly generated spaces.\oss
This means $G$ is a group and a compactly generated space,\oss
the inverse\qss $a\qff \longmapsto\qff a^{\fff -\dff 1}$\qss is a continuous map\qss
$G\ttoo G$\nnsp,\oss 
and the multiplication is continuous as a map\qss
$G\dff\times_k\dff G\ttoo G$\nnsp,\oss
but may be not continuous as a map\qss
$G\dff\times_c\dff G\ttoo G$\nnsp.\oss

For the rest of this section we assume that\qss
$p\dff \colon\dff E\ttoo B$\qss 
is a continuous map of compactly generated spaces,\oss
and that an action of $G$ on $E$\nnsp,\oss
continuous as a map\qss
$E\dff \times_k\dff G\ttoo E$\nnsp,\oss
is given.\oss
The map $p$ together with the action of $G$ on $E$\dss
is said to be a\qss
\emph{locally trivial principal bundle}\qss 
in the category of compactly generated spaces
if $B$ can be covered by open subsets $U$ 
such that $p$ is\qss \emph{trivial}\pss over $U$\dnsp,\oss i.e.\qss there exists
a homeomorphism $h$ such that the diagram\dss\vspace*{6pt}
\begin{equation*}
\quad
\begin{tikzcd}[column sep=normal, row sep=huge]\dis
p^{\fff -\dff 1}\fff(\fff U \fff) \arrow[rr, "\dis h"]
\arrow[rd,  
"\dis p"']
&
& 
U \times_k\dff G
\arrow[ld,  
"\dis p_{\dff U}"]
\\
&
U
& 
\end{tikzcd}
\end{equation*}

\vspace*{-3pt}
(where\dss $p_{\dff U}$\dss is the projection to the factor $U$\dnsp)\pss
is commutative
both as the diagram of topological spaces and as the diagram of $G$\dnsp-sets,\oss
where $G$ acts on the set $U\times G$ in the obvious manner\halfff.\oss

\mypar{Lemma.}{p-wp}
\emph{If\dss $p$\dss is a locally trivial principal bundle,\oss
then\dss $p$\dss is a weakly principal\dss bundle.\oss}

\proof\qss
Suppose that $p$ is a locally trivial principal bundle.\oss

Since\dss $U \times_k\dff G$\dss is the product in the category of compactly generated spaces,\oss
the standard proof of the fact that every locally trivial bundle is a Serre fibration
applies without any changes to this situation.\oss 
One can also use the fact that for a compact Hausdorff space $X$\dss
(in particular\halfff,\pss for a cube)\qss the map\qss $X\ttoo U \times_k\dff G$\qss
is continuous if an only if it is continuous as a map\qss $X\ttoo U \times_c\dff G$\nnsp.\oss
This allows to deduce Serre's property in the present framework from
the Serre's property for the usual locally trivial bundles.\oss
Therefore $p$ satisfies the condition\qss (i).

Since $G$ acts freely on the product\dss $U \times G$\dss
and\qss
$U\qff =\qff (\dff U \times G \dff)/G$\qss as a set\halfff,\pss
$p$ satisfies\qss (ii). 

It remains to prove\qss (iii).\oss
As we pointed out above,\oss
the map $g$ satisfying\qss (\ref{lifts})\qss exists and is unique.\oss
We need only to check that it is continuous.\oss
Since the continuity is a local property,\oss we may assume
that the image\qss $p\dff \circ\dff \tau\fff(\dff \Delta_{\fff n}\fff)$\qss
is contained in an open set $U$ such that $p$ is trivial over $U$\dnsp.\oss
By composing everything with the map $h$ from the above diagram,\oss
we see that it is sufficient to prove that\dss if\pss
$\tau\fff,\pff \sigma\dff \colon\dff \Delta_{\fff n}\ttoo U \times_k G$\pss
are continuous maps such that\qss\vspace*{3pt}
\[
\quad
p_{\dff U}\dff \circ\dff \tau
\off =\off 
p_{\dff U}\dff \circ\dff \sigma\dff,
\]

\vspace*{-9pt}
then the unique map\qss
$g\dff \colon\dff \Delta_{\fff n}\ttoo G$\qss
such that\qss
$\sigma\qff =\qff g\dff \cdot\dff \tau$\qss
is continuous.\oss

Let\oss
$f
\off =\off 
p_{\dff U}\dff \circ\dff \tau
\off =\off 
p_{\dff U}\dff \circ\dff \sigma$\nnsp.\oss
Since $\Delta_{\fff n}$ is a compact Hausdorff space
and\dss hence is a compactly generated space,\pss
there exist continuous maps\qss
$t\fff,\pff s\dff \colon\dff \Delta_{\fff n}\ttoo G$\qss
such that\vspace*{3pt}
\[
\quad
\tau\fff(x)\qff =\qff (\dff f\fff(x)\fff,\pff t\fff(x) \dff)
\hspace*{1em}\mbox{ and }\hspace*{1em}
\sigma\fff(x)\qff =\qff (\dff f\fff(x)\fff,\pff s\fff(x) \dff)
\]

\vspace*{-9pt}
for all\qss $x\qff \in\qff \Delta_{\fff n}$\nnsp.\oss
Obviously,\pss
$g\fff(x)\qff =\qff s\fff(x)\dff \cdot\dff t\fff(x)^{\fff -\dff 1}$\qss
for all\qss $x\qff \in\qff \Delta_{\fff n}$\nnsp.\oss
Since\dss $G\times_k G$\dss is a product in the category of compactly generated spaces,\oss
the map\qss 
$x\qff \longrightarrow (\dff s\fff(x)\fff,\pff t\fff(x) \dff)$\qss
is continuous as a map\qss
$\Delta_{\fff n}\ttoo G\times_k G$\nnsp.\oss
Since the inverse map\qss $G\ttoo G$\qss 
and the product map\qss $G\times_k G\ttoo G$\qss are continuous,\oss
it follows that $g$ is continuous.\oss  \eproof

\myuppar{Numerable bundles.}
Recall that a covering of a space is called\qss \emph{numerable}\qss
if there exists a partition of unity subordinated to this covering.\oss
The bundle\qss
$p\dff \colon\dff E\ttoo B$\qss 
is said to be\qss \emph{numerable}\qss
if\trs there exists an open numerable covering $\mathcal{U}$ of $B$ such
that $p$ is trivial over every element\qss $U\qff \in\qff \mathcal{U}$\dnsp.\oss
As\sss is\sss well known,\qss $p$\dss is\sss numerable\sss if\dss $B$ is\sss paracompact\halfff.\oss

The remaining part of\trs this section\sss 
is\sss devoted to the proof\dss of\trs Theorem\qss \ref{cw-bundle}\qss below.\oss
The author hoped to find it in textbooks,\oss but failed even with research papers.

\mypar{Lemma.}{induced-bundle}
\emph{Suppose that\dss
locally trivial principal bundle\qss
$p\dff \colon\dff E\ttoo B$\qss is numerable 
and\dss that\qss
$f\dff \colon\dff B'\ttoo B$\qss is a continuous map.\oss
Let\qss
$p'\dff \colon\dff E'\ttoo B'$\qss
be the bundle induced from $p$ by $f$\nnsp.\oss
If\qss $f$ is a homotopy equivalence,\oss then the canonical map\qss
$f^\sim\dff \colon\dff E'\ttoo E$\qss is a homotopy equivalence.\oss}

\proof\qss 
We refer to\qss \cite{td}\qss for the standard results of the bundle theory.\oss
Their proofs in\qss \cite{td}\qss apply without any changes to compactly generated spaces
and locally trivial bundles in the category of compactly generated spaces,\oss
if one takes into the account the fact that by Lemma\qss \ref{times-lc}\qss 
$X\dff \times_c\dff [\fff 0\fff,\qff 1\fff]$\qss
is compactly generated when $X$ is compactly generated.\oss

Let\qss $g\dff \colon\dff B\ttoo B'$\qss be a homotopy equivalence inverse to $f$
and\dss let\qss
$p''\dff \colon\dff E''\ttoo B$\qss
be the bundle induced from $p'$\dss by $g$\nnsp.\oss 
Let\qss
$g^\sim\dff \colon\dff E''\ttoo E'$\qss
be the canonical map.\oss
The bundle $p''$\dss is induced from the bundle $p$ by\qss $f\dff \circ\dff g$\nnsp.\oss
Since $p$ is numerable and\qss $f\dff \circ\dff g$\qss is homotopic to the identity\qss $\id_{\dff B}$\nnsp,\oss
the bundle $p''$\dss is isomorphic to $p$ over $B$\nnsp.\oss 
See\qss \cite{td},\oss Theorem\qss 14.3.2.\oss
By composing $g^\sim$ with an isomorphism\qss
$E\ttoo E''$\qss from $p$ to $p''$\dss we get a bundle map\vspace*{3pt}
\begin{equation*}
\quad
\begin{tikzcd}[column sep=huge, row sep=huge]\dis
E
\arrow[r,  
"\dis h"]
\arrow[d,  
"\dis p"]
& 
E'
\arrow[d,  
"\dis p'"]
\\
B
\arrow[r,  
"\dis g"]
& 
B'
\end{tikzcd}
\end{equation*}

\vspace*{-6pt}
The composition\qss $f^\sim \circ\dff h$\qss is a bundle map covering\qss
$f\dff \circ\dff g$\nnsp.\oss
Since\qss $f\dff \circ\dff g$\qss is homotopic to the identity\qss $\id_{\dff B}$\nnsp,\oss
the homotopy lifting theorem\qss
(see\qss \cite{td},\oss 14.3.4)\qss
implies that\qss $f^\sim \circ\dff h$\qss is homotopic to a bundle map\qss 
$j\dff \colon\dff E\ttoo E$\qss covering
the identity map\qss $\id_{\dff B}$\nnsp.\oss
Since the map $j$ covers a homeomorphism,\oss it is itself a homeomorphism\qss
(see \cite{td},\oss remarks after Proposition\qss 14.1.6).\oss
It follows that\qss 
$f^\sim \circ\dff h 
\dff \colon\dff
E\ttoo E$\qss is a homotopy equivalence.\oss

Since the bundle $p'$\dss is induced from a numerable bundle,\oss it is itself numerable.\oss
Therefore the homotopy lifting theorem applies to the bundle map\qss 
$h\dff \circ\dff f^\sim
\dff \colon\dff
E'\ttoo E'$\qss
covering\qss $g\dff \circ f$\qss and a homotopy connecting\qss $g\dff \circ f$\qss
with\qss $\id_{\dff B'}$\nnsp.\oss
It follows that\qss 
$h\dff \circ\dff f^\sim$\qss
is homotopic a bundle map\qss 
$j'\dff \colon\dff E'\ttoo E'$\qss covering
the identity map\qss $\id_{\dff B'}$\nnsp.\oss
As above,\oss the map $j'$ is a homeomorphism,\oss
and\dss hence\qss 
$h\dff \circ\dff f^\sim
\dff \colon\dff
E'\ttoo E'$\qss
is a homotopy equivalence.\oss
Since both maps\qss 
$f^\sim \circ\dff h$\qss
and\qss
$h\dff \circ\dff f^\sim$\dss
are homotopy equivalences,\oss
the map $f^\sim$\dss is a homotopy equivalence,\oss
as claimed.\oss  \eproof

\myuppar{CW-complexes.}
Suppose now that $B$ is a\dss CW-complex\halfff,\oss
and\dss let\dss $B_n$\dss be the $n$\dnsp-th skeleton of $B$\dnsp.\oss
Then\dss $B_n$\dss is obtained from\dss $B_{n\dff -\dff 1}$\dss by glueing a collection of
$n$\dnsp-dimensional discs.\oss
In more details,\oss
let us denote by\dss $\mathbb{D}^n$\dss the disjoint union of these discs,\oss
and\dss by\dss $\mathbb{S}^{n\dff -\dff 1}$\dss the disjoint union of their boundary spheres,\pss
$\mathbb{S}^{n\dff -\dff 1}\qff \subset\qff \mathbb{D}^n$\dnsp.\oss
Then\dss $B_n$\dss is obtained
by glueing to\dss $B_{n\dff -\dff 1}$\dss the space\dss $\mathbb{D}^n$\dss by a continuous map\qss
$\alpha_{\fff n}\qff \colon\qff \mathbb{S}^{n\dff -\dff 1}\ttoo B_{n\dff -\dff 1}$\dnsp.\oss
Let\qss
$B_{n\dff -\dff 1}\qff \sqcup\qff \mathbb{D}^n$\qss
be the disjoint union of\qss $B_{\fff n\dff -\dff 1}$\dss and\qss $\mathbb{D}^n$\nnsp.\oss
There is a continuous map\vspace*{3pt}
\[
\quad
\varphi_{n}
\qff \colon\qff
B_{\fff n\dff -\dff 1}\qff \sqcup\qff \mathbb{D}^n
\ttoo
B_{\fff n}
\]

\vspace*{-9pt}
equal to the inclusion
on\dss $B_{\fff n\dff -\dff 1}$\nnsp,\oss
to\dss $\alpha_{\fff n}$\dss on\dss $\mathbb{S}^{n\dff -\dff 1}\qff \subset\qff \mathbb{D}^n$\nnsp,\oss
and inducing a homeomorphism\vspace*{3pt}
\[
\quad
\mathbb{D}^n\qff \smallsetminus\qff \mathbb{S}^{n\dff -\dff 1}
\ttoo
B_n\qff \smallsetminus\qff B_{n\dff -\dff 1}\dff.
\]

\vspace*{-9pt}
Moreover\halfff,\pss the topology of\dss $B_n$\dss is the quotient topology of\qss
$B_{n\dff -\dff 1}\qff \sqcup\qff \mathbb{D}^n$\qss
induced\dss by\dss $\varphi_{n}$\nnsp.\oss

\myuppar{Locally trivial principal bundles over\dss CW-complexes.}
We continue to assume that $B$ is a\dss CW-complex
and keep the above notations related to $B$\nnsp.\oss
Let\pss $E_{\fff n}\off =\off p^{\fff -\dff 1}\fff(\dff B_n \fff)$\pss
and\dss let\qss\vspace*{3pt}
\[
\quad
\pi_{\fff n}\dff \colon\dff \mathbf{E}_{\fff n}\ttoo B_{\fff n\dff -\dff 1}\qff \sqcup\qff \mathbb{D}^n
\]

\vspace*{-9pt}
be the bundle induced from $p$\dss by\dss $\varphi_{n}$\nnsp.\oss 
Then\oss 
$\mathbf{E}_{\fff n}
\off =\off
E_{\fff n\dff -\dff 1}
\qff \sqcup\off 
\mathbb{E}_{\fff n}$\nnsp,\qff\oss where\qff\oss 
$\mathbb{E}_{\fff n}
\off =\off
\pi_{\fff n}^{\fff -\dff 1}\fff(\dff \mathbb{D}^n\fff)$\dnsp.\oss
Let\oss
$\Phi_n
\qff \colon\qff
E_{\fff n\dff -\dff 1}
\qff \sqcup\off 
\mathbb{E}_{\fff n}
\off =\off
\mathbf{E}_{\fff n}
\ttoo
E_{\fff n}$\oss
be the canonical map.\qff\oss
Then the diagram\vspace*{3pt}
\begin{equation*}
\quad
\begin{tikzcd}[column sep=huge, row sep=huge]\dis
E_{\fff n\dff -\dff 1}\qff \sqcup\qff \mathbb{E}_{\fff n}
\arrow[r,  
"\dis \Phi_n"]
\arrow[d,  
"\dis \pi_{\fff n}"]
& 
E_{\fff n}
\arrow[d,  
"\dis p"]
\\
B_{\fff n\dff -\dff 1}\qff \sqcup\qff \mathbb{D}^n
\arrow[r,  
"\dis \varphi_{n}"]
& 
B_n
\end{tikzcd}
\end{equation*}

\vspace*{-6pt}
is commutative.\oss

\mypar{Lemma.}{cw-bundle-topology}
\emph{\dnsp$\Phi_n$\qss
is a quotient map and\dss hence\qss
$E_{\fff n}$\qss can be obtained by glueing\pss
$\mathbb{E}_{\fff n}$\qss
to\qss $E_{n\dff -\dff 1}$\qss
along a continuous map\qff\oss
$\pi_{\fff n}^{\fff -\dff 1}\fff(\dff \mathbb{S}^{n\dff -\dff 1}\dff)
\qff\ttoo\qff
E_{n\dff -\dff 1}$\nsp.\oss}

\proof\qss
It is sufficient to prove that $B_n$ can be covered by open sets $U$
such that\dss $\Phi_n$\dss is a quotient map over\dss $U$\dnsp,\oss
i.e.\qss such that the map\vspace*{3pt}
\[
\quad
\chi
\qff \colon\qff
\pi_{n}^{\fff -\dff 1}\fff
\left(\dff
\varphi_{n}^{\fff -\dff 1}\fff(\dff U \dff)
\dff\right)
\ttoo
p^{\fff -\dff 1}\fff(\dff U \dff)
\]

\vspace*{-9pt}
induced by\dss $\Phi_n$\dss is a quotient map.\oss
Since $p$ is locally trivial,\oss it is sufficient to show that $\chi$ is a quotient map\dss
if the bundle\qss $E_{\fff n}\ttoo B_n$\qss is trivial over $U$\dnsp.\oss
Suppose that this is the case.\oss
Then\dss $\pi_n$\dss is trivial over\qss
$\varphi_{n}^{\fff -\dff 1}\fff(\dff U \dff)$\qss and\dss
there is a commutative diagram\vspace*{3pt}
\begin{equation*}
\quad
\begin{tikzcd}[column sep=huge, row sep=huge]\dis
\pi_{n}^{\fff -\dff 1}\fff
\left(\dff
\varphi_{n}^{\fff -\dff 1}\fff(\dff U \dff)
\dff\right)
\arrow[r, "\dis \chi"]
\arrow[d,  
"\dis h'"]
& 
p^{\fff -\dff 1}\fff(\dff U \dff)
\arrow[d,  
"\dis h"]
\\
\varphi_{n}^{\fff -\dff 1}\fff(\dff U \dff)\dff \times_k\dff G
\arrow[r]
& 
U\dff \times_k\dff G\qff,
\end{tikzcd}
\end{equation*}

\vspace*{-6pt}
where 
the vertical arrows are homeomorphisms
and the lower horizontal arrow is the map\qss
$\varphi_{n}^{\fff -\dff 1}\fff(\dff U \dff)\ttoo U$\qss
multiplied by\dss $\id_{\dff G}$\nnsp.\oss
Since\dss $B$\dss is a\dss CW-complex\halfff,\oss
the map\qss
$\varphi_{n}^{\fff -\dff 1}\fff(\dff U \dff)\ttoo U$\qss
is a quotient map.\oss
By Lemma\qss \ref{times-quot}\qss this implies that the lower horizontal arrow is a quotient map.\oss
Hence the top horizontal arrow $\chi$ is also a quotient map.\oss  \eproof

\mypar{Lemma.}{cw-bundle-colimit}
\emph{The topology of\pss $E$\qss is the weak topology defined by the subspaces\qss
$E_{\fff n}$\nnsp.\oss}

\proof\qss
The proof is similar to the proof of Lemma\qss \ref{cw-bundle-topology}.\oss
It is sufficient to prove this claim locally,\oss over open sets\qss
$U\qff \subset\qff B$\qss such that $p$ is trivial over $U$\dnsp.\oss
Therefore,\oss it is sufficient to prove that\qss
$U\dff \times_k\dff G$\qss has the weak topology defined by subspaces\qss
$(\dff U\dff \cap\dff B_n\fff)\dff \times_k\dff G$\nnsp.\oss
Since the space $B$ has the weak topology defined by the subspaces $B_n$\nnsp,\oss
the subspace $U$ has the weak topology defined by the subspaces\dss $U\dff \cap\dff B_n$\nnsp.\oss
It remains to use the fact that this property survives multiplication\qss
$\times_k Z$\qss by any compactly generated space $Z$\nnsp.\oss
See\qss \cite{st},\oss Theorem\qss 10.1.\oss  \eproof

\mypar{Lemma.}{cw-bundle-total}
\emph{If\qss the\dss fiber\qss $G$\qss of\qss the bundle\qss
$p\dff \colon\dff E\ttoo B$\qss is\dss a\dss CW-complex\halfff,\oss
then\dss $E$\dss is\dss homotopy equivalent\dss to\sss a\dss CW-complex.\oss}

\proof\qss
The space\dss $E_{\fff 0}$\dss is a\dss CW-complex\halfff,\oss
being a disjoint union of fibers which are assumed to be\dss CW-complexes.\oss
Let\qss $X_{\fff 0}\qff =\qff E_{\fff 0}$\nnsp.\oss
Suppose that  
CW-complexes\dss $X_{\fff m}$\dss and homotopy equivalences\qss 
$f_m\dff \colon\dff E_{\fff m}\ttoo X_{\fff m}$\qss
are already defined\dss for\qss $m\qff \leq\qff n\qff -\qff 1$\qss 
and\dss that\dss $X_{\fff l}$\dss
is a subcomplex of\dss $X_{\fff m}$\dss and $f_m$ extends\dss $f_l$\qss if\qss $l\qff \leq\qff m$\nnsp.\oss

By Lemma\qss \ref{cw-bundle-topology}\qss
$E_{\fff n}$\dss results from glueing\pss
$\mathbb{E}_{\fff n}$\qss
to\qss $E_{n\dff -\dff 1}$\qss
along a continuous map\qff\oss\vspace*{3pt}
\[
\quad
g_{\fff n}
\qff \colon\qff
\pi_{\fff n}^{\fff -\dff 1}\fff(\dff \mathbb{S}^{n\dff -\dff 1}\dff)
\qff\ttoo\qff
E_{n\dff -\dff 1}\qff.
\]

\vspace*{-9pt}
Since all components of\dss $\mathbb{D}^n$\dss are contractible,\oss
the induced bundle\dss $\pi_{\fff n}$\dss is trivial over\dss $\mathbb{D}^n$\dss
and\dss hence there exists a homeomorphism\dss\vspace*{3pt}
\[
\quad
\mathbb{E}_{\fff n}
\qff\ttoo\qff
\mathbb{D}^n\dff \times_k\dff G
\]

\vspace*{-9pt}
taking\qss
$\pi_{\fff n}^{\fff -\dff 1}\fff(\dff \mathbb{S}^{n\dff -\dff 1}\dff)$\qss
to\qss
$\mathbb{S}^{n\dff -\dff 1}\dff \times_k\dff G$\nnsp.\oss
We may treat this homeomorphism as an identification.\oss
Since $G$ is assumed to be a\dss CW-complex\halfff,\pss
the product\qss
$\mathbb{D}^n\dff \times_k\dff G$\qss admits a structure of a\dss CW-complex such that\qss
$\mathbb{S}^{n\dff -\dff 1}\dff \times_k\dff G$\qss
is a subcomplex\halfff.\oss
Therefore $E_{\fff n}$ is the result of glueing of the CW-complex\qss
$\mathbb{D}^n\dff \times_k\dff G$\qss
to\dss $E_{\fff n\dff -\dff 1}$\dss by the map\dss $g_{\fff n}$\dss defined on the subcomplex\qss
$\mathbb{S}^{n\dff -\dff 1}\dff \times_k\dff G$\nnsp.\oss

This implies,\oss in particular\halfff,\oss that\dss $E_{\fff n\dff -\dff 1}$\dss
is a neighborhood deformation retract of\dss $E_{\fff n}$\dss in the sense of\qss \cite{st},\oss
and hence satisfies the homotopy extension property\qss
(see\qss \cite{po},\oss Lecture\qss 2,\oss Prop\-o\-si\-tion\qss 3,\oss for example).\oss
Hence the homotopy equivalence\qss 
$f_{n\dff -\dff 1}\dff \colon\dff E_{\fff n\dff -\dff 1}\ttoo X_{\fff n\dff -\dff 1}$\qss 
extends to a homotopy equivalence between $E_{\fff n}$ and the result of glueing of\qss
$\mathbb{D}^n\dff \times_k\dff G$\qss to\dss $X_{\fff n\dff -\dff 1}$\dss by\vspace*{3pt}
\[
\quad
f_{n\dff -\dff 1}\dff \circ\dff g_{\fff n}
\qff \colon\qff
\mathbb{S}^{n\dff -\dff 1}\dff \times_k\dff G
\qff\ttoo\qff
X_{\fff n\dff -\dff 1}
\]

\vspace*{-9pt}
See\qss \cite{td},\oss Proposition\qss 5.1.10,\oss 
or\qss \cite{po},\oss Lecture\qss 2,\oss p.\qss 92,\oss
or\qss \cite{b},\oss Theorem\qss 7.5.7.\oss

Let $h_{\fff n}$ be a cellular map homotopic to\qss $f_{n\dff -\dff 1}\dff \circ\dff g_{\fff n}$\dnsp,\oss
and\dss let\dss $X_{\fff n}$\dss be the result of glueing of the\dss CW-complex\qss
$\mathbb{D}^n\dff \times_k\dff G$\qss to\dss $X_{\fff n\dff -\dff 1}$\dss by $h_{\fff n}$\nnsp.\oss 
Then\dss $X_{\fff n}$\dss is a\dss CW-complex containing\dss $X_{\fff n\dff -\dff 1}$\qss
as a subcomplex and\dss
$E_{\fff n}$\dss is homotopy equivalent to\dss $X_{\fff n}$\dss also.\oss
Moreover\halfff,\pss $f_{n\dff -\dff 1}$\dss extends to a homotopy equivalence\qss
$f_{n}\dff \colon\dff E_{\fff n}\ttoo X_{\fff n}$\dss
(see\qss \cite{b},\oss Corollary\qss 7.5.5,\oss for example).\oss

Let $X$ be the union of the\dss CW-complexes $X_{\fff n}$\nnsp,\oss
and let\qss 
$f\dff \colon\dff E\ttoo X$\qss
be the map equal to $f_n$ on $X_{\fff n}$\nnsp.\oss
Then $X$ has a unique structure of a\dss CW-complex such that every $X_{\fff n}$
is a subcomplex.\oss
The topology of $X$ is the weak topology defined by the subspaces $X_{\fff n}$\nnsp.\oss
Lemma\qss \ref{cw-bundle-colimit}\qss implies that $f$ is continuous.\oss
For every $n$ the subspace\dss $X_{\fff n\dff -\dff 1}$\nnsp,\oss being a subcomplex of\dss $X_{\fff n}$\nnsp,\oss
is a neighborhood deformation retract of\dss $X_{\fff n}$\nnsp.\oss
Also,\oss as we proved above,\oss for every $n$ the subspace\dss $E_{\fff n\dff -\dff 1}$\dss
is a neighborhood deformation retract of\dss $E_{\fff n}$\nnsp.\oss
Therefore,\oss the fact that each $f_n$ is a homotopy equivalence implies that $f$
is a homotopy equivalence.\oss
See,\oss for example,\oss \cite{td},\oss Proposition 5.2.9,\oss
or\qss \cite{po},\oss Lecture\qss 10,\oss Additional\dss Material,\oss Theorem\qss 1\qss 
and\dss Proposition\qss 2.\oss  \eproof

\mypar{Theorem.}{cw-bundle}
\emph{Suppose that\qss
$p\dff \colon\dff E\ttoo B$\qss is a numerable locally trivial principal bundle with the fiber\dss $G$\nnsp.\oss
If\qss $G$\dss is\dss a\dss CW-complex and\qss $B$\dss is homotopy equivalent to a\dss CW-complex,\oss
then\dss $E$\dss is homotopy equivalent to a\dss CW-complex.\oss}

\proof\qss
This follows from Lemmas\qss \ref{cw-bundle-total}\qss and\qss \ref{induced-bundle}.\oss  \eproof

\mysection{McCord\qss classifying\qss spaces\qss and\qss principal\qss bundles}{classifying-spaces}

\vspace*{6pt}
\myuppar{The spaces\qss $B\fff(\fff G\fff,\pff X\dff)$\dnsp.}
Let $X$ be a weakly Hausdorff compactly generated space 
and\dss let $G$ be a topological group in the category of
weakly Hausdorff compactly generated spaces.\oss
Suppose that a base point $*$ of $X$ is chosen.\oss
Let\dss $B\fff(\fff G\fff,\pff X\dff)$\dss be the set of functions\qss
$u\dff \colon\dff X\ttoo G$\qss 
such that\qss $u\fff(*)\qff =\qff e$\nnsp,\oss where $e$ is the unit of $G$\nnsp,\oss
and\qss $u\fff(x)\qff =\qff e$\qss for all\dss but finitely many points $x\qff \in\qff X$\nnsp.\oss
Then\dss $B\fff(\fff G\fff,\pff X\dff)$\dss is itself a group with respect to the
point-wise multiplication $\oplus$ of maps\qss $X\ttoo G$\nnsp,\oss
defined by the formula\qss 
$(\fff u \oplus v\fff)\fff(x)\qff =\qff u\fff(x)\dff \cdot\dff v\fff(x)$\dnsp.\oss
The unit of this group is the map\dss $\mathbf{e}$\dss
taking the value $e$ at all points of $X$\nnsp.\oss
Let\qss $g\qff \in\qff G$\nnsp.\oss 
For\qss $x\qff \in\qff X$\nnsp,\pss
$x\qff \neq\qff *$\nnsp,\oss
let\dss $g\fff x$\dss be the map\qss $X\ttoo G$\qss taking the value $g$ at $x$
and the value $e$ at all other points of $X$\nnsp.\oss
Let\dss $g\dff *\qff =\qff \mathbf{e}$\nnsp.\oss

The set\dss $B\fff(\fff G\fff,\pff X\dff)$\dss is equipped with a topology in the following manner\halfff.\oss
For each\qss $n\qff \geq\qff 0$\nnsp,\oss
let\dss be\dss $B_{\fff n}\fff(\fff G\fff,\pff X\dff)$\dss
the subset of functions $u$ taking the value of $e$ at all\dss but\qss $\leq\qff n$\qss points of $X$\nnsp.\oss
Clearly,\pss $B_{\fff 0}\fff(\fff G\fff,\pff X\dff)\qff =\qff \{\dff \mathbf{e} \dff\}$\qss
and\dss for\qss $n\qff \geq\qff 1$\qss the set\dss $B_{\fff n}\fff(\fff G\fff,\pff X\dff)$\dss
consists of elements of the form\qss\vspace*{3pt}
\[
\quad
g_{\fff 1}\fff x_{\fff 1}\dff \oplus\dff g_{\fff 2}\fff x_{\fff 2}\dff \oplus\dff 
\ldots\dff \oplus\dff g_{\fff n}\fff x_{\fff n}\dff.
\]

\vspace*{-9pt}
The space\dss $B_{\fff n}\fff(\fff G\fff,\pff X\dff)$\dss
is equipped with the topology of a quotient space induced by the map\vspace*{3pt}
\[
\quad
\mu_{\fff n}
\qff \colon\qff
(\dff G\dff \times\dff X\dff)^{\fff n}\ttoo B_{\fff n}\fff(\fff G\fff,\pff X\dff)\dff,
\]

\vspace*{-9pt}
defined by the formula\oss\vspace*{3pt}
\[
\quad
\mu_{\fff n}\fff
(\dff (\fff g_{\fff 1}\fff,\pff x_{\fff 1}\fff)\fff,\pff
(\fff g_{\fff 2}\fff,\pff x_{\fff 2}\fff)\fff,\pff \ldots\fff,\pff (\fff g_{\fff n}\fff,\pff x_{\fff n}\fff)\dff)
\off =\off
g_{\fff 1}\fff x_{\fff 1}\dff \oplus\dff g_{\fff 2}\fff x_{\fff 2}\dff \oplus\dff 
\ldots\dff \oplus\dff g_{\fff n}\fff x_{\fff n}\dff.
\]

\vspace*{-9pt}
Each\dss
$B_{\fff n}\fff(\fff G\fff,\pff X\dff)$\dss
is a closed subspace of\dss $B_{\fff n\dff +\dff 1}\fff(\fff G\fff,\pff X\dff)$\dss
by Lemma\qss 6.2\qss of\qss \cite{mcc}.\oss
Therefore\vspace*{3pt}
\[
\quad
B_{\fff 0}\fff(\fff G\fff,\pff X\dff)
\qff \subset\qff 
B_{\fff 1}\fff(\fff G\fff,\pff X\dff)
\qff \subset\qff 
B_{\fff 2}\fff(\fff G\fff,\pff X\dff)
\qff \subset\qff \ldots
\]

\vspace*{-9pt}
is a sequence of spaces, whose union is equal to\dss 
$B\fff(\fff G\fff,\pff X\dff)$\dnsp.\oss
The set\dss $B\fff(\fff G\fff,\pff X\dff)$\dss
is equipped with the weak topology defined by 
the spaces\dss $B_{\fff n}\fff(\fff G\fff,\pff X\dff)$\dnsp.\oss
If\dss $X$\dss is homeomorphic to\qss $S^{\fff 0}$\dnsp,\oss i.e.\qss
if\dss $X$\dss consists of\dss two points,\oss 
then\oss 
$B\fff(\fff G\fff,\pff X\dff)
\off =\off 
B_{\fff 1}\fff(\fff G\fff,\pff X\dff)
\off =\off
G$\nnsp.\oss

\mypar{Theorem.}{b-space}
\emph{The space\pss $B\fff(\fff G\fff,\pff X\dff)$\dss
is a weakly Hausdorff compactly generated space.\oss
If\qss $G$ is a topological abelian group,\oss
then\qss $B\fff(\fff G\fff,\pff X\dff)$\dss
is also a topological abelian group.\oss}

\proof\qss
See\qss \cite{mcc},\oss Lemma\qss 6.5\qss
and Proposition\qss 6.6.\oss  \eproof

\mypar{Theorem.}{cw}
\emph{Suppose that\qss $G$\dss is a discrete abelian group
and that\qss $X$\dss is a simplicial complex equipped with the weak topology.\oss
Then\pss $B\fff(\fff G\fff,\pff X\dff)$\dss admits structure of a CW-complex.\oss}

\proof\qss 
See\qss \cite{mcc},\oss Section\qss 7.\oss  \eproof

\myuppar{Induced maps.} Let $G$ be a topological group and\qss 
$\varphi\dff \colon\dff X\ttoo Y$\qss
be a continuous map of based weakly Hausdorff compactly generated spaces.\oss
If either $G$ is abelian or the map\dss $\varphi$\dss
is\dss injective on\qss $X\dff \smallsetminus\dff \{\dff * \dff\}$\nnsp,\oss 
then\dss $\varphi$\dss
induces a continuous map\vspace*{3pt}
\[
\quad
\varphi_{*}\dff \colon\dff 
B\fff(\fff G\fff,\pff X\dff)\ttoo B\fff(\fff G\fff,\pff Y\dff)
\]

\vspace*{-9pt}
acting by the formula\vspace*{3pt}
\[
\quad
\varphi_*\fff
\left(\dff
g_{\fff 1}\fff x_{\fff 1}\dff \oplus\dff g_{\fff 2}\fff x_{\fff 2}\dff \oplus\dff 
\ldots\dff \oplus\dff g_{\fff n}\fff x_{\fff n}
\dff\right)
\off =\off
g_{\fff 1}\fff \varphi\fff(x_{\fff 1})\dff \oplus\dff g_{\fff 2}\fff \varphi\fff(x_{\fff 2})\dff \oplus\dff 
\ldots\dff \oplus\dff g_{\fff n}\fff \varphi\fff(x_{\fff n})\dff.
\]

\vspace*{-9pt}
The assumptions about\dss $G$\dss and\dss $\varphi$\dss are needed\dss for $\varphi_*$ to be correctly defined.\oss
In order to see this,\oss it is better to describe $\varphi_*$ in terms of functions\qss
$X\fff,\pff Y\ttoo G$\nnsp.\oss
In these terms\vspace*{3pt}
\[
\quad
\varphi_*\fff(u)(y)
\off =\off
\sum_{x\qff \in\qff \varphi^{\fff -\dff 1}\fff(y)}\qff u\fff(x)\dff.
\]

\vspace*{-9pt}
There is no natural way to order the points in the preimage\dss $\varphi^{\fff -\dff 1}\fff(y)$\dnsp.\oss
In order for the sum to be correctly defined,\oss
it is sufficient to assume either that $G$ is abelian,\oss
or that every sum consists of $0$ or $1$ terms,\oss
except\halfff,\oss perhaps,\oss the sum corresponding to\qss $y\qff =\qff *$\nnsp.\oss
In other terms,\oss it is sufficient to  assume that either $G$ is abelian,\oss or $\varphi$ 
is\dss injective on\qss $X\dff \smallsetminus\dff \{\dff * \dff\}$\nnsp.\oss

\mypar{Theorem.}{induced-maps}
\emph{Let\qss $G$\dss be an abelian topological group and\qss
$\varphi\dff \colon\dff X\ttoo Y$\qss
be a continuous map of based weakly Hausdorff compactly generated spaces.\oss
If\qss $\varphi$\dss is a quotient map,\oss then}\vspace*{3pt}
\[
\quad
\varphi_{*}\dff \colon\dff 
B\fff(\fff G\fff,\pff X\dff)\ttoo B\fff(\fff G\fff,\pff Y\dff)
\] 

\vspace*{-9pt}
\emph{is also a quotient map.\oss
If\dss $\varphi$ is a closed injective map,\oss
then $\varphi_*$ also has this property.\oss}

\proof\qss
See\qss \cite{mcc},\oss Proposition\qss 6.7.\oss  \eproof

\mypar{Theorem.}{induced-homotopies}
\emph{Let\qss $G$\dss be an abelian topological group and\qss
$\varphi_t\dff \colon\dff X\ttoo Y$\dnsp,\pss $t\qff \in\qff [\fff 0\fff,\qff 1 \fff]$\qss
be a homotopy of continuous maps of based weakly Hausdorff compactly generated spaces.\oss
Then}\qss\vspace*{3pt}
\[
\quad
(\varphi_t)_{*}\dff \colon\dff 
B\fff(\fff G\fff,\pff X\dff)\ttoo B\fff(\fff G\fff,\pff Y\dff)
\] 

\vspace*{-9pt}
\emph{is a homotopy.\oss
If\pss $X\fff,\pff Y$\trs are homotopy equivalent\halfff,\oss
then\qss $B\fff(\fff G\fff,\pff X\dff)$\dss and\qss $B\fff(\fff G\fff,\pff Y\dff)$\dss
are homotopy equivalent\halfff.\oss
In particular\halfff,\oss if\pss $X$\trs is contractible,\oss
then so is\qss $B\fff(\fff G\fff,\pff X\dff)$\dnsp.\oss}

\proof\qss
See\qss \cite{mcc},\oss Proposition\qss 6.10.\oss  \eproof

\myuppar{Principal bundles.}
Let\dss $G$\dss be a discrete group and\dss
$(\dff X\fff,\pff A \dff)$\dss be a pair of\dss based
simplicial complexes equipped with the weak topology.\oss
Suppose\sss that\sss either\dss the group\dss $G$\dss is\dss abelian,\oss
or\qss
$(\dff X\fff,\pff A \dff)$\dss
is homeomorphic to\dss
$(\dff I\fff,\pff \partial\fff I \dff)$\nnsp,\oss 
where\qss $I$\trs is an interval and\qss $\partial\fff I$\qss is its boundary.\oss
Let\qss $i\dff \colon\dff A\ttoo X$\qss be the inclusion map and\qss 
$p\dff \colon\dff X\ttoo X/A$\qss
be the quotient map.\oss
The maps\dss $i$\dss and\dss $p$ induce maps\vspace*{-3pt}
\begin{equation*}
\quad
\begin{tikzcd}[column sep=normal, row sep=huge]\dis
B\fff(\fff G\fff,\pff A\dff) \arrow[r, "\dis i_*"]
&
B\fff(\fff G\fff,\pff X\dff) \arrow[r, "\dis p_*"]
& 
B\fff(\fff G\fff,\pff X/A\dff)
\end{tikzcd}
\end{equation*}

\vspace*{-9pt}
If\dss $G$\dss is abelian,\oss then
Theorem\qss \ref{induced-maps}\qss implies that $i_*$ is a closed embedding and $p_*$ is a quotient map.\oss
The map $i_*$ together with the group operation $\oplus$ on\dss $B\fff(\fff G\fff,\pff X\dff)$\dss
define an action of\dss $B\fff(\fff G\fff,\pff A\dff)$\dss on\dss $B\fff(\fff G\fff,\pff X\dff)$\dnsp.\oss
By\qss \cite{mcc},\oss Lemma\qss 8.3,\pss $p_*$ induces a canonical homeomorphism\vspace*{3pt}
\[
\quad
B\fff(\fff G\fff,\pff X \dff)\bigl/\dff B\fff(\fff G\fff,\pff A \dff)\bigr.
\ttoo
B\fff(\fff G\fff,\pff X/A \dff)\dff.
\]

\vspace*{-9pt}
If\dss
$(\dff X\fff,\pff A \dff)$\dss
is homeomorphic to\dss
$(\dff I\fff,\pff \partial\fff I \dff)$\nnsp,\oss
this is true if the topology of\dss
$B\fff(\fff G\fff,\pff I\dff)$\dss
and\dss
$B\fff(\fff G\fff,\pff \partial\fff I\dff)$\dss
is defined slightly differently.\oss
See\qss \cite{mcc},\oss Section\qss 9.\oss
With this topology\dss $B\fff(\fff G\fff,\pff I\dff)$\dss is contractible\qss
and since\dss $B\fff(\fff G\fff,\pff \partial\fff I\dff)$\dss is equal to $G$\nnsp,\oss
the map $p_*$ is the universal cover of a\dss 
$K\dff(\dff G\fff,\pff 1\dff)$\dnsp-space.

\mypar{Theorem.}{main-bundles}
\emph{Under the above assumptions,\oss
the map\qss
$p_*\qff \colon\qff
B\fff(\fff G\fff,\pff X \dff)
\ttoo
B\fff(\fff G\fff,\pff X/A \dff)$\qss
is a numerable locally trivial principal
bundle with the fiber\qss $B\fff(\fff G\fff,\pff A \dff)$\dnsp.\oss}

\proof\qss
See\qss \cite{mcc},\oss Theorems\qss 8.8\qss 
and\qss 9.17.\oss  \eproof

\myuppar{Eilenberg-MacLane spaces and universal bundles.}
Let $G$ be a discrete group and $n$ be a natural number\halfff.\oss
Suppose that either $G$ is abelian or\qss $n\qff =\qff 1$\nnsp.\oss
Let\dss $D^{n}$\dss be the standard $n$\dnsp-di\-men\-sional disc
and\dss $S^{n\dff -\dff 1}$\dss be its boundary sphere.\oss
Let us identify\qss $D_n/\fff S^{n\dff -\dff 1}$\dss
with the standard $n$\dnsp-sphere\dss $S^n$\dnsp.\oss
By\dss Theorem\qss \ref{main-bundles}\qss there is a principal numerable
bundle\vspace*{3pt}
\begin{equation}
\label{mccord-bundle}
\quad
p_{\fff G}^{\fff n}
\qff \colon\qff
B\fff(\fff G\fff,\pff D^{n} \dff)
\ttoo
B\fff(\fff G\fff,\pff S^n \dff)
\end{equation}

\vspace*{-9pt}
with the fiber\dss $B\fff(\fff G\fff,\pff S^{n\dff -\dff 1} \dff)$\dnsp.\oss
If\qss $n\qff =\qff 1$\nnsp,\oss
then the fiber is $G$ and the bundle is a covering space.\oss
By Theorem\qss \ref{induced-homotopies}\qss the space\dss
$B\fff(\fff G\fff,\pff D^{n} \dff)$\dss 
is contractible.\oss
By Section\qss \ref{w-principal}\qss the bundle\qss $p_{\fff G}^{\fff n}$\qss
is a Serre fibration.\oss
An induction by $n$\nnsp,\oss 
starting with\oss
$B\fff(\fff G\fff,\pff S^0 \dff)\off =\off G$\nnsp,\oss
and using the homotopy sequence of the Serre fibration\qss $p_{\fff G}^{\fff n}$\qss
to\dss go\dss from\qss $n\qff -\qff 1$\qss to\dss $n$\nnsp,\oss
shows that\dss
$B\fff(\fff G\fff,\pff S^n \dff)$\dss
is an Eilenberg--MacLane space of the type\dss $K\fff(\fff G\fff,\pff n\fff)$\dnsp.\oss
Cf.\qss \cite{mcc},\oss Corollary\qss 10.6.\oss
Since\dss $B\fff(\fff G\fff,\pff D^{n} \dff)$\dss is contractible,\oss
the bundle\dss $p_{\fff G}^{\fff n}$\dss is a\qss \emph{universal bundle}\pss 
in the category of compactly generated spaces.\oss
Its base\dss $B\fff(\fff G\fff,\pff S^n \dff)$\dss is a\dss
$K\fff(\fff G\fff,\pff n\fff)$\dnsp-space and its fiber\dss 
$B\fff(\fff G\fff,\pff S^{n\dff -\dff 1} \dff)$\dss is a\dss
$K\fff(\fff G\fff,\pff n\qff -\qff 1 \fff)$\dnsp-space.\oss
By Section\qss \ref{w-principal}\qss
$p_{\fff G}^{\fff n}$\dss 
is a weakly principal bundle with the fiber\dss
$B\fff(\fff G\fff,\pff S^{n\dff -\dff 1} \dff)$\dnsp.\oss

If\qss $n\qff \geq\qff 2$\nnsp,\oss
then the spaces\qss
$B\fff(\fff G\fff,\pff D^{n} \dff)$\nnsp,\pss
$B\fff(\fff G\fff,\pff S^n \dff)$\nnsp,\pss
$B\fff(\fff G\fff,\pff S^{n\dff -\dff 1} \dff)$\qss
are\dss CW-complexes by Theorem\qss \ref{cw}.\oss
For the further use,\oss let us replace the above\dss $p_{\fff G}^{\fff 1}$\dss
by the universal cover of a\dss CW-complex of\dss the type $K\dff(\dff G\fff,\pff 1\dff)$\dnsp.\oss
Obviously,\oss it is a weakly principal bundle with the fiber $G$\nnsp.

\mysection{Spaces\qss with\qss amenable\qss fundamental\qss group}{amenable}

\vspace*{6pt}
\myuppar{Means.}
For a set $S$ we will denote by $B\fff(\fff S\fff)$ the vector space of all bounded real functions on $S$\nnsp.\oss 
As is well known,\qss $B\fff(\fff S\fff)$\dss is a Banach space with the norm\qss 
\[
\quad
\|\dff f\trf\|
\off =\off\sup\nolimits_{\fff s\qff \in\qff S}\qff 
|\dff f\dff(s\fff)\trf|\qff.
\] 
A linear functional\pss $m\dff \colon\qff B\fff(\fff S\fff)\toto \rrr$\pss 
is called a\qss \emph{mean on}\qss $B\fff(\fff S\fff)$\trs if\oss
$|\qff m\fff(\dff f \dff)\qff|\qff \leq\qff \|\qff f \qff\|$\oss
for all\qss $f\qff \in\qff B\fff(\fff S\fff)$\qss and\qss
$m\fff(\fff \mathbf{1}\fff)
\off =\off
1$\nnsp,\oss
where\qss $\mathbf{1}\dff(s\fff)\off =\off 1$\qss for all\qss $s\qff \in\qff S$\dnsp.\oss
In other words,\pss $m$ is a bounded functional such that\qss
$\|\qff m \qff\|\qff \leq\qff 1$\qss 
and\qss
$m\fff(\fff \mathbf{1}\fff)
\off =\off
1$\nnsp.\oss
Usually means are defined by requiring that
\begin{equation}
\label{classical-mean-definition}
\quad
\inf\qff f
\off \leq\off
m\fff(\dff f \dff)
\off \leq\off
\sup\dff f
\end{equation}
for all\dss $f\qff \in\qff B\fff(\fff S\fff)$\dnsp,\oss 
where the infimum and the supremum are taken over the set $S$\dnsp.\oss
This definition motivates the term\qss \emph{mean}.\oss
The two definitions are equivalent by the following lemma.

\mypar{Lemma.}{definitions-means}
\emph{A linear functional\dss $m$\dss
is a mean if and only if\qss (\ref{classical-mean-definition})\qss holds for all\qss
$f\qff \in\qff B\fff(\fff S\fff)$\dnsp.\oss}

\proof\qss
Let\dss $f\qff \in\qff B\fff(\fff S\fff)$\dnsp,\oss
and\qss $f'\off =\off f\qff -\pff \inf\qff f$\dnsp.\oss
Then\qss  
$\|\qff f' \qff\|\off =\off \sup\dff f\qff -\pff \inf\qff f$\qss
and\dss hence\oss
\[
\quad
m\fff(\dff f \dff)
\off =\off
m\fff(\dff f' \dff)\qff +\qff \inf\qff f
\off \leq\off
\left(\dff
\sup\dff f\qff -\qff \inf\qff f
\dff\right)
\qff +\pff
\inf\dff f
\off =\off
\sup\dff f\dff.
\]
By applying this inequality to\dss $-\qff f$\qss in the role of $f$\dnsp,\oss
we see that\qss
$\inf\dff f\qff \leq\qff m\fff(\dff f \dff)$\dnsp.\oss
Conversely,\oss suppose that\qss (\ref{classical-mean-definition})\qss holds.\oss
This immediately implies that\qss
$m\fff(\fff \mathbf{1}\fff)
\off =\off
1$\nnsp.\oss
Since\qss $\|\qff f \qff\|$\qss is equal to the largest of\dss the numbers\qss
$\sup\dff f$\qss and\qss 
$-\qff \inf\qff f\off =\off \sup\dff(\dff -\qff f \dff)$\dnsp,\oss
the inequality\qss (\ref{classical-mean-definition})\qss applied to $f$ and\dss $-\qff f$\dss
implies that\qss $\|\qff m \qff\|\qff \leq\qff 1$\nnsp.\oss  \eproof

\myuppar{Invariant means and amenable groups.}
Let $G$ be a group acting on a set $S$ from the right\halfff.\oss 
Then $G$ acts on $B\fff(\fff S\fff)$ from the left  
by the formula\qss 
$g\cdot f\dff (\fff s\fff )\qff =\qff f\dff (\fff s\dff \cdot\dff g \dff)$\dnsp,\oss 
where\qss $g\qff \in\qff G$\nnsp,\qss $f\qff \in\qff B\fff(\fff S\fff)$\nnsp,\qss 
$s\qff \in\qff S$\nnsp.\oss 
A mean $m$ on\dss $B\fff(\fff S\fff)$\dss is called\qss \emph{right invariant}\pss if\qss 
$m\fff(\fff g\cdot f \dff)\qff =\qff m\fff(\fff f \dff)$\qss 
for all\qss 
$g\qff \in\qff G$\nnsp,\qss 
$f\qff \in\qff B\fff(\fff S\fff)$\dnsp.\oss 
Usually the right invariant means will be called simply\qss \emph{invariant means}.\oss
If there exists an invariant mean on\dss $B\fff(\fff G\fff)$\dss
with respect to the action of the group $G$ on itself by the right translations,\oss   
then $G$ is said to be\qss \emph{amenable}. 

One can define in an obvious way the notion of a\qss \emph{left invariant}\qss mean.\oss
It turns out that if\dss $G$ is amenable,\oss
then there exists a mean on\dss $B\fff(\fff G\fff)$\dss 
which is simultaneously right and left invariant\qss
(see\qss \cite{gr},\oss Lemmas\qss 1.1.1\qss and\qss 1.1.3),\oss
but we will not need this result\halfff.\oss

\mypar{Theorem.}{abelian}
\emph{If\dss a\sss group is abelian,\oss then it is amenable.}\oss

\proof\qss
See\qss \cite{gr},\oss Theorem\qss 1.2.1.\oss  \eproof

\myuppar{Free transitive actions.}
Suppose that a group $G$ acts on a set $S$ on the left
and\dss that this action is free and transitive.\oss
Then every\qss $s\qff \in\qff S$\qss defines a bijection\qss
$r_s\dff \colon\dff G\toto S$\qss by the rule\qss
$g\qff \longmapsto\qff g\cdot s$\nnsp,\oss where\qss $g\qff \in\qff G$\nnsp.\oss
If\qss $s\fff,\pff t\qff \in\qff S$\nnsp,\oss
then the bijections\dss $r_s$\dss and\dss $r_t$\dss differ by a right translation of $G$\nnsp.\oss 
Indeed,\oss since the action is transitive,\oss
there exists\qss $h\qff \in\qff G$\qss such that\qss $t\qff =\qff h\cdot s$\qss
and\dss hence\qss
\[
\quad
r_t\fff(\fff g\fff)
\off =\off
g\cdot t
\off =\off
g\cdot (\fff h\cdot s\fff)
\off =\off
(\fff g\cdot h\fff)\cdot s
\off =\off
r_s\fff(\fff g\cdot h\fff)\dff.
\]
It follows that\dss if\qss $f\qff \in\qff B\fff(\fff S\fff)$\dnsp,\oss
then\oss
\[
\quad
f\dff \circ\dff r_t\fff(\fff g\fff)
\off =\off
f\dff \circ\dff r_s\fff(\fff g\cdot h\fff)
\off =\off
h\cdot (f\dff \circ\dff r_s)\dff(\fff g\fff)
\]
for every\qss $g\qff \in\qff G$\nnsp.\oss
In other terms,\pss
$f\dff \circ\dff r_t
\off =\off
h\cdot (f\dff \circ\dff r_s)$\dnsp.\oss
If $m$ is a mean on\dss $B\fff(\fff G\fff)$\dss and\qss $s\qff \in\qff S$\nnsp,\oss
then\qss
$f\qff \longmapsto\qff m\dff(\fff f\dff \circ\dff r_s\fff)$\qss
is a mean on\dss $B\fff(\fff S\fff)$\dnsp.\oss
If\dss the mean $m$ is invariant and\qss $s\fff,\pff t\qff \in\qff S$\nnsp,\oss
then\oss
\[
\quad
m\dff(\fff f\dff \circ\dff r_t\fff)
\off =\off
m\dff(\fff h\cdot (f\dff \circ\dff r_s)\fff)
\off =\off
m\dff(\fff f\dff \circ\dff r_s\fff)\dff,
\]
and hence the mean\qss
$f\qff \longmapsto\qff m\dff(\fff f\dff \circ\dff r_s\fff)$\qss
on\dss $B\fff(\fff S\fff)$\dss
does not depend on the choice of $s$\nnsp.\oss
The mean\qss
$f\qff \longmapsto\qff m\dff(\fff f\dff \circ\dff r_s\fff)$\qss 
on\dss $B\fff(\fff S\fff)$\dss
is said to be\qss \emph{induced by $m$ and the action of\qss $G$ on $S$}\dss
(which should be free and transitive).\oss

\myuppar{Push-forwards.}
Let\qss $\pi\dff \colon\dff G\toto H$\qss be a surjective homomorphism.\oss
The map\qss
$f\qff \longmapsto\qff f\dff \circ\dff \pi$\qss 
is a bounded operator\qss
$B\fff(\fff H\fff)\toto B\fff(\fff G\fff)$\dnsp.\oss
In fact\halfff,\oss its norm is obviously equal to $1$\nnsp.\oss
If\dss $m$\dss is an invariant mean on\dss $B\fff(\fff G\fff)$\dnsp,\oss
then the map $\pi_{\fff *} m$ defined by the formula
\[
\quad
\pi_{\fff *} m
\qff \colon\qff
f\qff \longmapsto\qff m\fff(\dff f\dff \circ\dff \pi\dff)
\] 
is an invariant mean on\dss $B\fff(\fff H\fff)$\dss
called the\qss \emph{push-forward}\qss of $m$ by $\pi$\nnsp.\oss

\myuppar{Coherent sequences of\dss invariant means.}
Let $G$ be a topological group.\oss
Suppose that $G$ is either abelian,\oss or discrete and amenable.\oss
For a topological space $X$ let $G^{\dff X}$
be the group of continuous maps\qss $X\toto G$ considered as a discrete group.\oss
Recall that $\Delta_{\fff n}$ is the standard $n$\dnsp-dimensional simplex\halfff.\oss
If $G$ is abelian,\oss then the group\dss $G^{\dff \Delta_{\fff n}}$\dss is abelian,\oss
and\dss hence is amenable by Theorem\qss \ref{abelian}.\oss
If $G$ is discrete and amenable,\oss then\dss $G^{\dff \Delta_{\fff n}}$\dss
is equal to\dss $G$\dss and\dss hence is a\-me\-na\-ble.\oss
Let\dss $M_{\fff n}$\dss be the set of invariant means on\dss 
$B\fff(\fff G^{\dff \Delta_{\fff n}}\dff)$\dnsp,\oss
so\qss $M_{\fff n}\qff \subset\qff B\fff(\fff G^{\dff \Delta_{\fff n}}\dff)^*$\dnsp.\oss
Suppose that\qss
\[
\quad
m_{\fff 0}\fff,\off m_{\fff 1}\fff,\off \ldots\fff,\off m_{\fff n}\fff,\off \ldots
\]
is an either finite or infinite sequence of\dss invariant means\qss
$m_{\fff n}\qff \in\qff M_{\fff n}$\nnsp.\oss
If it is a finite sequence,\oss
then $n$ ranges over\oss
$0\fff,\pff 1\fff,\pff 2\fff,\pff \ldots\fff,\pff N$\oss
for some\qss $N\qff \geq\qff 0$\nnsp,\oss
if it is an infinite sequence,\oss
then $n$ ranges over all non-negative integers\oss
$0\fff,\pff 1\fff,\pff 2\fff,\pff \ldots\pff $\nnsp.\oss
If\qss 
$i\dff \colon\dff \Delta_{\fff n\trf -\dff 1}\toto \Delta_{\fff n}$\qss
is a simplicial embedding,\oss
then the map\qss
$\sigma\qff \longmapsto\qff \sigma\dff \circ\dff i$\qss
is a homomorphism 
\begin{equation*}
\quad
\pi^{\fff i}\dff \colon\dff G^{\dff \Delta_{\fff n}}\toto G^{\dff \Delta_{\fff n\trf -\dff 1}}\dff,
\end{equation*} 
which is obviously surjective.\oss
The sequence
is called\qss \emph{coherent}\oss if
\begin{equation}
\label{coherence}
\quad
\pi^{\fff i}_*\fff m_{\fff n}
\off =\off
m_{\fff n\trf -\dff 1}
\end{equation} 
for every $n$ 
such that\dss $m_{\fff n}$\dss is defined 
and every $i$ as above.\oss
Equivalently,\oss\vspace*{2pt} 
\begin{equation}
\label{coherence-one}
\quad
m_{\fff n}\fff(\dff f\dff \circ\dff \pi^{\fff i}\dff)
\off =\off
m_{\fff n\trf -\dff 1}\fff(\dff f\dff)
\end{equation}

\vspace*{-10pt}
for every simplicial embedding\qss 
$i\dff \colon\dff \Delta_{\fff n\trf -\dff 1}\toto \Delta_{\fff n}$\qss
and\qss 
$f\qff \in\qff B\fff(\fff G^{\dff \Delta_{\fff n\trf -\dff 1}}\dff)$\pss
if\dss $m_{\fff n}$\dss is defined.\oss

\myuppar{\dnsp$\Sigma_{\dff n\dff +\dff 1}$\nsp\dnsp-invariant means on\qss 
$B\fff(\fff G^{\dff \Delta_{\fff n}}\dff)$\dnsp.}
The natural action of the symmetric group\dss 
$\Sigma_{\dff n\dff +\dff 1}$\dss 
on\dss $\Delta_{\fff n}$\dss 
induces an action of the group\dss 
$\Sigma_{\dff n\dff +\dff 1}$\dss on the group\dss $G^{\dff \Delta_{\fff n}}$\nnsp.\oss
In turn,\oss 
this action induces an action of\trs 
$\Sigma_{\dff n\dff +\dff 1}$\sss 
on the dual space\dss 
$B\fff(\fff G^{\dff \Delta_{\fff n}}\dff)^*$\dnsp.\oss
Obviously,\oss this action leaves
the space of means on\dss 
$B\fff(\fff G^{\dff \Delta_{\fff n}}\dff)$\dss invariant\halfff.\oss
Moreover\halfff,\oss 
since the group\dss 
$\Sigma_{\dff n\dff +\dff 1}$\dss 
acts on\dss 
$G^{\dff \Delta_{\fff n}}$\dss by automorphisms,\oss 
this action leaves the set\dss $M_{\fff n}$\dss 
of\dss invariant means invariant\halfff.\oss
We will say that an invariant mean $m$ on\dss 
$B\fff(\fff G^{\dff \Delta_{\fff n}}\dff)$\dss
is\dss 
\emph{$\Sigma_{\dff n\dff +\dff 1}$\nsp\dnsp-invariant}\pss 
if $m$ is fixed by the action of\dss 
$\Sigma_{\dff n\dff +\dff 1}$\dss
on the dual space\dss 
$B\fff(\fff G^{\dff \Delta_{\fff n}}\dff)^*$\dnsp.\oss
If $m$ is an arbitrary invariant mean on\dss 
$B\fff(\fff G^{\dff \Delta_{\fff n}}\dff)$\dnsp,\oss 
then\vspace*{3pt}
\[
\quad
\frac{1}{\fff (n\qff +\qff 1)\fff!\fff}\off
\sum\nolimits_{\qff \sigma\qff \in\qff \Sigma_{\dff n\dff +\dff 1}}\dff 
m^\sigma
\]

\vspace*{-9pt}
is\dss a\dss $\Sigma_{\dff n\dff +\dff 1}$\nsp\dnsp-invariant mean,\oss 
where the action of\dss $\Sigma_{\dff n\dff +\dff 1}$\dss is written as 
$(m\fff,\pff \sigma)\qff\longmapsto\qff m^\sigma$\nnsp.\oss

Suppose that\dss $m_{\fff n}$\dss is\dss a\dss 
$\Sigma_{\fff n\dff +\dff 1}$\nsp\dnsp-invariant mean on\qss 
$B\fff(\fff G^{\dff \Delta_{\fff n}}\dff)$\qss
and that\oss
$i\dff \colon\dff \Delta_{\fff n\dff -\dff 1}\toto \Delta_{\fff n}$\oss
is a simplicial embedding.\oss
Let us consider the push-forward\vspace*{2pt} 
\[
\quad
m_{\fff n\dff -\dff 1}
\off =\qff\off 
\pi^{\fff i}_*\fff m_{\fff n}\dff.
\] 

\vspace*{-10pt}
Since\dss $m_{\fff n}$\dss is\dss 
$\Sigma_{\fff n\dff +\qff 1}$\dnsp-invariant\halfff,\oss
the push-forward 
$m_{\fff n\dff -\dff 1}$ is independent of the choice of the simplicial embedding\dss 
$i$\dss 
and\dss is\dss 
$\Sigma_{\fff n}$\dnsp-invariant\halfff.\oss 
Since $m_{\fff n\dff -\dff 1}$ is independent of the choice of $i$\nnsp,\oss
the coherence condition\qss (\ref{coherence})\qss holds for every
simplicial embedding\qss
$i\dff \colon\dff \Delta_{\fff n\dff -\dff 1}\toto \Delta_{\fff n}$\nnsp.\oss

\mypar{Theorem.}{coherent-exists-finite}
\emph{Let\qss $G$\dss be a topological group and\qss $N\qff \geq\qff 0$\nnsp.\oss
If\pss $G$\trs is either abelian,\pss or dis\-crete and amenable,\oss
then there exists a coherent sequence of\dss invariant means\qss
$m_{\fff 0}\fff,\pff m_{\fff 1}\fff,\pff m_{\fff 2}\fff,\pff 
\ldots\fff,\pff m_{\dff N}$\nnsp.}

\proof\qss
As we saw above,\oss
under assumptions of the theorem
the groups\dss $G^{\dff \Delta_{\fff n}}$\dss 
are amenable.\oss
In particular\halfff,\oss 
there exist an invariant mean $m_{\trf N}$ on\dss 
$B\fff(\fff G^{\dff \Delta_{\dff N}}\dff)$\dnsp.\oss 
Moreover\halfff,\oss $m_{\trf N}$\dss can be chosen to be\qss $\Sigma_{\fff n\dff +\qff 1}$\dnsp-invariant\halfff.\oss
Let\qss
$i\dff \colon\dff \Delta_{\dff N\qff -\qff 1}\ttoo \Delta_{\dff N}$\qss
be an arbitrary simplicial embedding\dss
and\dss let\vspace*{2pt}
\[
\quad
m_{\trf N\qff -\qff 1}
\off =\off
\pi^{\fff i}_*\fff m_{\trf N}\dff.
\]

\vspace*{-10pt} 
By repeating this construction 
we get a sequence\oss 
$m_{\trf N}\fff,\pff 
m_{\trf N\qff -\qff 1}\fff,\pff 
m_{\trf N\qff -\qff 2}\fff,\off 
\ldots \fff,\off 
m_{\dff 0}$\oss
of\dss invariant means satisfying\qss (\ref{coherence})\qss
for all\qss $n\qff \leq\qff N$\nnsp.\oss  \eproof

\myuppar{Coboundary maps.}
Let us review the basic notations related to the coboundary maps.
Let $Z$ be a topological space.\oss
For every\oss
$n\off =\off  1\fff,\pff 2\fff,\pff \ldots $\oss 
and\dss every\oss
$k\off =\off 0\fff,\pff 1\fff,\pff \ldots\fff,\pff n$\nnsp,\qff\oss 
let\oss\vspace*{3pt} 
\[
\quad
\partial_k
\qff \colon\qff
S_n\fff(\dff Z\dff)\ttoo S_{n\dff -\dff 1}\fff(\dff Z\dff)
\]

\vspace*{-9pt}
be the map taking a singular simplex\qss
$\tau\dff \colon\dff \Delta_{\fff n}\toto Z$\qss
into\dss its\dss $k$\dnsp-th\dss face\qss
$\partial_k\dff \tau\dff \colon\dff \Delta_{\fff n\dff -\qff 1}\toto Z$\nnsp.\oss
By the definition,\oss
$\partial_k\dff \tau
\off =\off
\tau\dff \circ\dff i_{\fff k}$\oss
for some simplicial embedding\qss
$i_{\fff k}\dff \colon\dff \Delta_{\fff n\dff -\qff 1}\toto \Delta_{\fff n}$\nnsp.\oss
The coboundary of a cochain 
$f\qff \in\qff C^{\fff n\dff -\dff 1}\fff(\dff Z\dff)$
is\dss the\dss cochain 
$\partial f\qff \in\qff C^{\fff n}\fff(\dff Z\dff)$ 
defined\dss by the formula\vspace*{3pt}
\[
\quad
\partial f\fff(\fff \tau\fff)
\off\off =\off\off
\sum_{k\qff =\qff 0}^{n}(-\qff 1)^k\dff f\dff(\dff \partial_k\dff \tau \dff)
\off\off =\off\off
\sum_{k\qff =\qff 0}^{n}(-\qff 1)^k\dff f\dff \circ\dff \partial_k\dff(\fff \tau\fff)\dff,
\]

\vspace*{-9pt}
where\qss
$\tau\qff \in\qff S_n\fff(\dff Z\dff)$\dnsp.\oss

\vspace{6pt}
\mypar{Theorem.}{chain-map}
\emph{Let\qss $G$\dss be a topological group.\oss
Suppose that\pss $G$\trs is either abelian,\oss or discrete and amenable.\oss
Let\qss 
$p\dff \colon\dff X\toto Y$\qss 
be a weakly principal bundle with the fiber\dss $G$\nnsp.\oss 
Then every coherent sequence\qss
$\{\dff m_{\fff n} \dff \}$\qss
leads to a sequence of homomorphisms}\qss\vspace*{3pt} 
\[
\quad
p_{\fff *}\dff \colon\dff B^{\fff n}\fff(\dff X \dff)\toto B^{\fff n}\fff(\dff Y \dff)
\]

\vspace*{-9pt}
\emph{defined for the same $n$ as\qss
$\{\dff m_{\fff n} \dff \}$\qss
and such that\oss
$p_{\fff *}\circ\fff p^{\fff *}\off =\dff\off \id$\nnsp,\oss
$\|\dff p_{\fff *} \dff\|\off =\off 1$\nnsp,\oss
and\qss $p_{\fff *}$\qss commutes with the coboundaries 
of the chain complexes\qss
$B^{\fff *}\fff(\dff X \dff)$\dnsp,\qss  
$B^{\fff *}\fff(\dff Y \dff)$\qss
whenever this makes sense.}

\vspace{3pt}
\proof\qss
In what follows the number $n$ is assumed to be such that\dss $m_{\fff n}$\dss
is defined.\oss

For each singular simplex\qss 
$\sigma\dff \colon\dff \Delta_{\fff n}\toto Y$\qss
let us denote by\dss $C_{\dff \sigma}$\dss 
the set of singular simplices\qss 
$\tau\dff \colon\dff \Delta_{\fff n}\toto X$\qss
such that\qss $p\dff \circ\dff \tau\qff =\qff \sigma$\nnsp.\oss
Since $p$ is a weakly principal bundle with the fiber\dss $G$\nnsp,\oss
the group\dss $G^{\dff \Delta_{\fff n}}$\dss
acts on\dss $C_{\dff \sigma}$\dss from the left
and this action is free and transitive.\oss
Therefore the mean $m_{\fff n}$ together with this action induces a mean on\dss
$B\fff(\fff C_{\dff \sigma}\fff)$\dnsp,\oss
which we denote by $m_{\dff \sigma}$\nnsp.\oss

Let us define a linear map\oss
$p_{\fff *}\dff \colon\dff B^{\fff n}\fff(\dff X\dff)\toto B^{\fff n}\fff(\dff Y\dff)$\oss
as follows.\oss
For a cochain\pss $f\qff \in\qff B^{\fff n}\fff(\fff X\fff)$\pss
let\pss
$p_{\fff *}\fff(\fff f\fff)\qff \in\qff B^{\fff n}\fff(\fff Y\fff)$\pss
be the cochain taking the value\vspace*{3pt}
\[
\quad
p_{\fff *}\fff(\fff f\fff)(\fff \sigma\fff)
\off =\off
m_{\dff \sigma}\fff(\fff f \mid C_{\dff \sigma}\fff)
\]

\vspace*{-9pt}
on each singular simplex\qss
$\sigma\dff \colon\dff \Delta_{\fff n}\toto Y$\dnsp,\oss
where\dss $f \mid C_{\dff \sigma}$\dss
is the restriction of $f$ to the set\dss $C_{\dff \sigma}$\nnsp.\oss
Since\qss
$m_{\dff \sigma}\dff \colon\dff B\fff(\fff C_{\dff \sigma}\fff)\toto \rrr$\qss
is a mean,\oss
$\|\dff p_{\fff *}\fff(\fff f\fff) \dff\|\qff \leq\qff \|\dff f \dff\|$\oss
for\dss all\qss
$f\qff \in\qff B^{\fff n}\fff(\fff Y\fff)$\dnsp.\oss
It follows that\qss
$\|\dff p_{\fff *} \dff\|\qff \leq\qff 1$\nnsp.\oss
At the same time\qss
$m_{\dff \sigma}\fff(\fff \mathbf{1}\fff)\qff =\qff 1$\qss
implies that\qss
$\|\dff p_{\fff *} \dff\|\qff \geq\qff 1$\qss
and\dss hence\qss
$\|\dff p_{\fff *} \dff\|\qff =\qff 1$\nnsp.\oss
In addition,\qss
$m_{\dff \sigma}\fff(\fff \mathbf{1}\fff)\qff =\qff 1$\qss
implies that\qss
$m_{\dff \sigma}\fff(\fff a\mathbf{1}\fff)\qff =\qff a$\qss
for every\qss
$a\qff \in\qff \rrr$\qss
and\dss hence\qss
$p_{\fff *}\dff\circ\fff p^*\off =\dff\off \id$\nnsp.\oss

It remains to show that\dss the maps\dss $p_{\fff *}$\dss
commute with the coboundaries of the chain complexes\qss
$B^*\fff(\dff X\dff)$\dnsp,\qss $B^*\fff(\dff Y\dff)$\dnsp.\oss
Let\qss $\sigma\dff \colon\dff \Delta_{\fff n}\toto Y$\qss 
be a singular $n$\dnsp-simplex\halfff,\qss
$i\dff \colon\dff \Delta_{\fff n\dff -\qff 1}\toto \Delta_{\fff n}$\qss
be  a simplicial embedding,\oss
and\oss
$\tau
\off =\off 
\sigma\dff \circ\dff i$\nnsp.\oss
Let\oss
$\dis
\delta
\qff \colon\dff
C_{\dff \sigma}
\ttoo
C_{\dff \tau}$\oss
be the map\qss
$\rho\qff \longmapsto\qff \rho\dff \circ\dff i$\nnsp.\oss
We claim that\vspace*{3pt}
\begin{equation}
\label{coherence-two}
\quad
m_{\dff \tau}\fff(\dff f\dff)
\off =\off 
m_{\dff \sigma}\fff(\dff f\dff \circ\dff \delta\dff)
\end{equation}

\vspace*{-9pt}
for all 
$f\qff \in\qff B\fff(\fff C_{\dff \tau}\fff)$\nnsp.\oss
In order to prove this,\oss
let us choose an arbitrary simplex\dss
$\overline{\sigma}\qff \in\qff C_{\dff \sigma}$\nnsp.\oss
Then the simplex\oss
$\overline{\tau}
\off =\off 
\overline{\sigma}\dff \circ\dff i$\oss
belongs to\dss $C_{\dff \tau}$\nnsp.\qss\oss
The diagram
\[
\quad
\begin{tikzcd}[column sep=huge, row sep=huge]\dis
G^{\dff \Delta_{\fff n}}  \arrow[r, "\dis \pi^{\fff i}"] \arrow[d, "\dis r_{\dff \overline{\sigma}}"]
& 
G^{\dff \Delta_{\fff n\dff -\dff 1}} \arrow[d, "\dis r_{\dff \overline{\tau}}"]
\\  
C_{\dff \sigma} \arrow[r, "\dis\delta"] 
&   
C_{\dff \tau}
\end{tikzcd}
\]

\vspace*{-6pt}
is obviously commutative.\oss
The 
maps\oss
$r_{\dff \overline{\sigma}}$\nnsp,\oss $r_{\dff \overline{\tau}}$\oss 
are bijections
and can be used to identify the top row of the diagram with the bottom row.\oss
By the definition of\dss  
$m_{\dff \sigma}$\dss and\dss $m_{\dff \tau}$\dss
this identification turns\dss  
$m_{\fff n}$\dss and\dss $m_{\fff n\trf -\dff 1}$\dss
into\dss
$m_{\dff \sigma}$\dss and\dss $m_{\dff \tau}$\dss
respectively.\oss
Hence\qss (\ref{coherence-two})\qss follows from\qss (\ref{coherence-one}).

Now we are ready to prove that the maps\dss $p_{\fff *}$\dss
commute with the differentials.\pss
Let\trs $f\qff \in\qff B^{\fff n\dff -\dff 1}\fff(\dff X\dff)$
and\dss let\qss
$\sigma\dff \colon\dff \Delta_{\fff n}\toto Y$\qss 
be a singular $n$\dnsp-simplex in $Y$\nnsp.\oss 
Then\vspace*{6pt} 
\[
\quad
p_{\fff *}\fff(\dff \partial f\dff)\fff(\fff \sigma \fff)
\off\off =\off\off
m_{\fff n}\fff(\dff \partial f \mid C_\sigma \dff)
\]

\vspace*{-33pt} 
\[
\quad
\phantom{p_{\fff *}\fff(\dff \partial f\dff)\fff(\fff \sigma \fff)
\off\off }
=\off\off
m_{\fff n}\qff \left(\qff
\sum_{k\qff =\qff 0}^{n}\qff
(-\qff 1)^k\qff    
f\circ\partial_k \off \Biggl|\off C_\sigma 
\qff\right)
\off\off =\off\off
\sum_{k\qff =\qff 0}^{n}\qff
(-\qff 1)^k\qff
m_{\fff n}\qff 
\bigl(\qff
f\circ\partial_k \pff\bigl|\pff C_\sigma 
\qff\bigr).
\]

\vspace{-3pt}
For every\dss $k$\dss the map\dss
$\partial_k$\dss 
takes\dss $C_{\fff \sigma}$\dss into\dss $C_{\fff \tau}$\nnsp,\oss
where\oss $\tau\off =\off \partial_k\dff\sigma$\nnsp.\oss
The map\qss
$C_{\fff \sigma}\toto C_{\fff \tau}$\qss
induced by\dss $\partial_k$\dss
is nothing else but the above map $\delta$
corresponding to\qss
$i\qff =\qff i_{\fff k}$\nnsp.\oss
Hence\qss (\ref{coherence-two})\qss implies that\vspace*{3pt}
\[
\quad
m_{\fff n}\fff 
\left(\qff
f\circ\partial_k \pff\bigl|\pff C_\sigma 
\qff\right)
\off =\off
m_{\fff n\dff -\dff 1}\fff 
\left(\qff
f \pff\bigl|\pff C_{\qff \partial_k\dff \sigma} 
\qff\right).
\]

\vspace*{-9pt}
It follows that\vspace*{-2pt}
\[
\quad
p_{\fff *}\fff(\dff \partial f\dff)\fff(\fff \sigma \fff)
\off\off =\off\off
\sum_{k\qff =\qff 0}^{n}\qff
(-\qff 1)^k\qff
m_{\fff n}\fff 
\left(\qff
f\circ\partial_k \pff\bigl|\pff C_\sigma 
\qff\right)
\]

\vspace*{-27pt} 
\[
\quad
\phantom{p_{\fff n}\fff(\partial f)\fff(\fff \sigma \fff)
\off\off }
=\off\off
\sum_{k\qff =\qff 0}^{n}\qff
(-\qff 1)^k\qff
m_{\fff n\dff -\dff 1}\fff 
\left(\qff
f \pff\bigl|\pff C_{\qff \partial_k\dff \sigma} 
\qff\right)
\]

\vspace*{-27pt} 
\[
\quad
\phantom{p_{\fff n}\fff(\partial f)\fff(\fff \sigma \fff)
\off\off }
=\off\off
\sum_{k\qff =\qff 0}^{n}\qff
(-\qff 1)^k\qff
p_{\fff *}\fff(\dff f\dff)(\dff \partial_k\dff \sigma\dff)
\off\off 
=\off\off
\partial\dff p_{\fff *}\dff (\dff f\dff)(\fff \sigma\fff)
\]

\vspace*{-6pt}
Therefore\oss
$p_{\fff *}\dff \circ\dff \partial 
\off =\off
\partial\dff \circ\dff p_{\fff *}$\nnsp.\oss
This completes the proof\halfff.\oss  \eproof

\mypar{Corollary.}{injectivity}
\emph{Under assumptions of Theorem\qss \ref{chain-map}\qss
the induced map}\vspace*{3pt}
\[
\quad
\quad
p^{\fff *}
\dff \colon\dff 
\widehat{H}^{\fff i}\fff(\dff Y \dff)
\ttoo 
\widehat{H}^{\fff i}\fff(\dff X \dff)
\]

\vspace*{-9pt}
\emph{of the bounded cohomology groups
is injective for every $i$\nnsp.\oss}

\proof\qss
Let $N$ be a natural number\qss $\geq\qff i\qff +\qff 1$\nnsp.\oss
By Theorem\qss \ref{coherent-exists-finite}\qss
there exists a coherent sequence\pss 
$m_{\fff 0}\fff,\pff m_{\fff 1}\fff,\pff m_{\fff 2}\fff,\pff 
\ldots\fff,\pff m_{\dff N}$\pss for\qss $G$\nnsp.\oss
By\dss Theorem\qss \ref{chain-map}\qss for each\qss
$i\qff \leq\qff N$\qss
there exists a homomorphism\vspace*{3pt}
\[
\quad
p_{\fff *}
\dff \colon\dff 
B^{\fff i}\fff(\dff X \dff)
\ttoo 
B^{\fff i}\fff(\dff Y \dff)
\]

\vspace*{-9pt}
such that\pss
$p_{\fff *}
\qff \circ\qff 
p^{\fff *}
\pff =\pff 
\id$\pss
and these homomorphisms commute with the coboundary maps.\oss
Hence for each\oss
$i\off =\off 0\fff,\pff 1\fff,\pff \ldots\fff,\pff N\qff -\qff 1$\oss
these homomorphisms induce a homomorphism\vspace*{3pt}
\begin{equation*}
\quad
p_{\fff *}
\dff \colon\dff 
\widehat{H}^{\fff i}\fff(\dff X \dff)
\ttoo 
\widehat{H}^{\fff i}\fff(\dff Y \dff)
\end{equation*}

\vspace*{-9pt}
which is a left inverse to\pss
$\dis
p^{\fff *}
\dff \colon\dff 
\widehat{H}^{\fff i}\fff(\dff Y \dff)
\ttoo 
\widehat{H}^{\fff i}\fff(\dff  X \dff)$\dnsp.\oss
It follows that $p^{\fff *}$ is injective.\oss  \eproof

\myuppar{Killing homotopy groups.}
Let $X$ be a path-connected space homotopy e\-quiv\-a\-lent to a\dss CW-complex\halfff.\oss
Suppose that\oss 
$\pi_{\fff i}\dff(\dff X\dff)
\off =\off 
0$\oss 
for all\qss $i\qff <\qff n$\qss
and\dss let\oss $\pi\off =\off \pi_{\fff n}\fff(\dff X \dff)$\dnsp.\oss 
For\qss $n\qff \geq\qff 2$\nnsp,\oss 
let\dss\vspace*{3pt}
\[
\quad
B_{\dff \pi}
\off =\off
B\fff(\dff \pi\fff,\pff S^n \dff)\dff,\hspace*{1.5em}
E_{\dff \pi}
\off =\off
B\fff(\dff \pi\fff,\pff D^n \dff)\fff,
\hspace*{1.5em}\mbox{ and }\hspace*{1.5em}
G_{\dff \pi}
\off =\off
B\fff(\dff \pi\fff,\pff S^{n\dff -\dff 1} \dff)\fff,\quad
\]

\vspace*{-9pt}
be the spaces from Section\qss \ref{classifying-spaces}.\oss
As we saw in Section\qss \ref{classifying-spaces},\oss
the spaces\qss
$B_{\dff \pi}$\nnsp,\pss $E_{\dff \pi}$\nnsp,\pss and\qss $G_{\dff \pi}$\qss
are\dss CW-com\-plexes,\oss
the spaces\qss $B_{\dff \pi}$\qss and\qss  $G_{\dff \pi}$\qss
are Eilenberg--MacLane spaces of the types\qss
$K\fff(\dff \pi\fff,\pff n \dff)$\dss
and\dss
$K\fff(\dff \pi\fff,\pff n\qff -\qff 1 \dff)$\dss
respectively\qss
and the space\qss $E_{\dff \pi}$\qss is contractible.\oss
By Theorem\qss \ref{b-space}\qss  $G_{\dff \pi}$\dss
is a topological group in the category of compactly generated spaces.\oss
With these notations,\oss the principal numerable bundle\qss 
(\ref{mccord-bundle})\qss takes the form\vspace*{6pt}
\begin{equation}
\label{pi-bundle}
\quad
p_{\fff \pi}^{\fff n}
\qff \colon\qff
E_{\dff \pi}
\ttoo
B_{\dff \pi}
\end{equation}

\vspace*{-6pt}
and has\qss $G_{\dff \pi}$\qss as the fiber\halfff.\oss
If\qss $n\qff =\qff 1$\nnsp,\oss then we take as\qss (\ref{pi-bundle})\qss
the universal cover of any CW-complex which is a\dss $K\fff(\dff \pi\fff,\pff 1 \dff)$\dnsp-space.\oss
By Lemma\qss \ref{p-wp}\dss $p_{\fff \pi}^{\fff n}$\dss
is a weakly principal bundle.\oss
In par\-tic\-u\-lar\halfff,\pss
$p_{\fff \pi}^{\fff n}$\dss
is a Serre fibration.\oss
The homotopy sequence of\dss $p_{\fff \pi}^{\fff n}$\dss has $0$\dnsp'{\halfff}s 
at all places except\vspace*{3pt}
\[
\quad
\pi_{\fff n}\fff(\dff B_{\dff \pi}\fff)
\ttoo
\pi_{\fff n\dff -\dff 1}\fff(\dff G_{\dff \pi}\fff)\dff.
\]

\vspace*{-9pt}
This boundary map is nothing else but\qss $\id_{\dff \pi}\dff \colon\dff \pi\toto \pi$\nnsp.\oss
Since $X$ is homotopy equivalent to a\dss CW-complex and\qss $B_{\dff \pi}$\qss
is\dss a\dss CW-complex of the type\qss
$K\fff(\dff \pi\fff,\pff n \dff)$\dnsp,\oss there exists a map\vspace*{3pt}
\[
\quad
j\dff \colon\dff X\ttoo 
B_{\dff \pi}
\]

\vspace{-9pt}
such the the induced map of homotopy groups $\pi_{\fff n}\fff(\fff j\dff)$ is an isomorphism.\oss
Let\qss\vspace*{3pt} 
\[
\quad
p\dff \colon\dff Y\ttoo X
\]

\vspace*{-9pt}
be the bundle induced from the bundle\qss $p_{\fff \pi}^{\fff n}$\qss
by the map\dss $j$\nnsp.\oss
It has the same fiber\dss $G_{\dff \pi}$\dss as the bundle\qss $p_{\fff \pi}^{\fff n}$\nnsp.\oss
Since\qss $p_{\fff \pi}^{\fff n}$\qss is a Serre fibration,\oss
the bundle\dss $p$\dss is a Serre fibration also.\oss
By comparing the homotopy sequences of the bundles\dss $p$\dss 
and\qss $p_{\fff \pi}^{\fff n}$\nsp,\oss
we see that the boundary map\vspace*{6pt}
\[
\quad
\pi_{\fff n}\fff(\dff X\fff)
\off =\off
\pi
\ttoo
\pi_{\fff n\dff -\dff 1}\fff(\dff G_{\dff \pi}\fff)
\]

\vspace*{-6pt}
is an isomorphism.\oss
Since\oss $\pi_{\fff m}\fff(\dff G_{\dff \pi}\fff)\off =\off 0$\oss 
for\qss $m\off \neq\off n\qff -\qff 1$\nnsp,\oss
it follows that\vspace*{6pt}
\[
\quad
\pi_{\fff n}\fff(\dff Y \dff)
\off =\off
0
\hspace*{1.5em}\mbox{ and }\hspace*{1.5em}
p_{\fff *}
\qff \colon\qff
\pi_{\fff m}\fff(\dff Y \dff)
\ttoo
\pi_{\fff m}\fff(\dff X \dff)
\]

\vspace*{-6pt}
is an isomorphism for\qss $m\off \neq\off n$\nnsp.\oss
One may say that the space $Y$ resulted from\qss \emph{killing the $n$\dnsp-th homotopy group of\pss $X$\nnsp.\oss}
The classical version of this construction is due to Cartan and Serre\qss \cite{cs}\qss 
and uses instead of the principal bundle\qss $p_{\fff \pi}^{\fff n}$\qss 
the path space Serre fibration\qss\vspace*{3pt}
\[
\quad
PK(\dff \pi\fff,\pff n \dff)\ttoo K(\dff \pi\fff,\pff n \dff)\dff,
\]

\vspace*{-9pt}
where\dss $PK(\dff \pi\fff,\pff n \dff)$\dss is the space of paths in a\dss
$K(\dff \pi\fff,\pff n \dff)$\dnsp-space starting at a fixed point\halfff.\oss

Since the bundle $p$ is induced from the bundle\dss $p_{\fff \pi}^{\fff n}$\dss
and\dss $p_{\fff \pi}^{\fff n}$\dss is a numerable locally trivial principal bundle,\pss
$p$ is also a numerable locally trivial principal bundle.\oss
By Theorem\qss \ref{cw-bundle}\qss this implies that $Y$ is homotopy equivalent to
a\dss CW-complex\halfff.\oss
This allows to apply the same construction to\dss $Y$\dss and\dss $n\qff +\qff 1$\dss
in the roles of\dss $X$\dss and\dss $n$\dss respectively,\oss
and continue in this way.

\myuppar{Iterated killing.}
Let $X$ be a path-connected space homotopy equivalent to a\dss CW-complex\halfff.\oss
One can start killing homotopy groups with the first non-zero 
group 
and\dss then iterate the construction.\oss
This procedure is also due to Cartan and Serre\qss \cite{cs}.\oss
It\dss leads to a sequence of maps\vspace*{6pt}
\begin{equation}
\label{tower}
\quad
\begin{tikzcd}[column sep=large, row sep=normal]\dis
\ldots\quad \arrow[r, "\dis p_{\fff n}"]
& 
X_{\fff n} \arrow[r, "\dis p_{\fff n\dff -\dff 1}\off"]
& 
X_{\fff n\dff -\dff 1} \arrow[r, "\dis p_{\fff n\dff -\dff 2}\off"]
&   
\quad \ldots\quad \arrow[r, "\dis p_{\fff 2}"]
&
X_{\fff 2} \arrow[r, "\dis p_{\fff 1}"]
&
X_{\fff 1}
\end{tikzcd}
\end{equation}

\vspace*{-3pt}
such that\oss $X_{\fff 1}\off =\off X$\nnsp,\oss
$\pi_{\fff i}\fff(\dff X_{\fff n}\fff)\off =\off 0$\oss
if\pss $i\qff <\qff n$\nnsp,\oss
$\pi_{\fff i}\fff(\dff X_{\fff n}\fff)\off =\off \pi_{\fff i}\fff(\dff X \dff)$\oss
if\oss $i\qff \geq\qff n$\nnsp,\oss
and each map\vspace*{6pt}
\begin{equation*}
\quad
p_{\fff n}\dff \colon\dff X_{\fff n\dff +\dff 1}\ttoo X_{\fff n}
\end{equation*}

\vspace*{-6pt}
is a weakly principal bundle having a topological group\dss $G_{\fff n}$\dss 
as a fiber\halfff.\oss
The group\dss $G_{\fff n}$\dss
is abelian for\qss $n\qff \geq\qff 2$\qss
and is discrete and isomorphic to\dss $\pi_{\fff 1}\fff(\dff X \dff)$\dss
for\qss $n\qff =\qff 1$\nnsp.\oss
Each space\dss $X_{\fff n}$\dss is homotopy equivalent to a\dss CW-complex\halfff,\pss
and each group\dss $G_{\fff n}$\dss is a\dss CW-complex\halfff.\oss

\myuppar{Partial contracting homotopies.}
Let $X$ be a topological space,\oss
and let $n$ be a natural num\-ber\halfff.\oss
A\qss \emph{partial contracting homotopy up to dimension}\dss $n$\dss 
is a sequence of homomorphisms\vspace*{6pt}
\[
\quad
\begin{tikzcd}[column sep=large, row sep=normal]\dis
\rrr
& 
B^{0}\fff(\dff X \dff) \arrow[l, "\dis \off K^{\dff 0}"]
& 
B^{1}\fff(\dff X \dff) \arrow[l, "\dis \off K^{\dff 1}"]
&   
\quad \ldots\quad \arrow[l, "\dis \off K^{\dff 2}"]
&
B^{n}\fff(\dff X \dff) \arrow[l, "\dis \off K^{\fff n}"]
\end{tikzcd}
\]

\vspace*{-3pt}
such that\oss
$\dis
\partial
\qff \circ\qff
K^{\dff i}
\off +\off 
K^{\dff i\qff +\qff 1}
\qff \circ\pff 
\partial
\off =\off
\id$\oss
for all\qss $0\qff \leq\qff i\qff \leq\qff n\qff -\qff 1$\nnsp,\oss
where\dss $\partial$\dss is,\pss
as usual,\pss
the coboundary operator\halfff.\oss
It is called\qss \emph{bounded}\pss if\dss all\dss $K^{\fff i}$\dss
are bounded operators,\oss
and\qss \emph{strictly bounded}\pss if\qss
$\|\dff K^{\fff i} \dff\|\qff \leq\qff 1$\qss 
for\dss all\qss 
$0\qff \leq\qff i\qff \leq\qff n$\nnsp.\oss 

Suppose that $X$ is path-connected\dss and\pss 
$\pi_{\fff i}\dff(\dff X\dff)
\off =\off 
0$\pss 
for all\qss $i\qff <\qff n$\nnsp.\oss
Then there exists a strictly bounded partial contracting homotopy up to dimension $n$\nnsp.\oss
A construction of such partial homotopy is
implicitly contained in the most proofs\dss of\dss the vanishing of the singular cohomology groups
of contractible spaces.\oss
Since we need to deal with non-contractible spaces
and need to control the norm of the involved maps,\oss
we present such a construction from the scratch.\oss

\vspace*{3pt}
Let us begin with constructing a sequence of maps\vspace*{6pt}
\[
\quad
\begin{tikzcd}[column sep=large, row sep=normal]\dis
\{\trf 1 \trf\} \arrow[r, "\dis L^{0}"]
& 
S_{\halfff 0}\fff(\dff X\fff) \arrow[r, "\dis L^{1}\off"]
& 
S_{1}\fff(\dff X\fff) \arrow[r, "\dis L^{2}\off"]
&   
\quad \ldots\quad \arrow[r, "\dis L^{n}\off"]
&
S_{n}\fff(\dff X\fff)
\end{tikzcd}
\]

\vspace*{-3pt}
such that
for every\qss $\sigma\qff \in\qff S_i\dff(\dff X \dff)$\dnsp,\qss
$i\qff \leq\qff n\qff -\qff 1$\nnsp,\oss
the equality\vspace*{3pt}
\begin{equation}
\label{l-homotopy}
\quad
\partial\fff\halfff L^{i\qff +\qff 1\dff}\dff(\fff \sigma \fff)
\off\off +\off\off
\sum_{k\qff =\qff 0}^i\qff
(-\qff 1)^{\fff k}\qff
L^{i}\dff(\trf \partial_{\fff k}\dff\sigma \trf)
\off\off =\off\off
\sigma\dff.
\end{equation}

\vspace*{-6pt}
holds in the chain group\qss $C_{\dff i}\dff(\dff X \dff)$\dnsp.\oss
Let us choose a base point\qss 
$b\qff \in\qff X$\qss 
and set\qss 
$L^{0}\dff(\fff 1\fff)\qff =\qff \beta$\nnsp,\oss
where $\beta$ 
is the $0$\dnsp-simplex with the image $\{\dff b\dff\}$\nnsp.\oss
Given a sin\-gu\-lar $0$\dnsp-simplex\qss
$\sigma\dff \colon\dff \Delta_{\fff 0}\toto X$\nnsp,\oss
let us connect its only vertex with the base point $b$ by a path in\qss $X$\nnsp.\qff\oss
Such a path leads to a singular\oss $1$\dnsp-simplex\qss
$L^{1}\dff(\fff \sigma\fff)\dff \colon\dff \Delta_{\fff 1}\toto X$\qss
such that\qss
$\partial\fff\halfff L^{1}\dff(\fff \sigma \fff)
\off =\off
\sigma\qff -\qff L^{0}\dff(\fff 1\fff)$\dnsp.\oss
Next\halfff,\oss 
if\qss $\sigma\dff \colon\dff \Delta_{\fff 1}\toto X$\qss
is a sin\-gu\-lar $1$\dnsp-simplex\fff,\pss
then the three $1$\dnsp-simplices\qss
$L^{1}\dff(\fff \partial_{\fff 0}\fff\sigma\fff)$\dnsp,\qss
$\sigma$\nnsp,\qss
$L^{1}\dff(\fff \partial_{\fff 1}\fff\sigma\fff)$\qss
form a loop.\oss
If\qss $n\qff \geq\qff 2$\nnsp,\oss
then $X$ is simply-connected and\dss hence
this loop is contractible in $X$\nnsp.\oss
In this case there exists singular $2$\dnsp-simplex\qss
$L^{2}\dff(\fff \sigma\fff)\dff \colon\dff \Delta_{\fff 2}\toto X$\qss
such that\oss
$\partial\fff\halfff L^{2}\dff(\fff \sigma \fff)
\off =\off
\sigma
\qff -\qff 
L^{1}\dff(\dff \partial_{\fff 0}\fff\sigma \fff)
\qff +\qff
L^{1}\dff(\dff \partial_{\fff 1}\fff\sigma \fff)$\dnsp.\oss
Since $X$ is\dss $(n\qff -\qff 1)$\nsp\dnsp-connected,\oss 
we can continue in this way until we get all the maps\dss $L_{\fff i}$\dss 
needed.\oss

Now we can define\oss\
$\dis
K^{\dff i\dff +\qff 1}\dff \colon\dff
B^{i\dff +\dff 1}\fff(\dff X \dff)
\ttoo
B^{i}\fff(\dff X \dff)$\oss
as the map induced by\oss\vspace*{3pt}
\[
\quad
L^{i\dff +\dff 1}\dff \colon\dff
S_{i}\fff(\dff X\fff)
\ttoo
S_{i\dff +\dff 1}\fff(\dff X\fff)\dff.
\]

\vspace*{-9pt}
Since $K^{\dff i}$ is induced by a map
between the sets of singular simplices,\oss
its norm is\qss
$\|\dff K^{\dff i} \dff\|\qff \leq\qff 1$\nnsp.\oss
It follows from\qss (\ref{l-homotopy})\qss that\qss
$K^{\dff i}$ with\qss $i\qff \leq\qff n$\qss
form a strictly bounded partial contracting homotopy.\oss

\vspace{3pt}
\mypar{Theorem.}{simply-connected-homology} 
\emph{Suppose that\pss $X$\dss 
is a path-connected space homotopy equivalent to a\dss CW-complex and\dss having amenable fundamental group.\oss
Then\oss 
$\widehat{H}^{\fff i}\fff(\dff X\dff)
\off =\off 
0$\oss 
for all\qss $i\qff \geq\qff 1$\nnsp.\oss}

\proof\qss 
Iterated killing of homotopy groups leads
to a sequence of maps\vspace*{6pt}
\[
\quad
\begin{tikzcd}[column sep=large, row sep=normal]\dis
\ldots\quad \arrow[r, "\dis p_{\fff n}"]
& 
X_{\fff n} \arrow[r, "\dis p_{\fff n\dff -\dff 1\off}"]
& 
X_{\fff n\dff -\dff 1} \arrow[r, "\dis p_{\fff n\dff -\dff 2}\off"]
&   
\quad \ldots\quad \arrow[r, "\dis p_{\fff 2}"]
&
X_{\fff 2} \arrow[r, "\dis p_{\fff 1}"]
&
X_{\fff 1}
\end{tikzcd}
\]

\vspace*{-3pt}
such that\oss $X_{\fff 1}\off =\off X$\nnsp,\oss
$\pi_{\fff i}\dff(\dff X_{\fff n}\fff)\qff =\qff 0$\qss
if\qss $i\qff <\qff n$\nnsp,\oss
$\pi_{\fff i}\dff(\dff X_{\fff n}\fff)\qff =\qff \pi_{\fff i}\dff(\dff X \dff)$\qss
if\qss $i\qff \geq\qff n$\nnsp,\oss
and each map\vspace*{3pt}
\[
\quad
p_{\fff n}\dff \colon\dff X_{\fff n\dff +\dff 1}\ttoo X_{\fff n}
\]

\vspace*{-9pt}
is a principal $G_{\fff n}$\nsp\dnsp-bundle
for a topological group $G_{\fff n}$
which is either abelian,\oss or discrete and amenable.\oss
Let\qss $i\qff \geq\qff 1$\nnsp,\oss and\dss let us choose some\qss $n\qff \geq\qff i$\nnsp.\oss
By Corollary\qss \ref{injectivity}\qss
the maps\vspace*{3pt}
\begin{equation*}
\quad
p_{\fff m}^{\fff *}
\qff \colon\qff 
\widehat{H}^{\fff i}\fff(\dff X_{\fff m} \fff)
\qff \ttoo\qff 
\widehat{H}^{\fff i}\fff(\dff  X_{\fff m\dff +\dff 1}\fff)
\end{equation*}

\vspace*{-9pt}
are injective for all\dss $m$\nnsp.\oss
This implies that the induced map\vspace*{3pt}
\begin{equation}
\label{induce-b}
\quad
\left(\dff
p_{\fff 1}\dff \circ\dff p_{\fff 2}\dff \circ\dff \ldots \dff \circ\dff p_{\fff n\dff +\dff 1}
\dff\right)^*
\colon\pff
\widehat{H}^{\fff i}\fff(\dff X \dff)
\qff \ttoo\qff 
\widehat{H}^{\fff i}\fff(\dff  X_{\fff n\dff +\dff 1}\fff)
\end{equation}

\vspace*{-9pt}
is injective.\oss
Since $X_{\fff n\dff +\dff 1}$ is path-connected\dss and\pss 
$\pi_{\fff i}\dff(\dff X_{\fff n\dff +\dff 1} \dff)
\off =\off 
0$\pss 
for all\qss $i\qff <\qff n\qff +\qff 1$\nnsp,\oss
there is a bounded partial contracting homotopy up to dimension\qss $n\qff +\qff 1$\qss for\qss
$X_{\fff n\dff +\dff 1}$\nnsp.\oss 
Since\qss $n\qff \geq\qff i$\nnsp,\oss
it follows that\oss
$\widehat{H}^{\fff i}\fff(\dff X_{\fff n\dff +\dff 1} \dff)
\off =\off
0$\dnsp,\oss
and\dss the injectivity of\qss (\ref{induce-b})\qss 
implies that\oss
$\widehat{H}^{\fff i}\fff(\dff X \dff)
\off =\off
0$\nnsp.\oss  \eproof

\myuppar{Countable abelian groups.}
If $\pi$ is a countable abelian group,\oss 
then for every\qss $n\qff \geq\qff 1$\qss 
there is an Eilenberg--MacLane space\dss $G_{\dff \pi}$\dss of the type\dss 
$K\fff(\dff \pi\fff,\pff n\qff -\qff 1 \dff)$\dss
which is a topological group in usual sense.\oss
This was proved by Milnor\qss \cite{m2},\oss 
Corollary to Theorem\qss 3.\oss 
The space\dss $G_{\dff \pi}$\dss is the geometric realization
of a semi-simplicial complex
and\dss hence is a\dss CW-complex\halfff.\oss
Alternatively,\oss one can construct such a topological group\dss $G_{\dff \pi}$\dss
by applying the Dold--Thom\qss \cite{dt}\qss construction\qss $AG\fff(\dff \bullet\dff)$\qss 
to a countable simplicial complex which is 
a Moore space of the type\dss $(\dff \pi\fff,\pff n\qff -\qff 1 \dff)$\dnsp.\oss

Milnor's universal bundle\qss \cite{m2}\qss
is a map of the form\qss (\ref{pi-bundle})\qss
which is a principal\dss $G_{\dff \pi}$\dnsp-bundle in the classical sense,\oss
and\dss hence is a weakly principal bundle with the fiber\dss $G_{\dff \pi}$\nnsp.\oss
If the base of the Milnor's universal bundle
is not a\dss CW-complex\halfff,\oss then one can replace this universal bundle 
by another one induced from it by a weak homotopy equivalence.\oss
Even better\halfff,\oss one can start with a countable\dss CW-complex\dss $B_{\dff \pi}$\dss which is
an Eilenberg--MacLane space of the type\dss 
$K\fff(\dff \pi\fff,\pff n \dff)$\dnsp.\oss
By\qss \cite{m1},\oss Theorem\qss 5.2\fff(3),\oss there exists a principal bundle
with all the desired properties.\oss

It follows that if all homotopy groups\dss $\pi_{\fff n}\fff(\dff X \dff)$\dss are countable,\oss
for example if\dss $X$\dss is homotopy equivalent to countable\dss CW-complex\halfff,\oss
then one can prove Theorem\qss \ref{simply-connected-homology}\qss without using 
compactly generated spaces and the McCord classifying spaces from Section\qss
\ref{classifying-spaces}.

\myuppar{Contracting\qss chain\qss homotopies.}
The rest of this section is devoted to the proof 
of the existence of strictly bounded contracting chain homotopies\dss $K^{\fff i}$\dss 
defined for all\qss $i\qff \geq\qff 0$\qss at once.\oss
This result will not be be needed till Section\qss \ref{spaces}.

\mypar{Theorem.}{coherent-exists-infinite}
\emph{Suppose that\qss $G$\dss is either an abelian topological group,\oss
or an amenable discrete group.\oss
Then there exists a coherent sequence of\dss invariant means\qss
$\{\dff m_{\fff n} \dff \}$\nnsp,\oss
$n
\off =\off
0\fff,\pff 1\fff,\pff 2\fff,\pff \ldots\off $\dnsp. }

\proof\qss
Let\dss $\mathcal{M}_{\fff n}$\dss be the set of all
$\Sigma_{\dff n\dff +\dff 1}$\nsp\dnsp-invariant means on\qss 
$B\fff(\fff G^{\dff \Delta_{\fff n}}\dff)$\dnsp.\oss
By remarks preceding Theorem\qss \ref{coherent-exists-finite}\dss 
$\mathcal{M}_{\fff n}\off \neq\off \emptyset$\dss  
and\dss if\qss
$i\dff \colon\dff \Delta_{\fff n\dff -\qff 1}\toto \Delta_{\fff n}$\qss
is a simplicial embedding,\oss
then
the map\qss
$m\qff \longmapsto\qff \pi^{\fff i}_*\fff m$\qss
is independent on $i$ and
takes\dss $\mathcal{M}_{\fff n}$\trs to\dss $\mathcal{M}_{\fff n\dff -\dff 1}$\nnsp.\oss
Thus we have a projective system\vspace*{1.5pt}
\begin{equation}
\label{projective}
\quad
\mathcal{M}_{\fff 0}\off \longleftarrow\off \mathcal{M}_{\fff 1}\off \longleftarrow\off \mathcal{M}_{\fff 2}
\off \longleftarrow\off \ldots\off\off.
\end{equation}

\vspace*{-10.5pt}
The set\dss $\mathcal{M}_{\fff n}$\trs 
is a subset of the set\dss $M_{\fff n}$\dss 
of all invariant means and hence is contained 
in the unit ball of\dss 
the dual Banach space\dss 
$B\fff(\fff G^{\dff \Delta_{\fff n}}\dff)^*$\dnsp.\pss
The set\dss $\mathcal{M}_{\fff n}$\trs 
is obviously closed in the weak\dnsp$^*$ to\-pol\-o\-gy.\oss
By the Banach-Alaoglu theorem\qss
(see\qss \cite{r},\oss Theorem\qss 3.15)\qss 
the unit ball is compact in this topology,\oss 
and\dss hence\dss $\mathcal{M}_{\fff n}$\trs
is also compact\halfff.\oss 
The map\qss 
$f\qff \longmapsto f\dff \circ\dff \pi^{\fff i}$\qss 
is a bounded operator\vspace*{1.5pt}
\[
\quad
B\fff(\fff G^{\dff \Delta_{\fff n\dff -\dff 1}}\dff)
\ttoo
B\fff(\fff G^{\dff \Delta_{\fff n}}\dff)
\]

\vspace*{-10.5pt}
It follows that the dual map\oss 
$\dis
B\fff(\fff G^{\dff \Delta_{\fff n}}\dff)^*
\ttoo
B\fff(\fff G^{\dff \Delta_{\fff n\dff -\dff 1}}\dff)^*$
is also bounded,\oss
and hence is continuous in the weak\dnsp$^*$ topology.\oss
The push-forward map\qss
$m\qff \longmapsto\qff \pi^{\fff i}_*\fff m$\qss
is the restriction of\dss this dual map to the set of invariant means 
and\dss hence is continuous in the weak\dnsp$^*$\nsp topology.\oss
It follows that\qss (\ref{projective})\qss is a projective system
of continuous in the  weak\dnsp$^*$\nsp topology  maps.\oss
Since the sets\dss $\mathcal{M}_{\fff n}$\trs are compact in this topology,\oss
the limit of this projective system is nonempty\qss 
(cf.\qss \cite{bo},\oss 
Chap.\qss 1,\oss Sec.\qss 9,\oss n\dnsp$^\circ$\dnsp\ 6,\oss Proposition\qss 8).\oss
Any point of this limit is a coherent sequence.\oss  \eproof

\myuppar{Compatible partial contracting homotopies.}
Suppose\dss that in the sequence\vspace*{2pt}
\begin{equation}
\label{tower-principal}
\quad
\begin{tikzcd}[column sep=large, row sep=normal]\dis
\ldots\quad \arrow[r, "\dis p_{\fff n}"]
& 
X_{\fff n} \arrow[r, "\dis p_{\fff n\dff -\dff 1}\off"]
& 
X_{\fff n\dff -\dff 1} \arrow[r, "\dis p_{\fff n\dff -\dff 2}\off"]
&   
\quad \ldots\quad \arrow[r, "\dis p_{\fff 2}"]
&
X_{\fff 2} \arrow[r, "\dis p_{\fff 1}"]
&
X_{\fff 1}
\end{tikzcd}
\end{equation}

\vspace*{-9pt}
the groups\oss
$\pi_{\fff i}\fff(\dff X_{\fff n}\fff)\off =\off 0$\oss
for\oss $i\qff <\qff n$\nnsp,\oss
and each\dss
$p_{\fff n}$\dss
is a principal $G_{\fff n}$\nsp\dnsp-bundle for a
topological group $G_{\fff n}$ which is either abelian or\dss discrete and amenable.\oss
Suppose that for each\qss $n\qff \geq\qff 1$\vspace*{3pt}
\[
\quad
\begin{tikzcd}[column sep=large, row sep=normal]\dis
\rrr
& 
B^{0}\fff(\dff X_{\fff n} \dff) \arrow[l, "\dis K^{\dff 0}_{\fff n}"]
& 
B^{1}\fff(\dff X_{\fff n} \dff) \arrow[l, "\dis K^{\dff 1}_{\fff n}"]
&   
\quad \ldots\quad \arrow[l, "\dis K^{\fff 2}_{\dff n}"]
&
B^{n\dff -\dff 1}\fff(\dff X_{\fff n} \dff) \arrow[l, "\dis K^{\dff n\dff -\dff 1}_{\fff n}"]
\end{tikzcd}
\]

\vspace*{-9pt}
is a partial contracting homotopy.\oss
These partial homotopies are said to be\qss
\emph{compatible}\pss if\qss\vspace*{2pt} 
\begin{equation}
\label{homotopy-compatibility}
\quad
p_{\fff n\dff *}
\qff \circ\pff  
K^{\dff i}_{\fff n\dff +\dff 1} 
\qff \circ\qff 
p_{\fff n}^{\fff *}
\off\qff =\off\qff 
K^{\fff i}_{\fff n}
\end{equation}

\vspace*{-10pt}
for every\dss $i\qff \leq\qff n\qff -\qff 1$\nnsp,\oss
where\qss 
$p_{\fff n\dff *}
\dff \colon\dff
B^{i}\fff(\dff X_{\fff n\dff +\dff 1} \dff)
\toto
B^{i}\fff(\dff X_{\fff n} \dff)$\dss are the maps from Theorem\qss \ref{chain-map}.\oss
By Theorem\qss \ref{coherent-exists-finite}\qss
one can assume that they are defined for\qss $i\qff \leq\qff n\qff -\qff 1$\nnsp.\oss
Suppose that for each $n$\vspace*{2pt}
\begin{equation*}
\quad
\begin{tikzcd}[column sep=large, row sep=normal]\dis
\{\dff 1\dff\} \arrow[r, "\dis L^n_{\dff 0}"]
& 
S_{\halfff 0}\fff(\dff X_{\fff n}\fff) \arrow[r, "\dis L^n_{\dff 1}\off"]
& 
S_{1}\fff(\dff X_{\fff n}\fff) \arrow[r, "\dis L^n_{\dff 2}\off"]
&   
\quad \ldots\quad \arrow[r, "\dis L^n_{\dff n\dff -\dff 1}"]
&
S_{n\dff -\dff 1}\fff(\dff X_{\fff n}\fff)
\end{tikzcd}
\end{equation*}

\vspace*{-6pt}
is a sequence of maps
such that for every\qss 
$\sigma\qff \in\qff S_i\dff(\dff X_{\fff n} \dff)$\dnsp,\qss
$i\qff \leq\qff n\qff -\qff 2$\nnsp,\oss
the equality\vspace*{3pt}
\begin{equation}
\label{l-homotopy-level}
\quad
\partial\fff\halfff L_{\dff i\dff +\dff 1\dff}^{n}\dff(\fff \sigma \fff)
\off\off +\off\off
\sum_{k\qff =\qff 0}^i\qff
(-\qff 1)^{\fff k}\qff
L_{\dff i}^{n}\qff(\dff \partial_k\fff\sigma \dff)
\off\off =\off\off
\sigma
\end{equation}

\vspace*{-6pt}
holds in the chain group\qss $C_{\dff i\qff +\qff 1}\dff(\dff X_{\fff n} \fff)$\dnsp.\oss
Sequences\dss $L^{\fff n}$\dss are said to be\qss
\emph{compatible}\qss if\vspace*{3pt}
\begin{equation}
\label{l-compatibility}
\quad
p_{\fff n}\pff \circ\pff L_{\dff i}^{n\dff +\dff 1}
\off =\off 
L_{\dff i}^{n}\pff \circ\qff p_{\fff n}
\end{equation}

\vspace{-9pt}
for every\qss $i\qff \leq\qff n\qff -\qff 2$\nnsp,\oss
where we denote by the same symbol $p_{\fff n}$ the map\qss
$S_\bullet\fff(\dff X_{\fff n\dff +\dff 1}\fff)
\ttoo
S_\bullet\fff(\dff X_{\fff n}\fff)$\qss
induced by\qss
$p_{\fff n}
\dff \colon\dff
X_{\fff n\dff +\dff 1}\ttoo X_{\fff n}$\nnsp.\oss

\mypar{Lemma.}{l-to-h}
\emph{If\dss there exist compatible sequences\oss 
$L^{1},\off\off L^2,\off\off L^3,\off\off \ldots \off\off$\dnsp,\oss
then there exist compatible strictly bounded partial homotopies\oss
$K_{\fff 1}\dff,\off\off K_{\fff 2}\dff,\off\off K_{\fff 3}\dff,\off\off \ldots \off\off$\dnsp.\oss}

\proof\qss
Suppose that\dss $L^{\fff n}$\dss are compatible sequences.\oss 
For\qss $i\qff \leq\qff n$\qss let\oss\vspace*{3pt}
\[
\quad
K^{\fff i}_{\fff n}\dff \colon\dff
B^{i}\fff(\dff X_{\fff n} \dff)
\ttoo
B^{i\dff -\dff 1}\fff(\dff X_{\fff n} \dff)
\]

\vspace*{-9pt}
be the map induced by\oss 
$L_{\fff i}^{\fff n}\dff \colon\dff
S_{i\dff -\dff 1}\fff(\dff X_{\fff n}\fff)
\ttoo
S_{i}\fff(\dff X_{\fff n}\fff)
$\dnsp.\oss 

The identity\qss (\ref{l-homotopy})\qss implies that the maps
$K^{\fff i}_{\fff n}$ with\qss $i\qff \leq\qff n$\qss
form a partial contracting homotopy.\oss
Since $K^{\fff i}_{\fff n}$ is induced by a map
between the sets of singular simplices,\oss
its norm is\qss
$\|\dff K^{\fff i}_{\fff n} \dff\|\qff \leq\qff 1$\nnsp.\oss
Therefore this partial contracting homotopy is strictly bounded.\oss

It remains to check 
the compatibility condition\qss (\ref{homotopy-compatibility}).\oss
If\dss 
$f\qff \in\qff B^i\fff(\dff X_{\fff n}\fff)$\dnsp,\oss then\vspace*{6pt}
\[
\quad
p_{\fff n\dff *}
\dff \circ\pff  
K^{\fff i}_{\fff n\dff +\dff 1} 
\qff \circ\qff 
p_{\fff n}^{\fff *}\pff 
(\dff f\dff)
\off =\off
p_{\fff n\dff *}
\dff \circ\pff  
K^{\fff i}_{\fff n\dff +\dff 1}\pff
(\dff f\dff \circ\dff p_{\fff n}\dff)
\]

\vspace*{-33pt}
\[
\quad
\phantom{p_{\fff n\dff *}
\dff \circ\pff  
K^{\fff i}_{\fff n\dff +\dff 1} 
\qff \circ\qff 
p_{\fff n}^{\fff *}\pff 
(\dff f\dff)
\off }
=\off
p_{\fff n\dff *}\dff
(\dff f\dff \circ\dff p_{\fff n}\dff \circ\dff L_{\fff i}^{\fff n\dff +\dff 1} \dff)
\off =\off
p_{\fff n\dff *}\pff
(\dff f\dff \circ\dff L_{\fff i}^{\fff n}\dff \circ\dff p_{\fff n} \dff)
\]

\vspace*{-33pt}
\[
\quad
\phantom{p_{\fff n\dff *}
\dff \circ\pff  
K^{\fff i}_{\fff n\dff +\dff 1} 
\qff \circ\qff 
p_{\fff n}^{\fff *}\pff 
(\dff f\dff)
\off =\off
p_{\fff n\dff *}\dff
(\dff f\dff \circ\dff p_{\fff n}\dff \circ\dff L_{\fff i}^{\fff n\dff +\dff 1} \dff)
\off}
=\off
p_{\fff n\dff *}
\dff \circ\dff
p_{\fff n}^{\fff *}\pff
(\dff f\dff \circ\dff L_{\fff i}^{\fff n} \dff)
\off =\off
f\dff \circ\dff L_{\fff i}^{\fff n}
\]

\vspace*{-33pt}
\[
\quad
\phantom{p_{\fff n\dff *}
\dff \circ\pff  
K^{\fff i}_{\fff n\dff +\dff 1} 
\qff \circ\qff 
p_{\fff n}^{\fff *}\pff 
(\dff f\dff)
\off =\off
p_{\fff n\dff *}\dff
(\dff f\dff \circ\dff p_{\fff n}\dff \circ\dff L_{\fff i}^{\fff n\dff +\dff 1} \dff)
\off =\off
p_{\fff n\dff *}
\dff \circ\dff
p_{\fff n}^{\fff *}\pff
(\dff f\dff \circ\dff L_{\fff i}^{\fff n} \dff)
\off }
=\off
K^{\fff i}_{\fff n}\pff (\dff f\dff)\dff,
\]

\vspace*{-6pt}
where on the second line we used the compatibility condition\qss (\ref{l-compatibility}).\qff\oss  
This\dss proves\qss (\ref{homotopy-compatibility}).\oss  \eproof

\mypar{Lemma.}{l-exist}
\emph{If\dss for every $n$ the topological group\dss $G_{\fff n}$\dss
is\dss $(n\qff -\qff 2)$\nsp\dnsp-connected,\oss
then there exist compatible sequences\oss 
$L^{1},\off\off L^2,\off\off L^3,\off\off \ldots \off\off$\dnsp.\oss}

\proof\qss
Let us choose the base points\qss
$b_n\qff \in\qff X_{\fff n}$\qss
in such a way that\pss
$p_{\fff n}\fff (\fff b_{n\dff +\dff 1}\fff)
\off =\off 
b_{n}$\qss
for all $n$\nnsp.\oss
Let\qss $m\qff \geq\qff 1$\nnsp.\oss
Suppose that we already constructed the sequences\oss 
$L^{1},\off\off L^2,\off\off L^3,\off\off \ldots, \off\off L^m$\oss
such that\oss
$L^n_{\dff 0}\trf(\dff 1\dff)
\off =\off 
b_n$\oss 
for all\qss $n\qff \leq\qff m$\qss
and\qss (\ref{l-compatibility})\qss holds
for every\qss $n\qff \leq\qff m\qff -\qff 1$\qss
and\qss $i\qff \leq\qff n\qff -\qff 2$\nnsp.\oss

Let us construct\qss $L^{m\dff +\dff 1}$\qss in such a way that\qss
$L^{m\dff +\dff 1}_{\dff 0}\trf(\dff 1\dff)
\off =\off 
b_{m\dff +\dff 1}$\oss
and\qss (\ref{l-compatibility})\qss holds
for\qss $n\qff =\qff m$\qss
and every\qss $i\qff \leq\qff n\qff -\qff 2$\nnsp.\oss
Suppose that\qss $1\qff \leq\qff k\qff \leq\qff m$\qss and the maps\vspace*{1.5pt}
\[
\quad
L_{\fff i}^{m\dff +\dff 1}\dff \colon\dff
S_{i\dff -\dff 1}\fff(\dff X_{\fff m\dff +\dff 1}\fff)
\ttoo
S_{i}\fff(\dff X_{\fff m\dff +\dff 1}\fff)
\]

\vspace*{-10.5pt}
are already constructed for\qss $i\qff \leq\qff k\qff -\qff 1$\nnsp.\oss
Arguing by induction,\oss
we may assume that\qss
(\ref{l-compatibility})\qss
holds if\qss
$n\qff =\qff m$\qss and\qss $i\qff \leq\qff k\qff -\qff 1$\nnsp.\oss
Let us construct the map\qss $L_{\fff k}^{\fff m\dff +\dff 1}$\nnsp.\qff\oss
If\qss $k\qff =\qff m$\nnsp,\oss
then\qss
(\ref{l-compatibility})\qss
with\dss $n\qff =\qff m$\dss
imposes no restrictions on\dss $L_{\fff k}^{\fff m\dff +\dff 1}$\dss
and we can construct\dss $L_{\fff k}^{\fff m\dff +\dff 1}$\dss as in Section\qss \ref{amenable}.

Hence we may assume that\qss $k\qff \leq\qff m\qff -\qff 1$\nnsp.\oss
In the rest of the proof
we will omit the subscripts and superscripts of\dss $L$\nnsp.\pss
Let\oss 
$\overline{\sigma}
\qff \colon\qff
\Delta_{\fff k\dff -\dff 1}\ttoo X_{\fff m}$\oss
be an arbitrary singular $(k\qff -\qff 1)$\nsp\dnsp-simplex\fff,\oss 
and\dss let\oss\vspace*{1.5pt}
\[
\quad
\sigma\off =\off p_{\fff m}\dff(\dff \overline{\sigma} \dff)
\hspace*{1em}\mbox{ and }\hspace*{1em}
\tau\off =\off L\trf (\dff \sigma \dff)\dff.
\]

\vspace*{-10.5pt}
We need to define\dss 
$L\fff\halfff (\dff \overline{\sigma} \dff)$\dss 
in such a way that\oss 
$p_{\fff m}\dff \circ\qff L\qff (\trf \overline{\sigma} \trf)
\off =\off  
\tau$\nnsp.\oss 
The boundary\dss 
$\partial\trf L\trf (\trf \overline{\sigma} \trf)$\dss 
is already defined and we can take as\dss 
$L\trf (\trf \overline{\sigma} \trf)$\dss 
any singular $k$\nsp\dnsp-simplex\qss 
$\overline{\tau}\dff \colon\dff \Delta_{\fff k}\ttoo X_{\fff m\dff +\dff 1}$\qss 
such that\oss\vspace*{1.5pt}
\[
\quad
\partial\trf \overline{\tau}
\off =\off
\partial\trf L\trf (\trf \overline{\sigma} \trf)
\hspace*{1em}\mbox{ and }\hspace*{1em}
p_{\fff m}\dff(\dff \overline{\tau} \dff)
\off =\off 
\tau\dff.
\]

\vspace*{-10.5pt}
The boundary\dss 
$\partial\trf L\trf (\trf \overline{\sigma} \trf)$\dss
defines a continuous map\oss  
$\rho 
\qff \colon\qff 
b\fff\Delta_{\fff k}
\ttoo 
X_{\fff m\dff +\dff 1}$\nnsp,\oss
where\dss 
$b\fff\Delta_{\fff k}$\dss
is the geometric boundary of the simplex 
$\Delta_{\fff k}$\nnsp.\oss
The inductive assumption implies that\qss\vspace*{3pt}
\[
\quad
p_{\fff m}\dff \circ\dff \rho
\off =\off
\tau\dff | \dff b\fff\Delta_{\fff k}\dff,
\]

\vspace{-9pt}
where\qss
$\tau\dff | \dff b\fff\Delta_{\fff k}$\qss
is the restriction of\dss $\tau$\dss
to\dss
$b\fff\Delta_{\fff k}$\nnsp.\oss 
Let\oss 
$\widetilde{\tau}
\dff \colon\dff \Delta_{\fff k}\ttoo X_{\fff m\dff +\dff 1}$\oss 
be 
an arbitrary singular $k$\nsp\dnsp-simplex 
such that\oss 
$p_{\fff m}\pff \circ\pff \widetilde{\tau}
\off =\off
\tau$\nnsp.\oss 
Since\dss $p_{\fff m}$\dss is a principal $G_{\fff m}$\nsp\dnsp-bundle,\oss 
the map\dss 
$\rho$\dss differs from the restriction\qss
$\widetilde{\tau}\dff | \dff b\fff\Delta_{\fff j}$\qss
by a continuous map\oss 
\[
\quad
d
\dff \colon\dff
b\fff\Delta_{\fff k}
\ttoo
G_{\fff m}\dff. 
\]
Since\dss $b\fff\Delta_{\fff k}$\dss
is homeomorphic to the sphere of dimension\qss 
$k\qff -\qff 1\qff \leq\qff (\fff m\qff -\qff 1\fff)\qff -\qff 1\qff =\qff m\qff -\qff 2$\qss
and\dss $G_{\fff m}$\dss is\dss $(m\qff -\qff 2)$\nsp\dnsp-connected,\oss
one can extend\dss $d$\dss
to a continuous map\oss 
\[
\quad
D
\dff \colon\dff
\Delta_{\fff k}
\ttoo
G_{\fff m}\oss
\]
and take as\dss $\overline{\tau}$\dss the map differing from\qss $\widetilde{\tau}$\qss
by\qss $D$\nnsp.\oss
Since the singular simplex\oss 
$\overline{\sigma}
\qff \colon\qff
\Delta_{\fff k\dff -\dff 1}\ttoo X_{\fff m}$\oss
was arbitrary,\oss
this completes the induction step and hence the proof of the lemma.\oss  \eproof

\mypar{Theorem.}{simply-connected-homotopy} 
\emph{Let\qss $X$\dss 
be a path-connected space homotopy equivalent to a\dss CW-com\-plex\halfff.\oss 
If\qss $\pi_{\fff 1}\fff(\dff X \dff)$\qss is amenable,\oss
then there exists a strictly bounded contracting homotopy}\vspace*{3pt}
\[
\quad
\begin{tikzcd}[column sep=large, row sep=normal]\dis
\rrr
& 
B^{0}\fff(\dff X \dff) \arrow[l, "\dis K^{\fff 0}"]
& 
B^{1}\fff(\dff X \dff) \arrow[l, "\dis K^{\fff 1}"]
&   
\quad \ldots\off, \arrow[l, "\dis K^{\fff 2}"]
\end{tikzcd}
\]

\vspace*{-9pt}
\emph{i.e.\qss a chain homotopy
between\pss $\id$\qss and\pss $0$\dss 
such that\pss 
$\|\dff K^{\fff n} \dff\|\qff \leq\qff 1$\qss for all $n$\nnsp.\oss}

\proof\qss
The iterated killing of homotopy 
groups leads
to a sequence of maps and spaces of the form\qss (\ref{tower-principal})\qss 
such that\oss $X_{\fff 1}\off =\off X$\nnsp,\oss
$\pi_{\fff i}\dff(\dff X_{\fff n}\fff)\qff =\qff 0$\qss
if\qss $i\qff <\qff n$\nnsp,\oss
and each map\dss 
$p_{\fff n}$\dss 
is a principal $G_{\fff n}$\nsp\dnsp-bundle
for a topological group $G_{\fff n}$
which is either abelian or\dss is discrete and amenable.\oss
Each\dss
$G_{\fff n}$\dss is
an Eilenberg--MacLane space of type\qss 
$K\fff(\dff \pi_{\fff n}\fff(\dff X \dff)\fff,\pff n\dff -\dff 1\dff)$\dnsp,\oss
and\dss hence is\dss
$(n\qff -\qff 2)$\nsp\dnsp-connected for every $n$\nnsp.\oss
By
Lemma\qss \ref{l-exist}\qss 
there exist compatible sequences\oss 
$L^{1},\off\off L^2,\off\off L^3,\off\off \ldots \off\off$\dnsp,\oss
and\dss by\dss Lemma\qss \ref{l-to-h}\qss 
there exist compatible strictly bounded partial homotopies\oss
$K_{\fff 1}^{\fff *}\dff,\off\off 
K_{\fff 2}^{\fff *}\dff,\off\off 
K_{\fff 3}^{\fff *}\dff,\off\off \ldots \off\off$\dnsp.\oss

In order to speak about compatible partial homotopies,\oss
it was sufficient to know that
the maps\vspace*{3.5pt}
\[
\quad
p_{\fff n\dff *}
\dff \colon\dff 
B^{\fff i}\fff(\dff X_{\fff n\dff +\dff 1} \dff)
\toto 
B^{\fff i}\fff(\dff X_{\fff n} \dff)\qss
\]

\vspace*{-8.5pt}
from Theorem\qss \ref{chain-map}\qss
are defined for all\qss
$i\qff \leq\qff n\qff -\qff 1$\nnsp.\oss
But in view of Theorem\qss \ref{coherent-exists-infinite}\qss
one may assume that they are defined for all\qss $i\qff \geq\qff 0$\nnsp.\oss
Let\oss
$K^{\fff i}
\dff \colon\dff 
B^{i}\fff(\dff X \dff)
\ttoo 
B^{i\dff -\dff 1}\fff(\dff X \dff)$\oss 
be the map\vspace*{5pt}
\[
\quad
K^{\fff i}
\off =\off 
p_{\fff 1\dff *}
\qff \circ\qff 
\off \ldots\off
\qff \circ\pff
p_{\fff m\dff *}
\qff \circ\off
K^{\fff i}_{\fff m}
\pff \circ\pff
p_{\fff m}^{\fff *}
\pff \circ\pff
\off \ldots\off
\qff \circ\pff 
p_{\fff 1}^{\fff *}\qff,
\]

\vspace*{-7pt}
where $m$ is any integer\dss $\geq\qff i\qff +\qff 1$\nnsp.\oss
Condition\dss (\ref{homotopy-compatibility})\dss 
implies that this definition does not depend on the choice of $m$\nnsp.\oss
The sequence $K^{\fff \bullet}$ turns out to be a contracting homotopy.\oss
Indeed,\vspace*{6pt}
\[
\quad
\partial
\pff \circ\pff 
K^{\fff i}
\off +\off 
K^{\fff i\dff +\dff 1}
\qff \circ\pff 
\partial
\off\off =\off\off
\partial\pff \circ\qff p_{\fff 1\dff *}\qff \circ\qff
\off \ldots\quad 
+
\quad \ldots\off
\qff \circ\pff 
p_{\fff 1}^{\fff *}
\qff \circ\pff 
\partial 
\]

\vspace*{-30pt}
\[
\quad
=\off\off
p_{\fff 1\dff *}
\qff \circ\qff 
\off \ldots\off
\qff \circ\qff
p_{\fff m\dff *}
\qff \circ\pff
\partial
\qff \circ\qff
K^{\fff i}_{\fff m}
\pff \circ\pff
p_{\fff m}^{\fff *}
\qff \circ\qff
\off \ldots\off
\qff \circ\qff 
p_{\fff 1}^{\fff *}
\]

\vspace*{-30pt}
\[
\quad
\hspace*{6.52em}
\off +\off\off
p_{\fff 1\dff *}
\qff \circ\qff 
\off \ldots\off
\qff \circ\qff
p_{\fff m\dff *}
\qff \circ\pff
K^{\fff i\dff +\dff 1}_{\fff m}
\pff \circ\pff
\partial
\pff \circ\pff
p_{\fff m}^{\fff *}
\qff \circ\qff
\off \ldots\off
\qff \circ\qff 
p_{\fff 1}^{\fff *}
\]

\vspace*{-30pt}
\[
\quad
=\off\off
p_{\fff 1\dff *}
\qff \circ\qff 
\off \ldots\off
\qff \circ\qff
p_{\fff m\dff *}
\qff \circ\off
\left(\qff \partial
\pff \circ\pff
K^{\fff i}_{\fff m}
\off +\off
K^{\fff i\dff +\dff 1}_{\fff m}
\pff \circ\pff
\partial
\qff\right)
\off \circ\qff
p_{\fff m}^{\fff *}
\qff \circ\qff
\off \ldots\off
\qff \circ\qff 
p_{\fff 1}^{\fff *}
\]

\vspace*{-30pt}
\[
\quad
=\off\off
p_{\fff 1\dff *}
\qff \circ\qff 
\off \ldots\off
\qff \circ\qff
p_{\fff m\dff *}
\qff \circ\qff
p_{\fff m}^{\fff *}
\qff \circ\qff
\off \ldots\off
\qff \circ\qff 
p_{\fff 1}^{\fff *}
\off\off =\off\off
\id\
\]

\vspace*{-3pt}
for every\qss $m\qff \geq\qff i\qff +\qff 1$\nnsp,\oss
where we used the facts that\pss
$\partial
\pff \circ\pff
K^{\fff i}_{\fff m}
\off +\off
K^{\fff i\dff +\dff 1}_{\fff m}
\pff \circ\pff
\partial
\off =\off
\id$\pss
and\dss that\qss 
$p_{\fff n\dff *}$\nnsp,\pss $p_{\fff n}^{\fff *}$\pss
commute with the differentials.\oss
Moreover\halfff,\qss
$\|\dff K^{\fff i} \dff\|\qff \leq\qff 1$\qss
because all the norms\qss
$\|\dff p_{\fff n\dff *} \dff\|$\nnsp,\qss
$\|\dff K^{\fff i}_{\fff m} \dff\|$\nnsp,\qss
$\|\dff p_{\fff n}^{\fff *} \dff\|$\qss
are\qss
$\leq\qff 1$\nnsp.\oss
Hence\qss $K^{\fff *}$\qss
is the promised chain homotopy.\oss  \eproof

\mysection{Weak\qss equivalences}{weak-equivalences}

\vspace*{9pt}
\myuppar{The\dss $l_{\fff 1}$\dnsp-norm of chains.}
The\dss $l_{\fff 1}$\dnsp-norm\qss $\|\qff  c \qff\|_{\dff 1}$\qss
of a singular chain\qss\vspace*{3pt} 
\[
\quad 
c
\off\off =\off\off
\sum_{\sigma} \qff c_{\dff \sigma}\dff \sigma
\qff \in\qff C_{\fff m}\fff(\dff X \dff)\dff,
\]

\vspace*{-12pt} 
where\qss $c_{\dff \sigma}\qff \in\qff \rrr$\qss 
and the sum is taken over all singular simplices\qss 
$\sigma\qff \in\qff S_m\fff(\dff X \dff)$\dnsp,\oss
is defined as\oss\vspace*{3pt} 
\[
\quad 
\|\qff  c \qff\|_{\dff 1}
\off\off =\off\off
\sum_{\sigma} \qff |\qff c_{\dff \sigma}\qff |
\]

\vspace*{-12pt}
By the definition of singular chains,\oss the sums above involve only finite number
of non-zero coefficients\dss $c_{\dff \sigma}$\nnsp,\oss
and\dss hence\qss $\|\qff  c \qff\|_{\dff 1}$\qss is a well defined real number\qss
({\fff}i.e.\qss is\qss $<\qff \infty$\nnsp).\oss
The\dss $l_{\fff 1}$\dnsp-norm is one of the main notions of\dss Gromov's theory,\oss
but in this section it plays only a technical role.

\vspace*{3pt}
\myuppar{Natural chain homotopies.}
Let\qss $\iota_n\dff \colon\dff
\Delta_{\fff n}\ttoo \Delta_{\fff n}$\qss
be the identity map of\dss $\Delta_{\fff n}$\dss considered
as a singular simplex\qss 
$\iota_n\qff \in\qff S_n\fff(\dff \Delta_{\fff n}\dff)$\dnsp.\oss
There exist singular chains\qss
$h_{\fff n}
\qff \in\qff 
C_{\fff n\dff +\dff 1}\fff(\dff \Delta_{\fff n}\dff \times\qff  I \dff)$\qss
such that\vspace*{3pt}
\begin{equation}
\label{universal-homotopy}
\quad  
\partial\fff h_{\fff n}
\off\off =\off\off
\iota_n\dff \times\dff 1 
\off -\off
\iota_n\dff \times\dff 0
\off\off -\off\off
\sum_{i\qff =\qff 0}^n\off
(\dff -\qff 1\fff)^{\fff i}\pff
\bigl(\dff
\partial_{\fff i}\dff\iota_n\dff \times\qff I
\trf\bigr)_*\fff (\dff h_{\fff n\dff -\dff 1} \fff)\dff,
\end{equation}

\vspace{-9pt}
for all\qss $n\qff \geq\qff 0$\nnsp,\oss
where\qss 
$\partial_{\fff i}\dff\iota_n
\dff \colon\dff
\Delta_{\fff n\dff -\dff 1}\ttoo \Delta_{\fff n}$\qss
is the\dss $i$\dnsp-th face of the singular simplex\dss $\iota_n$\nnsp.\oss
Such chains are provided by 
natural chain homotopy between 
the morphisms\qss\vspace*{3pt} 
\[
\quad
C_{\fff *}\fff(\dff \Delta_{\fff n} \dff)
\qff \ttoo\qff
C_{\fff *}\fff(\dff \Delta_{\fff n}\dff \times\qff I \dff)
\]

\vspace*{-9pt}
induced by the maps\oss  
$\iota_n\dff \times\dff 0\fff,\pff 
\iota_n\dff \times\dff 1
\qff \colon\qff
\Delta_{\fff n}
\qff \ttoo\qff
\Delta_{\fff n}\dff \times\qff I\qff, 
$\oss 
or\halfff,\oss better\halfff,\oss
such chains are the main ingredient in the construction of 
natural chain homotopy
between the morphisms\qss\vspace*{3pt} 
\[
\quad
C_{\fff *}\fff(\dff Y \dff)
\qff \ttoo\qff
C_{\fff *}\fff(\dff Y\dff \times\qff I \dff)
\]

\vspace*{-9pt}
induced by the maps\oss 
$\id_{\dff Y}\dff \times\dff 0\fff,\pff 
\id_{\dff Y}\dff \times\dff 1
\qff \colon\qff
Y
\qff \ttoo\qff
Y\dff \times\qff I$\oss
for an arbitrary space $Y$\dnsp.\oss

\vspace*{3pt}
\myuppar{\dnsp$k$\dnsp-connected pairs of spaces.}
Let $X$ be a topological space and $A$ be a subspace of $X$\nnsp.\oss
As usual,\oss let
$D^m$ be the standard $m$\dnsp-dimensional disc and
$S^{\fff m\dff -\dff 1}$ be its boundary sphere.\oss
The pair $(\dff X\fff,\pff A \fff)$ is said to be\dss 
\emph{$k$\dnsp-connected}\oss if\dss for\qss 
$m\qff \leq\qff k$\qss every map\qss\vspace*{3pt}
\[
\quad
f\dff \colon\dff (\dff D^m\fff,\pff S^{\fff m\dff -\dff 1}\fff)
\ttoo
(\dff X\fff,\pff A \fff)
\]

\vspace*{-9pt}
is homotopic relatively to $A$ to a map having image in $A$\nnsp.\oss

\myuppar{Bounded Eilenberg complexes.}
We need a bounded cohomology version of a classical construction from the singular homology theory
going back to Eilenberg\qss \cite{e}.\oss
For natural number $k$
let\dss $(\fff \Delta_{\fff n}\fff)_k$\dss be the union of all\qss $\leq\qff k$\dnsp-dimensional 
faces of the standard $n$\dnsp-simplex 
$\Delta_{\fff n}$\nnsp.\oss
Let\qss $S_n\fff(\dff X\fff,\pff A \fff)_k$\qss
be the set of all singular simplices\qss
$\sigma\dff \colon\dff \Delta_{\fff n}\ttoo X$\qss
such that\vspace*{3pt}
\[
\quad
\sigma\dff\bigl(\qff (\fff \Delta_{\fff n}\fff)_k \qff\bigr)
\off \subset\off
A\qff.
\]

\vspace{-9pt}
If\qss
$\sigma\qff \in\qff S_n\fff(\dff X\fff,\pff A \fff)_k$\dnsp,\oss
then\qss
$\partial_i\fff\sigma\qff \in\qff S_{n\dff -\dff 1}\fff(\dff X\fff,\pff A \fff)_k$\qss
for every face\dss $\partial_i\fff\sigma$\dss of\dss $\sigma$\nnsp.\oss
Let\dss $B^{\fff n}\fff(\dff X\fff,\pff A \dff)_k$\dss
be the space of all bounded functions\qss
$S_n\fff(\dff X\fff,\pff A \fff)_k
\ttoo
\rrr$\nnsp.\oss
The space\dss $B^{\fff n}\fff(\dff X\fff,\pff A \dff)_k$\dss
is a Banach space with respect to the supremum norm.\oss
The coboundary maps\qss\vspace*{3pt}
\[
\quad
d_{\dff n}
\qff \colon\qff
B^{\fff n}\fff(\dff X\fff,\pff A \dff)_k
\qff \ttoo \qff
B^{\fff n\dff +\dff 1}\fff(\dff X\fff,\pff A \dff)_k
\]

\vspace*{-9pt}
are defined
by the same formula as for $B^{\fff \bullet}\dff(\dff X \dff)$\dnsp.\oss
The spaces\qss $B^{\fff n}\fff(\dff X\fff,\pff A \dff)_k$\qss
together with the co\-bound\-a\-ry maps\dss $d_{\dff n}$\dss form a cochain complex
denoted by\qss
$B^{\fff \bullet}\fff(\dff X\fff,\pff A \dff)_k$\dnsp.\oss
Let\qss
$\widehat{H}^{\fff n}\fff(\dff X\fff,\pff A \dff)_k$\qss
be the cohomology groups of this complex\halfff.\oss
These cohomology groups
inherit a semi-norm from the Banach norm of\dss the spaces\dss
$B^{\fff n}\fff(\dff X\fff,\pff A \dff)_k$\nnsp.\oss
The obvious restriction map\vspace*{3pt}
\[
\quad 
\rho
\qff \colon\qff
B^{\fff \bullet}\dff(\dff X \dff)
\qff \ttoo\qff
B^{\fff \bullet}\fff(\dff X\fff,\pff A \dff)_k
\]

\vspace*{-9pt}
commutes with the boundary operators and hence induces homomorphisms\vspace*{3pt}
\begin{equation*}
\quad
\rho_{\fff *}
\qff \colon\qff
\widehat{H}^{\fff n}\fff(\dff X \dff)
\qff \ttoo\qff
\widehat{H}^{\fff n}\fff(\dff X\fff,\pff A \dff)_k
\end{equation*}

\vspace*{-9pt}
\mypar{Lemma.}{k-equivalence-lemma}
\emph{If\qss $(\dff X\fff,\pff A \fff)$\qss
is $k$\dnsp-connected,\oss
then\qss $\rho_{\fff *}$\dss is an isometric isomorphism for
all\qss $n$\nnsp.}

\proof\qss
If\qss $(\dff X\fff,\pff A \fff)$\qss
is $k$\dnsp-connected,\oss
then one can assign to each\qss 
$\sigma\qff \in\qff S_n\fff(\dff X \fff)$\qss
a homotopy\qss\vspace*{3pt}
\[
\quad  
P\fff(\fff \sigma\fff)
\dff \colon\dff 
\Delta_{\fff n}\dff \times\qff I 
\qff \ttoo\qff
X\dff,
\]

\vspace{-9pt}
where\qss $I\qff =\qff [\fff 0\fff,\pff 1\fff]$\nnsp,\oss
in such a way that the following conditions hold.\vspace{-9pt}
\begin{itemize}
\item[({\fff}i{\fff})] \dnsp$P\fff(\fff \sigma\fff)_0\off =\off \sigma$\nnsp.\oss
\item[({\fff}ii{\fff})] \dnsp$P\fff(\fff \sigma\fff)_1\qff \in\qff S_n\fff(\dff X\fff,\pff A \fff)_k$\nnsp.\oss
\item[({\fff}iii{\fff})] If\pss $\sigma\qff \in\qff S_n\fff(\dff X\fff,\pff A \fff)_k$\nnsp,\oss
then\qss $P\fff(\fff \sigma\fff)_t$\qss is a constant homotopy.
\item[({\fff}iv{\fff})] \dnsp$P\fff(\fff \sigma\fff)
\dff \circ\dff 
(\trf \partial_i\fff\sigma\dff \times\dff \id_{\dff I} \trf)
\off =\off
P\fff(\trf \partial_i\fff\sigma \dff)$\oss
for each face\dss $\partial_i\fff\sigma$\dss of\dss $\sigma$\dnsp.
\end{itemize}

\vspace{-6pt}
The construction is similar to the one of the maps\dss $L^{n}$\dss
in Section\qss \ref{amenable}.\oss
See\qss \cite{sp},\oss Section\qss 7.4,\oss Lemma\qss 7\qss and Theorem\qss 8,\oss
or\qss \cite{td},\oss Section\qss 9.5.

Let us define\qss
$\alpha
\dff \colon\dff
B^{\fff n}\fff(\dff X\fff,\pff A \dff)_k
\ttoo
B^{\fff n}\dff(\dff X \dff)$\qss by\vspace*{3pt} 
\[
\quad
\alpha\fff(\dff f\dff)(\fff \sigma\fff)
\off =\off
f\fff(\dff P\fff(\fff \sigma\fff)_1 \dff)\dff.
\]

\vspace*{-9pt}
By the property\qss ({\fff}ii{\fff})\qss the map $\alpha$ is well defined,\oss
by\qss ({\fff}iv{\fff})\qss the map $\alpha$ is a chain map,\oss
and\qss ({\fff}iii{\fff})\qss implies that\qss
$\rho\dff \circ\dff \alpha$\qss is the identity map.\oss

Let us now define\qss
$k_{\fff n\dff +\dff 1}
\dff \colon\dff
B^{\fff n\dff +\dff 1}\dff(\dff X \dff)
\ttoo
B^{\fff n}\dff(\dff X \dff)$\qss by\vspace*{3pt} 
\[
\quad
k_{\fff n\dff +\dff 1}\fff(\dff f\dff)(\fff \sigma\fff)
\off =\off
f\dff\bigl(\dff P\fff(\fff \sigma\fff)_* (\dff h_{\fff n} \dff)\dff\bigr)\dff,
\]

\vspace*{-9pt}
where in the right hand side the cochain\qss
$f\qff \in\qff B^{\fff n\dff +\dff 1}\dff(\dff X \dff)$\qss
is interpreted as a linear functional\qss
$C_{\fff n\dff +\dff 1}\dff(\dff X \dff)
\ttoo 
\rrr$\nnsp.\oss
A priori\dss $k_{\fff n\dff +\dff 1}\fff(\dff f\dff)$\dss is only a singular cochain,\oss
and we need to check that this cochain is bounded.\oss
But\halfff,\oss obviously,\oss\vspace*{3pt} 
\[
\quad
\left|\qff k_{\fff n\dff +\dff 1}\fff(\dff f\dff)(\fff \sigma\fff) \qff\right|
\off \leq\off
\|\qff f \qff\|\cdot \|\qff P\fff(\fff \sigma\fff)_* (\dff h_{\fff n} \dff) \qff\|_{\dff 1}
\off \leq\off
\|\qff f \qff\|\cdot \|\qff  h_{\fff n} \qff\|_{\dff 1}\dff,
\]

\vspace{-9pt}
and\dss hence the cochain\dss
$k_{\fff n\dff +\dff 1}\fff(\dff f\dff)$\dss is indeed bounded.\oss
Moreover\halfff,\oss the above inequalities imply that\dss $k_{\fff n\dff +\dff 1}$\dss
is a bounded operator and\qss
$\|\qff k_{\fff n\dff +\dff 1} \qff\|
\qff \leq\qff
\|\qff  h_{\fff n} \qff\|_{\dff 1}$\nnsp,\oss
but this is not needed for the proof\halfff.

The properties\qss ({\fff}i{\fff})\qss and\qss ({\fff}iv{\fff})\qss
together with\qss (\ref{universal-homotopy})\qss
imply that\dss $k_{\fff \bullet}$\dss is a chain homotopy between\qss
$\alpha\dff \circ\dff \rho$\qss and the identity morphism of\dss
$B^{\fff *}\dff(\dff X \dff)$\dnsp.\oss
Cf.\qss \cite{sp},\oss Section\qss 7.4,\oss the proof of Lemma\qss 7,\oss
or\qss \cite{td},\oss the end of the proof of Theorem\qss 9.5.1.\oss
It follows that\qss
$\alpha\dff \circ\dff \rho$\qss
induces the identity map of the cohomology of\dss
$B^{\fff *}\dff(\dff X \dff)$\dnsp,\oss
i.e.\qss of\dss the bounded cohomology groups\dss
$\widehat{H}^{\fff n}\fff(\dff X \dff)$\dnsp.\oss
Since\qss
$\rho\dff \circ\dff \alpha$\qss is the identity map,\oss
it follows that the chain maps $\alpha$ and $\rho$ induce 
the mutually inverse isomorphisms between\dss
$\widehat{H}^{\fff n}\fff(\dff X \dff)$\dss
and\dss
$\widehat{H}^{\fff n}\fff(\dff X\fff,\pff A \dff)_k$\nnsp.\oss

It remains to prove that these isomorphisms are isometric.\oss
It is sufficient to prove that both induced isomorphisms in cohomology
are bounded with the norm\qss $\leq\qff 1$\nnsp.\oss
In order to prove this,\oss
it is sufficient to prove that
$\alpha$ and $\rho$ are
bounded operators with the norm\qss $\leq\qff 1$\nnsp.\oss
The operator $\rho$ is defined as a restriction and\dss hence its norm is\qss $\leq\qff 1$\nnsp.\oss 
The operator $\alpha$ is dual to a map between sets of simplices,\oss
and\dss hence\qss $\|\qff \alpha \qff\|\qff \leq\qff 1$\nnsp.\oss
In more details,\oss\vspace*{3pt}
\[
\quad
\left|\qff \alpha\fff(\dff f\dff)(\fff \sigma\fff) \qff\right|
\off =\off
\left|\qff f\fff(\dff P\fff(\fff \sigma\fff)_1 \dff) \qff\right|\qff 
\off \leq\off
\|\qff f \qff\|
\]

\vspace*{-9pt}
because\dss $P\fff(\fff \sigma\fff)_1$\dss is a singular simplex\qss
(compare with the inequalities for the norm\dss $\|\qff k_{\fff n\dff +\dff 1} \qff\|$\nnsp).\oss
It follows that\qss $\|\qff \alpha \qff\|\qff \leq\qff 1$\nnsp.\oss
This completes the proof\halfff.\oss  \eproof

\mypar{Theorem.}{k-equivalence-theorem}
\emph{If\dss the pair\qss $(\dff X\fff,\pff A \fff)$\qss
is $k$\dnsp-connected,\oss
then the homomorphism}\vspace*{2pt}
\begin{equation*}
\quad
i^{\fff *}
\qff \colon\qff
\widehat{H}^{\fff n}\fff(\dff X \dff)
\qff \ttoo\qff
\widehat{H}^{\fff n}\fff(\dff A \dff)
\end{equation*}

\vspace*{-10pt}
\emph{induced by the inclusion\qss $i\dff \colon\dff A\ttoo X$\qss is an isometric isomorphism
for\qss $n\qff \leq\qff k\qff -\qff 1$\nnsp.\oss}

\proof\qss
If\qss $n\qff \leq\qff k$\nnsp,\oss
then,\oss obviously,\oss
$S_n\fff(\dff X\fff,\pff A \fff)_k
\off =\off
S_n\fff(\fff A \fff)$\oss
and\dss hence\oss
$B^{\fff n}\fff(\dff X\fff,\pff A \fff)_k
\off =\off
B^{\fff n}\fff(\fff  A \fff)$\dnsp.\oss
Moreover\halfff,\oss
the restriction map $\rho$ is nothing else but the restriction map\qss
$B^{\fff n}\fff(\dff  X \dff)
\ttoo
B^{\fff n}\fff(\fff  A \fff)$\dnsp.\oss
In view of this,\oss the theorem follows from Lemma\qss \ref{k-equivalence-lemma}.\oss
({\fff}In order to ensure the isomorphism property for\qss $n\qff =\qff k$\nnsp,\oss
an isomorphism of cochains in dimension\dss $k\qff +\qff 1$\qss is needed.)  \eproof

\mypar{Theorem.}{k-equivalence-maps}
\emph{Suppose that the spaces\qss $X\fff,\pff A$\qss are path-connected and the map\qss
$\varphi\dff \colon\dff
A\ttoo X$\qss is\dss a\dss $k$\dnsp-equivalence,\oss 
i.e.\qss the induced map of\dss homotopy groups}\vspace*{3pt}
\[
\quad
\varphi_{*}
\qff \colon\qff
\pi_{\fff n}\fff(\dff A\fff,\pff a \dff )
\qff \ttoo\qff
\pi_{\fff n}\fff(\dff X\fff,\pff \varphi(\fff a\fff) \dff )
\]

\vspace*{-9pt}
\emph{is an isomorphisms for\qss $n\qff \leq\qff k\qff -\qff 1$\qss
and an epimorphisms for\qss $n\qff =\qff k$\nnsp.\oss
Then}\vspace*{3pt}
\[
\quad
\varphi^{*}
\qff \colon\qff
\widehat{H}^{\fff n}\fff(\dff X \dff)
\qff \ttoo\qff
\widehat{H}^{\fff n}\fff(\dff A \dff)
\]

\vspace*{-9pt}
\emph{is an isometric isomorphism for\qss $n\qff \leq\qff k\qff -\qff 1$\nnsp.\oss}

\proof\qss
The mapping cylinder $Z$ of $\varphi$
contains $X$ as a deformation retract and $A$ as a subspace.\oss
Moreover\halfff,\oss the inclusion\qss $A\ttoo Z$\qss is homotopic in $Z$
to the composition of $\varphi$ with the inclusion\qss $X\ttoo Z$\nnsp.\oss
As in the usual cohomology theory,\oss this reduces the theorem to the case
when the map\qss $A\ttoo X$\qss is the inclusion of a subspace.\oss
But if $A$ is a subspace of $X$\nnsp,\oss
then the map\qss $A\ttoo X$\qss 
is a\dss $k$\dnsp-equivalence if and only if the pair\dss
$(\dff X\fff,\pff A \fff)$\dss is\dss $k$\dnsp-connected.\oss
Therefore,\oss the theorem follows from Theorem\qss \ref{k-equivalence-theorem}.\oss  \eproof

\mypar{Corollary.}{weak-equivalence-maps}
\emph{Suppose that the spaces\qss $X\fff,\pff A$\qss are path-connected and the map\qss
$\varphi\dff \colon\dff
A\ttoo X$\qss is\dss a weak equivalence,\oss 
i.e.\qss the induced map of\dss homotopy groups}\vspace*{3pt}
\[
\quad
\varphi_{*}
\qff \colon\qff
\pi_{\fff n}\fff(\dff A\fff,\pff a \dff )
\qff \ttoo\qff
\pi_{\fff n}\fff(\dff X\fff,\pff \varphi(\fff a\fff) \dff )
\]

\vspace*{-9pt}
\emph{is an isomorphisms for all\dss $n$\nnsp.\oss
Then}\vspace*{3pt}
\[
\quad
\varphi^{*}
\qff \colon\qff
\widehat{H}^{\fff n}\fff(\dff X \dff)
\qff \ttoo\qff
\widehat{H}^{\fff n}\fff(\dff A \dff)
\]

\vspace*{-9pt}
\emph{is an isometric isomorphism for all\dss $n$\nnsp.\oss}  \eproof

\myuppar{Extending the results to arbitrary topological spaces.}
As is well known,\oss every path-con\-nected space is weakly homotopy equivalent to a\dss CW-complex\halfff.\oss
In view of this,\oss
Corollary\qss \ref{weak-equivalence-maps}\qss allows to extend the main results of the bounded
cohomology theory from spaces homotopy equivalent to CW-complexes to arbitrary spaces.\oss
For example,\oss the conclusions of Theorems\qss \ref{simply-connected-homology}\qss
and\qss \ref{simply-connected-homotopy}\qss hold for arbitrary spaces.

\mysection{Elements\qss of\qss homological\qss algebra}{algebra}

\vspace*{6pt}
\myuppar{Bounded $G$\dnsp-modules.}
Let $G$ be a discrete group.\oss
A\qss \emph{bounded left\dss $G$-module}\qss 
is defined as a real semi-normed space $V$ 
together with a left action of $G$ on $V$ 
such that\qss 
$\|\dff g\cdot v \dff\|\qff \leq\qff \|\dff v \dff\|$\qss 
for all\qss 
$g\qff \in\qff G$\qss and\qss $v\qff \in\qff V$\nnsp.\oss 
We will call the bounded left $G$-modules simply \emph{$G$\dnsp-modules}.\oss
If $V$ and $W$ are two $G$\dnsp-modules,\oss
then a \emph{$G$\dnsp-morphism}\qss from $V$ to $W$ 
is defined as a bounded linear operator\qss $V\toto W$\qss 
commuting with the action of $G$\nnsp.\oss

For every semi-normed space $V$ there is an action of $G$ on $V$
defined by\qss
$g\cdot v\off =\off v$\qss for all\qss $g\qff \in\qff G$\nnsp,\pss $v\qff \in\qff V$\dnsp.\oss
This action and the corresponding structure of a\dss $G$\dnsp-module on $V$ are
called\qss \emph{trivial}.\oss
The simplest semi-normed space is $\rrr$ with the absolute value function being the semi-norm.\oss
The corresponding trivial\dss $G$\dnsp-module $\rrr$ is the simplest\dss $G$\dnsp-module.\oss

For a bounded left\dss $G$\dnsp-module $V$ let\dss $B\dff(\fff G\fff,\pff V\dff)$\dss 
be the space of\dss functions\qss 
$f\dff \colon\dff G\toto V$\qss  
such that\vspace*{3pt} 
\[
\quad
\|\dff f \dff\|
\off =\off 
\sup\qff \left\{\off \|\dff f\dff(\fff g\fff) \trf\| \off  \left| \off g\qff \in\qff G \off\right\} \right.
\off \leq\off \infty\qff. 
\]

\vspace*{-9pt}
The space $B\dff(\fff G\fff,\pff V \dff)$
is a Banach space with the norm $\|\dff f \trf\|$\nnsp,\oss 
and the action of $G$ defined by\vspace*{3pt} 
\[
\quad
(\fff h\cdot f \dff)\dff(\fff g\fff)
\off =\off
h\cdot\left(\dff
f\dff(\fff g\halfff h\fff)
\dff\right)
\]

\vspace*{-9pt} 
turns it into a bounded $G$\dnsp-module.\oss 
If\dss $V$ is a trivial $G$\dnsp-module,\oss 
then this action takes the form\vspace*{3pt} 
\[
\quad
(\fff h\cdot f \dff)\dff(\fff g\fff)
\off =\off
f\dff(\fff g\halfff h\fff)\dff.
\]

\vspace*{-9pt} 
\myuppar{Relatively injective $G$-modules.}
A $G$\dnsp-morphism of $G$\dnsp-modules\qss 
$i\dff \colon\dff V_1\toto V_2$\qss
is said to be\qss \emph{strong\-ly injective}\qss 
if there exists a bounded linear map\qss 
$\sigma\dff \colon\dff V_2\toto V_1$\qss such that\qss 
$\sigma\dff \circ\dff i\qff =\qff \id$\qss 
and\qss 
$\|\dff \sigma \dff\|\qff \leq\qff 1$\nnsp.\oss 
Here the map $\sigma$ is not assumed to be a morphism of $G$\dnsp-modules.\oss
Obviously,\oss a strongly injective $G$\dnsp-morphism is injective.\oss

A $G$-module $U$ is said to be\qss \emph{relatively injective}\qss 
if for every strongly injective $G$\dnsp-morphism of $G$\dnsp-modules\qss 
$i\dff \colon\dff V_1\toto V_2$\qss and any $G$\dnsp-morphism of $G$\dnsp-modules\qss 
$\alpha\dff \colon\dff V_1\toto U$\qss 
there exists a $G$\dnsp-morphism\qss 
$\beta\dff \colon\dff V_2\toto U$\qss 
such that\qss 
$\beta\dff \circ\dff i
\off =\off 
\alpha$\qss 
and\qss 
$\|\dff \beta \dff\|\qff \leq\qff \|\dff \alpha \dff\|$\nnsp.\oss
See the following diagram.\vspace*{0pt}
\begin{equation*}
\quad
\begin{tikzcd}[column sep=huge, row sep=huge]\dis
V_1 
\arrow[r, shift left=3pt, "\dis i"]
\arrow[d, "\dis \alpha"']
& 
V_2 
\arrow[l, shift left=3pt, "\dis \sigma"]
\arrow[ld, dashed, "\dis \beta"] 
\\ 
U
& 
&  
\end{tikzcd}
\end{equation*}

\vspace{-9pt}
The following lemma provides us with all relatively injective modules we will need.

\mypar{Lemma.}{bgv}
\emph{For every\dss $G$\dnsp-module\dss $V$\dss the\dss $G$\dnsp-module\qss 
$B\dff(\fff G\fff,\pff V\dff)$\qss is relatively injective.\oss}

\proof\qss
Suppose that we are in the situation of the above diagram with\qss
$U\qff =\qff B\dff(G\fff,\pff V\fff)$\dnsp.\oss
Given $i$ and $\alpha$\nnsp,\oss we need to construct $\beta$\nnsp.\oss
Let us define $\beta$ by the formula\vspace*{3pt}
\[
\quad
\beta\dff(\fff w\dff)(\fff g\fff)
\off =\off
g^{\fff -\dff 1}\nsp\cdot\dff
\left(\dff 
\alpha\dff
\circ\dff \sigma\dff(\fff g\cdot w \dff)
(\fff 1\fff)
\dff\right)\dff,
\]

\vspace*{-9pt}
where\qss $w\qff \in\qff V_2$\qss and\qss $g\qff \in\qff G$\nnsp.\oss
A calculation,\oss which we,\oss contrary to the tradition,\oss
do not omit\halfff,\oss shows 
that $\beta$ commutes with $G$ and\qss 
$\beta\dff \circ\dff i
\qff =\qff
\alpha$\nnsp.\oss
Namely,\oss if\qss
$w\qff \in\qff V_2$\qss and\qss $g\fff,\pff h\qff \in\qff G$\nnsp,\oss
then\vspace*{6pt}
\[
\quad 
\beta\dff(\fff h\cdot w\fff)(\fff g\fff)
\off =\off
g^{\fff -\dff 1}\nsp\cdot\dff
\left(\dff 
\alpha\dff
\circ\dff \sigma\dff(\fff g\cdot(\fff h\cdot w \dff) \dff) 
(\fff 1 \fff)
\dff\right)
\]

\vspace*{-30pt}
\[
\quad
\phantom{\beta\dff(\fff h\cdot w\fff)(\fff g\fff)
\off }
=\off
g^{\fff -\dff 1}\nsp\cdot\dff
\left(\dff 
\alpha\dff
\circ\dff \sigma\dff(\fff (\fff g\fff h \fff)\cdot w \dff) \dff) 
(\fff 1 \fff)
\dff\right)
\]

\vspace*{-30pt}
\[
\quad
\phantom{\beta\dff(\fff h\cdot w\fff)(\fff g\fff)
\off }
=\off
\left(\dff
h\dff
(\fff g\fff h \fff)^{\fff -\dff 1}
\dff\right)
\dff \cdot\dff
\left(\dff 
\alpha\dff
\circ\dff \sigma\dff(\fff (\fff g\fff h \fff)\cdot w \dff) \dff) 
(\fff 1 \fff)
\dff\right)
\]

\vspace*{-30pt}
\[
\quad
\phantom{\beta\dff(\fff h\cdot w\fff)(\fff g\fff)
\off }
=\off
h\cdot
\Bigl(
(\fff g\fff h \fff)^{\fff -\dff 1}
\nsp\cdot\dff 
\left(\dff
\alpha\dff
\circ\dff \sigma\dff(\fff (\fff g\fff h \fff)\cdot w \dff) \dff) 
(\fff 1 \fff)
\dff\right)
\qff\Bigr)
\]

\vspace*{-30pt}
\[
\quad
\phantom{\beta\dff(\fff h\cdot w\fff)(\fff g\fff)
\off }
=\off
h\cdot
\left(\dff 
\beta\dff(\fff w\dff)(\fff g\fff h\fff)
\dff\right)
\off =\off
\left(\dff
h\cdot\beta\dff(\fff w \dff)
\dff\right)
(\fff g \fff)
\]

\vspace*{-6pt}
and\dss hence $\beta$ is a $G$\dnsp-morphism.\off\oss
Also,\oss if\qss $v\qff \in\qff V_1$\nnsp,\oss then\vspace*{6pt}
\[
\quad 
\beta\dff(\dff i\dff(\fff v \dff) \dff)(\fff g \fff)
\off =\off
g^{\fff -\dff 1}\nsp\cdot\dff
\left(\dff 
\alpha\dff
\circ\dff \sigma\dff(\fff g\cdot i\dff(\fff v \dff) \dff)
(\fff 1\fff)
\dff\right)
\]

\vspace*{-30pt}
\[
\quad
\phantom{\beta\dff(\dff i\dff(\fff v \dff) \dff)(\fff g \fff)
\off }
=\off
g^{\fff -\dff 1}\nsp\cdot\dff
\left(\dff 
\alpha\dff
\circ\dff \sigma\dff \circ\dff  i\dff(\dff g\cdot v \dff) 
(\fff 1\fff)
\dff\right)
\]

\vspace*{-30pt}
\[
\quad
\phantom{\beta\dff(\dff i\dff(\fff v \dff) \dff)(\fff g \fff)
\off }
=\off
g^{\fff -\dff 1}\nsp\cdot\dff
\left(\dff 
\alpha\dff
(\dff g\cdot   v \dff) 
(\fff 1\fff)
\dff\right)
\off =\off
g^{\fff -\dff 1}\nsp\cdot\dff
\Bigl(\dff
\left(\dff 
g\cdot \alpha\dff
(\dff  v \dff)
\dff\right) 
(\fff 1\fff)
\dff\Bigr)
\]

\vspace*{-30pt}
\[
\quad
\phantom{\beta\dff(\dff i\dff(\fff v \dff) \dff)(\fff g \fff)
\off =\off
g^{\fff -\dff 1}\nsp\cdot\dff
\left(\dff 
\alpha\dff
(\dff g\cdot   v \dff) 
(\fff 1\fff)
\dff\right)
\off }
=\off
g^{\fff -\dff 1}\nsp\cdot\dff
\Bigl(\dff 
g\cdot 
\left(\dff
\alpha\dff
(\dff  v \dff) 
(\fff g\fff)
\dff\right)
\dff\Bigr)
\]

\vspace*{-30pt}
\[
\quad
\phantom{\beta\dff(\dff i\dff(\fff v \dff) \dff)(\fff g \fff)
\off =\off
g^{\fff -\dff 1}\nsp\cdot\dff
\left(\dff 
\alpha\dff
(\dff g\cdot   v \dff) 
(\fff 1\fff)
\dff\right)
\off }
=\off
g^{\fff -\dff 1}\nsp\cdot\dff
g\cdot 
\left(\dff
\alpha\dff
(\dff  v \dff) 
(\fff g\fff)
\dff\right)
\off =\off
\alpha\dff(\fff v \dff)(\fff g\fff)
\]

\vspace*{-6pt}
and\dss hence\qss $\beta\dff \circ\dff i\qff =\qff \alpha$\nnsp.\off\oss
Finally,\pss
$\|\qff \beta \qff\|\qff \leq\qff \|\qff \alpha \qff\|$\oss
because for every\qss $w\qff \in\qff V_2$ \vspace*{6pt}
\[
\quad
\|\qff \beta(\fff w \dff)(\fff g \fff) \qff\|
\off =\off
\|\qff g^{\fff -\dff 1}\nsp\cdot\dff
\left(\dff 
\alpha\dff
\circ\dff \sigma\dff(\fff g\cdot w \dff)
(\fff 1\fff)
\dff\right)
\qff\|
\off \leq\off 
\|\qff \alpha\dff \circ\dff \sigma\dff(\fff g\cdot w \dff)(\fff 1 \fff) \qff\|
\]

\vspace*{-30pt}
\[
\quad
\phantom{\|\qff \beta(\fff w \dff)(\fff g \fff) \qff\|
\off }
\leq\off 
\|\dff \alpha \trf\| \cdot \|\dff \sigma \dff\| \cdot \|\dff g\cdot w \trf\|
\off \leq\off 
\|\dff \alpha \trf\| \cdot \|\dff w \trf\|
\hspace*{1.5em}\mbox{  \eproof  }
\]

\vspace{-6pt}
\myuppar{Resolutions.}
A\qss \emph{resolution}\pss or\halfff,\oss more precisely,\oss 
a\dss \emph{$G$\dnsp-resolution}\qss 
of a $G$\dnsp-module $U$ is defined as an exact sequence of $G$\dnsp-modules and 
$G$\dnsp-morphisms of the form\vspace*{3pt}
\begin{equation}
\label{resolution}
\quad
\begin{tikzcd}[column sep=large, row sep=normal]\dis
0 \arrow[r]
& 
U \arrow[r, "\dis d_{\dff -\dff 1}\off"]
& 
U_{0} \arrow[r, "\dis d_{\dff 0}\off"]
&   
U_{1} \arrow[r, "\dis d_{\dff 1}\off"]
&
U_{2} \arrow[r, "\dis d_{\dff 2}\off"]
&
\off \ldots \off.
\end{tikzcd}
\end{equation}

\vspace*{-6pt}
A\qss \emph{contracting homotopy}\pss for the
resolution\qss (\ref{resolution})\qss 
is a sequence of\dss bounded operators\vspace*{6pt}
\begin{equation*}
\quad
\begin{tikzcd}[column sep=large, row sep=normal]\dis
U 
& 
U_{0} \arrow[l, "\dis \off K_{\dff 0}"]
&   
U_{1} \arrow[l, "\dis \off K_{\dff 1}"]
&
U_{2} \arrow[l, "\dis \off K_{\dff 2}"]
&
\off \ldots \off \arrow[l]
\end{tikzcd}
\end{equation*}

\vspace{-9pt}
such that\oss\vspace*{3pt}
\begin{equation}
\label{homotopy-identity}
\quad
d_{\dff n\dff -\dff 1}\dff \circ\dff K_{\dff n}
\qff +\qff
K_{\dff n\dff +\dff 1}\dff \circ\dff d_{\dff n}
\off =\off
\id
\end{equation}

\vspace*{-6pt}
for\qss $n\qff \geq\qff 1$\nnsp,\oss
$K_{\dff 0}\dff \circ\dff d_{\dff -\dff 1}\off =\off \id_{\dff U}$\nnsp,\oss
and\qss 
$\|\qff K_{\dff n} \qff\|\qff \leq\qff 1$\qss for all $n$\nnsp.\oss

The resolution\qss (\ref{resolution})\qss is said to be\qss
\emph{relatively injective}\pss if all $G$\dnsp-modules $U_n$ are relatively injective,\oss
and\trs \emph{strong}\qss if it admits a\qss
\emph{contracting homotopy}.\oss 
A\qss \emph{split resolution}\pss\vspace*{6pt}
\begin{equation}
\label{resolution-u}
\quad
\begin{tikzcd}[column sep=large, row sep=normal]\dis
0 \arrow[r]
& 
U \arrow[r, shift left=3pt, "\dis d_{\dff -\dff 1}\off"]
& 
U_{0} \arrow[r, shift left=3pt, "\dis d_{\dff 0}\off"] 
\arrow[l, shift left=3pt, "\dis \off K_{\dff 0}"]
&   
U_{1} \arrow[r, shift left=3pt, "\dis d_{\dff 1}\off"] 
\arrow[l, shift left=3pt, "\dis \off K_{\dff 1}"]
&
U_{2} \arrow[r, shift left=3pt, "\dis d_{\dff 2}\off"] 
\arrow[l, shift left=3pt, "\dis \off K_{\dff 2}"]
&
\off \ldots \off,
\arrow[l, shift left=3pt, "\dis \off K_{\dff 3}"]
\end{tikzcd}
\end{equation}

\vspace*{-3pt}
is defined as a strong resolution together with a contracting homotopy.\oss

\vspace*{3pt}
\mypar{Lemma.}{split-strongly-injective}
\emph{If\qss (\ref{resolution-u})\qss is a split resolution,\oss
then\sss $d_{\dff -\dff 1}$\sss is strongly injective,\oss
as also the morphisms}\vspace*{3pt}
\[
\quad
d_{\dff n}'
\qff \colon\qff
U_n\left/\trf\kernel d_{\dff n}
\ttoo
U_{n\dff +\dff 1}\right.\dff.
\]

\vspace*{-9pt}
\emph{induced by morphisms\sss $d_{\dff n}$\sss
for all\oss
$n\off =\off 0\fff,\pff 1\fff,\pff 2\fff,\pff \ldots \off$\dnsp.\oss}

\proof\qss
Since\qss 
$K_{\dff 0}\dff \circ\dff d_{\dff -\dff 1}\off =\off \id$\qss
and\qss
$\|\qff K_{\dff n} \qff\|\qff \leq\qff 1$\nnsp,\oss
the morphism\qss
$d_{\dff -\dff 1}$\qss
is strongly injective.\oss
Let\vspace*{3pt}
\[
\quad
q_{\fff n}\dff \colon\dff U_n\ttoo U_n\left/\trf\kernel d_{\dff n}\right.
\]

\vspace*{-9pt}
be the canonical projection and\dss let\oss 
$\sigma_{\fff n}
\off =\off 
q_{\fff n}\dff \circ\dff K_{\dff n\dff +\dff 1}
\qff \colon\qff
U_{n\dff +\dff 1}
\ttoo
U_n\left/\trf\kernel d_{\dff n}\right. 
$\nnsp.\oss 
The homotopy identity\qss (\ref{homotopy-identity})\qss
implies that\vspace*{3.5pt}
\[
\quad
q_{\fff n}\dff \circ\dff d_{\dff n\dff -\dff 1}\dff \circ\dff K_{\dff n}
\qff +\qff
q_{\fff n}\dff \circ\dff K_{\dff n\dff +\dff 1}\dff \circ\dff d_{\dff n}
\off =\off
q_{\fff n}\dff.
\]

\vspace*{-8.5pt}
The exactness of\qss (\ref{resolution})\qss implies that\qss
$\image d_{\dff n\dff -\dff 1}\off =\off \kernel d_{\dff n}$\qss
and\dss hence\qss
$q_{\fff n}\dff \circ\dff d_{\dff n\dff -\dff 1}
\off =\off
0$\nnsp.\oss
Therefore\oss
$\sigma_{\fff n}\dff \circ\dff d_{\dff n}
\off =\off
q_{\fff n}\dff \circ\dff K_{\dff n\dff +\dff 1}\dff \circ\dff d_{\dff n}
\off =\off
q_{\fff n}
$\nnsp.\oss 
It follows that\qss
$\sigma_{\fff n}\dff \circ\dff d_{\dff n}'\off =\off \id$\nnsp.\oss

Obviously,\pss
$\|\dff q_{\fff n} \dff\|\qff \leq\qff 1$\qss
and\dss therefore\qss\vspace*{3pt}
\[
\quad
\|\dff \sigma_{\fff n} \dff\|
\off =\off
\|\dff q_{\fff n}\dff \circ\dff K_{\dff n\dff +\dff 1}\dff\|
\off \leq\off 
\|\dff q_{\fff n}\dff\|\dff \cdot\dff \|\dff K_{\dff n\dff +\dff 1}\dff\|
\off \leq\off
1\cdot 1
\off =\off
1\dff.
\]

\vspace*{-9pt}
These properties of\dss $\sigma_{\fff n}$\dss imply that\dss $d_{\dff n}'$\dss
is strongly injective.\oss  \eproof

\vspace{6pt}
\mypar{Lemma.}{main-homology-lemma}
\emph{Suppose that\qss (\ref{resolution-u})\qss 
is a split resolution of\pss $U$\qss and}\vspace*{5.5pt}
\begin{equation}
\label{complex-v}
\quad
\begin{tikzcd}[column sep=large, row sep=normal]\dis
0 \arrow[r]
& 
V \arrow[r, "\dis d_{\dff -\dff 1}\off"]
& 
V_{0} \arrow[r, "\dis d_{\dff 0}\off"]
&   
V_{1} \arrow[r, "\dis d_{\dff 1}\off"]
&
V_{2} \arrow[r, "\dis d_{\dff 2}\off"]
&
\off \ldots \off,
\end{tikzcd}
\end{equation}

\vspace*{-3.5pt}
\emph{is a complex of $G$\dnsp-modules,\oss 
i.e.\qss $d_{\dff n\dff +\dff 1}\dff \circ\dff d_{\dff n}\qff =\qff 0$\qss for\qss 
$n\qff \geq\qff -\qff 1$\nnsp.\oss
If\dss all\dss $G$\dnsp-modules\qss $V_n$\nnsp,\oss $n\qff \geq\qff 0$\nnsp,\oss
are relatively injective,\oss
then any\dss $G$\dnsp-morphism\qss 
$u\dff \colon\dff U\toto V$\qss 
can be extended to a $G$\dnsp-morphism from the resolution\qss (\ref{resolution-u})\qss
to the complex\qss (\ref{complex-v}),\oss 
i.e. to a commutative diagram}\vspace*{9pt}
\begin{equation}
\label{morphism}
\quad
\begin{tikzcd}[column sep=large, row sep=huge]\dis
0 \arrow[r]
& 
U \arrow[r, "\dis d_{\dff -\dff 1}\off"]
\arrow[d, "\dis u"']
& 
U_{0} \arrow[r, "\dis d_{\dff 0}\off"] 
\arrow[d, "\dis u_{\dff 0}"']
&   
U_{1} \arrow[r, "\dis d_{\dff 1}\off"] 
\arrow[d, "\dis u_{\dff 1}"']
&
U_{2} \arrow[r, "\dis d_{\dff 2}\off"] 
\arrow[d, "\dis u_{\dff 2}"']
&
\off \ldots \off
\\
0 \arrow[r]
& 
V \arrow[r, "\dis d_{\dff -\dff 1}\off"]
& 
V_{0} \arrow[r, "\dis d_{\dff 0}\off"]
&   
V_{1} \arrow[r, "\dis d_{\dff 1}\off"]
&
V_{2} \arrow[r, "\dis d_{\dff 2}\off"]
&
\off \ldots \off,
\end{tikzcd}
\end{equation}

\vspace*{0pt}
\emph{in which all maps\qss $u_{\dff i}$\nnsp,\qss $i\qff \geq\qff 0$\nnsp,\oss 
are $G$\dnsp-morphisms.\oss}

\vspace*{6pt}
\proof\qss
Since\dss $d_{\dff -\dff 1}$\dss is strongly injective by Lemma\qss \ref{split-strongly-injective}\qss
and the\dss $G$\dnsp-module $V_0$ is relatively injective,\oss 
there exists\dss a\dss $G$\dnsp-morphism\qss
$u_{\dff 0}\dff \colon\dff U_0\ttoo V_0$\qss
such that\qss\vspace*{4pt}
\[
\quad
u_{\dff 0}\dff \circ\dff d_{\dff -\dff 1}
\off =\off 
d_{\dff -\dff 1}\dff \circ\dff u\dff,
\]

\vspace*{-8pt}
i.e.\qss
the leftmost square of\qss (\ref{morphism})\qss is commutative.\oss
Suppose that $G$\dnsp-morphisms\oss
$u_{\dff 0}\fff,\pff u_{\dff 1}\fff,\pff \ldots\fff,\pff u_{\dff n}$\oss
are already constructed and that 
involving them squares of\qss (\ref{morphism})\qss are commutative.\oss
Then\vspace*{4pt}
\[
\quad 
d_{\dff n}\dff \circ\dff u_{\dff n}\dff \circ\dff d_{\dff n\dff -\dff 1}
\off =\off
d_{\dff n}\dff \circ\dff d_{\dff n\dff -\dff 1}\dff \circ\dff u_{\dff n\dff -\dff 1}
\off =\off
0
\]

\vspace{-8pt}
and\dss hence\qss $d_{\dff n}\dff \circ\dff u_{\dff n}$\qss 
is equal to zero on the image\dss $\image d_{\dff n\dff -\dff 1}$\nnsp.\oss
Since\qss
$\image d_{\dff n\dff -\dff 1}\off =\off \kernel d_{\dff n}$\nnsp,\oss
the\dss $G$\dnsp-morphism\qss $d_{\dff n}\dff \circ\dff u_{\dff n}$\qss
induces a\dss $G$\dnsp-morphism\vspace*{3pt}
\[
\quad
u_{\dff n}'
\qff \colon\qff
U_n\left/\trf\kernel d_{\dff n}
\ttoo
V_{n\dff +\dff 1}\right.\dff.
\]

\vspace{-9pt}
Since the morphism\dss $d_{\dff n}'$\dss
from Lemma\qss \ref{split-strongly-injective}\qss is strongly injective
and $V_{n\dff +\dff 1}$ is relatively injective,\oss
there exists a\dss $G$\dnsp-morphism\qss
$u_{\dff n\dff +\dff 1}\dff \colon\dff U_{n\dff +\dff 1}\ttoo V_{n\dff +\dff 1}$\qss
such that\qss
$u_{\dff n}'\off =\off u_{\dff n\dff +\dff 1}\dff \circ\dff d_{\dff n}'$\qss
and\dss hence\vspace*{3pt}
\[
\quad
d_{\dff n}\dff \circ\dff u_{\dff n}
\off =\off
u_{\dff n\dff +\dff 1}\dff \circ\dff d_{\dff n}\qff.
\]

\vspace*{-9pt}
The induction completes the proof\dss of\dss the existence of\dss
$G$\dnsp-morphisms\oss
$u_{\dff 0}\fff,\pff u_{\dff 1}\fff,\pff u_{\dff 2}\fff,\pff \ldots \off$\nnsp.\oss  \eproof

\mypar{Lemma.}{homotopy-uniqueness}
\emph{Under the assumptions of Lemma\qss \ref{main-homology-lemma},\oss
every two extensions\qss $u_{\dff \bullet}$\qss of\qss $u$\qss are chain homotopic
by a chain homotopy consisting of\qss $G$\dnsp-morphisms.\oss}

\proof\qss
It is sufficient to prove that if
in the diagram\qss (\ref{morphism})\qss $u\qff =\qff 0$\nnsp,\oss
then the chain map\dss $u_{\dff \bullet}$\dss is chain homotopic to zero.\oss
Suppose that\qss $u\qff =\qff 0$\nnsp.\oss 
Then\qss
$u_{\dff 0}\dff \circ\dff d_{\dff -\dff 1}
\off =\off
d_{\dff -\dff 1}\dff \circ\dff u
\off =\off
0$\nnsp.\oss
Since\qss
$\image d_{\dff -\dff 1}\off =\off \kernel d_{\dff 0}$\nnsp,\oss
it follows that $u_{\dff 0}$ defines a morphism\vspace*{3pt}
\[
\quad
u_{\dff 0}'
\qff \colon\qff
U_0\left/\trf\kernel d_{\dff 0}
\ttoo
V_{0}\right.\dff.
\]

\vspace{-9pt}
Since\dss $d_{\dff 0}'$\dss is a strongly injective\dss $G$\dnsp-morphism
and $V_0$ is a relatively injective\dss $G$\dnsp-module,\oss
there exists a\dss $G$\dnsp-morphism\qss
$k_{\dff 1}\dff \colon\dff U_1\ttoo V_0$\qss 
such that\qss
$k_{\dff 1}\dff \circ\dff d_{\dff 0}'\off =\off u_{\dff 0}'$\qss
and\dss hence\qss
$k_{\dff 1}\dff \circ\dff d_{\dff 0}\off =\off u_{\dff 0}$\nnsp.\oss
Suppose that for\oss
$m\off =\off 1\fff,\pff 2\fff,\pff \ldots\fff,\pff n$\oss morphisms\qss\vspace*{4pt}
\[
\quad
k_{\dff m}\dff \colon\dff U_m\ttoo V_{m\dff -\dff 1}
\]

\vspace{-8pt}
are already constructed
in such a way that\oss\vspace*{4pt}
\[
\quad
d_{\dff m\dff -\dff 1}\dff \circ\dff k_{\dff m}
\qff +\qff
k_{\dff m\dff +\dff 1}\dff \circ\dff d_{\dff m}
\off =\off
u_{\dff m}
\]

\vspace{-8pt}
if\qss $1\qff \leq\qff m\qff \leq\qff n\qff -\qff 1$\nnsp.\oss
Then\oss
$\dis
k_{\dff n}\dff \circ\dff d_{\dff n\dff -\dff 1}
\off =\off
u_{\dff n\dff -\dff 1}
\qff -\qff
d_{\dff n\dff -\dff 2}\dff \circ\dff k_{\dff n\dff -\dff 1}$\oss
and\dss therefore\vspace*{6pt}
\[
\quad
\bigl(\dff
u_{\dff n}
\qff -\qff
d_{\dff n\dff -\dff 1}\dff \circ\dff k_{\dff n}
\dff\bigr)
\dff \circ\dff
d_{\dff n\dff -\dff 1}
\off =\off
u_{\dff n}
\dff \circ\dff
d_{\dff n\dff -\dff 1}
\qff -\qff
d_{\dff n\dff -\dff 1}\dff \circ\dff k_{\dff n}\dff \circ\dff
d_{\dff n\dff -\dff 1}
\]

\vspace*{-30pt}
\[
\quad
\phantom{\bigl(\dff
u_{\dff n}
\qff -\qff
d_{\dff n\dff -\dff 1}\dff \circ\dff k_{\dff n}
\dff\bigr)
\dff \circ\dff
d_{\dff n\dff -\dff 1}
\off }
=\off
u_{\dff n}
\dff \circ\dff
d_{\dff n\dff -\dff 1}
\qff -\qff
d_{\dff n\dff -\dff 1}\dff \circ\dff
\left(
u_{\dff n\dff -\dff 1}
\qff -\qff 
d_{\dff n\dff -\dff 2}\dff \circ\dff k_{\dff n\dff -\dff 1}
\dff\right)
\]

\vspace*{-30pt}
\[
\quad
\phantom{\bigl(\dff
u_{\dff n}
\qff -\qff
d_{\dff n\dff -\dff 1}\dff \circ\dff k_{\dff n}
\dff\bigr)
\dff \circ\dff
d_{\dff n\dff -\dff 1}
\off }
=\off
\bigl(\dff
u_{\dff n}
\dff \circ\dff
d_{\dff n\dff -\dff 1}
\qff -\qff
d_{\dff n\dff -\dff 1}
\dff \circ\dff 
u_{\dff n\dff -\dff 1}
\bigr)\dff
\qff +\qff
d_{\dff n\dff -\dff 1}
\dff \circ\dff
d_{\dff n\dff -\dff 2}\dff \circ\dff k_{\dff n\dff -\dff 1}
\]

\vspace*{-30pt}
\[
\quad
\phantom{\bigl(\dff
u_{\dff n}
\qff -\qff
d_{\dff n\dff -\dff 1}\dff \circ\dff k_{\dff n}
\dff\bigr)
\dff \circ\dff
d_{\dff n\dff -\dff 1}
\off }
=\off
0\qff +\qff 0
\off =\off
0\dff.
\]

\vspace{-6pt}
It follows that the\dss $G$\dnsp-morphism\qss
$u_{\dff n}
\qff -\qff
d_{\dff n\dff -\dff 1}\dff \circ\dff k_{\dff n}$\qss
defines a\dss $G$\dnsp-morphism\vspace*{6pt}
\[
\quad
u_{\dff n}'
\qff \colon\qff
U_n\left/\trf\kernel d_{\dff n}
\ttoo
V_{n}\right.\dff.
\]

\vspace{-6pt}
Since\dss $d_{\dff n}'$\dss is strongly injective
and $V_n$ is a relatively injective,\oss
there exists a\dss $G$\dnsp-morphism\qss\vspace*{4pt}
\[
\quad
k_{\dff n\dff +\dff 1}\dff \colon\dff U_{n\dff +\dff 1}\ttoo V_n
\]

\vspace*{-8pt} 
such that\qss
$k_{\dff n\dff +\dff 1}\dff \circ\dff d_{\dff n}'\off =\off u_{\dff n}'$\nnsp.\oss
Then\vspace*{4pt}
\[
\quad
k_{\dff n\dff +\dff 1}\dff \circ\dff d_{\dff n}
\off =\off
u_{\dff n}
\qff -\qff
d_{\dff n\dff -\dff 1}\dff \circ\dff k_{\dff n}
\]

\vspace{-8pt}
and therefore\qss
$d_{\dff n\dff -\dff 1}\dff \circ\dff k_{\dff n}
\qff +\qff
k_{\dff n\dff +\dff 1}\dff \circ\dff d_{\dff n}
\off =\off
u_{\dff n}$\nnsp.\oss
An induction completes the proof\dss of\dss the existence of
a chain homotopy between\dss $u_{\dff \bullet}$\dss and zero.\oss  \eproof

\vspace*{6pt}
\myuppar{The norm of morphisms in an extension.}
It would be nice if the extension\dss $u_{\dff \bullet}$\dss
in Lemma\qss \ref{main-homology-lemma}\qss
could be chosen in such a way that\qss
$\|\qff u_{\dff i} \qff\|\qff \leq\qff \|\qff u \qff\|$\qss
for all\qss $i\qff \geq\qff 0$\nnsp.\oss
It is instructive to see what estimate of the norms is implicit in the proof\dss
of\dss Lemma\qss \ref{main-homology-lemma}.\oss

\vspace*{5.25pt}
For every vector\dss $x\qff \in\qff U_n\left/\trf\kernel d_{\dff n}\right.$\dss
and\dss every\qss $\varepsilon\qff >\qff 0$\qss there exists a vector\qss
$y\qff \in\qff U_n$\qss such that\qss
$q_{\fff n}\fff(\fff y \fff)\qff =\qff x$\qss
and\qss
$\|\qff y \qff\|\qff \leq\qff \|\qff x \qff\|\qff +\qff \varepsilon$\nnsp.\oss
It follows that\qss
$u_{\dff n}'\fff(\fff x \fff)
\off =\off
d_{\dff n}\dff \circ\dff u_{\dff n}\fff(\fff y \fff)$\qss 
and\dss hence\vspace*{4pt}
\[
\quad
\|\qff u_{\dff n}'\fff(\fff x \fff) \qff\|
\off =\off
\|\qff d_{\dff n}\dff \circ\dff u_{\dff n}\fff(\fff y \fff) \qff\|
\off \leq\off
\|\qff d_{\dff n}\dff \circ\dff u_{\dff n} \qff\|\dff \cdot\dff \|\qff y \qff\|
\]

\vspace*{-33pt}
\[
\quad
\phantom{\|\qff u_{\dff n}'\fff(\fff x \fff) \qff\|
\off =\off
\|\qff d_{\dff n}\dff \circ\dff u_{\dff n}\fff(\fff y \fff) \qff\|
\off }
\off \leq\off
\|\qff d_{\dff n}\dff \circ\dff u_{\dff n} \qff\|\dff \cdot\dff \|\qff x \qff\|
\qff +\qff
\|\qff d_{\dff n}\dff \circ\dff u_{\dff n} \qff\|\dff \cdot\dff \varepsilon\dff.
\]

\vspace{-5pt}
Since\qss $\varepsilon\qff >\qff 0$\qss is arbitrary,\pss 
$\|\qff u_{\dff n}'\fff(\fff x \fff) \qff\|
\qff \leq\qff
\|\qff  d_{\dff n}\dff \circ\dff u_{\dff n} \qff\|\dff \cdot\dff \|\qff x \qff\|$\nnsp.\oss
Therefore\qss
$\|\qff u_{\dff n}' \qff\|
\qff \leq\qff
\|\qff  d_{\dff n}\dff \circ\dff u_{\dff n} \qff\|$\nnsp.\oss
Since the morphism\dss $d_{\dff n}'$\dss is strongly injective
and the\dss $G$\dnsp-module $V_{n\dff +\dff 1}$ is relatively injective,\oss
there exists a\dss $G$\dnsp-morphism\qss
$u_{\dff n\dff +\dff 1}\dff \colon\dff U_{n\dff +\dff 1}\ttoo V_{n\dff +\dff 1}$\qss
such that\qss\vspace*{6pt}
\[
\quad 
\|\qff u_{\dff n\dff +\dff 1} \qff\|
\qff \leq\qff
\|\qff u_{\dff n}' \qff\|
\qff \leq\qff
\|\qff d_{\dff n}\dff \circ\dff  u_{\dff n} \qff\|
\qff \leq\qff
\|\qff d_{\dff n} \qff\|
\dff \cdot\dff
\|\qff u_{\dff n} \qff\|
\]

\vspace{-6pt}
and\oss
$u_{\dff n}'\off =\off u_{\dff n\dff +\dff 1}\dff \circ\dff d_{\dff n}'$\nnsp.\off\oss
It\dss follows\dss that\qss\vspace*{3pt}
\[
\quad
\|\qff u_{\dff n\dff +\dff 1} \qff\|
\qff \leq\qff
\|\qff u \qff\|\dff
\prod_{i\qff =\qff -\qff 1}^n\qff \|\qff d_{\dff n} \qff\|\qff.
\]

\vspace{-9pt}
Suppose that\qss
$\|\qff d_{\dff n} \qff\|\qff \leq\qff n\qff +\qff 2$\nnsp,\oss
as is the case for the so-called standard resolution of\dss $V$\dss
and for resolutions defined\dss by topological spaces.\oss
Then the above proof\dss 
leads to the estimate\qss\vspace*{6pt}
\begin{equation*}
\quad
\|\qff u_{\dff n} \qff\|\qff \leq\qff (n\qff +\qff 1)\dff!\dff \cdot\dff \|\qff u \qff\|\qff,
\end{equation*}

\vspace*{-6pt}
implicitly contained in the work of\qss R.\dss Brooks\qss \cite{br}.\oss
But by Theorem\qss \ref{comparing-to-standard}\qss below,\oss
for the standard resolution there is indeed 
an extension\dss 
$u_{\dff \bullet}$\dss
such that\qss
$\|\qff u_{\dff i} \qff\|\qff \leq\qff \|\qff u \qff\|$\qss
for all\qss $i\qff \geq\qff 0$\nnsp.\oss

\myuppar{The\dss $G$\dnsp-modules\dss $B\dff(\fff G^{\fff n\dff +\dff 1},\pff U\dff)$\dnsp.}
For a bounded left\dss $G$\dnsp-module $U$ and\qss $n\qff \geq\qff 0$\qss 
let\dss $B\dff(\fff G^{\fff n\dff +\dff 1},\pff U\dff)$\dss 
be the space of\dss functions\qss 
$f\dff \colon\dff G^{\fff n\dff +\dff 1}\toto U$\qss  
such that\vspace*{4pt} 
\[
\quad
\|\dff f \dff\|
\off =\off 
\sup\qff \left\{\off \|\dff f\dff(\fff g_{\fff 0}\fff,\pff g_{\fff 1}\fff,\pff \ldots\fff,\pff g_{\fff n}\fff) \trf\| 
\off  \left| \off 
(\fff g_{\fff 0}\fff,\pff g_{\fff 1}\fff,\pff \ldots\fff,\pff g_{\fff n}\fff)\qff \in\qff G^{\fff n\dff +\dff 1} \off\right\} \right.
\off \leq\off \infty\qff. 
\]

\vspace*{-8pt}
\dnsp$B\dff(\fff G^{\fff n\dff +\dff 1},\pff U \dff)$\dss
is a Banach space with the norm $\|\dff f \dff\|$\nnsp,\oss 
and the action\dss $\bullet$\dss of $G$ defined by\vspace*{4pt} 
\[
\quad
(\fff h\bullet\nsp f \dff)\dff(\fff g_{\fff 0}\fff,\pff g_{\fff 1}\fff,\pff \ldots\fff,\pff g_{\fff n}\fff)
\off =\off
h\cdot\left(\dff
f\dff(\fff g_{\fff 0}\fff,\pff g_{\fff 1}\fff,\pff \ldots\fff,\pff g_{\fff n}\halfff h\fff)
\dff\right)
\]

\vspace{-8pt} 
turns it into a bounded $G$\dnsp-module.\oss 
Let the\dss $G$\dnsp-module\dss $b\dff(\fff G^{\fff n\dff +\dff 1},\pff U \dff)$\dss be equal to\dss 
$B\dff(\fff G^{\fff n\dff +\dff 1},\pff U \dff)$\dss as a Banach space,\oss
but with the action\dss $*$ of\dss $G$\dss defined by\vspace*{3pt} 
\[
\quad
(\fff h *\nsp f \dff)\dff(\fff g_{\fff 0}\fff,\pff g_{\fff 1}\fff,\pff \ldots\fff,\pff g_{\fff n}\fff)
\off =\off
h\cdot\left(\dff
f\dff(\fff g_{\fff 0}\fff,\pff g_{\fff 1}\fff,\pff \ldots\fff,\pff g_{\fff n}\fff)\dff.
\dff\right)
\]

\vspace*{-9pt} 
Obviously,\pss 
$B\dff(\fff G^{\fff n\dff +\dff 1},\pff U \dff)
\off =\off
B\dff(\fff G\fff,\pff b\dff(\fff G^{\fff n},\pff U \dff) \dff)\
$\qss 
and\dss hence\dss $B\dff(\fff G^{\fff n\dff +\dff 1},\pff U \dff)$\dss 
is a relatively injective\dss $G$\dnsp-module by Lemma\qss \ref{bgv}.\oss
If $U$ is a trivial\dss $G$\dnsp-module,\oss
then\dss $b\dff(\fff G^{\fff n\dff +\dff 1},\pff U \dff)$\dss is also trivial.

\myuppar{The standard resolution.}
For a bounded left\dss $G$\dnsp-module $U$ let
us consider the sequence\vspace*{5pt}
\begin{equation}
\label{st-res}
\quad
\begin{tikzcd}[column sep=large, row sep=normal]\dis
0 \arrow[r]
& 
U \arrow[r, "\dis d_{\dff -\dff 1}\off"]
& 
B\dff(\fff G\fff,\pff U \dff) \arrow[r, "\dis d_{\dff 0}\off"]
&   
B\dff(\fff G^{\dff 2},\pff U \dff) \arrow[r, "\dis d_{\dff 1}\off"]
&
B\dff(\fff G^{\dff 3},\pff U \dff) \arrow[r, "\dis d_{\dff 2}\off"]
&
\off \ldots \off,
\end{tikzcd}
\end{equation}

\vspace*{-3pt}
where\oss  
$d_{\dff -\dff 1}\dff(\fff v \fff)(\fff g\fff)\off =\off v$\oss 
for\dss all\qss
$v\qff \in\qff U$\dnsp,\qss $g\qff \in\qff G$\nnsp,\off\oss
and\dss\vspace*{6pt}
\[
\quad
d_{\dff n}\dff(\fff f \fff)
(\fff g_{\fff 0}\fff,\pff g_{\fff 1}\fff,\pff \ldots\fff,\pff g_{\fff n\dff +\dff 1}\fff)
\off =\off
(\dff -\qff 1\dff)^{n\dff +\dff 1}\fff
f\dff(\fff g_{\fff 1}\fff,\pff g_{\fff 2}\fff,\pff \ldots\fff,\pff g_{\fff n\dff +\dff 1}\fff)
\]

\vspace*{-30pt}
\[
\quad
\phantom{d_{\dff n}\dff(\fff f \fff)
(\fff g_{\fff 0}\fff,\pff g_{\fff 1}\fff,\pff \ldots\fff,\pff g_{\fff n\dff +\dff 1}\fff)
\off =\off}
+\off
\sum_{i\qff =\qff 0}^{n}\qff (-\qff 1)^{\dff n\dff -\dff i}\dff 
f\dff(\dff g_{\fff 0}\fff,\pff \ldots\fff,\pff
g_{\dff i}\fff g_{\fff i\dff +\dff 1}\dff,\pff \ldots\fff,\pff   g_{\dff n\dff +\dff 1} \dff)\dff,
\]

\vspace*{-6pt}
for all\qss $n\qff \geq\qff 0$\qss
and\oss $g_{\fff 0}\fff,\pff g_{\fff 1}\fff,\pff \ldots\fff,\pff g_{\fff n\dff +\dff 1}
\qff \in\qff G$\nnsp.\off\oss
Obviously,\oss the maps $d_{\dff n}$ are $G$\dnsp-morphisms,\oss 
and a standard calculation shows that\oss 
$\dis
d_{\dff n\dff +\dff 1}\dff \circ\trf  d_{\dff n}\off =\off 0$\oss 
for all\qss $n\qff \geq\qff -\qff 1$\nnsp,\oss
i.e.\qss that\qss (\ref{st-res})\qss is a complex\halfff.\oss
See Appendix\qss \ref{categories-classifying}\qss for a non-calculational proof\halfff.\oss
Let us consider also the sequence\vspace*{6pt}
\begin{equation}
\label{st-homotopy}
\quad
\begin{tikzcd}[column sep=large, row sep=normal]\dis
U 
& 
B\dff(\fff G\fff,\pff U \dff) \arrow[l, "\dis \off K_{\dff 0}"]
&   
B\dff(\fff G^{\dff 2}\fff,\pff U \dff) \arrow[l, "\dis \off K_{\dff 1}"]
&
B\dff(\fff G^{\dff 3}\fff,\pff U \dff) \arrow[l, "\dis \off K_{\dff 2}"]
&
\off \ldots \off\off, \arrow[l]
\end{tikzcd}
\end{equation}

\vspace*{-3pt}
where\oss 
$
K_{\dff n}\dff(\fff f \fff)(\fff g_{\dff 0}\fff,\pff \ldots\fff,\pff g_{\dff n\dff -\dff 1} \fff)
\off =\off
f\dff(\fff  g_{\dff 0}\fff,\pff \ldots\fff,\pff g_{\dff n\dff -\dff 1}\fff,\pff 1 \fff)$\oss
for\dss all\oss
$g_{\fff 0}\fff,\pff g_{\fff 1}\fff,\pff \ldots\fff,\pff g_{\fff n\dff -\dff 1}
\qff \in\qff G$\nnsp.\off\oss
Obviously,\pss 
$\|\qff K_{\dff n} \qff\|\qff \leq\qff 1$\qss
for all\qss $n\qff \geq\qff 0$\nnsp.\oss
A much easier standard calculation shows that\qss (\ref{st-homotopy})\qss 
is a contracting homotopy for\qss (\ref{st-res}).\oss 
It follows that the sequence\qss (\ref{st-res})\qss is exact\halfff.\oss 
Therefore the chain complex\qss (\ref{st-res})\qss 
together with the contracting homotopy\qss (\ref{st-homotopy})\qss 
is a strong resolution of\dss $U$\nnsp.\oss
This resolution is called the\qss \emph{standard resolution}.\oss 
It is relatively injective by Lemma\qss \ref{bgv}.\oss

\myuppar{Comparing a split resolution with the standard one.}
Suppose that\vspace*{6pt}
\begin{equation*}
\quad
\begin{tikzcd}[column sep=large, row sep=normal]\dis
0 \arrow[r]
& 
U \arrow[r, shift left=3pt, "\dis d_{\dff -\dff 1}\off"]
& 
U_{0} \arrow[r, shift left=3pt, "\dis d_{\dff 0}\off"] 
\arrow[l, shift left=3pt, "\dis \off K_{\dff 0}"]
&   
U_{1} \arrow[r, shift left=3pt, "\dis d_{\dff 1}\off"] 
\arrow[l, shift left=3pt, "\dis \off K_{\dff 1}"]
&
U_{2} \arrow[r, shift left=3pt, "\dis d_{\dff 2}\off"] 
\arrow[l, shift left=3pt, "\dis \off K_{\dff 2}"]
&
\off \ldots \off,
\arrow[l, shift left=3pt, "\dis \off K_{\dff 3}"]
\end{tikzcd}
\end{equation*}

\vspace*{-3pt}
is a split resolution of\pss $U$\dnsp.\oss
Let us define the maps\vspace*{3pt}
\[
\quad
k_{\dff n}\fff,\pff u_{\dff n}
\qff \colon\qff
U_n
\qff \ttoo\qff
B\dff(\fff G^{\dff n\dff +\dff 1}\fff,\pff U \dff)\dff,
\]

\vspace*{-9pt}
where\qss $n\qff \geq\qff 0$\nnsp,\oss
by the formulas\vspace*{3pt}
\[
\quad
k_{\dff n}\dff(\dff f \dff)(\dff g_{\dff 0}\fff,\pff \ldots\fff,\pff g_{\dff n} \dff)
\off =\off
K_{\dff 0}\dff(\dff g_{\dff 0}\cdot\dff K_{\dff 1}\dff(\off \ldots\off K_{\dff n\dff -\dff 1}\dff
(\dff g_{\dff n\dff -\dff 1}\cdot\dff K_{\dff n}\dff(\dff g_{\dff n}\cdot f \dff) \dff)\off \ldots \off)\dff)\dff,
\]

\vspace{-33pt}
\begin{equation*}
\quad
u_{\dff n}\dff(\dff f \dff)(\dff g_{\dff 0}\fff,\pff \ldots\fff,\pff g_{\dff n} \dff)
\off =\off
\left(\dff g_{\dff 0}\fff g_{\dff  1}\fff \ldots\fff g_{\fff n}\dff\right)^{\dff -\dff 1}
\cdot\dff
\Bigl(\dff
k_{\dff n}\dff(\dff f \dff)(\dff g_{\dff 0}\fff,\pff \ldots\fff,\pff g_{\dff n} \dff)
\dff\Bigr)\dff.
\end{equation*}

\vspace{-6pt}
An equivalent way to define\dss $u_{\dff n}$\dss 
is to set\qss $u_{\dff -\dff 1}\qff =\qff \id_{\dff U}$\qss
and use the recursive relation\vspace*{3pt}
\begin{equation}
\label{u-recursion-twisted}
\quad
u_{\dff m}\dff(\dff f \dff)(\dff g_{\dff 0}\fff,\pff \ldots\fff,\pff g_{\dff m} \dff)
\off =\off
g_{\dff m}^{\dff -\dff 1}\cdot\dff
\Bigl(\dff
u_{\dff m\dff -\dff 1}
\left(\dff 
K_{\dff m}\dff(\dff g_{\dff m}\cdot f \dff) 
\dff\right)
(\dff g_{\dff 0}\fff,\pff \ldots\fff,\pff g_{\dff m\dff -\dff 1} \dff)
\dff\Bigr)\dff
\end{equation}

\vspace*{-9pt}
in order to define\dss $u_{\dff m}$\qss for\qss $m\qff \geq\qff 0$\nnsp.\oss

\vspace*{0pt}
\mypar{Lemma.}{formula-is-morhphism}
\emph{The maps\qss $u_{\dff n}$\qss are\dss $G$\dnsp-morphisms.}

\proof\qss
Obviosly,\oss
$k_{\dff n}\dff(\dff h\cdot\nsp f \dff)(\dff g_{\dff 0}\fff,\pff \ldots\fff,\pff g_{\dff n} \dff)
\off =\off
k_{\dff n}\dff(\dff f \dff)(\dff g_{\dff 0}\fff,\pff \ldots\fff,\pff g_{\dff n}\fff h \dff)$\dnsp.\off\oss
Therefore\qss\vspace*{3pt} 
\[
\quad
u_{\dff n}\dff(\dff h\cdot f \dff)(\dff g_{\dff 0}\fff,\pff \ldots\fff,\pff g_{\dff n} \dff)
\phantom{\Bigl(\dff\Bigr)}
\]

\vspace*{-33pt}
\[
\quad
=\off
\left(\dff g_{\dff 0}\fff g_{\dff  1}\fff \ldots\fff g_{\fff n}\dff\right)^{\dff -\dff 1}
\cdot\dff
k_{\dff n}\dff(\dff h\cdot f \dff)(\dff g_{\dff 0}\fff,\pff \ldots\fff,\pff g_{\dff n} \dff)
\]

\vspace*{-30pt}
\[
\quad
=\off
\left(\dff g_{\dff 0}\fff g_{\dff  1}\fff \ldots\fff g_{\fff n}\dff\right)^{\dff -\dff 1}
\cdot\dff
k_{\dff n}\dff(\dff f \dff)(\dff g_{\dff 0}\fff,\pff \ldots\fff,\pff g_{\dff n}\fff h \dff)
\]

\vspace*{-30pt}
\[
\quad
=\off
\Bigl(\dff
h
\cdot
\left(\dff g_{\dff 0}\fff g_{\dff  1}\fff \ldots\fff g_{\fff n}\fff h \dff\right)^{\dff -\dff 1}
\dff\Bigr)
\cdot\dff
k_{\dff n}\dff(\dff f \dff)(\dff g_{\dff 0}\fff,\pff \ldots\fff,\pff g_{\dff n}\fff h \dff)
\]

\vspace*{-30pt}
\[
\quad
=\off
h
\dff \cdot\dff
\Bigl(\dff
\left(\dff g_{\dff 0}\fff g_{\dff  1}\fff \ldots\fff g_{\fff n}\fff h \dff\right)^{\dff -\dff 1}
\cdot\dff
k_{\dff n}\dff(\dff f \dff)(\dff g_{\dff 0}\fff,\pff \ldots\fff,\pff g_{\dff n}\fff h \dff)
\dff\Bigr)
\]

\vspace*{-30pt}
\[
\quad
=\off
h
\dff \cdot\dff
\left(\dff
u_{\dff n}\dff(\dff f \dff)(\dff g_{\dff 0}\fff,\pff \ldots\fff,\pff g_{\dff n}\fff h \dff)
\dff\right)
\off =\off
\left(\dff
h
\dff\bullet\dff
u_{\dff n}\dff(\dff f \dff)
\dff\right)
(\dff g_{\dff 0}\fff,\pff \ldots\fff,\pff g_{\dff n} \dff)\dff.
\]

\vspace{-6pt}
It follows that\oss 
$u_{\dff n}\dff(\dff h\cdot\nsp f \dff)
\off =\off
h
\dff \bullet\dff
u_{\dff n}\dff(\dff f \dff)$\oss and\dss hence\dss
$u_{\dff n}$\dss is a $G$\dnsp-morphism.\oss  \eproof

\vspace{3pt}
\mypar{Theorem.}{comparing-to-standard}
\emph{The sequence of maps\qss 
$u_{\dff \bullet}
\off =\off
\{\dff u_{\dff n} \dff\}$\qss 
is a morphism of resolutions}\vspace*{9pt}
\begin{equation*}
\quad
\begin{tikzcd}[column sep=large, row sep=huge]\dis
0 \arrow[r]
& 
U \arrow[r]
\arrow[d, "\dis u_{\dff -\dff 1}"']
& 
U_{0} \arrow[r] 
\arrow[d, "\dis u_{\dff 0}"']
&   
U_{1} \arrow[r] 
\arrow[d, "\dis u_{\dff 1}"']
&
U_{2} \arrow[r] 
\arrow[d, "\dis u_{\dff 2}"']
&
\off \ldots \off
\\
0 \arrow[r]
& 
U \arrow[r]
& 
B\dff(\fff G\fff,\pff U \dff) \arrow[r]
&   
B\dff(\fff G^{\dff 2}\fff,\pff U \dff) \arrow[r]
&
B\dff(\fff G^{\dff 3}\fff,\pff U \dff) \arrow[r]
&
\off \ldots \off,
\end{tikzcd}
\end{equation*}

\vspace*{0pt}
\emph{extending\qss $u_{\dff -\dff 1}\qff =\qff \id_{\trf U}$\qss
and such that\qss $\|\qff u_{\dff n} \qff\|\qff \leq\qff 1$\qss
for all\qss $n\qff \geq\qff 0$\nnsp.\oss}

\vspace{3pt}
\proof\qss
Since\qss $\|\qff K_{\dff m} \qff\|\qff \leq\qff 1$\qss and\qss 
$\|\qff g\cdot f \qff\|\qff \leq \qff \|\qff f \qff\|$\qss 
for all\qss 
$g\qff \in\qff G$\nnsp,\qss
$f\qff \in\qff U_{\qff m}$\nnsp,\qss
$m\qff \geq\qff 0$\nnsp,\oss\vspace*{3pt}
\[
\quad
u_{\dff n}\dff(\dff f \dff)(\dff g_{\dff 0}\fff,\pff \ldots\fff,\pff g_{\dff n} \dff)
\off \leq\off
\|\qff f \qff\|
\]

\vspace*{-9pt}
for every\oss
$g_{\dff 0}\fff,\pff \ldots\fff,\pff g_{\dff n}\qff \in\qff G$\oss
and\qss $f\qff \in\qff U_{\dff n}$\nnsp,\oss
and\dss hence\oss
$
\|\dff u_{\dff n}\dff(\dff f \dff) \qff\|
\off \leq\off
\|\qff f \qff\|
$\oss
for every\qss $f\qff \in\qff U_{\dff n}$\nnsp.\oss
It follows that\qss
$\|\dff u_{\dff n} \qff\|\qff \leq\qff 1$\qss
for every $n$\nnsp.\oss

It remains to check that\oss
$d_{\dff n}\dff \circ\dff u_{\dff n}
\off =\off
u_{\dff n\dff +\dff 1}\dff \circ\dff d_{\dff n}$\oss
for all\qss $n\qff \geq\qff -\qff 1$\nnsp.\oss
We will prove this using induction by $n$\nnsp.\oss
Since\dss $d_{\dff -\dff 1}$\dss is a\dss $G$\dnsp-morphism and\oss
$K_{\dff 0}\dff \circ\dff d_{\dff -\dff 1}\off =\off \id$\nnsp,\oss\vspace*{6pt} 
\[
\quad
u_{\dff 0}\dff(\dff d_{\dff -\dff 1}\dff(\fff v \fff) \dff)(\dff g \dff)
\off =\off
g^{\fff -\dff 1}\cdot\dff
K_{\dff 0}\dff(\dff g\dff \cdot\dff d_{\dff -\dff 1}\dff (\fff v \fff) \dff)
\]

\vspace*{-30pt}
\[
\quad
\phantom{u_{\dff 0}\dff(\dff d_{\dff -\dff 1}\dff(\fff v \fff) \dff)(\dff g \dff)
\off }
=\off
g^{\fff -\dff 1}\cdot\dff
K_{\dff 0}\dff \left(\dff d_{\dff -\dff 1}\dff(\fff g\cdot v \fff)\dff\right)
\off =\off
g^{\fff -\dff 1}\cdot\dff
g\dff \cdot\dff v
\off =\off
v\dff.
\]

\vspace{-6pt}
At the same time\oss $u_{\dff -\dff 1}\off =\off \id_{\dff U}$\oss and\dss hence\oss\vspace*{3pt}
\[
\quad
d_{\dff -\dff 1}\dff(\dff u_{\dff -\dff 1}\dff (\fff v \fff) \dff)(\dff g \dff)
\off =\off
d_{\dff -\dff 1}\dff (\fff v \fff)(\dff g \dff)
\off =\off
v 
\]

\vspace{-9pt}
It\dss follows\dss that\oss
$d_{\dff -\dff 1}\dff \circ\dff u_{\dff -\dff 1}
\off =\off
u_{\dff 0}\dff \circ\dff d_{\dff -\dff 1}$\nnsp.\off\oss
Suppose now that\oss
$d_{\dff n\dff -\dff 1}\dff \circ\dff  u_{\dff n\dff -\dff 1}
\off =\off
u_{\dff n}\dff \circ\dff  d_{\dff n\dff -\dff 1}$\oss
and prove that\oss 
$d_{\dff n}\dff \circ\dff  u_{\dff n}
\off =\off
u_{\dff n\dff +\dff 1}\dff \circ\dff  d_{\dff n}$\nnsp.\oss
By the relation\qss (\ref{u-recursion-twisted})\qss 
with\qss $m\qff =\qff n\qff +\qff 1$\qss we have\vspace*{6pt}
\[
\quad
u_{\dff n\dff +\dff 1}\dff(\dff d_{\dff n}\dff(\dff f \dff) \dff)
(\dff g_{\dff 0}\fff,\pff \ldots\fff,\pff g_{\dff n\dff +\dff 1} \dff)
\]

\vspace*{-27pt}
\[
\quad
=\off
g_{\dff n\dff +\dff 1}^{\dff -\dff 1}\cdot\dff
\Bigl(\dff
u_{\dff n}
\left(\dff 
K_{\dff n\dff +\dff 1}\dff(\dff g_{\dff n\dff +\dff 1}\cdot d_{\dff n}\dff(\dff f \dff) \dff) 
\dff\right)
(\dff g_{\dff 0}\fff,\pff \ldots\fff,\pff g_{\dff n} \dff)
\dff\Bigr)
\] 

\vspace*{-27pt}
\[
\quad
=\off
g_{\dff n\dff +\dff 1}^{\dff -\dff 1}\cdot\dff
\Bigl(\dff
u_{\dff n}
\left(\dff 
K_{\dff n\dff +\dff 1}\dff(\dff d_{\dff n}\dff(\dff g_{\dff n\dff +\dff 1}\cdot f \dff) \dff) 
\dff\right)
(\dff g_{\dff 0}\fff,\pff \ldots\fff,\pff g_{\dff n} \dff)
\dff\Bigr)
\]

\vspace*{-27pt}
\[
\quad
=\off
g_{\dff n\dff +\dff 1}^{\dff -\dff 1}\cdot\dff
\Bigl(\dff
u_{\dff n}
\left(\dff 
K_{\dff n\dff +\dff 1}\dff \circ\dff d_{\dff n}\dff(\dff g_{\dff n\dff +\dff 1}\cdot f \dff) 
\dff\right)
(\dff g_{\dff 0}\fff,\pff \ldots\fff,\pff g_{\dff n} \dff)
\dff\Bigr)\dff.
\]

\vspace{-3pt}
The homotopy identity\qss (\ref{homotopy-identity})\qss implies that\vspace*{6pt}
\[
\quad
u_{\dff n}
\left(\dff 
K_{\dff n\dff +\dff 1}\dff \circ\dff d_{\dff n}\dff(\dff g_{\dff n\dff +\dff 1}\cdot f \dff) 
\dff\right)
(\dff g_{\dff 0}\fff,\pff \ldots\fff,\pff g_{\dff n} \dff)\dff.
\]

\vspace*{-27pt}
\[
\quad
=\off
u_{\dff n}
\left(\dff 
g_{\dff n\dff +\dff 1}\cdot f
\off -\off
d_{\dff n\dff -\dff 1}\dff \circ\dff K_{\dff n}\dff(\dff g_{\dff n\dff +\dff 1}\cdot f \dff) 
\dff\right)
(\dff g_{\dff 0}\fff,\pff \ldots\fff,\pff g_{\dff n} \dff)
\]

\vspace*{-27pt}
\[
\quad
=\off
u_{\dff n}
\left(\dff 
g_{\dff n\dff +\dff 1}\cdot f
\dff\right)
(\dff g_{\dff 0}\fff,\pff \ldots\fff,\pff g_{\dff n} \dff)
\off -\off
u_{\dff n}
\left(\dff 
d_{\dff n\dff -\dff 1}\dff \circ\dff K_{\dff n}\dff(\dff g_{\dff n\dff +\dff 1}\cdot f \dff) 
\dff\right)
(\dff g_{\dff 0}\fff,\pff \ldots\fff,\pff g_{\dff n} \dff)
\]

\vspace*{-27pt}
\[
\quad
=\off
u_{\dff n}
\left(\dff 
g_{\dff n\dff +\dff 1}\cdot f
\dff\right)
(\dff g_{\dff 0}\fff,\pff \ldots\fff,\pff g_{\dff n} \dff)
\off -\off
u_{\dff n}
\dff \circ\dff 
d_{\dff n\dff -\dff 1}\dff 
\left(\dff 
K_{\dff n}\dff(\dff g_{\dff n\dff +\dff 1}\cdot f \dff) 
\dff\right)
(\dff g_{\dff 0}\fff,\pff \ldots\fff,\pff g_{\dff n} \dff)
\]

\vspace*{-27pt}
\[
\quad
=\off
u_{\dff n}
\left(\dff 
g_{\dff n\dff +\dff 1}\cdot f
\dff\right)
(\dff g_{\dff 0}\fff,\pff \ldots\fff,\pff g_{\dff n} \dff)
\off -\off
d_{\dff n\dff -\dff 1}
\dff
\left(\dff  
u_{\dff n\dff -\dff 1}\dff 
(\dff 
K_{\dff n}\dff(\dff g_{\dff n\dff +\dff 1}\cdot f \dff) 
\dff)
\dff\right)
(\dff g_{\dff 0}\fff,\pff \ldots\fff,\pff g_{\dff n} \dff)\dff,
\]

\vspace*{-3pt}
where at the last step we used the 
inductive assumption\oss
$u_{\dff n}\dff \circ\dff  d_{\dff n\dff -\dff 1}
\off =\off
d_{\dff n\dff -\dff 1}\dff \circ\dff  u_{\dff n\dff -\dff 1}$\nnsp.\oss 
The definition of\qss $d_{\dff n\dff -\dff 1}$\qss
together with the recursive relation\qss (\ref{u-recursion-twisted})\qss for\qss $m\qff =\qff n$\qss
imply that the last dispayed expression is equal to\oss\vspace*{6pt}
\[
\quad
u_{\dff n}
\left(\dff 
g_{\dff n\dff +\dff 1}\cdot f
\dff\right)
(\dff g_{\dff 0}\fff,\pff \ldots\fff,\pff g_{\dff n} \dff)
\off\off -\off\off 
(\dff -\qff 1\dff)^{n}\dff 
u_{\dff n\dff -\dff 1}\dff (\dff K_{\dff n}\dff (g_{\dff n\dff +\dff 1}\cdot f \dff)\dff)
(\dff g_{\dff 0}\fff,\pff \ldots\fff,\pff g_{\dff n\dff -\dff 1} \dff)
\]

\vspace*{-27pt}
\[
\quad\phantom{=\off }
\hspace*{7.7em}
-\off\off
\sum_{i\qff =\qff 0}^{n\qff -\qff 1}(\dff -\dff 1)^{\dff n\dff -\dff 1\dff -\dff i}\dff 
u_{\dff n\dff -\dff 1}\dff (\dff K_{\dff n}\dff (\dff g_{\dff n\dff +\dff 1}\cdot f \dff) \dff)
(\dff g_{\dff 0}\fff,\pff \ldots\fff,\pff 
g_{\dff i}\fff g_{\dff i\dff +\dff 1}\fff,\pff \ldots\fff,\pff g_{\dff n} \dff)
\]

\vspace*{-27pt}
\[
\quad
=\off
u_{\dff n}
\left(\dff 
g_{\dff n\dff +\dff 1}\cdot f
\dff\right)
(\dff g_{\dff 0}\fff,\pff \ldots\fff,\pff g_{\dff n} \dff)
\off\off -\off\off
(\dff -\qff 1\dff)^{n}\dff 
g_{\dff n\dff +\dff 1}\dff\cdot\dff\Bigl(\dff 
u_{\dff n}\dff(\dff f \dff)(\dff g_{\dff 1}\fff,\pff 
\ldots\fff,\pff
g_{\dff n\dff +\dff 1} \dff)
\dff\Bigr)
\]

\vspace*{-27pt}
\[
\quad
\phantom{=\off }
\hspace*{7.7em}
-\off\off
\sum_{i\qff =\qff 0}^{n\qff -\qff 1}(\dff -\dff 1)^{\dff n\dff -\dff 1\dff -\dff i}\dff 
g_{\dff n\dff +\dff 1}\dff\cdot\dff\Bigl(\dff
u_{\dff n}\dff
(\dff 
f
\dff) 
(\dff g_{\dff 0}\fff,\pff \ldots\fff,\pff
g_{\dff i}\fff g_{\dff i\dff +\dff 1}\fff,\pff \ldots\fff,\pff g_{\dff n\dff +\dff 1} \dff)
\dff\Bigr)
\]

\vspace*{-27pt}
\[
\quad
=\off
g_{\dff n\dff +\dff 1}\dff \cdot\dff
\Bigl(\dff 
u_{\dff n}\dff(\dff f \dff)(\dff g_{\dff 1}\fff,\pff \ldots\fff,\pff g_{\dff n}\fff g_{\dff n\dff +\dff 1} \dff)
\dff\Bigr)
\off\off +\off\off
(\dff -\qff 1\dff)^{n\dff +\dff 1}\dff 
g_{\dff n\dff +\dff 1}\dff\cdot\dff\Bigl(\dff 
u_{\dff n}\dff(\dff f \dff)(\dff g_{\dff 1}\fff,\pff 
\ldots\fff,\pff
g_{\dff n\dff +\dff 1} \dff)
\dff\Bigr)
\]

\vspace*{-27pt}
\[
\quad
\phantom{=\off }
\hspace*{7.7em}
+\off\off
\sum_{i\qff =\qff 0}^{n\qff -\qff 1}(\dff -\dff 1)^{\dff n\dff -\dff i}\dff 
g_{\dff n\dff +\dff 1}\dff \cdot\dff
\Bigl(\dff 
u_{\dff n}\dff
(\dff 
f
\dff) 
(\dff g_{\dff 0}\fff,\pff \ldots\fff,\pff
g_{\dff i}\fff g_{\dff i\dff +\dff 1}\fff,\pff \ldots\fff,\pff g_{\dff n\dff +\dff 1} \dff)\
\dff\Bigr)
\]

\vspace*{-27pt}
\[
\quad
=\off
(\dff -\qff 1\dff)^{n\dff +\dff 1}\dff 
g_{\dff n\dff +\dff 1}\dff\cdot\dff\Bigl(\dff 
u_{\dff n}\dff(\dff f \dff)(\dff g_{\dff 1}\fff,\pff 
\ldots\fff,\pff
g_{\dff n\dff +\dff 1} \dff)
\dff\Bigr)
\]

\vspace*{-27pt}
\[
\quad
\phantom{=\off }
\hspace*{7.7em}
+\off\off
\sum_{i\qff =\qff 0}^{n}(\dff -\dff 1)^{\dff n\dff -\dff i}\dff 
g_{\dff n\dff +\dff 1}\dff \cdot\dff
\Bigl(\dff 
u_{\dff n}\dff
(\dff 
f
\dff) 
(\dff g_{\dff 0}\fff,\pff \ldots\fff,\pff
g_{\dff i}\fff g_{\dff i\dff +\dff 1}\fff,\pff \ldots\fff,\pff g_{\dff n\dff +\dff 1} \dff)\
\dff\Bigr)
\]

\vspace{-27pt}
\[
\quad
=\off
g_{\dff n\dff +\dff 1}\dff \cdot\dff
\Bigl(\dff 
d_{\dff n}\dff(\dff u_{\dff n}\dff(\dff f \dff) \dff)
(\dff g_{\dff 0}\fff,\pff \ldots\fff,\pff g_{\dff n\dff +\dff 1} \dff)
\dff\Bigr)\dff.
\]

\vspace{-3pt}
By combining the above calculations,\oss
we conclude that\vspace*{6pt}
\[
\quad
u_{\dff n\dff +\dff 1}\dff(\dff d_{\dff n}\dff(\dff f \dff) \dff)
(\dff g_{\dff 0}\fff,\pff \ldots\fff,\pff g_{\dff n\dff +\dff 1} \dff)
\off =\off
d_{\dff n}\dff(\dff u_{\dff n}\dff(\dff f \dff) \dff)
(\dff g_{\dff 0}\fff,\pff \ldots\fff,\pff g_{\dff n\dff +\dff 1} \dff)\
\]

\vspace*{-6pt}
for all\qss $f\qff \in\qff U_{n}$\qss 
and all\qss
$g_{\dff 0}\fff,\pff g_{\dff 1}\fff,\pff \ldots\fff,\pff g_{\dff n\dff +\dff 1}
\qff \in\qff G$\nnsp.\oss
Therefore\oss
$u_{\dff n\dff +\dff 1}\dff \circ\dff d_{\dff n}
\off =\off
d_{\dff n}\dff \circ\dff u_{\dff n}$\nnsp.\oss  \eproof

\vspace*{1pt}
\myuppar{Bounded cohomology of\dss discrete groups.}
For any $G$\dnsp-module\dss $W$\dss let\vspace*{4pt}   
\[
\quad
W^{\dff G}
\off =\off
\left\{\off 
v\qff \in\qff W \off \left| \off gv\off =\off v
\off\off\mbox{ for all}\off\off 
g\qff \in\qff G 
\off\right\}\dff. \right.
\]

\vspace*{-8pt}
be the subspace of $G$\dnsp-invariant elements in $W$\nnsp.\oss
Given a\dss $G$\dnsp-module $U$ and a\dss $G$\dnsp-resolution\qss (\ref{resolution})\qss of $U$\dnsp,\oss
we can form the complex of subspaces of invariant elements\vspace*{3pt}
\begin{equation}
\label{u-inv}
\quad
\begin{tikzcd}[column sep=large, row sep=normal]\dis
0 \arrow[r]
& 
U_{0}^{\dff G} \arrow[r, "\dis d_{\dff 0}\off"]
&   
U_{1}^{\dff G} \arrow[r, "\dis d_{\dff 1}\off"]
&
U_{2}^{\dff G} \arrow[r, "\dis d_{\dff 2}\off"]
&
\off \ldots \off.
\end{tikzcd}
\end{equation}

\vspace*{-6pt}
From now on,\oss
we will often shorten the notations for resolutions and complexes
by omitting explicit references to the maps $d_{\dff n}$\nnsp,\pss
where\qss $n\qff \geq\qff -\qff 1$\qss or\qss $\geq\qff 0$\dnsp.\oss
So,\oss the resolution\qss (\ref{resolution})\qss can be denoted
simply by\dss $U_\bullet$\nnsp,\oss
and the complex\qss (\ref{u-inv})\qss by\dss $U_\bullet^{\dff G}$\dnsp.\oss
Let\vspace*{4pt}
\[
\quad 
\mathcal{H}^{\fff n}(\dff G\fff,\qff U_\bullet\fff)
\off\off =\off\off
H^{\fff n}\left(\dff U_\bullet^{\dff G} \trf\right)
\off\off =\off\off
\kernel\dff 
\left(\dff d_{\dff n}\qff \left|\qff U_n^{\dff G}\right.
\dff\right)
\Bigl/\qff
\image\dff 
\left(\dff d_{\dff n\dff -\dff 1}\qff \left|\qff U_{n\dff -\dff 1}^{\dff G}\right.
\dff\right)
\Bigr.
\]

\vspace*{-8pt}
be the $n${\dnsp}-th cohomology
space of the complex\qss (\ref{u-inv}).\oss
Being a sub-quotient of\dss the semi-normed space $U_n$\nnsp,\oss
it inherits from $U_n$ a semi-norm.\oss
This semi-norm is not a norm if the image\qss 
$\dis
d_{\dff n\dff -\dff 1}\dff(\qff U_{n\dff -\dff 1}^{\dff G}\dff)$\qss
is not closed.\oss
But if $U_n$ is a Banach space,\oss this semi-norm is\qss \emph{complete}\qss
in the sense that the quotient by the subspace of elements with the norm
$0$ is a Banach space.\oss

There is a preferred resolution of any $G$\dnsp-module $U$\dnsp,\oss
namely,\oss the standard resolution\qss (\ref{st-res}).\oss
The\qss \emph{bounded cohomology of\qss $G$\dss with coefficients in}\dss
$U$ are defined as the cohomology spaces\vspace*{4pt}
\[
\quad
\widehat{H}^{\fff n}\fff(\dff G\fff,\qff U \dff)
\off\off =\off\off
\mathcal{H}^{\fff n}\left(\dff 
G\fff,\pff B\dff(\fff G^{\dff \bullet\dff +\dff 1}\dnsp,\pff U \dff) 
\fff\right)
\off\off =\off\off
H^{\fff n}\left(\dff B\dff(\fff G^{\dff \bullet\dff +\dff 1}\dnsp,\pff U \dff)^{\dff G} \qff \right)
\]

\vspace*{-8pt}
of the complex of invariant subspaces of the standard resolution.\oss
See Appendix\qss \ref{invariants}\qss for a more explicit description of\dss this subcomplex\halfff.\oss
The spaces\qss
$\widehat{H}^{\fff n}\fff(\dff G\fff,\qff U \dff)$\qss
are semi-normed real vector spaces.\oss 
The bounded cohomology with coefficients in $\rrr$ considered as a trivial\dss $G$\dnsp-module
will\dss be denoted simply by\dss $\widehat{H}^{\fff n}\fff(\dff G \dff)$\dnsp.\oss
The space\dss $\widehat{H}^{\fff 2}\fff(\dff G \dff)$\dss is always Hausdorff\halfff.\oss
See Appendix\qss \ref{second}.\oss
If\qss $n\qff \geq\qff 3$\dnsp,\oss
then $\widehat{H}^{\fff n}\fff(\dff G \dff)$
may happen to be non-Hausdorff\halfff,\oss
as was shown by\dss T.\dss Soma\dss \cite{so1},\qss \cite{so2}.

\myuppar{Bounded cohomology and other resolutions.}
Suppose that\qss $U\fff,\pff V$\qss are two $G$\dnsp-modules and\qss
$u
\dff \colon\dff
U\ttoo V$\qss
be a $G$\dnsp-morphism.\oss
Suppose that\dss $U_\bullet$\dss is a strong relatively injective resolution of $U$
and\dss $V_\bullet$\dss is a strong relatively injective resolution of\dss $V$\dnsp.\oss

By Lemma\qss \ref{main-homology-lemma}\qss there exists a $G$\dnsp-morphism\oss 
$u_{\dff \bullet}
\qff \colon\qff
U_\bullet
\qff \ttoo\qff
V_\bullet
$\oss 
extending $u$\nnsp,\oss and\dss by\dss Lemma\qss \ref{homotopy-uniqueness}\qss such an extension $u_{\dff \bullet}$
is unique up to chain homotopies consisting of\dss $G$\dnsp-morphisms.\oss
Being\dss a $G$\dnsp-morphism,\pss $u_{\dff \bullet}$ defines a homomorphism\vspace*{3pt}
\[
\quad
u_{\dff *}
\qff \colon\qff
\mathcal{H}^{\fff n}(\dff G\fff,\qff U_\bullet\fff)
\qff \ttoo\qff
\mathcal{H}^{\fff n}(\dff G\fff,\qff V_\bullet\fff)
\]

\vspace*{-9pt}
for every\qss $n\qff \geq\qff 0$\dnsp.\oss
Since $u_{\dff \bullet}$ is unique up to chain homotopies consisting of\dss $G$\dnsp-morphisms,\pss
$u_{\dff *}$ depends only on $u$\nnsp,\oss and
since $u_{\dff \bullet}$ consists of\dss bounded maps,\oss
all maps $u_{\dff *}$ are bounded.\oss

Suppose now that\qss $U_\bullet\fff,\pff U_\bullet'$\qss are two
resolutions of the same $G$\dnsp-module $U$\dnsp,\oss
which are both strong and relatively injective.\oss
Then the identity map\dss $\id_{\dff U}$\dss extends to\dss $G$\dnsp-morphisms\vspace*{2pt}
\begin{equation*}
\quad
\begin{tikzcd}[column sep=normal, row sep=normal]\dis
i_{\dff \bullet}
\qff \colon\qff
U_\bullet
\off 
\arrow[r, shift left=3pt]
&
\off 
U'_\bullet
\off\dff \colon\dff
i'_{\fff \bullet}
\arrow[l, shift left=3pt]
\end{tikzcd}
\end{equation*}

\vspace*{-10pt}
which,\oss in\dss turn,\oss leads\dss to canonical\dss 
homomorphisms\vspace*{2pt}
\begin{equation*}
\quad
\begin{tikzcd}[column sep=normal, row sep=normal]\dis
i_{\dff *}
\qff \colon\qff
\mathcal{H}^{\fff n}(\dff G\fff,\qff U_\bullet\fff)
\off 
\arrow[r, shift left=3pt]
&
\off 
\mathcal{H}^{\fff n}(\dff G\fff,\qff U'_\bullet\fff) 
\off\dff \colon\dff
i'_{\fff *}
\arrow[l, shift left=3pt]
\end{tikzcd}
\end{equation*}

\vspace*{-10pt}
for every\qss $n\qff \geq\qff 0$\dnsp.\oss
The composition\qss 
$i'_{\dff \bullet}\dff \circ\trf i_{\dff \bullet}
\dff \colon\dff
U_\bullet
\ttoo
U_\bullet$\qss
extends the identity $G$\dnsp-morphism $\id_{\dff U}$\nnsp.\oss
But the identity morphism of\dss
$U_\bullet$\dss also extends\dss $\id_{\dff U}$\nnsp.\oss
By\dss Lemma\qss \ref{homotopy-uniqueness}\qss this implies that\dss 
$i'_{\dff \bullet}\dff \circ\trf i_{\dff \bullet}$\qss is chain homotopic to
the identity morphism of\dss
$U_\bullet$\dss and hence the map\vspace*{3pt}
\[
\quad
\mathcal{H}^{\fff n}(\dff G\fff,\qff U_\bullet\fff)
\qff \ttoo\qff
\mathcal{H}^{\fff n}(\dff G\fff,\qff U_\bullet\fff)
\]

\vspace*{-9pt}
induced\dss by\qss
$i'_{\dff \bullet}\dff \circ\trf i_{\dff \bullet}$\qss
is equal to the identity
for every\qss $n\qff \geq\qff 0$\dnsp.\oss
It follows that\dss the composition\qss
$i'_{\dff *}\dff \circ\trf i_{\dff *}$\qss
is equal to the identity.\oss
By the same argument\qss
$i_{\dff *}\dff \circ\trf i'_{\dff *}$\qss
is equal to the identity.\oss
It follows that\qss $i_{\dff *}\fff,\pff i'_{\dff *}$\qss
are mutually inverse isomorphism of\dss vector spaces.\oss
Since\qss $i_{\dff *}\fff,\pff i'_{\dff *}$\qss
are bounded,\oss
they are isomorphisms of\dss topological vector spaces.\oss

It follows that up to a canonical isomorphism\qss 
$\mathcal{H}^{\fff n}(\dff G\fff,\qff U_\bullet\fff)$\qss 
does not depend on the choice of\dss the resolution\dss $U_\bullet$\dss
as a topological vector space.\oss
In particular\halfff,\oss it is isomorphic to\qss
$\widehat{H}^{\fff n}\fff(\dff G\fff,\qff U \dff)$\qss
as a topological vector space.\oss
But its semi-norm depends on\dss $U_\bullet$\nnsp.\oss

\mypar{Theorem.}{norm-minimality}
\emph{Suppose that\dss $U_\bullet$\dss is a strong relatively injective resolution
of\dss a\dss $G$\dnsp-module\dss $U$\dnsp.\oss
Then for every\qss
$n\qff \geq\qff 0$\qss 
there exists a canonical\dss isomorphism of\dss topological vector spaces}\vspace*{3pt}
\[
\quad
\mathcal{H}^{\fff n}(\dff G\fff,\qff U_\bullet\fff)
\qff \ttoo\qff
\widehat{H}^{\fff n}\fff(\dff G\fff,\qff U \dff)
\]

\vspace*{-9pt}
\emph{which is a bounded operator of\dss norm\qss $\leq\qff 1$\nnsp.\oss}

\proof\qss
The existence of such an isomorphism follows from Theorem\qss \ref{comparing-to-standard},\oss
and its independence on any choices follows from\dss Lemma\qss \ref{homotopy-uniqueness}.\oss  \eproof

\tikzcdset{column sep/mycolumnsize6/.initial=6em}

\mysection{Bounded\qss cohomology\qss and\qss the\qss fundamental\qss group}{spaces}

\vspace*{6pt}
\myuppar{Discrete principal bundles.}
Let $G$ be a discrete group.\oss 
Suppose that\oss 
$p\dff \colon\dff \mathcal{X}\ttoo X$\oss
is a locally trivial principal right\dss $G$\dnsp-bundle.\oss
Then the group $G$ acts freely on $\mathcal{X}$ from the right\halfff,\pss
the quotients space $\mathcal{X}/G$ is equal to $X$\nnsp,\oss
and $p$ is a covering space projection.\oss
The action of $G$ on $\mathcal{X}$ induces a left action of $G$ 
on the vector spaces\qss $B^{\fff n}\fff(\dff \mathcal{X} \dff)$\qss 
and thus turns them into left $G$\dnsp-modules.\oss 
The projection of the bundle $p$ induces isometric isomorphisms\oss\vspace*{2.5pt}
\[
\quad 
p^*
\qff \colon\qff 
B^{\fff n}(\dff X \dff)
\qff \ttoo\qff
B^{\fff n}(\dff \mathcal{X} \dff)^{\fff G}
\]

\vspace*{-9.5pt} 
commuting with the differentials,\oss 
and\dss hence 
an isometric isomorphism of the complexes\vspace*{2.5pt}
\begin{equation*}
\quad
p^*
\qff \colon\qff 
B^{\fff \bullet}(\dff X \dff)
\qff \ttoo\qff
B^{\fff \bullet}(\dff \mathcal{X} \dff)^{\fff G}
\end{equation*}

\vspace{-9.5pt}
Of course,\oss the cohomology of\dss $B^{\fff \bullet}(\dff X \dff)$\dss is the bounded cohomology\dss
$\widehat{H}^{\dff *}\fff(\dff X \dff)$\dnsp.\oss

\myuppar{Strictly acyclic bundles.}
The bundle\oss 
$p\dff \colon\dff \mathcal{X}\ttoo X$\oss
is said to be\qss \emph{strictly acyclic}\qss if
$\mathcal{X}$ admits strictly bounded
contracting homotopy in the sense of\dss Theorem\qss \ref{simply-connected-homotopy}.\oss
By\dss Theorem\qss \ref{simply-connected-homotopy}\qss this is the case if\dss
the fundamental group\qss $\pi_{\fff 1}\fff(\dff \mathcal{X} \dff)$\qss
is amenable.\oss
If $p$ is strictly acyclic,\oss then\dss 
$B^{\fff \bullet}(\dff \mathcal{X} \dff)$\dss 
is a strong resolution
of the trivial\dss $G$\dnsp-module $\rrr$
and the cohomology of the complex\dss
$B^{\fff \bullet}(\dff \mathcal{X} \dff)^{\fff G}$\dss
are the cohomology denoted in Section\qss \ref{algebra}\qss by\qss
$\mathcal{H}^{\fff *}\fff(\dff G\fff,\qff B^{\fff \bullet}\fff(\dff \mathcal{X} \dff) \fff)$\dnsp.\oss
It follows that in this case\vspace*{3pt}
\begin{equation}
\label{p-identification}
\quad
\widehat{H}^{\dff *}\fff(\dff X \dff)
\off =\off
\mathcal{H}^{\fff *}\fff(\dff G\fff,\qff B^{\fff \bullet}\fff(\dff \mathcal{X} \dff) \dff)\dff.
\end{equation}

\vspace*{-9pt}
On the other hand,\oss if $p$ is strictly acyclic,\oss then by
Theorem\qss \ref{comparing-to-standard}\qss 
there exists a\dss $G$\dnsp-morphism\vspace*{2.5pt}
\[
\quad
u_{\dff \bullet}
\qff \colon\qff
B^{\fff \bullet}\dff(\dff \mathcal{X} \dff)
\qff \ttoo\qff 
B\dff(\fff G^{\dff \bullet\dff +\dff 1} \dff)
\]

\vspace{-9.5pt}
from\dss $B^{\fff \bullet}(\dff \mathcal{X} \dff)$\dss to the standard\dss $G$\dnsp-resolution
of $\rrr$
extending\dss $\id_{\dff \rrr}$\dss and
consisting of maps $u_{\dff n}$ of the norm\qss $\leq\qff 1$\nnsp.\oss
Since\dss $B\dff(\fff G^{\dff \bullet\dff +\dff 1} \dff)$\dss is relatively injective,\oss
by Lemma\qss \ref{homotopy-uniqueness}\qss 
$u_{\dff \bullet}$ is unique up to chain homotopies.\oss
By passing to\dss $G$\dnsp-invariants and then to the cohomology\halfff,\pss
$u_{\dff \bullet}$ leads to a map\vspace*{3pt}
\[
\quad
u\fff(\fff p \fff)_{*}
\qff \colon\qff
\mathcal{H}^{\fff *}\fff(\dff G\fff,\qff B^{\fff \bullet}\fff(\dff \mathcal{X} \dff) \dff)
\qff \ttoo\qff
\mathcal{H}^{\fff *}\fff(\dff G\fff,\qff B\dff(\fff G^{\dff \bullet\dff +\dff 1} \dff) \dff)
\off\qff =\off\qff
\widehat{H}^{\dff *}\fff(\dff G \dff)
\]

\vspace*{-9pt}
depending only on the action of $G$ on $\mathcal{X}$\nnsp,\oss
i.e.\oss only on the principal\dss $G$\dnsp-bundle $p$\nnsp.\oss
In view of\qss (\ref{p-identification})\qss and 
the definition of\dss $\widehat{H}^{\dff *}\fff(\dff G \dff)$\nnsp,\oss
the map\dss $u\fff(\fff p \fff)_{*}$\dss can be interpreted as a map\vspace*{3pt}
\begin{equation}
\label{p-map}
\quad
u\fff(\fff p \fff)_{*}
\qff \colon\qff
\widehat{H}^{\dff *}\fff(\dff X \dff)
\qff \ttoo\qff
\widehat{H}^{\dff *}\fff(\dff G \dff)\dff.
\end{equation}

\vspace{-9pt}
It depends only on the principal\dss $G$\dnsp-bundle $p$\nnsp.\oss
Since $u_{\dff \bullet}$ consists on the maps of the norm\qss $\leq\qff 1$\nnsp,\oss
the norm of\dss $u\fff(\fff p \fff)_{*}$\dss is also\qss $\leq\qff 1$\nnsp.\oss

\myuppar{A morphism\qss
$B\dff(\fff G^{\dff \bullet\dff +\dff 1} \dff)\ttoo B^{\fff \bullet}\dff(\dff \mathcal{X} \dff)$\dnsp.}
Let us construct a morphism of resolutions\vspace*{2pt}
\[
\quad
r_{\fff \bullet}
\qff \colon\qff
B\dff(\fff G^{\dff \bullet\dff +\dff 1} \dff)\ttoo B^{\fff \bullet}\dff(\dff \mathcal{X} \dff)
\]

\vspace{-10pt}
extending\dss $\id_{\dff \rrr}$\dss
and consisting of maps $r_{\fff n}$ of the norm\qss $\leq\qff 1$\nnsp.\oss

Let $F$ be a fundamental set for 
the action of $G$ on $\mathcal{X}$\dnsp.\oss
i.e.\qss a\dss subset\qss $F\off \subset\off \mathcal{X}$\qss
such that $F$ intersects each $G$\dnsp-orbit in exactly one point\halfff.\oss
Let\qss $n\qff \geq\qff 0$\qss
and\dss let\dss $v_{\fff i}$\dss be the $i$\dnsp-th vertex 
of\dss $\Delta_{\dff n}$\nnsp.\qff\oss
Let\vspace*{3pt}
\[
\quad
s_{\fff n}
\qff \colon\qff
S_n\fff(\dff \mathcal{X} \dff)
\qff \ttoo\qff
G^{\fff n\dff +\dff 1}
\] 

\vspace*{-9pt}
be the map defined as follows.\oss
For a singular simplex\qss 
$\sigma\dff \colon\dff \Delta_{\dff n}\ttoo \mathcal{X}$\oss 
let\vspace*{3pt}
\[
\quad
s_{\fff n}\fff(\dff \sigma \dff)
\off =\off
(\dff g_{\dff 0}\fff,\pff g_{\dff 1}\fff,\pff \ldots\fff,\pff g_{\dff n} \dff)\dff,
\]

\vspace*{-9pt} 
where\oss 
$g_{\dff 0}\fff,\pff g_{\dff 1}\fff,\pff \ldots\fff,\pff g_{\dff n}$\oss 
are the unique elements of $G$ such that\vspace*{9pt}
\[
\quad
\sigma (\fff v_{\fff 0} \dff)
\hspace*{1.2em}
\off\qff \in\off\off 
F\qff
g_{\dff n}\qff,
\]

\vspace{-33pt}
\[
\quad
\sigma (\fff v_{\fff 1} \dff)
\hspace*{1.2em}
\qff\off \in\off\off
F\qff 
g_{\dff n\dff -\dff 1}\dff g_{\dff n}\qff,
\]

\vspace{-33pt}
\[
\quad
\sigma (\fff v_{\fff 2} \dff)
\hspace*{1.25em}
\qff\off \in\off\off
F\qff 
g_{\dff n\dff -\dff 2}\dff g_{\dff n\dff -\dff 1}\dff g_{\dff n}\qff,
\]

\vspace{-33pt}
\[
\quad
\ldots\ldots\ldots\ldots\qff,
\]

\vspace{-30pt}
\[
\quad
\sigma (\fff v_{\fff n\dff -\dff 1} \dff)
\qff\off \in\off\off 
F\qff
g_{\dff 1}\off \ldots\off g_{\dff n\dff -\dff 1}\dff g_{\dff n}\qff,
\]

\vspace{-33pt}
\[
\quad
\sigma (\fff v_{\fff n} \dff)
\hspace*{1.2em}
\qff\off \in\off\off 
F\qff
g_{\dff 0}\dff g_{\dff 1}\off \ldots\off g_{\dff n\dff -\dff 1}\dff g_{\dff n}\qff.
\]

\vspace{-3pt}
Let\oss 
$r_{\fff n}
\qff \colon\qff B\fff(\dff G^{\fff n\dff +\dff 1} \fff)
\qff \ttoo\qff 
B^{\fff n}\fff(\dff \mathcal{X} \dff)$\oss 
be the map defined
by the formula\oss\vspace*{2pt}
\[
\quad
r_{\fff n}\dff(\fff f \fff)(\dff \sigma \dff)
\off =\off
f\fff(\dff s_{\fff n}\fff(\dff \sigma \dff)\dff) 
\]

\vspace*{-8pt}
A direct verification shows that $r_{\fff \bullet}$ 
commutes with the differentials and\dss hence is a morphism of resolutions.\oss
Obviously,\oss  
$r_{\fff \bullet}$ 
extends\dss $\id_{\dff \rrr}$\dss and consists
of maps $r_{\fff n}$ with the norm\qss $\leq\qff 1$\nnsp.\oss

A better way to see that $r_{\fff \bullet}$ 
commutes with the differentials is to note that $s_{\fff \bullet}$ commutes with
the face operators of the singular complex\dss
$S_{\bullet}\fff(\dff \mathcal{X} \dff)$\dss
and of\dss the nerve $\mathcal{NG}$ of the category $\mathcal{G}$ associated with the group $G$\dss
(see Appendix\qss \ref{categories-classifying}).\oss 
This immediately implies that $r_{\fff \bullet}$ 
commutes with the differentials.\oss

\mypar{Lemma.}{acyclic-bundles-ri}
\emph{The\qss $G$\dnsp-modules\qss
$B^{\fff n}\dff(\dff \mathcal{X} \dff)$\qss
are relatively injective for all\dss $n$\nnsp.\oss}

\proof\qss
Let $F$ be a fundamental set for 
the action of $G$ on $\mathcal{X}$\dnsp.\oss
Let\qss $S_n\fff(\dff \mathcal{X}\nsp,\off F \dff)$\qss 
be the set of singular simplices\oss 
$\Delta_{\dff n}\ttoo \mathcal{X}$\oss
taking the first vertex $v_{\fff 0}$ of $\Delta_{\dff n}$ into $F$\nnsp.\oss
If\qss
$\sigma\dff \colon\dff \Delta_{\dff n}\ttoo \mathcal{X}$\qss
is a singular simplex,\oss
then\qss
$\sigma\off =\off \tau\dff g$\qss
for unique\qss $\tau\qff \in\qff S_n\fff(\dff \mathcal{X}\nsp,\off F \dff)$\qss
and\qss $g\qff \in\qff G$\nnsp.\oss
Hence the map\vspace*{2pt}
\[
\quad
S_n\fff(\dff \mathcal{X}\nsp,\off F \dff)\dff \times\qff G
\qff \ttoo\qff
S_n\fff(\dff \mathcal{X}\dff)
\]

\vspace*{-10pt}
defined by\qss $(\tau\fff,\pff g)\qff \longmapsto\qff \tau\dff g$\qss 
is a bijection.\oss
This bijection is equivariant with respect to the obvious right action of $G$ on\qss 
$S_n\fff(\dff \mathcal{X}\nsp,\pff F \dff)\dff \times\qff G$\qss
and\dss the action of\dss $G$ on\dss $S_n\fff(\dff \mathcal{X} \dff)$\dss
induced by the action of\dss $G$\sss on\sss $\mathcal{X}$\dnsp.\oss
Therefore this bijection leads to an isometric isomorphism of\dss $G$\dnsp-modules\vspace*{2pt}
\[
\quad
B\fff(\dff G\fff,\pff B\fff(\dff S_n\fff(\dff \mathcal{X}\nsp,\off F \dff) \dff) \dff)
\qff \ttoo\qff
B^{\fff n}\fff(\dff \mathcal{X} \dff)
\]

\vspace*{-10pt}
where the Banach space $B\fff(\dff S_n\fff(\dff \mathcal{X}\nsp,\off F \dff)\dff)$
is considered as a trivial $G$\dnsp-module.\oss
In view of this isomorphism Lemma\qss \ref{bgv}\qss implies that the $G$\dnsp-module\qss
$B^{\fff n}\fff(\dff \mathcal{X} \dff)$\qss is relatively injective.\oss  \eproof

\mypar{Theorem.}{acyclic-bundles-isomorphism}
\emph{If\oss 
$p\dff \colon\dff \mathcal{X}\ttoo X$\oss
is strictly acyclic,\oss then\qss $u\fff(\fff p \fff)_{*}$\qss
is an isometric isomorphism.\oss
The inverse map is induced by\oss
$r_{\fff \bullet}
\qff \colon\qff
B\dff(\fff G^{\dff \bullet\dff +\dff 1} \dff)\ttoo B^{\fff \bullet}\dff(\dff \mathcal{X} \dff)$\dnsp.}

\proof\qss
The morphisms of resolutions $u_{\dff \bullet}$ and $r_{\fff \bullet}$ induce maps\vspace*{1pt}
\begin{equation*}
\quad
\begin{tikzcd}[column sep=normal, row sep=normal]\dis
u_{\dff *}
\qff \colon\qff
H^{\fff n}\fff(\dff G\fff,\pff B^{\fff \bullet}\dff(\dff \mathcal{X} \dff) \fff)
\off 
\arrow[r, shift left=3pt]
&
\off 
H^{\fff n}\fff(\dff G\fff,\pff B\dff(\fff G^{\dff \bullet\dff +\dff 1} \dff) \fff) 
\off\dff \colon\dff
r_{\fff *}
\arrow[l, shift left=3pt]
\end{tikzcd}
\end{equation*}

\vspace*{-10.5pt}
which have the norm\qss $\leq\qff 1$\qss together with $u_{\dff n}$ and $r_{\fff n}$\nnsp.\oss
If\dss the compositions\qss $r_{\fff *}\dff \circ\dff u_{\dff *}$\qss
and\qss $u_{\dff *}\dff \circ\dff r_{\fff *}$\qss 
are equal to the identity,\oss then $r_{\fff *}$ and $u_{\dff *}$
are mutually inverse isomorphisms,\oss
and since they both have the norm\qss $\leq\qff 1$\nnsp,\oss
even isometric isomorphisms.\oss
In order to prove that\qss $r_{\fff *}\dff \circ\dff u_{\dff *}$\qss
and\qss $u_{\dff *}\dff \circ\dff r_{\fff *}$\qss 
are equal to the identity,\oss
it is sufficient to prove that the compositions\vspace*{3pt}
\[
\quad
u_{\dff \bullet}\dff \circ\dff r_{\fff \bullet}
\qff \colon\qff
B\dff(\fff G^{\dff \bullet\dff +\dff 1} \dff)
\qff \ttoo\qff 
B\dff(\fff G^{\dff \bullet\dff +\dff 1} \dff)
\hspace*{1.5em}\mbox{ and }\hspace*{1.5em}
r_{\fff \bullet}\dff \circ\dff u_{\dff \bullet}
\qff \colon\qff
B^{\fff \bullet}\dff(\dff \mathcal{X} \dff)
\qff \ttoo\qff 
B^{\fff \bullet}\dff(\dff \mathcal{X} \dff)
\]

\vspace*{-9.5pt}
are chain homotopic to the identity by chain homotopies consisting of\dss $G$\dnsp-morphisms.\oss
Both\qss $u_{\dff \bullet}\dff \circ\dff r_{\fff \bullet}$\qss
and\qss
$r_{\fff \bullet}\dff \circ\dff u_{\dff \bullet}$\qss
extend\dss  
$\id_{\dff \rrr}$\nnsp,\oss
as also do the identity morphisms of\dss 
$B\dff(\fff G^{\dff \bullet\dff +\dff 1} \dff)$\dss
and\dss
$B^{\fff \bullet}\dff(\dff \mathcal{X} \dff)$\dnsp.\oss
But\dss 
$B\dff(\fff G^{\dff \bullet\dff +\dff 1} \dff)$\dss 
is relatively injective by Lemma\qss \ref{bgv}\qss
and\dss
$B^{\fff \bullet}\dff(\dff \mathcal{X} \dff)$\dss
is relatively injective by Lemma\qss \ref{acyclic-bundles-ri}.\oss
Hence\qss
$u_{\dff \bullet}\dff \circ\dff r_{\fff \bullet}$\qss
and\qss
$r_{\fff \bullet}\dff \circ\dff u_{\dff \bullet}$\qss
are chain homotopic to the identity
by\dss Lemma\qss \ref{homotopy-uniqueness}.\oss  \eproof

\mypar{Theorem.}{spaces-groups-quotient}
\emph{Let\qss $X$\qss be a path-connected space and\dss let\oss 
$\Gamma\off =\off \pi_{\fff 1}\halfff(\dff X \dff)$\dnsp.\oss
Suppose that\qss $A$\qss is a normal amenable subgroup of\pss $\Gamma$\dnsp.\oss
Then there exists a canonical isometric isomorphism}\vspace*{2pt}
\[
\quad
\widehat{H}^{\dff *}\fff(\dff X \dff) 
\qff \ttoo \qff
\widehat{H}^{\dff *}\fff(\dff \Gamma/A \dff)\dff.
\]

\vspace{-10pt}
\emph{In particular\halfff,\oss there exists a canonical isometric isomorphism}\oss
$\widehat{H}^{\dff *}\fff(\dff X \dff) 
\qff \ttoo \qff
\widehat{H}^{\dff *}\fff(\dff \pi_{\fff 1}\halfff(\dff X \dff) \dff)$\dnsp.

\proof\qss
In view of Section\qss \ref{weak-equivalences}\qss we can assume that $X$ is 
a\dss CW-complex\halfff.\oss
Let\oss 
$p\dff \colon\dff \mathcal{X}\ttoo X$\oss
be the covering space of\dss $X$ corresponding to the subgroup $A$ of\qss 
$\Gamma\off =\off \pi_{\fff 1}\halfff(\dff X \dff)$\dnsp,\oss
and\dss let\qss $G\off =\off \Gamma/A$\nnsp.\oss 
Then $G$ acts freely on $\mathcal{X}$ and\qss 
$\mathcal{X}/G
\off =\off X$\nnsp.\oss
Hence\oss
$p\dff \colon\dff \mathcal{X}\ttoo X$\oss
is a locally trivial principal\dss $G$\dnsp-bundle.\oss
Since\qss $\pi_{\fff 1}\fff(\dff \mathcal{X} \dff)$\qss is isomorphic to $A$
and hence is amenable,\oss
Theorem\qss \ref{simply-connected-homotopy}\qss implies that the bundle $p$
is strictly acyclic.\oss
Therefore the theorem follows from
Theorem\qss \ref{acyclic-bundles-isomorphism}.\oss  \eproof

\mypar{Theorem\qss ({\fff}Mapping\dss theorem).}{mapping-theorem}
\emph{Let\qss $X\fff,\pff Y$\qss 
be two path-connected spaces and\oss
$\varphi\dff \colon\dff Y\ttoo X$\oss be a continuous map.\oss
If\dss the induced\dss homomorphism of the fundamental groups}\oss\vspace*{2pt} 
\[
\quad
\varphi_*
\dff \colon\dff 
\pi_{\dff 1}\fff(\dff Y \dff)\ttoo \pi_{\dff 1}\fff(\dff X \dff)
\]

\vspace*{-10pt} 
\emph{is surjective with amenable kernel,\oss 
then\oss 
$\varphi^*
\dff \colon\dff
\widehat{H}^{\dff *}\fff(\dff X \dff)\ttoo \widehat{H}^{\dff *}\fff(\dff Y \dff)$\oss
is an isometric isomorphism.}

\proof\qss
In view of Section\qss \ref{weak-equivalences}\qss we can assume that\qss $X\fff,\pff Y$\qss 
are\dss CW-complexes.\oss
Let\qss 
$\Gamma\off =\off \pi_{\dff 1}\fff(\dff Y \dff)$\dnsp,\pss
$G\off =\off \pi_{\dff 1}\fff(\dff X \dff)$\dnsp,\oss
and\pss
$A\off =\off \ker\qff \varphi_*$\nnsp.\oss
Then\qss $G\off =\off \Gamma/A$\nnsp.\oss
Let\oss 
$p\dff \colon\dff \mathcal{X}\ttoo X$\oss
be the universal covering space of\dss $X$\dss
and\oss 
$q\dff \colon\dff \mathcal{Y}\ttoo Y$\oss
be the covering space of $Y$ induced from $p$ by $\varphi$\nnsp.\oss
Then $q$ corresponds to the subgroup
$A
\qff \subset\qff 
\Gamma
\off =\off
\pi_{\fff 1}\halfff(\dff Y \dff)$\dnsp.\oss
The group $G$ acts on both $\mathcal{X}$ and $\mathcal{Y}$
and the canonical map\qss
$\Phi
\dff \colon\dff
\mathcal{Y}
\ttoo
\mathcal{X}$\qss
is\dss $G$\dnsp-equivariant\halfff,\oss
i.e.\dss 
$\Phi$\trs is a morphism of\dss $G$\dnsp-bundles.\oss
By Theorem\qss \ref{simply-connected-homotopy}\qss both $p$ and $q$ are strictly acyclic.\oss
Hence there exist morphisms\qss\vspace*{2pt}
\[
\quad
u_{\dff \bullet}
\qff \colon\qff
B^{\fff \bullet}\dff(\dff \mathcal{X} \dff)
\qff \ttoo\qff 
B\dff(\fff G^{\dff \bullet\dff +\dff 1} \dff)
\hspace*{1.5em}\mbox{ and }\hspace*{1.5em}
v_{\dff \bullet}
\qff \colon\qff
B^{\fff \bullet}\dff(\dff \mathcal{Y} \dff)
\qff \ttoo\qff 
B\dff(\fff G^{\dff \bullet\dff +\dff 1} \dff)
\]

\vspace{-10pt}
extending\dss $\id_{\dff \rrr}$\nnsp.\oss
The triangle\vspace*{-0.5pt}
\begin{equation*}
\quad
\begin{tikzcd}[column sep=normal, row sep=huge]\dis
B^{\fff \bullet}\dff(\dff \mathcal{X} \dff) 
\arrow[rr,  "\dis \Phi^{\fff *}"]
\arrow[dr, "\dis u_{\dff \bullet}"']
&
& 
B^{\fff \bullet}\dff(\dff \mathcal{Y} \dff)
\arrow[ld, "\dis v_{\dff \bullet}"] 
\\ 
& 
B\dff(\fff G^{\dff \bullet\dff +\dff 1} \dff)
&  
\end{tikzcd}
\end{equation*}

\vspace{-12pt}
does not need to be commutative,\oss
but it is commutative up to a chain homotopy consisting of $G$\dnsp-morphisms 
by Lemma\qss \ref{homotopy-uniqueness}.\oss
It follows that the triangle\vspace*{0pt}
\begin{equation*}
\quad
\begin{tikzcd}[column sep=normal, row sep=huge]\dis
\widehat{H}^{\dff *}\fff(\dff X \dff)
\arrow[rr,  "\dis \varphi^{\fff *}"]
\arrow[dr, "\dis u\fff(\fff p \fff)_{*}"']
&
& 
\widehat{H}^{\dff *}\fff(\dff Y \dff)
\arrow[ld, "\dis u\fff(\fff q \fff)_{*}"] 
\\ 
& 
\widehat{H}^{\dff *}\fff(\dff G \dff)\dff,
&  
\end{tikzcd}
\end{equation*}

\vspace{-12pt}
which one gets
after passing to the subspaces of\dss $G$\dnsp-invariants and then to cohomology,\oss
is commutative.\oss
The slanted arrows of this triangle are isometric isomorphisms
by Theorem\qss \ref{acyclic-bundles-isomorphism}.\oss
Hence the theorem follows from the commutativity of this triangle.\oss  \eproof

\myuppar{Functoriality.}
Let\qss 
$\alpha\dff \colon\dff \Gamma\ttoo G$\qss 
be a homomorphism of discrete groups.\oss
If\dss $U$ is a $G$\dnsp-module,\oss then $\alpha$ induces an action of\dss $\Gamma$
on $U$ by the rule\qss 
$h\cdot f\off =\off \alpha\fff(\fff h\dff)\cdot f
$\qss 
for\qss $h\qff \in\qff \Gamma$\dnsp,\qss $f\qff \in\qff U$\dnsp.\oss
The resulting structure of a $\Gamma$\nsp\dnsp-module on $U$ is said to be\qss
\emph{induced}\qss by $\alpha$\nnsp.\oss 
If $U$ is a $G$\dnsp-module,\oss $V$ is a $\Gamma$\nsp\dnsp-module,\pss 
and\qss $u\dff \colon\dff U\ttoo V$\qss is a $\Gamma$\nsp\dnsp-morphism 
with respect to the induced by $\alpha$ structure of a $\Gamma$\nsp\dnsp-module on $U$
then\qss $U^{\fff G}\qff \subset\qff U^{\fff \Gamma}$\qss and\dss hence 
$u$ induces a map\qss\vspace*{3pt}
\[
\quad
u\dff(\fff G\fff,\dff \Gamma\dff)
\qff \colon\qff
U^{\fff G}
\qff \ttoo\qff 
V^{\fff \Gamma}\dff.
\]

\vspace*{-9pt}
The obvious maps\oss 
$\alpha^{\fff *}
\qff \colon\qff
B\dff(\fff G^{\fff n +\dff 1} \dff)
\qff \ttoo\qff
B\dff(\fff \Gamma^{\fff n +\dff 1} \dff)
$\oss 
are $\Gamma$\dnsp-morphisms and
commute with the differentials of the standard resolutions.\oss
Together all these maps define an\dss
$\Gamma$\nsp\dnsp-morphism of\dss resolutions\qss
$B\dff(\fff G^{\fff \bullet\dff +\dff 1} \dff)
\ttoo
B\dff(\fff \Gamma^{\fff \bullet\dff +\dff 1} \dff)$\qss
and\dss hence the maps\dss
$\alpha^{\fff *}(\fff G\fff,\dff \Gamma\dff)$\dss induce a map\vspace*{3pt}
\[
\quad
\alpha^{\fff *}
\qff \colon\qff
\widehat{H}^{\dff *}\fff(\dff G \dff)
\qff \ttoo\qff 
\widehat{H}^{\dff *}\fff(\dff \Gamma \dff)\dff.
\]

\vspace*{-9pt}
in the bounded cohomology.\oss
See Appendix\qss \ref{functoriality-coefficients}\qss for a more general version of\dss
functoriality.

\tikzcdset{row sep/mysize6/.initial=4.8em}

\mypar{Theorem.}{functoriality-spaces-to-groups}
\emph{Let\qss $X\fff,\pff Y$\qss 
be path-connected spaces and\qss
$\varphi\dff \colon\dff Y\ttoo X$\qss is a continuous map.\oss
Let\oss
$\varphi_*
\dff \colon\dff 
\pi_{\dff 1}\fff(\dff Y \dff)\ttoo \pi_{\dff 1}\fff(\dff X \dff)
$\oss 
be the induced map of the fundamental groups.\oss
Then the square}\vspace*{2pt}
\begin{equation*}
\quad
\begin{tikzcd}[column sep=huge, row sep=mysize6]\dis
\widehat{H}^{\dff *}\fff(\dff X \dff) 
\arrow[r,  "\dis \varphi^*"]
\arrow[d]
& 
\widehat{H}^{\dff *}\fff(\dff Y \dff)
\arrow[d] 
\\ 
\widehat{H}^{\dff *}\fff(\dff \pi_{\dff 1}\fff(\dff X \dff) \dff)
\arrow[r, "\dis (\fff \varphi_* \fff)^{\fff *}"]
& 
\widehat{H}^{\dff *}\fff(\dff \pi_{\dff 1}\fff(\dff Y \dff) \dff)\dff,
\end{tikzcd}
\end{equation*}

\vspace{-6pt}
\emph{in which the vertical arrows are canonical isomorphisms 
from Theorem\qss \ref{spaces-groups-quotient},\oss
is commutative.}

\proof\qss
In view of Section\qss \ref{weak-equivalences}\qss we
can assume that\qss $X\fff,\pff Y$\qss are\dss CW-complexes.\oss
Let\oss
$\Gamma\off =\off \pi_{\dff 1}\fff(\dff Y \dff)$\dnsp,\oss
$G\off =\off \pi_{\dff 1}\fff(\dff X \dff)$\dnsp,\oss
and\oss
$\alpha\off =\off \varphi_*$\nnsp.\oss
Let\oss 
$p\dff \colon\dff \mathcal{X}\ttoo X$\oss
and\oss 
$q\dff \colon\dff \mathcal{Y}\ttoo Y$\oss
be the universal coverings of\dss $X$\dss
and\dss $Y$\dss respectively.\oss
There exists a map\oss
$\Phi\dff \colon\dff \mathcal{Y}\ttoo \mathcal{X}$\oss
such that the square\vspace*{2pt}
\begin{equation*}
\quad
\begin{tikzcd}[column sep=huge, row sep=huge]\dis
\mathcal{Y} 
\arrow[r,  "\dis \Phi"]
\arrow[d]
& 
\mathcal{X}
\arrow[d] 
\\ 
Y
\arrow[r, "\dis \varphi"]
& 
X
\end{tikzcd}
\end{equation*}

\vspace{-6pt}
is commutative\dss and\dss $\Phi$\dss is\dss $\alpha$\dnsp-equivariant\halfff,\oss
i.e.\qss if\pss $h\qff \in\qff \Gamma$\qss and\qss $y\qff \in\qff \mathcal{Y}$\dnsp,\oss
then\qss\vspace*{3pt} 
\[
\Phi\fff(\fff h\cdot y\fff)
\off =\off 
\alpha\fff(\fff h\dff)\cdot \Phi\fff(\fff y\fff)\dff.
\]

\vspace{-9pt}
The coverings $p$ and $q$ are locally trivial right\dss 
$G$\dnsp-bundle and\dss $\Gamma$\dnsp-bundle
respectively.\oss
By Theorem\qss \ref{simply-connected-homotopy}\qss both $p$ and $q$ are strictly acyclic
and\dss hence there exist morphisms\vspace*{3pt}
\[
\quad
u_{\dff \bullet}
\qff \colon\qff
B^{\fff \bullet}\dff(\dff \mathcal{X} \dff)
\qff \ttoo\qff 
B\dff(\fff G^{\dff \bullet\dff +\dff 1} \dff)
\hspace*{1.5em}\mbox{ and }\hspace*{1.5em}
v_{\dff \bullet}
\qff \colon\qff
B^{\fff \bullet}\dff(\dff \mathcal{Y} \dff)
\qff \ttoo\qff 
B\dff(\fff \Gamma^{\dff \bullet\dff +\dff 1} \dff)
\]

\vspace{-9pt}
extending\dss $\id_{\dff \rrr}$\nnsp.\oss
The square\vspace*{2pt}
\begin{equation*}
\quad
\begin{tikzcd}[column sep=huge, row sep=huge]\dis
B^{\fff \bullet}\dff(\dff \mathcal{X} \dff) 
\arrow[r,  "\dis \Phi^{\fff *}"]
\arrow[d, "\dis u_{\dff \bullet}"']
&
B^{\fff \bullet}\dff(\dff \mathcal{Y} \dff)
\arrow[d, "\dis v_{\dff \bullet}"] 
\\
B\dff(\fff G^{\dff \bullet\dff +\dff 1} \dff) 
\arrow[r, "\dis \alpha^*"]
& 
B\dff(\fff \Gamma^{\dff \bullet\dff +\dff 1} \dff)  
\end{tikzcd}
\end{equation*}

\vspace{-7pt}
does not need to be commutative,\oss
but it is commutative up to a chain homotopy consisting of $G$\dnsp-morphisms 
by Lemma\qss \ref{homotopy-uniqueness}.\oss
It follows that the square\vspace*{2pt}
\begin{equation*}
\quad
\begin{tikzcd}[column sep=huge, row sep=huge]\dis
\widehat{H}^{\dff *}\fff(\dff X \dff)
\arrow[r,  "\dis \varphi^{\fff *}"]
\arrow[d, "\dis u\fff(\fff p \fff)_{*}"']
&
\widehat{H}^{\dff *}\fff(\dff Y \dff)
\arrow[d, "\dis u\fff(\fff q \fff)_{*}"] 
\\ 
\widehat{H}^{\dff *}\fff(\dff G \dff)
\arrow[r, "\dis \alpha^*"]
& 
\widehat{H}^{\dff *}\fff(\dff \Gamma \dff) 
\end{tikzcd}
\end{equation*}

\vspace{-7pt}
which one gets
after passing to the subspaces of\dss $G$\dnsp-invariants and then to cohomology,\oss
is commutative.\oss
The vertical arrows of this square are the canonical isomorphisms
from Theorem\qss \ref{spaces-groups-quotient}.\oss
Hence the theorem follows from the commutativity of this square.\oss  \eproof

\myuppar{Bounded vs. classical cohomology.}
For any space $X$ the inclusions\qss 
$B^{\fff n}\dff(\dff X \dff)
\ttoo
C^{\fff n}\dff(\dff X \dff)$\qss
commute with the differentials and hence induce a canonical map\vspace*{3pt}
\[
\quad
h_{\trf X}
\qff \colon\qff
\widehat{H}^{\dff *}\fff(\dff X \dff)
\qff \ttoo\qff
H^{\dff *}\fff(\dff X \dff)
\]

\vspace*{-9pt}
from bounded cohomology to the usual singular cohomology with real coefficients.\oss
In order to define an analogue of this map for groups,\oss
one needs to define first the classical cohomology of a group $G$\nnsp.\oss
Let\dss 
$C\dff(\dff G^{\fff n} \dff)$\dss be the vector space of all functions
$G^{\fff n}\ttoo \rrr$\qss and\dss let the sequence\vspace*{3pt}
\begin{equation*}
\quad
\begin{tikzcd}[column sep=large, row sep=normal]\dis
0 \arrow[r]
& 
\rrr \arrow[r, "\dis d_{\dff -\dff 1}\off"]
& 
C\dff(\fff G \dff) \arrow[r, "\dis d_{\dff 0}\off"]
&   
C\dff(\fff G^{\dff 2} \dff) \arrow[r, "\dis d_{\dff 1}\off"]
&
C\dff(\fff G^{\dff 3} \dff) \arrow[r, "\dis d_{\dff 2}\off"]
&
\off \ldots \off,
\end{tikzcd}
\end{equation*}

\vspace*{-6pt}
be defined by the same formulas as for\dss $B\dff(\fff G^{\dff \bullet\dff +\dff 1} \dff)$\dnsp.\oss
This sequence is a complex of the classical\dss $G$\dnsp-modules\qss
(without any boundedness conditions).\oss
One may consider the subcomplex of its $G$\dnsp-invariant elements
and defined the classical cohomology\dss $H^{\dff *}\fff(\dff G \dff)$\dss
of the group $G$ with coefficients in the trivial\dss $G$\dnsp-module $\rrr$
as the cohomology of this subcomplex\halfff.\oss
The inclusions\qss 
$B\dff(\fff G^{\dff n +\dff 1} \dff)
\ttoo
C\dff(\fff G^{\dff n +\dff 1} \dff)$\qss
commute with the differentials and hence induce a canonical map\vspace*{4pt}
\[
\quad
h_{\trf G}
\qff \colon\qff
\widehat{H}^{\dff *}\fff(\dff G \dff)
\qff \ttoo\qff
H^{\dff *}\fff(\dff G \dff)
\]

\vspace*{-8pt}
from the bounded cohomology of\dss $G$\dss to the classical cohomology of\dss $G$\dss with real coefficients.\oss

\myuppar{The classical cohomology and the fundamental group.}
Let $G$ be a discrete group 
and\dss let\oss 
$p\dff \colon\dff \mathcal{X}\ttoo X$\oss
be a principal right\dss $G$\dnsp-bundle.\oss
As before,\pss $p$ induces an isomorphism\vspace*{4pt}
\begin{equation*}
\quad
p^*
\qff \colon\qff 
C^{\fff \bullet}(\dff X \dff)
\qff \ttoo\qff
C^{\fff \bullet}(\dff \mathcal{X} \dff)^{\fff G}\dff.
\end{equation*}

\vspace{-8pt}
The spaces\dss $C^{\fff n}(\dff \mathcal{X} \dff)$\dss
are relatively injective $G$\dnsp-modules
in the classical sense.\oss
The proof is similar to the proof of\dss Lemma\qss \ref{acyclic-bundles-ri}.\oss
The classical version of Lemma\qss \ref{main-homology-lemma}\qss
leads to a $G$\dnsp-morphism\vspace*{4pt}
\begin{equation}
\label{class-chains}
\quad
C\dff(\fff G^{\dff \bullet\dff +\dff 1} \dff)
\qff \ttoo\qff
C^{\fff \bullet}(\dff \mathcal{X} \dff)
\end{equation}

\vspace*{-8pt}
unique up to a chain homotopy by Lemma\qss \ref{homotopy-uniqueness}.\oss
By taking the spaces of\dss $G$\dnsp-invariant elements and passing
to cohomology one gets a canonical map\qss\vspace*{4pt}
\[
\quad
H^{\fff *}\fff(\dff G \dff)
\qff \ttoo\qff
H^{\fff *}\fff(\dff X \dff)\dff.
\]

\vspace*{-8pt}
Suppose now that\qss $X$\qss is a path-connected space
and\dss let\oss
$G\off =\off \pi_{\dff 1}\fff(\dff X \dff)$\dnsp.\oss
If\dss $X$\dss admits a universal covering space,\oss
then taking 
$p$
to be the universal covering leads to a canonical map\qss\vspace*{4pt}
\begin{equation}
\label{class}
\quad
H^{\fff *}\fff(\dff G \dff)
\qff \ttoo\qff
H^{\fff *}\fff(\dff X \dff)\dff.
\end{equation}

\vspace{-8pt}
In general,\oss one needs to begin by replacing $X$ by a weak equivalent\dss CW-complex\halfff.\oss

\vspace*{3pt}
\mypar{Theorem.}{bounded-vs-classical}
\emph{Let\qss $X$\qss be a path-connected space
and let\oss
$G\off =\off \pi_{\dff 1}\fff(\dff X \dff)$\dnsp.\oss
Then the square}\vspace*{6pt}
\begin{equation*}
\quad
\begin{tikzcd}[column sep=huge, row sep=mysize6]\dis
\widehat{H}^{\dff *}\fff(\dff G \dff) 
\arrow[r, "\dis h_{\trf G}"]
\arrow[d]
& 
H^{\dff *}\fff(\dff G \dff)
\arrow[d] 
\\ 
\widehat{H}^{\dff *}\fff(\dff X \dff)
\arrow[r, "\dis h_{\trf X}"]
& 
H^{\dff *}\fff(\dff X \dff)\dff,
\end{tikzcd}
\end{equation*}

\vspace{-3pt}
\emph{where the left vertical arrow is the inverse of the canonical isomorphism 
from Theorem\qss \ref{spaces-groups-quotient}\qss
and the right vertical arrow is the map\qss (\ref{class}),\oss
is commutative.}

\proof\qss
As usual,\oss we may assume that $X$ is a\dss CW-complex\halfff.\oss
Let\oss 
$p\dff \colon\dff \mathcal{X}\ttoo X$\oss
be the universal covering of\dss $X$\nnsp.\oss
The left vertical arrow is an instance of the inverse of the canonical map 
from Theorem\qss \ref{acyclic-bundles-isomorphism}\qss and hence
is induced by a morphism\oss
$r_{\fff \bullet}
\qff \colon\qff
B\dff(\fff G^{\dff \bullet\dff +\dff 1} \dff)\ttoo B^{\fff \bullet}\dff(\dff \mathcal{X} \dff)$\dnsp.\oss
The square\vspace*{6pt}
\begin{equation*}
\quad
\begin{tikzcd}[column sep=huge, row sep=mysize6]\dis
B\dff(\fff G^{\dff \bullet\dff +\dff 1} \dff) 
\arrow[r]
\arrow[d, "\dis r_{\fff \bullet}"']
& 
C\dff(\fff G^{\dff \bullet\dff +\dff 1} \dff)
\arrow[d] 
\\ 
B^{\fff \bullet}(\dff \mathcal{X} \dff)
\arrow[r]
& 
C^{\fff \bullet}(\dff \mathcal{X} \dff)\dff,
\end{tikzcd}
\end{equation*}

\vspace{-3pt}
where the horizontal arrows are inclusions and the right vertical arrow
is the map\qss (\ref{class-chains}),\oss
does not need to be commutative,\oss
but it is commutative up to a chain homotopy consisting of classical\dss $G$\dnsp-morphisms 
by the classical analogue of\dss Lemma\qss \ref{homotopy-uniqueness}.\oss
After passing to\dss $G$\dnsp-invariants and then to cohomology,\oss
this implies that the square from the theorem is commutative.\oss  \eproof

\myuppar{The Eilenberg--MacLane spaces.}
Perhaps,\oss the easiest way to explain to a topologist what is the cohomology of groups
is to define the cohomology of a group $G$
as the cohomology of an Eilenberg--MacLane space\dss $K\dff(\dff G\fff,\pff 1 \fff)$\dnsp.\oss
For the bounded cohomology 
such a definition is equivalent to the definition from Section\qss \ref{algebra}\qss
by Theorem\qss \ref{spaces-groups-quotient}.\oss
Moreover\halfff,\oss the two definitions are equivalent in a functorial
manner\halfff,\oss as follows from Theorem\qss \ref{functoriality-spaces-to-groups}.\oss
For the classical cohomology the two definitions are equivalent by 
the classical analogue of\dss Theorem\qss \ref{spaces-groups-quotient},\oss
and this equivalence is functorial by 
the classical analogue of\dss Theorem\qss \ref{functoriality-spaces-to-groups}.\oss
These classical analogues are well\dss known and are easier than their bounded cohomology counterparts,\oss
but apply only to Eilenberg--MacLane spaces.\oss

\tikzcdset{column sep/my-size-em/.initial=8em}
\tikzcdset{row sep/my-size-em/.initial=12ex}

Let $X$ be a path-connected space,\pss 
$G\off =\off \pi_{\dff 1}\fff(\dff X \dff)$\dnsp,\oss
and\dss let\qss
$\varphi
\dff \colon\dff
X
\ttoo 
K\dff(\dff G\fff,\pff 1 \fff)$\qss
be a map
inducing the identity map of the fundamental groups\qss
$G\off =\off \pi_{\dff 1}\fff(\dff X \dff)
\qff \ttoo\qff
\pi_{\dff 1}\fff(\dff K\dff(\dff G\fff,\pff 1 \fff) \dff)$\nnsp.\oss
Such a map $\varphi$ always exists if\dss $X$ is homotopy equivalent to a\dss CW-complex\halfff.\oss
Then the diagram\vspace*{6pt}
\begin{equation*}
\quad
\begin{tikzcd}[column sep=my-size-em, row sep=my-size-em]\dis
\widehat{H}^{\fff *}\fff(\dff K\dff(\dff G\fff,\pff 1 \fff) \dff) 
\arrow[r, "\dis h_{\trf K\dff(\dff G\fff,\pff 1 \fff)}"]
\arrow[d, "\dis \varphi^{*}"]
& 
H^{\fff *}\fff(\dff K\dff(\dff G\fff,\pff 1 \fff) \dff)  
\arrow[d, "\dis \varphi^{*}"']
\\
\widehat{H}^{\fff *}\fff(\dff X \dff)
\arrow[r, "\dis h_{\trf X}"]
& 
{H}^{\fff *}\fff(\dff X \dff)\dff, 
\end{tikzcd}
\end{equation*}

\vspace*{-3pt}
is commutative and,\oss moreover\halfff,\oss 
is canonically isomorphic to the diagram of\dss Theorem\qss \ref{bounded-vs-classical}.\oss
The proof\dss is left to the reader as an exercise.\oss

\myuppar{The Hirsch--Thurston class of groups.}
Hirsch and Thurston\qss \cite{ht}\qss defined a class of groups
$\mathcal{C}$ 
as the smallest class containing all amenable groups
and such that\dss
if\qss $G\fff,\pff H\qff \in\qff \mathcal{C}$\nnsp,\oss then the free product\qss
$G\dff *\dff H\qff \in\qff \mathcal{C}$\nnsp,\oss
and\dss if\qss $G\qff \in\qff \mathcal{C}$\qss and $G$ is a subgroup of finite index in $K$\nnsp,\oss
then\qss $K\qff \in\qff \mathcal{C}$\nnsp.\oss

The class $\mathcal{C}$ is closed under passing to subgroups of finite index.\oss
One can prove this by induction using Kurosh's theorem about
subgroups of free products of groups\qss
(see,\oss for example,\oss \cite{s},\oss Chapter\qss I,\oss Section\qss 5.5)\qss
and the fact that every subgroup of an amenable group is amenable\qss
(see\qss \cite{gr},\oss Theorem\qss 1.2.5).\oss
We leave details to the reader\halfff.\oss

\tikzcdset{column sep/my-size-ht/.initial=4em}
\tikzcdset{row sep/my-size-ht/.initial=10ex}
\tikzcdset{row sep/my-size-hts/.initial=8ex}

\tikzcdset{column sep/my-size-htm/.initial=4em}
\tikzcdset{row sep/my-size-htm/.initial=10ex}
\tikzcdset{row sep/my-size-htsm/.initial=10ex}

\vspace*{3pt}
\mypar{Theorem.}{class-c}
\emph{If\pss $G\qff \in\qff \mathcal{C}$\nnsp,\oss then\oss
$h_{\trf G}
\dff \colon\dff
\widehat{H}^{\fff *}\fff(\dff G \dff)\ttoo {H}^{\fff *}\fff(\dff G \dff)$\qss
is the zero homomorphism.\oss}

\proof\qss 
To begin with,\qss Theorem\qss \ref{simply-connected-homology}\qss
together with Theorem\qss \ref{spaces-groups-quotient}\qss 
implies that\qss
$\widehat{H}^{\fff *}\fff(\dff G \dff)\off =\off 0$\qss
if\dss $G$ is an amenable group.\oss
Therefore in this case\qss
$\widehat{H}^{\fff *}\fff(\dff G \dff)\ttoo {H}^{\fff *}\fff(\dff G \dff)$\qss
is a zero homomorphism.\oss
Suppose that\qss
$\widehat{H}^{\fff *}\fff(\dff G_{\dff i} \dff)\ttoo {H}^{\fff *}\fff(\dff G_{\dff i} \dff)$\qss
are zero homomorphisms for\qss $i\qff =\qff 1\fff,\pff 2$\nnsp.\oss
The canonical embeddings\qss $G_{\dff i}\ttoo G_{\dff 1} *\fff G_{\dff 2}$\nnsp,\oss
where\qss $i\qff =\qff 1\fff,\pff 2$\nnsp,\oss
induce the horizontal arrows of the diagram\vspace*{3pt}
\begin{equation*}
\quad
\begin{tikzcd}[column sep=my-size-htm, row sep=my-size-htm]\dis
\widehat{H}^{\fff *}(\dff G_{\dff 1} *\fff G_{\dff 2} \dff) \arrow[r]
\arrow[d]
& 
\widehat{H}^{\fff *}(\dff G_{\dff 1} \dff)
\dff \oplus\dff
\widehat{H}^{\fff *}(\dff G_{\dff 2} \dff)  
\arrow[d]
\\
{H}^{\fff *}(\dff G_{\dff 1} *\fff G_{\dff 2} \dff) \arrow[r]
& 
{H}^{\fff *}(\dff G_{\dff 1} \dff)
\dff \oplus\dff 
{H}^{\fff *}(\dff G_{\dff 2} \dff)\dff.
\end{tikzcd}
\end{equation*}

\vspace{-6pt}
Obviously,\oss this diagram is commutative.\oss
It is well known that the lower horizontal arrow is an isomorphism.\oss
Therefore the commutativity of the diagram implies that the homomorphism\qss
$\widehat{H}^{\fff *}(\dff G_{\dff 1} *\fff G_{\dff 2} \dff)
\ttoo
{H}^{\fff *}(\dff G_{\dff 1} *\fff G_{\dff 2} \dff)$\qss
is equal to zero.\oss
Finally,\oss suppose that $G$ is a subgroup of finite index of\dss a group $K$
and\qss
$\widehat{H}^{\fff *}\fff(\dff G \dff)\ttoo {H}^{\fff *}\fff(\dff G \dff)$\qss
is a zero homomorphism.\oss
The inclusion\qss $G\ttoo K$ induces the horizontal arrows of the commutative diagram\vspace*{3pt}
\begin{equation*}
\quad
\begin{tikzcd}[column sep=my-size-htm, row sep=my-size-htsm]\dis
\widehat{H}^{\fff *}(\dff K \dff) \arrow[r]
\arrow[d]
& 
\widehat{H}^{\fff *}(\dff G \dff) 
\arrow[d]
\\
{H}^{\fff *}(\dff K \dff) \arrow[r]
& 
{H}^{\fff *}(\dff G \dff)\dff.
\end{tikzcd}
\end{equation*}

\vspace*{-6pt}
Since the module of coefficients is the trivial module $\rrr$\nnsp,\oss
the lower horizontal arrow is injective.\oss 
It follows that\qss
$\widehat{H}^{\fff *}(\dff K \dff)\ttoo {H}^{\fff *}(\dff K \dff)$\qss
is a zero homomorphism.\oss
The theorem follows.\oss  \eproof

\myuppar{Remark.} It may happen that\qss $G\qff \in\qff \mathcal{C}$\nnsp,\oss
but\qss $\widehat{H}^{\fff *}(\dff G \dff)\qff \neq\qff 0$\nnsp.\oss
For example,\oss this is the case if $G$ is a free group with\qss $\geq\qff 2$\qss free generators.\oss
See\qss \cite{br}\qss or\qss \cite{gro}.\oss

\myuppar{Closed manifolds of negative curvature.}
Hirsch and Thurston\qss \cite{ht}\qss 
suggested that if\qss $M$\dss 
is a closed riemannian manifold of negative curvature,\oss
then\qss $\pi_{\dff 1}\fff(\dff M \dff)\qff \not\in\qff \mathcal{C}$\dnsp.\oss 
See\qss \cite{ht},\oss remarks after Corollary\qss 1.3.\oss
Theorem\qss \ref{class-c}\qss leads to a natural proof of this conjecture.\oss

If $M$ is orientable,\oss
then\qss
$\widehat{H}^{\fff *}(\dff M \dff)\ttoo {H}^{\fff *}(\dff M \dff)$\qss
is non-zero by a theorem of Thurston.\oss 
In fact\halfff,\oss the image of this homomorphism
contains the fundamental class of $M$\nnsp.\oss
See\qss \cite{gro},\oss \cite{t},\oss or\qss \cite{iy}\qss for a proof\halfff.\oss
It follows that in this case\qss
$\widehat{H}^{\fff *}(\dff \pi_{\dff 1}\fff(\dff M \dff) \dff)
\ttoo 
{H}^{\fff *}(\dff \pi_{\dff 1}\fff(\dff M \dff) \dff)$\qss
is non-zero and\dss therefore\qss
$\pi_{\dff 1}\fff(\dff M \dff)\qff \not\in\qff \mathcal{C}$\qss
by Theorem\qss \ref{class-c}.\oss
If $M$ is not orientable,\oss
one can pass to the orientation covering
of $M$ and conclude that\dss
$\pi_{\dff 1}\fff(\dff M \dff)$\dss contains a subgroup
of index $2$ not belonging to $\mathcal{C}$\dnsp.\oss
Since $\mathcal{C}$ is closed under passing to subgroups of finite index\halfff,\oss
this completes the proof\halfff.\oss

The first proof of this conjecture is due to\qss
N.\qss Gusevskii\qss \cite{gus}\qss and\dss based on completely different ideas.\oss
See\qss \cite{i1},\oss Section\qss (5.4)\qss for some remarks about his proof\halfff.\oss
I refer to them mostly in order to correct a blatant translation mistake in\dss the 
English translation of\qss \cite{i1}.\oss
The last phrase of the footnote in Section\qss (5.4)\qss
should be\qss \emph{``It is amazing that this argument went unnoticed''}.

\vspace*{3pt}
\myuppar{The role of the homological algebra.}
The proofs in this section rely
in an essential manner on Lemma\qss \ref{homotopy-uniqueness},\qss
i.e. on the homotopy uniqueness of morphisms between the resolutions.\oss
At the same time they do not rely on
Lemma\qss \ref{main-homology-lemma}\qss about the existence of such morphisms,\qss
despite the fact that the involved resolutions such as
$B\dff(\fff G^{\dff \bullet\dff +\dff 1} \dff)$
and
$B^{\fff \bullet}\dff(\dff \mathcal{X} \dff)$
are strong and relatively injective and\dss hence Lemma\qss \ref{main-homology-lemma}\qss 
applies.\oss
Instead of using Lemma\qss \ref{main-homology-lemma}\qss all\dss the needed morphisms are
constructed explicitly.\oss
The reason is that Lemma\qss \ref{main-homology-lemma}\qss does not provide the estimate\qss
$\leq\qff 1$\qss for the norm,\oss
which is needed to prove that induced maps in cohomology are isometries.\oss

In contrast with the bounded case,\oss the singular cochain complex\dss
$C^{\fff \bullet}(\dff \mathcal{X} \dff)$\dss of the universal covering space
$\mathcal{X}$ of a space $X$ is only rarely acyclic.\oss
By this reason in the classical situation one cannot expect to have
an analogue of\dss the map\dss $u\fff(\fff p \fff)_{*}$\dss
and one has to resort\dss to the classical analogue of\dss
Lemma\qss \ref{main-homology-lemma}\qss
which leads to a map going in the opposite direction.\oss 
The use of\dss this classical analogue causes no problem
since there is no norm to take care of\sss anyhow.

\vspace*{3pt}
\myuppar{The definitions of strongly injective morphisms and relatively injective modules.}
Somewhat surprizingly,\oss the details of these definitions do
not play any essential role.\oss
The requirements\dss 
$\|\dff \sigma \dff\|\qff \leq\qff 1$\dss 
in the definition of strongly injective morphisms,\oss
and\qss
$\|\dff \beta \dff\|\qff \leq\qff \|\dff \alpha \dff\|$\qss
in the definition of relatively injective modules
could be relaxed to the boundedness of the operators\dss 
$\sigma$\dss and\dss $\beta$\dss
without affecting any of the above results
except the parenthetical discussion of the norm of morphisms of an extension
after Lemma\qss \ref{homotopy-uniqueness}.\oss
These requirements are imposed because they are met every time
one needs these notions and potentially they are useful.
In contrast\halfff,\oss the assumption\qss 
$\|\qff K_{\dff n} \qff\|\qff \leq\qff 1$\qss
in the definition of strong resolutions
plays a key role in\dss Theorem\qss \ref{comparing-to-standard}\qss
and\dss hence in most of the results about isomorphisms being isometries.

\myuppar{Other situations.}
The technique developed in this section is quite flexible.\oss
In Appendix\qss \ref{borel}\qss this is illustrated by applying
it to straight Borel cochains,\oss playing an important role in\qss \cite{gro}.

\mysection{The\qss covering\qss theorem}{covering}

\vspace*{6pt}
\myuppar{Coverings.}
Let $\mathcal{U}$ be a covering of a paracompact\sss space $X$\nnsp.\oss
We will assume the following.\sss\vspace*{-3pt}
\begin{itemize}
\item[({\fff}C{\fff})] \dnsp$\mathcal{U}$ is either open,\qss or\sss 
closed and\dss locally finite.\pss 
Every finite intersection of\dss elements of $\mathcal{U}$  
is\dss path connected\sss and,\oss if\dss $\mathcal{U}$ is closed,\oss 
also homologically\sss locally connected.\oss
\end{itemize}

\vspace*{-3pt}
Recall\dss that\sss a space\sss $Z$\sss is\dss
\emph{homologically\dss locally\sss connected}\pss
if\trs for every\dss $n\qff \geq\qff 0$\nnsp,\oss 
every\sss $z\qff \in\qff Z$\nnsp,\oss
and every neighborhood\sss $U$\sss of\dss $z$\sss
there exists another neighborhood 
$V\qff \subset\qff U$ of\dss $z$\sss
such\dss that\dss the inclusion\sss homomorphism\sss
$H_{\fff n}\fff(\trf V,\qff \{\dff z\trf\}\dff;\qff \zzz\trf)
\qff \ttoo\qff
H_{\fff n}\fff(\trf U,\qff \{\dff z\trf\}\dff;\qff \zzz\trf)$\sss
is\dss equal\dss to zero.\oss
For example,\oss every manifold and every CW-complex\sss is\sss homologically\sss
locally connected.\oss

The condition\qss ({\fff}C{\fff})\qss is needed to apply the sheaf theory.\oss
It ensures that $X$ and\sss finite intersections of\dss elements of\dss $\mathcal{U}$\sss
behave nicely\sss with respect\sss to the singular cohomology\sss theory.\oss

A path-connected subset $Y$ of a space $X$ is said to be\qss \emph{amenable}\qss if
the image of the inclusion homomorphism\qss
$\pi_{\fff 1}\fff(\dff Y \dff)\ttoo \pi_{\fff 1}\fff(\dff X \dff)$\qss
is an amenable group.\oss
For example,\oss if the fundamental group\dss
$\pi_{\fff 1}\fff(\dff Y \dff)$\dss is amenable,\oss the $Y$ is amenable.\oss
This follows from the fact that any quotient group of an amenable group is amenable.\oss
A covering\sss $\mathcal{U}$ of\dss $X$\sss is\sss said\sss to be\qss \emph{amenable}\qss
if\dss every\sss its element\sss is\sss amenable,\oss
and\qss 
\emph{amenable in\sss the sense of\qss Gromov}\qss
if\dss the covering\sss by\dss the path-components of\dss elements of\dss 
$\mathcal{U}$\sss is\sss amenable.\oss

\myuppar{The nerve of a covering.}
Recall that the nerve\dss $N$\dss of a covering\dss $\mathcal{U}$\dss 
is a simplicial complex in the sense,\oss for example,\oss 
of\qss \cite{sp},\oss Section\qss 3.1.\oss
The vertices of\dss $N$\dss are in a one-to-one correspondence with the set\dss $\mathcal{U}$\dss
and the simplices of\dss $N$\dss are finite sets of vertices 
such that the intersection of the corresponding
elements of\dss $\mathcal{U}$\dss is non-empty.\oss
For each simplex\dss $\sigma$\dss of\dss $N$\dss we denote by\dss $|\fff \sigma \fff|$\dss 
the intersection of the elements of the covering\dss $\mathcal{U}$\dss 
corresponding to the vertices of\dss $\sigma$\nnsp.\oss
We will  assume  that the set of the vertices of\dss $N$\dss is ordered
by a linear order\dss $<$\nnsp.

\vspace*{3pt}
\mypar{Theorem\qss (Covering\dss theorem).}{covering-theorem}
\emph{Let\dss
$\mathcal{U}$\sss be an amenable covering of a paracompact\sss space\qss $X$ 
satisfying\sss the assumption\dss \textup{({\fff}C{\fff})}.\oss
Let\dss $N$\dss be the nerve of\dss $\mathcal{U}$ and let\dss $|\fff N \fff|$\dss
be the geometric realization of\dss $N$\nnsp.\oss
Then the canonical homomorphism\qss
$\widehat{H}^{\fff *}\fff(\dff X \dff)\ttoo {H}^{\fff *}\fff(\dff X \dff)$\qss
can be factored through the canonical homomorphism\qss
${H}^{\fff *}\fff(\dff |\fff N \fff | \dff)\ttoo {H}^{\fff *}\fff(\dff X \dff)$\dnsp,\oss
i.e.\qss the diagram of the solid arrows}\vspace{3pt}
\begin{equation}
\label{covering-triangle}
\quad
\begin{tikzcd}[column sep=normal, row sep=huge]\dis
\widehat{H}^{\fff *}\fff(\dff X \dff) \arrow[rr]
\arrow[rd, dashed]
&
& 
{H}^{\fff *}\fff(\dff X \dff) 
\\
&
{H}^{\fff *}\fff(\dff |\fff N \fff| \dff) \arrow[ru]
& 
\end{tikzcd}
\end{equation}

\vspace*{-9pt}
\emph{can be completed to a commutative diagram by a dashed arrow.\oss}

\newpage
\proof\qss
We start by reducing the theorem to its special case when the fundamental groups of
all components of elements of $\mathcal{U}$ and 
of finite intersections of elements of $\mathcal{U}$
are amenable.\oss 
Let $U$ be such a component\halfff,\oss
and let $\alpha$ be a loop in $U$ contractible in $X$\nnsp.\oss
Let $X'$ be the result of attaching a two-dimensional disc along the loop $\alpha$ to $X$\nnsp.\oss
Since $\alpha$ is contractible in $X$\nnsp,\oss
the fundamental groups of\dss $X'$ and $X$ are the same,\oss
and hence\qss
({\fff}by\dss Theorem\qss \ref{mapping-theorem})\qss
the bounded cohomology of\dss $X'$ and $X$ are the same.\oss
If we include the attached disc into all elements of the covering $\mathcal{U}$
containing $U$\nnsp,\oss we will a get a covering $\mathcal{U}'$ of the $X'$\dnsp.\oss
Clearly,\pss $\mathcal{U}'$ has the same nerve as $\mathcal{U}$\dnsp.\oss
Moreover\halfff,\oss there is a canonical homomorphism\oss 
${H}^{\fff *}\fff(\dff X' \dff)
\ttoo
{H}^{\fff *}\fff(\dff X \dff)$\dnsp.\oss

The glueing operation preserves the paracompactness.\oss
See\qss \cite{fp},\oss Proposition\qss A.5.1\dff({\fff}v{\fff})\qss and\qss
Exercise\qss 5\qss in the Appendix\qss (pp.\qss 273\qss and\qss 305).\oss 
If\dss  $\mathcal{U}$ is closed and locally finite,\oss then the same is true for $\mathcal{U}'$\dnsp.\oss
If\dss  $\mathcal{U}$ is open,\oss then $\mathcal{U}'$ is not open in general.\oss
In this case one needs to replace every $V\qff \in\qff \mathcal{U}'$ by an approptiate open neighborhood 
of $V$ in $X'$ having $V$ as a deformation retract\halfff.\oss
It follows that we can assume that $\mathcal{U}'$ satisfies the condition\qss ({\fff}C{\fff}).\oss

Therefore,\oss if the theorem is true for $X'$
and $\mathcal{U}'$\dnsp,\oss
then it is true for $X$ and $\mathcal{U}$\dnsp.\oss
By attaching discs in this way one can kill 
the kernel of\qss
$\pi_{\fff 1}\fff(\dff U \dff)\ttoo \pi_{\fff 1}\fff(\dff X \dff)$\qss
and hence turn $U$ into a subspace with amenable fundamental group\qss
(since the image of\qss
$\pi_{\fff 1}\fff(\dff U \dff)\ttoo \pi_{\fff 1}\fff(\dff X \dff)$\qss
is amenable by the assumption).\oss
Attaching such discs to all components $U$ as above
reduces the theorem to the case where all\dss $\pi_{\fff 1}\fff(\dff U \dff)$\dss are amenable.\oss
In the rest of the proof we will consider only this case.\oss

For\qss $p\qff \geq\qff 0$\qss let\dss $N_{p}$\dss be the set of $p$\dnsp-dimensional simplices of
the nerve $N$\nnsp.\oss
For\qss $p\fff,\pff q\qff \geq\qff 0$\qss let\vspace*{3pt}
\[
\quad
C^{\fff p}\fff(\dff N\fff,\pff B^{\fff q} \dff)
\off\off =\off\off
\prod_{\sigma\qff \in\qff N_p} B^{\fff q}\fff(\dff |\fff \sigma \fff| \dff)\dff.
\]

\vspace*{-9pt}
Let us consider for every\qss $\sigma\qff \in\qff N_p$\qss the complex\vspace*{3pt}
\begin{equation}
\label{b-sigma}
\quad
\begin{tikzcd}[column sep=large, row sep=huge]\dis
0 \arrow[r]
& 
\rrr \arrow[r]
& 
B^{\fff 0}\fff(\dff |\fff \sigma \fff| \dff) \arrow[r, "\dis d_{\dff 0}"] 
&   
B^{\fff 1}\fff(\dff |\fff \sigma \fff| \dff) \arrow[r, "\dis d_{\dff 1}"] 
&
\off \ldots \off.
\end{tikzcd}
\end{equation}

\vspace{-3pt}
For each\qss $p\qff \geq\qff 0$\qss
the product\sss of these complexes over all\qss
$\sigma\qff \in\qff N_p$\qss is the complex\vspace*{3pt}
\begin{equation}
\label{c-nerve}
\quad
\begin{tikzcd}[column sep=large, row sep=huge]\dis
0 \arrow[r]
& 
C^{\fff p}\fff(\dff N \dff) \arrow[r]
& 
C^{\fff p}\fff(\dff N\fff,\pff B^{\fff 0} \dff) \arrow[r, "\dis d_{\dff 0}"] 
&   
C^{\fff p}\fff(\dff N\fff,\pff B^{\fff 1} \dff) \arrow[r, "\dis d_{\dff 1}"] 
&
\off \ldots \off,
\end{tikzcd}
\end{equation}

\vspace*{-3pt}
where\dss $C^{\fff p}\fff(\dff N \dff)$\dss is the space of real simplicial
$p$\dnsp-cochains of $N$\nnsp.\oss

Let\dss $p\qff \geq\qff 0$\dss and\dss let $\sigma$ be
a $p$\dnsp-simplex $\sigma$ of\dss $N$\nnsp.\oss 
Let\qss\vspace*{1.5pt}
\[
\quad
v_{\fff 0}\off <\off v_{\fff 1}\off <\off \ldots\off <\off v_{p}
\]

\vspace*{-10.5pt}
be the vertices of $\sigma$ listed in the increasing order\halfff.\oss
For\oss
$i\qff =\qff 0\fff,\pff 1\fff,\pff \ldots\fff,\pff p$\oss
let\oss 
$\partial_{\fff i}\fff\sigma
\off =\off 
\sigma\qff \smallsetminus\qff \{\qff v_{\fff i} \qff\}$\nnsp.\oss
In other terms,\oss
$\partial_{\fff i}\fff\sigma$\qss
is the simplex with the vertices\oss
$v_{\fff 0}\fff,\pff \ldots\fff,\pff \widehat{v_{\fff i}}\dff,\pff \ldots\fff,\pff v_{p}$\nnsp.\qff\oss
Let\oss \vspace*{6pt}
\[
\quad
\Delta_{\trf \sigma\fff,\trf i}
\pff \colon\pff
B^{\fff q}\fff(\dff |\dff \partial_{\fff i}\dff\sigma \dff| \dff)
\ttoo
B^{\fff q}\fff(\dff |\fff \sigma \fff| \dff)
\oss \]

\vspace{-6pt}
be the map induced by the inclusion\qss 
$|\dff \sigma \fff|
\ttoo
|\fff \partial_{\fff i}\dff\sigma \fff|$\nnsp,\qff\oss
and\vspace*{3pt}
\[
\quad
i_{\fff \sigma}
\qff \colon\qff
B^{\fff q}\fff(\dff |\fff \sigma \fff| \dff)
\ttoo
C^{\fff p}\fff(\dff N\fff,\pff B^{\fff q} \dff)\dff,
\]

\vspace*{-33pt}
\[
\quad
p_{\fff \sigma}
\qff \colon\qff
C^{\fff p}\fff(\dff N\fff,\pff B^{\fff q} \dff)
\ttoo
B^{\fff q}\fff(\dff |\fff \sigma \fff| \dff)
\]

\vspace{-6pt}
be the inclusion of\dss and the projection
to a factor respectively.\oss
Let us define maps\vspace*{6pt}
\[
\quad
\delta_{\fff i}\dff,\off \delta
\qff \colon\qff
C^{\fff p\dff -\dff 1}\fff(\dff N\fff,\pff B^{\fff q} \dff)
\ttoo
C^{\fff p}\fff(\dff N\fff,\pff B^{\fff q} \dff)
\]

\vspace*{-6pt}
by the formulas\vspace*{6pt}
\[
\quad
\delta_{\fff i}
\off\off =\off\off
\prod_{\sigma\qff \in\qff N_p}\pff
i_{\fff \sigma}
\qff \circ\qff 
\Delta_{\trf \sigma\fff,\trf i}
\qff \circ\qff
p_{\trf \partial_{\fff i}\fff \sigma}
\hspace*{1.5em}\mbox{ and }\hspace*{1.5em}
\delta
\off\off =\off\off
\sum_{i\qff =\qff 0}^p\qff (\dff -\qff 1 \dff)^{\fff i}\qff \delta_{\fff i}\qff.
\]

\vspace{-6pt}
Each of the maps $\delta_{\fff i}$\nnsp,\oss as also the map $\delta$\nnsp,\oss
defines a chain map\vspace*{6pt}
\begin{equation*}
\quad
\begin{tikzcd}[column sep=large, row sep=huge]\dis
0 \arrow[r]
& 
C^{\fff p\dff -\dff 1}\fff(\dff N \dff) \arrow[r]
\arrow[d]
& 
C^{\fff p\dff -\dff 1}\fff(\dff N\fff,\pff B^{\fff 0} \dff)\arrow[r] 
\arrow[d]
&   
C^{\fff p\dff -\dff 1}\fff(\dff N\fff,\pff B^{\fff 1} \dff) \arrow[r] 
\arrow[d]
&
\off \ldots \off
\\
0 \arrow[r]
& 
C^{\fff p}\fff(\dff N \dff) \arrow[r]
& 
C^{\fff p}\fff(\dff N\fff,\pff B^{\fff 0} \dff) \arrow[r]
&   
C^{\fff p}\fff(\dff N\fff,\pff B^{\fff 1} \dff) \arrow[r]
&
\off \ldots \off,
\end{tikzcd}
\end{equation*}

\vspace{-3pt}
Let\oss 
$\Delta_{\trf \sigma}
\pff \colon\pff
B^{\fff q}\fff(\dff X \dff)
\ttoo
B^{\fff q}\fff(\dff |\dff \sigma \fff| \dff)
$\oss 
be the map induced by the inclusion\qss 
$|\fff \sigma \fff|
\ttoo
X$\nnsp,\oss
and\dss let\vspace*{6pt}
\[
\quad
\Delta
\off\off =\off\off
\prod_{\sigma\qff \in\qff N_0}\qff
\Delta_{\trf \sigma}
\]

\vspace{-6pt}
The map $\Delta$ defines a chain map\vspace*{6pt}
\begin{equation*}
\quad
\begin{tikzcd}[column sep=large, row sep=huge]\dis
0 \arrow[r]
& 
\rrr \arrow[r]
\arrow[d]
& 
B^{\fff 0}\fff(\dff X \dff) \arrow[r] 
\arrow[d]
&   
B^{\fff 1}\fff(\dff X \dff) \arrow[r] 
\arrow[d]
&
\off \ldots \off
\\
0 \arrow[r]
& 
C^{\fff 0}\fff(\dff N \dff) \arrow[r]
& 
C^{\fff 0}\fff(\dff N\fff,\pff B^{\fff 0} \dff) \arrow[r]
&   
C^{\fff 0}\fff(\dff N\fff,\pff B^{\fff 1} \dff) \arrow[r]
&
\off \ldots \off,
\end{tikzcd}
\end{equation*}

\vspace{-3pt}
All these chain maps together form a commutative diagram\vspace*{12pt}
\begin{equation}
\label{big-diagram}
\qquad
\begin{tikzcd}[column sep=large, row sep=huge]\dis
0 \arrow[r]
& 
\rrr \arrow[r]
\arrow[d]
& 
B^{\fff 0}\fff(\dff X \dff) \arrow[r] 
\arrow[d]
&   
B^{\fff 1}\fff(\dff X \dff) \arrow[r] 
\arrow[d]
&
\off \ldots \off
\\
0 \arrow[r]
& 
C^{\fff 0}\fff(\dff N \dff) \arrow[r]
\arrow[d]
& 
C^{\fff 0}\fff(\dff N\fff,\pff B^{\fff 0} \dff) \arrow[r]
\arrow[d]
&   
C^{\fff 0}\fff(\dff N\fff,\pff B^{\fff 1} \dff) \arrow[r]
\arrow[d]
&
\off \ldots \off
\\
0 \arrow[r]
& 
C^{\fff 1}\fff(\dff N \dff) \arrow[r]
\arrow[d]
& 
C^{\fff 1}\fff(\dff N\fff,\pff B^{\fff 0} \dff) \arrow[r]
\arrow[d]
&   
C^{\fff 1}\fff(\dff N\fff,\pff B^{\fff 1} \dff) \arrow[r]
\arrow[d]
&
\off \ldots \off
\\
0 \arrow[r]
& 
C^{\fff 2}\fff(\dff N \dff) \arrow[r]
\arrow[d]
& 
C^{\fff 2}\fff(\dff N\fff,\pff B^{\fff 0} \dff) \arrow[r]
\arrow[d]
&   
C^{\fff 2}\fff(\dff N\fff,\pff B^{\fff 1} \dff) \arrow[r]
\arrow[d]
&
\off \ldots \off
\\
\ &
\ldots &
\ldots &
\ldots &
\end{tikzcd}
\end{equation}

\vspace*{-3pt}
A standard computation shows that the columns of this diagram are complexes.\oss
In particular\halfff,\oss the family\qss
$C^{\fff p}\fff(\dff N\fff,\pff B^{\fff q} \dff)\dff_{p,\trf q\qff \geq\qff 0}$\qss
together with the differentials\qss $d_{\dff *}\fff,\pff \delta$\qss is a double complex.\oss

Let us construct the homomorphism\qss
$\widehat{H}^{\fff *}\fff(\dff X \dff) 
\ttoo
{H}^{\fff *}\fff(\dff |\fff N \fff| \dff)$\dnsp.\oss
Let $T$ be the total complex of the above double complex.\oss
By the definition,\oss\vspace*{5pt}
\[
\quad
T^{\dff n}
\off\off =\off\off
\bigoplus_{p\qff +\qff q\qff =\qff n}
C^{\fff p}\fff(\dff N\fff,\pff B^{\fff q} \dff)
\]

\vspace*{-8pt}
and the differential\qss
$d\dff \colon\dff T^{\dff n}\ttoo T^{\dff n\dff +\dff 1}$\qss
is defined by the formula\oss\vspace*{4pt}
\[
\quad
d\qff |\qff C^{\fff p}\fff(\dff N\fff,\pff B^{\fff q} \dff)
\off =\off 
d_{\dff q}\qff +\qff (\dff -\qff 1 \fff)^p\dff \delta \dff.\oss
\]

\vspace*{-8pt}
The diagram\qss (\ref{big-diagram})\qss leads to the maps\qss\vspace*{3pt}
\[
\quad
B^{\fff *}\fff(\dff X \dff)\ttoo T^{\dff *}
\hspace*{1.5em}\mbox{ and }\hspace*{1.5em}
C^{\fff *}\fff(\dff N \dff)\ttoo T^{\dff *},
\]

\vspace*{-9pt}
where\dss $B^{\fff *}\fff(\dff X \dff)$\dss and\dss $C^{\fff *}\fff(\dff N \dff)$\dss
are the top row and the left column of the diagram\qss (\ref{big-diagram})\qss
without\dss $\rrr$\nnsp.
The second of these two maps commutes with the differentials only up to a sign,\oss
but still induces a map in cohomology.\oss
The cohomology of the complex\dss $B^{\fff *}\fff(\dff X \dff)$\dss are nothing else
but the bounded cohomology\dss $\widehat{H}^{\fff *}\fff(\dff X \dff)$\dnsp.\oss
The cohomology of the complex\dss $C^{\fff *}\fff(\dff N \dff)$\dss are,\pss
by definition,\pss the cohomology of the simplicial complex $N$\nnsp.\oss
The latter are known to be equal to\dss ${H}^{\fff *}\fff(\dff |\fff N \fff| \dff)$\dnsp.\oss
Therefore,\oss there are two canonical maps\vspace*{4pt}
\[
\quad
\widehat{H}^{\fff *}\fff(\dff X \dff)\ttoo {H}^{\fff *}\fff(\dff T \dff)
\hspace*{1.5em}\mbox{ and }\hspace*{1.5em}
{H}^{\fff *}\fff(\dff |\fff N \fff| \dff)\ttoo {H}^{\fff *}\fff(\dff T \dff)\dff.
\]

\vspace{-8pt}
As we will see,\oss the second of these two maps is an isomorphism.\oss
Assuming this for a moment\halfff,\oss one can get the promised homomorphism\qss
$\widehat{H}^{\fff *}\fff(\dff X \dff)\ttoo {H}^{\fff *}\fff(\dff |\fff N \fff| \dff)$\qss
as the composition of the first map with the inverse of the second.\oss
By a well known theorem about double complexes\qss
(which is stated in Appendix\qss \ref{double-complexes}\qss for the convenience of the reader),\oss
in order to prove that the homomorphism\qss
${H}^{\fff *}\fff(\dff |\fff N \fff| \dff)\ttoo {H}^{\fff *}\fff(\dff T \dff)$\qss
is an isomorphism,\oss
it\dss is\sss sufficient\dss to\sss prove  
that\sss all\dss rows\sss of (\ref{big-diagram}),\oss
starting with the second one,\oss are exact\halfff.\oss
Each of these rows is the complex\qss (\ref{c-nerve})\qss for some\qss
$p\qff \geq\qff 0$\qss and\dss hence is a product\sss of several complexes of the form\qss (\ref{b-sigma}).\oss
For every $\sigma$ the fundamental group\sss $\pi_{\dff 1}\fff(\dff |\fff \sigma \fff| \dff)$\sss
is amenable and\dss hence\dss
$\widehat{H}^{\fff n}\fff(\dff |\fff \sigma \fff| \dff)
\qff =\qff
0$\qss
for\qss $n\qff >\qff 0$\dss
by Theorem\qss \ref{simply-connected-homology}.\oss
It follows that the complexes\qss (\ref{b-sigma})\qss are exact\halfff,\oss
and hence the rows of\qss (\ref{big-diagram})\qss are also exact\halfff.\oss

This completes the construction of the homomorphism\qss
$\widehat{H}^{\fff *}\fff(\dff X \dff)\ttoo {H}^{\fff *}\fff(\dff |\fff N \fff| \dff)$\dnsp.\oss
It remains to prove the commutativity of the diagram\qss (\ref{covering-triangle}).\oss
As the first step towards the proof of the commutativity,\oss
let us replace in the construction of the diagram\qss (\ref{big-diagram})\qss
the spaces\dss $B^{\fff q}\fff(\dff |\fff \sigma \fff| \dff)$\dss
by the spaces\dss $C^{\fff q}\fff(\dff |\fff \sigma \fff| \dff)$\dnsp,\oss
and the spaces\dss $B^{\fff q}\fff(\dff X \dff)$\dss
by the spaces\dss $C^{\fff q}\fff(\dff X \dff)$\dnsp.\oss
The new diagram has the same left column as the old one,\oss
its top row is the usual complex of the singular cochains of $X$\nnsp,\oss
and every term\dss
$C^{\fff p}\fff(\dff N\fff,\pff B^{\fff q} \dff)$\dss 
is replaced by the term\dss
$C^{\fff p}\fff(\dff N\fff,\pff C^{\fff q} \dff)$\dnsp.\oss 
There is an obvious canonical map from the diagram\qss (\ref{big-diagram})\qss
to the new diagram.\oss
In order to save some space,\oss
let us represent these two diagrams in the following 
self-explanatory manner\fff:\vspace*{6pt}
\begin{equation}
\label{two-diagrams}
\quad
\begin{tikzcd}[column sep=large, row sep=huge]\dis
& 
B^{\fff *}\fff(\dff X \dff)
\arrow[d]
\\
C^{\fff *}\fff(\dff N \dff) \arrow[r]
& 
C^{\fff *}\fff(\dff N\fff,\pff B^{\fff *} \dff)
\end{tikzcd}
\hspace*{1.5em}\mbox{ and }\hspace*{1.5em}
\begin{tikzcd}[column sep=large, row sep=huge]\dis
& 
C^{\fff *}\fff(\dff X \dff)
\arrow[d]
\\
C^{\fff *}\fff(\dff N \dff) \arrow[r]
& 
C^{\fff *}\fff(\dff N\fff,\pff C^{\fff *} \dff)\dff.
\end{tikzcd}
\end{equation}

\vspace*{-3pt}
Of course,\oss the cohomology of the complex\dss $C^{\fff *}\fff(\dff X \dff)$\dnsp,\oss
i.e.\qss of the top row of the second diagram,\oss
is nothing else but the usual singular cohomology\dss $H^{\fff *}\fff(\dff X \dff)$\dnsp.\oss
The other rows\sss $C^{\fff p}\fff(\dff N\fff,\pff C^{\fff *} \dff)$\nnsp,\qss $p\qff >\qff 0$\nnsp,\oss
of\dss the double complex\sss $C^{\fff *}\fff(\dff N\fff,\pff C^{\fff *} \dff)$\sss
are exact\sss only\dss if\trs the sets\sss $|\fff \sigma \fff|$\sss 
are acyclic with respect\dss to singular cohomology.\oss
But\dss would\sss be\sss the columns\sss 
$C^{\fff *}\fff(\dff N\fff,\pff C^{\fff q} \dff)$ 
of\dss $C^{\fff *}\fff(\dff N\fff,\pff C^{\fff *} \dff)$\sss
exact,\oss we could easily complete the proof\dss by referring to the same theorem
about double complexes with rows and columns interchanged.\oss
These columns are close\sss to being exact,\oss
but\sss in\sss general\sss are not\sss exact.\oss
Our next goal\sss is\sss to find a similar diagram with exact\sss columns.\oss

We will use\sss the sheaf\trs theory\dss to\sss this end.\oss
For every\qss $q\qff \geq\qff 0$\qss let\dss $\mathcal{C}^q$\dss be the presheaf on $X$
assigning to every open subset\qss $U\qff \subset\qff X$\qss the real vector space
of the real-valued cochains\dss $C^{\fff q}\fff(\dff U \dff)$\dss
and to every inclusion\qss $U\qff \subset\qff V$\qss
the restriction homomorphism\qss
$C^{\fff q}\fff(\dff V \dff)\ttoo C^{\fff q}\fff(\dff U \dff)$\dnsp.\oss
Let\dss $\Gamma^{\fff q}$\dss be the sheaf associated with the presheaf\dss $\mathcal{C}^q$\dnsp.\oss
The maps\qss
$d_{\dff q}
\dff \colon\dff
C^{\fff q}\fff(\dff U \dff)
\ttoo
C^{\fff q\dff +\dff 1}\fff(\dff V \dff)$\qss
define a morphism of\dss presheaves\qss
$d_{\dff q}
\dff \colon\dff
\mathcal{C}^q
\ttoo
\mathcal{C}^{q\dff +\dff 1}$\dnsp.\oss
Let\qss
$d_{\dff q}
\dff \colon\dff
\Gamma^{\fff q}
\ttoo
\Gamma^{\fff q\dff +\dff 1}$\qss
be the associated morphism of sheaves.\oss
Let us replace in the construction of the diagram\qss (\ref{big-diagram})\qss
the spaces\dss $B^{\fff q}\fff(\dff |\fff \sigma \fff| \dff)$\dss
by\dss $\Gamma^{\fff q}\fff(\dff |\fff \sigma \fff| \dff)$\dnsp,\oss
and\dss $B^{\fff q}\fff(\dff X \dff)$\dss
by\dss $\Gamma^{\fff q}\fff(\dff X \dff)$\dnsp.\oss
We need also use\dss
$d_{\dff q}
\dff \colon\dff
\Gamma^{\fff q}\dff(\dff |\fff \sigma \fff| \dff)
\ttoo
\Gamma^{\fff q\dff +\dff 1}\dff(\dff |\fff \sigma \fff| \dff)$\qss
instead of\trs
$d_{\dff q}
\dff \colon\dff
B^{\fff q}\fff(\dff |\fff \sigma \fff| \dff)
\ttoo
B^{\fff q\dff +\dff 1}\fff(\dff |\fff \sigma \fff| \dff)$\sss
and\dss the sheaf\dss restriction maps\qss\vspace{3pt}
\[
\quad
\Gamma^{\fff q}\dff(\dff |\dff \partial_{\fff i}\dff\sigma \dff| \dff)
\ttoo
\Gamma^{\fff q}\dff(\dff |\fff \sigma \fff| \dff)
\]

\vspace{-9pt}
instead of\dss
$\Delta_{\trf \sigma\fff,\trf i}$\nnsp.\oss
This leads to a new analogue\vspace{3pt}\vspace{-0.25pt}
\begin{equation}
\label{diagram-gamma}
\quad
\begin{tikzcd}[column sep=large, row sep=huge]\dis
& 
\Gamma^{\fff *}\fff(\dff X \dff)
\arrow[d]
\\
C^{\fff *}\fff(\dff N \dff) \arrow[r]
& 
C^{\fff *}\fff(\dff N\fff,\pff \Gamma^{\fff *} \dff)
\end{tikzcd}
\end{equation}

\vspace*{-6pt}
of the diagram\qss (\ref{big-diagram}).\oss
Let\dss $C^{\fff *}\fff(\dff N\fff,\pff \Gamma^{\fff q} \dff)$\dss
be the $q$\dnsp-th column of the this diagram without the term\dss 
$\Gamma^{\fff q}\fff(\dff X \dff)$\dnsp.\oss
By definition,\oss the cohomology of the complex\dss
$C^{\fff *}\fff(\dff N\fff,\pff \Gamma^{\fff q} \dff)$\dss
is nothing else but the cohomology\dss
$H^{\fff *}\fff(\trf \mathcal{U}\fff,\qff \Gamma^{\fff q}\trf)$\dss
of\trs the covering $\mathcal{U}$ with coefficients
in the sheaf $\Gamma^{\fff q}$\dnsp.\oss\vspace{6pt}

\textsc{Claim.}\off\oss
\emph{$H^{\trf 0}\dff(\trf \mathcal{U}\fff,\qff \Gamma^{\fff q}\trf)
\off =\off \Gamma^{\fff q}\fff(\trf X \trf)$
and\trs
$H^{\dff p}\dff(\trf \mathcal{U}\fff,\qff \Gamma^{\fff q}\trf)
\off =\off 0$\dss
for\dss $p\qff >\qff 0$\nnsp.\oss}

\subproof
The proof consists of references to the sheaf theory.\oss
The first\sss equality\sss follows,\oss
for example,\oss
from\qss \cite{go},\oss Chapter\qss II,\oss Theorem\qss 5.2.2,\oss
which applies because $\mathcal{U}$\sss is\sss either open,\pss or closed and\dss
locally\sss finite.\oss

Let\dss $\Phi$\sss be\sss the family of\dss all\sss closed subsets of\dss $X$\nnsp.\oss
Since\sss $X$\sss is\sss paracompact,\pss $\Phi$\sss is\dss a\sss paracompactifying\sss family.\oss
See\qss \cite{go},\qss Chapter\qss II,\oss Section\qss 3.2.\oss
Since\sss $\Phi$\sss includes all\sss closed subsets,\oss
the cohomology\sss groups\sss
$H^{\dff p}\dff(\trf \mathcal{U}\fff,\qff \Gamma^{\fff q}\trf)$
are equal\dss to\sss the cohomology\sss groups
$H^{\dff p}_{\dff \Phi}\dff(\trf \mathcal{U}\fff,\qff \Gamma^{\fff q}\trf)$
with\sss $\Phi$\sss being\dss the family\sss of\dss supports,\oss
and\dss the\sss results about\sss
$H^{\dff p}_{\dff \Phi}\dff(\trf \mathcal{U}\fff,\qff \Gamma^{\fff q}\trf)$\sss
apply\sss to\sss
$H^{\dff p}\dff(\trf \mathcal{U}\fff,\qff \Gamma^{\fff q}\trf)$\nnsp.\oss

Since $X$\sss is\sss paracompact,\oss
the sheaf\dss
$\Gamma^{\fff q}$ is\sss a soft\sss sheaf\halfff.\oss
See\qss \cite{go},\oss Chapter\qss II,\oss Example\qss 3.9.1.\oss
Moreover\halfff,\pss $\Gamma^{\fff q}$ is\sss a fine\sss sheaf\halfff.\oss
See\qss \cite{go},\oss Chapter\qss II,\oss Example\qss 3.9.1,\oss
the footnote on\dss p.\qss 161.\oss
Now\sss we can apply\qss Theorem\qss 5.2.3\qss from\qss \cite{go},\oss Chapter\qss II.\oss
If $\mathcal{U}$ is open,\oss then\sss the part\qss (b)\qss of\trs this theorem applies
because\sss $\Gamma^{\fff q}$ is\sss a\sss fine sheaf\halfff.\oss 
If $\mathcal{U}$ is closed and locally finite,\oss 
then\sss the part\qss (c)\qss of\trs this theorem applies
because\sss $\Gamma^{\fff q}$ is a soft\sss sheaf\halfff.\oss  
In both cases\sss this theorem implies that\sss the higher cohomology\sss groups,\pss
$H^{\dff p}\dff(\trf \mathcal{U}\fff,\qff \Gamma^{\fff q}\trf)$\dss
with\dss $p\qff >\qff 0$\nnsp,\oss vanish.\oss  \esubproof\vspace{6pt}

The above claim means,\oss in\dss particular\halfff,\oss
that\dss the columns\sss
$C^{\fff *}\fff(\trf N\fff,\pff \Gamma^{\fff q} \trf)$
of\trs the diagram\qss (\ref{diagram-gamma})\qss 
without\dss the term\sss 
$\Gamma^{\fff q}\fff(\trf X \trf)$
are indeed exact.\oss
In order\sss to relate\sss the double complex\sss
$C^{\fff *}\fff(\trf N\fff,\pff \Gamma^{\fff *} \trf)$\sss
with\sss the double complex\sss
$C^{\fff *}\fff(\trf N\fff,\pff C^{\fff *} \trf)$\nnsp,\oss
we need one more analogue of\dss
$C^{\fff *}\fff(\trf N\fff,\pff B^{\fff *} \trf)$\nnsp.\oss
Given a topological space $Y$\dnsp,\oss we can apply\sss the construction\sss
of\trs the sheaves\dss $\Gamma^{\dff q}$\sss
to $Y$ in the role of\dss $X$\sss and get\sss a sheaves 
$\Gamma^{\dff q}\off =\off \Gamma^{\dff q}_{\dff Y}$\sss on $Y$\dnsp.\oss
Let\sss 
$\gamma^{\dff q}\dff(\trf Y\trf)
\off =\off 
\Gamma^{\dff q}_{\dff Y}\dff(\trf Y\trf)$\sss
be\sss the space of\dss global sections of\dss the sheaf\dss
$\Gamma^{\dff q}_{\dff Y}$\nsp.\oss
If\trs $Y\qff \subset\qff X$\nnsp,\oss
then $\gamma^{\dff q}\dff(\trf Y\trf)$\sss is\sss in\sss general different\dss
from\sss $\Gamma^{\dff q}\dff(\trf Y\trf)$ because\sss
$\gamma^{\dff q}\dff(\trf Y\trf)$\sss is\sss defined\sss inside of\dss $Y$\dnsp,\oss
but\sss the construction of\dss $\Gamma^{\dff q}\dff(\trf Y\trf)$\sss involves
singular cochains on sets open in $X$\nnsp.\oss
By\dss the same reason\sss there are canonical\dss maps\sss
$C^{\dff q}\dff(\trf Y\trf )
\qff \ttoo\qff
\gamma^{\dff q}\dff(\trf Y\trf)$\nnsp,\oss
but\sss no natural\dss maps\sss
$C^{\dff q}\dff(\trf Y\trf )
\qff \ttoo\qff
\Gamma^{\dff q}\dff(\trf Y\trf)$
in\sss general.\oss

Let\dss us\sss replace in\sss the construction of\trs the diagram\qss (\ref{big-diagram})\qss
the spaces\sss $B^{\fff q}\fff(\dff |\fff \sigma \fff| \dff)$\sss
by the spaces\dss $\gamma^{\fff q}\fff(\dff |\fff \sigma \fff| \dff)$\dnsp,\oss
and\sss $B^{\fff q}\fff(\dff X \dff)$
by\dss
$\gamma^{\dff q}\dff(\trf X\trf)
\off =\off
\Gamma^{\fff q}\fff(\dff X \dff)$\dnsp.\oss
The maps\dss
$d_{\dff q}
\dff \colon\dff
\gamma^{\fff q}\dff(\dff |\fff \sigma \fff| \dff)
\ttoo
\gamma^{\fff q\dff +\dff 1}\dff(\dff |\fff \sigma \fff| \dff)$\dss
are defined as\sss before and\dss the role of\dss $\Delta_{\trf \sigma\fff,\trf i}$\sss 
is\dss played\sss by\dss the restriction maps\sss
$\gamma^{\fff q}\dff(\dff |\dff \partial_{\fff i}\dff\sigma \dff| \dff)
\ttoo
\gamma^{\fff q}\dff(\dff |\fff \sigma \fff| \dff)$\nnsp.\oss
This leads to the following\sss analogue 
of\trs the diagram\qss (\ref{big-diagram}).\oss\vspace{3pt}
\begin{equation}
\label{diagram-other-gamma}
\quad
\begin{tikzcd}[column sep=large, row sep=huge]\dis
& 
\gamma^{\fff *}\fff(\dff X \dff)
\arrow[d]
\\
C^{\fff *}\fff(\dff N \dff) \arrow[r]
& 
C^{\fff *}\fff(\dff N\fff,\pff \gamma^{\fff *} \dff)\qff.
\end{tikzcd}
\end{equation}

\vspace{-6pt}
The canonical\dss maps\sss
$C^{\dff q}\dff(\trf |\fff \sigma \fff| \trf)
\ttoo
\gamma^{\fff q}\dff(\trf |\fff \sigma \fff| \trf)$\sss
lead\sss to a morphism of\dss double complexes\vspace{3pt}
\[
\quad
C^{\fff *}\fff(\dff N\fff,\pff C^{\fff *} \dff)
\qff \ttoo\qff
C^{\fff *}\fff(\dff N\fff,\pff \gamma^{\fff *} \dff)
\pff
\]

\vspace{-9pt}
and\dss hence\sss to a map\sss from\sss the second diagram from\qss (\ref{two-diagrams})\qss
to\sss the diagram\qss (\ref{diagram-other-gamma}).\oss
If\trs $U\qff \subset\qff X$\sss is\sss an open subset,\oss 
then\sss
$\gamma^{\dff q}\dff(\trf U\trf)
\off =\off
\Gamma^{\dff q}\dff(\trf U\trf)$\nnsp.\oss
In\sss particular,\pss
$\gamma^{\dff q}\dff(\trf X\trf)
\off =\off
\Gamma^{\dff q}\dff(\trf X\trf)$\nnsp.\oss
Hence\sss
if\dss $\mathcal{U}$\sss is\sss an open covering,\oss
then\sss the diagrams\qss (\ref{diagram-other-gamma})\qss and\qss
(\ref{diagram-gamma})\qss are identical.\oss
If\trs $F\qff \subset\qff X$\sss is\sss a closed subset,\oss 
then\sss there are canonical\dss maps\sss
$\Gamma^{\dff q}\dff(\trf F\trf)
\qff \ttoo\qff
\gamma^{\dff q}\dff(\trf F\trf)$
commuting\sss with\sss the coboundary\sss operators\sss $d_{\dff q}$\nsp.\oss
These maps\sss lead\sss to a morphism of\dss double complexes\vspace{3pt}
\[
\quad
C^{\fff *}\fff(\dff N\fff,\pff \Gamma^{\fff *} \dff)
\qff \ttoo\qff
C^{\fff *}\fff(\dff N\fff,\pff \gamma^{\fff *} \dff)
\pff
\]

\vspace{-9pt}
and\dss hence\sss to a map\sss from\sss the diagram\qss (\ref{diagram-gamma})\qss
to\sss the diagram\qss (\ref{diagram-other-gamma}).\oss\vspace{-0pt}

Let\dss us combine\sss the diagrams\qss (\ref{two-diagrams}),\oss (\ref{diagram-gamma})\qss
and\qss (\ref{diagram-other-gamma})\qss into one as follows.\oss\vspace{3pt}
\begin{equation*}
\quad
\begin{tikzcd}[column sep=ttiny, row sep=hhuge]\dis
&
B^{\fff *}\fff(\dff X \dff) \arrow[rr]
\arrow[d]
&
&
C^{\fff *}\fff(\dff X \dff) \arrow[rr]
\arrow[d]
&
& 
\gamma^{\fff *}\fff(\dff X \dff) 
\arrow[d]
&
& 
\Gamma^{\fff *}\fff(\dff X \dff) \arrow[ll, "\dis ="']
\arrow[d]
\\
&
C^{\fff *}\fff(\dff N\fff,\pff B^{\fff *} \dff) \arrow[rr]
&
&
C^{\fff *}\fff(\dff N\fff,\pff C^{\fff *} \dff) \arrow[rr]
&
&
C^{\fff *}\fff(\dff N\fff,\pff \gamma^{\fff *} \dff) 
&
&
C^{\fff *}\fff(\dff N\fff,\pff \Gamma^{\fff *} \dff) \arrow[ll]
\\
C^{\fff *}\fff(\dff N \dff) \arrow[rr, "\dis ="] 
\arrow[ru]
& 
&
C^{\fff *}\fff(\dff N \dff) \arrow[rr, "\dis ="] 
\arrow[ru]
&
&
C^{\fff *}\fff(\dff N \dff)  
\arrow[ru]
&
&
C^{\fff *}\fff(\dff N \dff)\qff. \arrow[ll, "\dis ="']
\arrow[ru]
&
\end{tikzcd}
\end{equation*}

\vspace{-6pt}
Let\dss
$T^{\dff *}_{\dff C}$\nsp,\pss
$T^{\dff *}_{\dff \gamma}$\nsp,\oss 
and\dss $T^{\dff *}_{\dff \Gamma}$\dss
be\sss the\sss total\sss complexes of\trs the double complexes\vspace{1.5pt}
\[
\quad
C^{\dff *}\dff(\trf N\fff,\pff C^{\fff *} \trf)
\qff,\hspace*{1.5em}
C^{\dff *}\dff(\trf N\fff,\pff \gamma^{\fff *} \trf)\qff,
\hspace*{1.5em}\mbox{ and }\hspace*{1.5em}
C^{\dff *}\dff(\trf N\fff,\pff \gamma^{\fff *} \trf)
\]

\vspace{-10.5pt}
respectively.\oss
The previous diagram leads to the diagram\vspace{0pt}
\begin{equation*}
\quad
\begin{tikzcd}[column sep=ttiny, row sep=hhuge]\dis
&
B^{\fff *}\fff(\dff X \dff) \arrow[rr]
\arrow[d]
&
&
C^{\fff *}\fff(\dff X \dff) \arrow[rr]
\arrow[d]
&
& 
\gamma^{\fff *}\fff(\dff X \dff) 
\arrow[d]
&
& 
\Gamma^{\fff *}\fff(\dff X \dff) \arrow[ll, "\dis ="']
\arrow[d]
\\
&
T^{\dff *} \arrow[rr]
&
&
T^{\dff *}_{\dff C} \arrow[rr]
&
&
T^{\dff *}_{\dff \gamma} 
&
&
T^{\dff *}_{\dff \Gamma} \arrow[ll]
\\
C^{\fff *}\fff(\dff N \dff) \arrow[rr, "\dis ="] 
\arrow[ru]
& 
&
C^{\fff *}\fff(\dff N \dff) \arrow[rr, "\dis ="] 
\arrow[ru]
&
&
C^{\fff *}\fff(\dff N \dff)  
\arrow[ru]
&
&
C^{\fff *}\fff(\dff N \dff)\qff. \arrow[ll, "\dis ="']
\arrow[ru]
&
\end{tikzcd}
\end{equation*}

\vspace{-9pt}\vspace{6pt}

\textsc{Claim.}\off\oss
\emph{The map\dss
$T^{\dff *}_{\dff \Gamma}
\qff \ttoo\qff
T^{\dff *}_{\dff \gamma}$\dss
induces an\sss isomorphism\sss in cohomology.\oss}\vspace{-0.25pt}

\subproof
If\dss $\mathcal{U}$\sss is\sss an open,\pss
then\dss this map\dss is\dss the identity\sss map.\oss
Suppose\sss that\sss $\mathcal{U}$\sss is\sss closed and\dss locally\sss finite.\oss
By our assumptions,\oss in\sss this case\sss finite intersections\sss $F$\sss of\dss
elements of\dss $\mathcal{U}$ are homologically\dss locally\sss connected.\oss
This implies\sss that\sss 
$\Gamma^{\dff *}\dff(\trf F\trf)
\qff \ttoo\qff
\gamma^{\dff *}\dff(\trf F\trf)$\sss
induces isomorphisms in cohomology.\oss
See\qss \cite{bre},\oss Section\qss III.1,\pss
especially\dss the big\sss diagram\sss on\dss p.\qss 183.\oss
It\dss follows\sss that\dss the map of\dss double complexes\sss
$C^{\fff *}\fff(\dff N\fff,\pff \Gamma^{\fff *} \dff)
\qff \ttoo\qff
C^{\fff *}\fff(\dff N\fff,\pff \gamma^{\fff *} \dff)$\sss
induces an\sss isomorphism\sss in cohomology\sss groups of\dss each row.\oss
By\sss a well\dss known\sss theorem about\sss comparison of\dss double complexes\qss
(see\dss Appendix\qss \ref{double-complexes})\qss
this implies\sss that\dss
$T^{\dff *}_{\dff \Gamma}
\qff \ttoo\qff
T^{\dff *}_{\dff \gamma}$\dss
induces an\sss isomorphism\sss in cohomology.\oss  \esubproof\vspace{6pt}\vspace{-0.25pt}

In view of\trs this claim\sss the\sss last\sss diagram\dss leads\sss to\sss
the following\sss diagram\vspace{0pt}
\begin{equation*}
\quad
\begin{tikzcd}[column sep=ttiny, row sep=hhuge]\dis
&
\widehat{H}^{\fff *}\fff(\trf X \trf) \arrow[rr]
\arrow[d, "\dis j"']
&
&
H^{\fff *}\fff(\trf X \trf) \arrow[rr, "\dis k"]
\arrow[d]
&
& 
H^{\fff *}\fff(\trf \Gamma^{\fff *}\fff(\trf X \trf) \trf)
\arrow[d, "\dis j_{\dff \Gamma}"']
\\
&
H^{\fff *}\fff(\trf T^{\dff *} \trf) \arrow[rr]
&
&
H^{\fff *}\fff(\trf T^{\dff *}_{\dff C} \trf) \arrow[rr]
&
&
H^{\fff *}\fff(\trf T^{\dff *}_{\dff \Gamma} \trf)
\\
{H}^{\fff *}\fff(\trf |\fff N \fff| \trf) \arrow[rr, "\dis ="] 
\arrow[ru, "\dis i"']
& 
&
{H}^{\fff *}\fff(\trf |\fff N \fff| \trf) \arrow[rr, "\dis ="] 
\arrow[ru]
&
&
{H}^{\fff *}\fff(\trf |\fff N \fff| \trf) 
\arrow[ru, "\dis i_{\trf \Gamma}"']
&
\end{tikzcd}
\end{equation*}

\vspace{-6pt}
of cohomology groups,\pss
which is,\pss like all the previous ones,\pss obviously commutative.\oss
As we saw above,\oss the homomorphism\qss
${H}^{\fff *}\fff(\trf |\fff N \fff| \trf)
\qff \ttoo\qff 
{H}^{\fff *}\fff(\trf T^{\fff *} \trf)$\qss
is an isomorphism.\oss
Similarly,\oss since the columns of\qss (\ref{diagram-gamma})\qss are exact\halfff,\oss
a\sss well\dss known\dss theorem about\sss double complexes\qss (see Appendix\qss \ref{double-complexes})\qss
implies that\qss
$j_{\dff \Gamma}
\qff \colon\dff
H^{\fff *}\fff(\trf \Gamma^{\fff *}\fff(\trf X \trf) \trf)
\ttoo 
H^{\fff *}\fff(\trf T^{\dff *}_{\dff \Gamma} \trf)$\qss
is an isomorphism.\oss
Moreover\halfff,\oss\vspace{3pt}
\[
\quad
k
\qff \colon\qff
H^{\fff *}\fff(\trf X \trf)
\ttoo
H^{\fff *}\fff(\trf \Gamma^{\fff *}\fff(\trf X \trf) \trf)
\]

\vspace*{-9pt}
is also an isomorphism.\oss
See\qss \cite{sp},\oss Example\qss 6.7.9.\oss
Recall that the map\qss\vspace{3pt}
\[
\quad
\widehat{H}^{\fff *}\fff(\trf X \trf)
\qff \ttoo\qff 
{H}^{\fff *}\fff(\trf |\fff N \fff| \trf)
\]

\vspace{-9pt}
was defined as\qss
$i^{\dff -\dff 1}\fff \circ\dff j$\nnsp.\oss
Now we see that the map\oss\vspace{3pt}
\[
\quad
i^{\dff -\dff 1}\fff \circ\dff j
\qff \colon\qff
\widehat{H}^{\fff *}\fff(\trf X \trf)
\qff \ttoo\qff 
{H}^{\fff *}\fff(\trf |\fff N \fff| \trf) 
\]

\vspace{-9pt}
can be included into the commutative diagram\vspace{1.5pt}
\begin{equation*}
\quad
\begin{tikzcd}[column sep=normal, row sep=huge]\dis
\widehat{H}^{\fff *}\fff(\trf X \trf) \arrow[rr]
\arrow[rd, near start, "\dis i^{\dff -\dff 1}\fff \circ\dff j"']
&
& 
{H}^{\fff *}\fff(\trf X \trf) 
\\
&
{H}^{\fff *}\fff(\trf |\fff N \fff| \trf) 
\arrow[ru, near end, 
"\dis k^{\dff -\dff 1}\dff \circ\dff j_{\dff \gamma}^{\dff -\dff 1}\dff \circ\qff i_{\dff \gamma}"']
& 
\end{tikzcd}
\end{equation*}

\vspace{-10.5pt}
having the canonical map\qss
$\widehat{H}^{\fff *}\fff(\dff X \dff)
\ttoo
{H}^{\fff *}\fff(\dff X \dff)$\qss
as the horizontal arrow.\oss
This diagram has the promised form\qss (\ref{covering-triangle}).\oss
It remains only to settle the question if\trs the map\dss
$k^{\dff -\dff 1}\dff \circ\dff j_{\dff \gamma}^{\dff -\dff 1}\dff \circ\qff i_{\dff \gamma}$\dss
is the same as the canonical map\qss
${H}^{\fff *}\fff(\trf |\fff N \fff| \trf)
\ttoo
{H}^{\fff *}\fff(\trf X \trf)$\dnsp.\oss
From the point of\dss view of\trs the sheaf theory,\oss
the most natural way to define the canonical map\dss
${H}^{\fff *}\fff(\trf |\fff N \fff| \trf)
\ttoo
{H}^{\fff *}\fff(\trf X \trf)$\dss
is to define it as the map\dss
$k^{\dff -\dff 1}\dff \circ\dff j_{\dff \gamma}^{\dff -\dff 1}\dff \circ\qff i_{\dff \gamma}$\nnsp.\oss
If\trs this approach is adopted,\oss the question disappears.\oss
The readers who nevertheless prefer some other definition of\trs this map
may compare their preferred definition with this one.\oss
In any case,\oss this issue does not\dss belong to the bounded cohomology theory.\oss  \eproof

\mypar{Theorem\qss ({\dff}Vanishing\dss theorem).}{vanishing}
\emph{If\pss $X$\dss admits an open covering amenable in\dss the sense of\pss Gromov\qss
and\sss such that every point of\qss $X$\dss
is contained in no more that $m$ elements of this covering,\oss
then the canonical homomorphism\oss
$\widehat{H}^{\fff i}\fff(\dff X \dff)
\ttoo
{H}^{\fff i}\fff(\dff X \dff)$\oss
vanishes for\qss $i\qff \geq\qff m$\nnsp.\oss}

\proof\qss
The proof\dss is\dss similar\sss to\sss the proof\dss of\qss Theorem\qss \ref{covering-theorem}.\oss
One needs\sss to\sss replace\sss the\sss term\sss $\rrr$\sss in\sss the complex\qss
(\ref{b-sigma})\qss by\sss the product\sss of\dss copies of\dss $\rrr$\sss corresponding\sss
to path components of\dss $|\fff \sigma \fff|$\sss
and\dss then\sss replace\sss the spaces\sss $C^{\fff i}\fff(\dff N \dff)$\sss
in\qss (\ref{c-nerve})\qss by\sss some other spaces.\oss 
After\sss this\sss the cohomology\sss groups\dss $H^{\dff i}$\sss 
of\trs the\sss left\sss column of\trs the diagram\qss
(\ref{big-diagram})\qss may\sss be not\sss equal\dss to\sss the cohomology\sss groups
of\trs the nerve\sss $N$\nnsp,\oss but\sss still\sss vanish\sss in dimensions\dss
$>\qff \dim\dff N$\nnsp.\oss
The rest\sss of\trs the proof\dss applies and shows\sss that\dss the map\sss
$\widehat{H}^{\fff i}\fff(\dff X \dff)
\ttoo
{H}^{\fff i}\fff(\dff X \dff)$\sss
factors\sss throug\sss $H^{\fff i}$\nnsp.\oss
By\sss the assumption of\trs the\sss theorem\sss
$\dim\dff N\qff \leq\qff m\qff -\qff 1$\sss
and\dss hence\sss
$H^{\dff i}\off =\off 0$\sss
for\qss $i\qff \geq\qff m$\nnsp.\oss
The\sss theorem\dss follows.\oss \eproof

\myuppar{Remark.}
In\qss \cite{i1}\qss the author suggested that it may be
interesting\sss to\sss look for a pure group-the\-o\-ret\-ic version of\trs
Theorems\qss \ref{covering-theorem}\qss and\qss \ref{vanishing},\pss
without realizing\sss that\sss a part\sss of\dss proof\dss of\trs Theorem\qss \ref{class-c}\qss is\dss
such a version of\trs the special case of\dss Theorem\qss \ref{vanishing}\qss
corresponding\sss to amenable coverings by\dss two subsets.\oss
But\dss now\sss it\sss seems\sss that\trs
Theorems\qss \ref{covering-theorem}\qss and\qss \ref{vanishing}\qss show\sss that\dss 
the bounded cohomology\dss theory\sss does not\dss reduce\sss to\sss the
group\sss theory,\oss the results of\qss Section\qss \ref{spaces}\qss notwithstanding.\oss

\myuppar{Remark.}
It may happen that
the assumptions of\qss Vanishing Theorem hold\dss but\dss
$\widehat{H}^{\fff i}\fff(\dff X \dff)$\dss for some\sss $i\qff \geq\qff m$\sss
is\dss non-zero.\oss
For example,\oss if\dss $X$\sss is\dss a wedge of\trs two circles,\oss
then\dss
$\widehat{H}^{\fff 2}\fff(\dff X \dff)\qff \neq\qff 0$\nnsp.\oss
See\qss \cite{br}\qss or\qss \cite{gro}.\oss
At\dss the same\sss time\sss the obvious covering of $X$ by\sss two circles\sss
is\sss amenable.\oss

\myuppar{Remark.}
In\dss the case of\dss open coverings\sss 
the paracompactness assumptions\dss is\dss superfluous.\oss
See\qss \cite{i3}.\oss
This assumption\dss is\dss caused\dss by\sss references\sss 
to very\dss general\sss results of\trs the sheaf\trs theory.\oss\vspace{-0.125pt}

\myuppar{Small simplices.}
At the end of the proof of\dss Theorem\qss \ref{covering-theorem}\qss
we used the fact that the homomorphism\qss
$k
\dff \colon\dff
{H}^{\fff *}\fff(\dff X \dff)\ttoo H^{\fff *}\fff(\dff \gamma^{\fff *}\fff(\dff X \dff) \dff)$\qss
is an isomorphism.\oss
This reflects\sss the fact\dss that\dss the singular cohomology can
be computed\dss by\dss taking only small simplices into account\halfff.\oss
See\qss \cite{sp},\oss Theorem\qss 4.4.14\qss for a precise form of the latter claim.\oss
This is not the case for the bounded cohomology,\oss
and\dss the ideas of\trs the proof\dss of\trs Theorem\qss \ref{covering-theorem}\qss
allow\sss to make\sss this claim\sss precise.\oss
Let $\mathcal{U}$ be an amenable open covering of\dss $X$ and\sss let $N$ be its nerve.\oss
As in the proof of\dss Theorem\qss \ref{covering-theorem},\oss
we may assume that all finite intersections of the elements of $\mathcal{U}$ have
amenable fundamental groups.\oss
Let\dss $\mathcal{B}^q$\dss be the analogue for bounded cochains 
of\trs the presheaf\dss $\mathcal{C}^q$\dss
and let\dss $\beta^{\fff q}$\dss be the sheaf 
associated with the presheaf\dss $\mathcal{B}^q$\nnsp.\oss
Let\sss $C^{\fff *}\fff(\trf N\fff,\pff \beta^{\fff *} \trf)$
be the bounded analogue of\qss 
$C^{\fff *}\fff(\trf N\fff,\pff \Gamma^{\fff *} \trf)$\nnsp,\oss 
and\dss let\dss $T_{\dff \beta}$\dss be\sss the\sss total complex of\trs 
$C^{\dff *}\dff(\trf N\fff,\pff \beta^{\fff *} \trf)$\nnsp.\oss
Since\sss $X$\sss is\dss paracompact,\pss
the sheaf\trs $\beta^{\fff q}$\sss is\dss fine\dss 
by the same reason as $\Gamma^{\fff q}$\dnsp.\oss
It follows that the canonical map\qss\vspace{3pt}\vspace{-1pt}
\[
\quad
j_{\dff \beta}
\qff \colon\dff
H^{\fff *}\fff(\dff \beta^{\fff *}\fff(\dff X \dff) \dff)
\ttoo 
H^{\fff *}\fff(\dff T^{\dff *}_{\dff \beta} \dff)
\]

\vspace{-9pt}\vspace{-1pt}
is an isomorphism.\oss
Suppose now that the bounded cohomology can be computed by using only
the small simplices.\oss
Then for every open\qss $Y\qff \subset\qff X$\qss
the canonical map\qss\vspace{3pt}\vspace{-1pt}
\[
\quad
\widehat{H}^{\fff *}\fff(\dff Y \dff)
\ttoo
H^{\fff *}\fff(\dff \beta^{\fff *}\fff(\dff Y \dff) \dff)
\]

\vspace{-9pt}\vspace{-1pt}
is an isomorphism and\dss therefore\dss 
$H^{\fff i}\fff(\dff \beta^{\fff *}\fff(\dff |\fff \sigma \fff| \dff) \dff)
\off =\off
\widehat{H}^{\fff i}\fff(\dff |\fff \sigma \fff| \dff)
\off =\off
0$\sss
for every\qss $i\qff \geq\qff 1$\qss and every simplex $\sigma$\nnsp.\oss
It follows that\sss all\dss rows of\dss 
$C^{\fff *}\fff(\trf N\fff,\pff \beta^{\fff *} \trf)$\nnsp,\oss
starting with the second one,\oss are exact
and hence the canonical map\qss
$H^{\fff *}\fff(\dff |\fff N \fff| \dff)
\ttoo
H^{\fff *}\fff(\dff T^{\fff *}_{\dff \beta} \dff)$\qss
is an isomorphism.\oss
Since the maps\qss\vspace{3pt}\vspace{-1pt}
\[
\quad
H^{\fff i}\fff(\dff \beta^{\fff *}\fff(\dff X \dff) \dff)
\ttoo
H^{\fff *}\fff(\dff T^{\fff *}_{\dff \beta} \dff)
\hspace*{1.5em}\mbox{ and }\hspace*{2em}
\widehat{H}^{\fff *}\fff(\dff X \dff)
\ttoo
H^{\fff *}\fff(\dff \beta^{\fff *}\fff(\dff X \dff) \dff)
\]

\vspace{-9pt}\vspace{-1pt}
are also isomorphisms,\oss
it follows that\qss
$H^{\fff *}\fff(\dff |\fff N \fff| \dff)$\qss
is canonically isomorphic to\qss
$\widehat{H}^{\fff *}\fff(\dff X \dff)$\dnsp.\oss
Since the condition\qss({\fff}C{\fff})\qss is very weak,\oss
it is easy to see that this cannot be true.\oss
It\sss is\sss also not\sss difficult\dss to give an example of\dss a covering $\mathcal{U}$
for which\sss
$H^{\fff *}\fff(\trf |\fff N \fff| \trf)$\sss
is\sss not\dss isomorphic to\qss
$\widehat{H}^{\fff *}\fff(\trf X \trf)$\dnsp.\oss
In fact\halfff,\oss the covering of the wedge of two circles by these two circles
is such an example.\oss

\newpage
\mysection{An\qss algebraic\qss analogue\qss of\qss the\qss mapping\qss theorem}{a-mapping-theorem}

\vspace*{6pt}
\myuppar{The $G$\dnsp-modules $F\dff(\fff G^{\fff n\dff +\dff 1},\pff U \dff)$\dnsp.}
Let $U$ be a\dss $G$\dnsp-module and\qss $n\qff \geq\qff 0$\nnsp.\oss 
Let\dss $F\dff(\fff G^{\fff n\dff +\dff 1},\pff U\dff)$\dss
be equal to\qss $B\dff(\fff G^{\fff n\dff +\dff 1},\pff U\dff)$\dss
as a Banach space,\oss but with the action $\bullet$ of $G$ defined by\vspace*{3pt} 
\[
\quad
(\fff h\bullet\nsp c \dff)\dff(\fff g_{\fff 0}\fff,\pff g_{\fff 1}\fff,\pff \ldots\fff,\pff g_{\fff n}\fff)
\off =\off
h\cdot\left(\dff
c\dff(\fff g_{\fff 0}\halfff h\fff,\pff g_{\fff 1}\halfff h\fff,\pff \ldots\fff,\pff g_{\fff n}\halfff h\fff)
\dff\right).
\]

\vspace{-9pt}
The map\qss 
$s
\dff \colon\dff
G^{\fff n\dff +\dff 1}\ttoo G^{\fff n\dff +\dff 1}$\qss 
defined by\vspace*{3pt} 
\[
\quad
s\dff(\dff g_{\dff 0}\fff,\pff 
g_{\dff 1}\fff,\pff 
\ldots\fff,\pff 
g_{\dff n}\dff)
\off = \off
(\dff g_{\dff n}\fff,\pff 
g_{\dff n\dff -\dff 1}\fff g_{\dff n}\fff,\pff 
\ldots\fff,\pff 
g_{\dff 0}\fff g_{\dff 1}\fff \ldots\fff g_{\dff n}\dff)
\]

\vspace{-9pt}
is an equivariant bijection from $G^{\fff n\dff +\dff 1}$ with the right action\vspace*{3pt} 
\[
\quad
(\dff g_{\dff 0}\fff,\pff 
g_{\dff 1}\fff,\pff 
\ldots\fff,\pff 
g_{\dff n}\dff)\cdot h
\off =\off
(\dff g_{\dff 0}\fff,\pff 
g_{\dff 1}\fff,\pff 
\ldots\fff,\pff 
g_{\dff n}\halfff h\dff)
\]

\vspace{-9pt}
of\dss $G$ to $G^{\fff n\dff +\dff 1}$ with the right action\vspace*{3pt} 
\[
\quad
(\dff g_{\dff 0}\fff,\pff 
g_{\dff 1}\fff,\pff 
\ldots\fff,\pff 
g_{\dff n}\dff)\cdot h
\off =\off
(\dff g_{\dff 0}\halfff h\fff,\pff 
g_{\dff 1}\halfff h\fff,\pff 
\ldots\fff,\pff 
g_{\dff n}\halfff h\dff)
\]

\vspace{-9pt}
of\dss $G$\nnsp.\oss
It follows that the map\vspace*{3pt} 
\[
\quad
s^{\fff *}
\qff \colon\qff
F\dff(\fff G^{\fff n\dff +\dff 1},\pff U\dff)
\qff \ttoo
B\dff(\fff G^{\fff n\dff +\dff 1},\pff U\dff)
\]

\vspace{-9pt}
defined by\qss $s^{\fff *}(\dff c \dff)\off =\off c\dff \circ\dff s$\qss
is an isomorphism of\dss $G$\dnsp-modules.\oss
Obviously,\pss $s^{\fff *}(\dff c \dff)\off =\off f$\nnsp,\oss where\vspace*{3pt} 
\[
\quad
f\dff (\dff g_{\dff 0}\fff,\pff 
g_{\dff 1}\fff,\pff 
\ldots\fff,\pff 
g_{\dff n}\dff)
\off =\off
c\dff
(\dff g_{\dff n}\fff,\pff 
g_{\dff n\dff -\dff 1}\fff g_{\dff n}\fff,\pff 
\ldots\fff,\pff 
g_{\dff 0}\fff g_{\dff 1}\fff \ldots\fff g_{\dff n}\dff)\dff,
\]

\vspace{-6pt}
\myuppar{The homogeneous form of the standard resolution.}
Let us consider the sequence\vspace*{6pt}
\begin{equation*}
\quad
\begin{tikzcd}[column sep=large, row sep=normal]\dis
0 \arrow[r]
& 
U \arrow[r, "\dis d_{\dff -\dff 1}\off"]
& 
F\dff(\fff G\fff,\pff U \dff) \arrow[r, "\dis d_{\dff 0}\off"]
&   
F\dff(\fff G^{\dff 2},\pff U \dff) \arrow[r, "\dis d_{\dff 1}\off"]
&
F\dff(\fff G^{\dff 3},\pff U \dff) \arrow[r, "\dis d_{\dff 2}\off"]
&
\off \ldots \off,
\end{tikzcd}
\end{equation*}

\vspace*{-3pt}
where\oss  
$d_{\dff -\dff 1}\dff(\fff v \fff)(\fff g\fff)\off =\off v$\oss 
for\dss all\qss
$v\qff \in\qff U$\dnsp,\qss $g\qff \in\qff G$\nnsp,\off\oss
and\dss\vspace*{6pt}
\[
\quad
d_{\dff n}\dff(\fff c \fff)
(\fff g_{\fff 0}\fff,\pff g_{\fff 1}\fff,\pff \ldots\fff,\pff g_{\fff n\dff +\dff 1}\fff)
\off =\off
\sum_{i\qff =\qff 0}^{n\dff +\dff 1}\qff (-\qff 1)^{\dff i}\dff 
c\dff(\dff g_{\fff 0}\fff,\pff \ldots\fff,\pff
\widehat{g_{\dff i}}\dff,\pff \ldots\fff,\pff   g_{\dff n\dff +\dff 1} \dff)\dff,
\]

\vspace*{-6pt}
for all\qss $n\qff \geq\qff 0$\qss
and\oss $g_{\fff 0}\fff,\pff g_{\fff 1}\fff,\pff \ldots\fff,\pff g_{\fff n\dff +\dff 1}
\qff \in\qff G$\nnsp.\off\oss
A direct verification shows that the maps $s^{\fff*}$ define an isomorphism
of the above sequence with the standard resolution\qss (\ref{st-res})\qss
of $U$\dnsp.\oss

Given\qss $n\qff \geq\qff 0$\qss and\qss $i\off =\off 1\fff,\pff 2\fff,\pff \ldots\fff,\pff n$\nnsp,\oss
let
\[
\quad
\pr_{\fff i\fff,\qff G}
\qff \colon\qff 
G^n
\qff \ttoo\qff  
G^{n\dff -\dff 1}
\]
be the projection of $G^n$ onto the product of all factors of\oss
$G^n\off =\off G\dff \times\dff G\dff \times\dff \ldots\dff \times\dff G$\oss
except the $i$\dnsp-th one.\oss
In terms of these projections the differential $d_{\dff n}$ takes the form\vspace*{6pt}
\[
\quad
d_{\dff n}\dff(\fff c \fff)\dff
(\dff g\dff)
\off =\off
\sum_{i\qff =\qff 0}^{n\dff +\dff 1}\qff (-\qff 1)^{\dff i}\dff 
c\dff \circ\trf \pr_{\fff i\fff,\qff G}\dff
(\dff g\dff)\dff,
\]

\vspace*{-6pt}
where\qss $n\qff \geq\qff 0$\qss
and\oss 
$g
\qff \in\qff 
G^{\fff n\dff +\dff 2}$\nnsp.\off\oss

\myuppar{Means of vector-valued functions.}
Let $U$ be a Banach space.\oss
Given a set $S$\dnsp,\oss we denote by $B\fff(\fff S\fff, U\fff)$ the vector space of 
all maps\qss $f\dff \colon\dff S\ttoo U$\qss such that the function\qss
$s\qff \longmapsto\qff \|\qff f\dff(s \fff) \qff\|$\qss
from $S$ to $\rrr$ is bounded.\oss
The vector space\dss $B\fff(\fff S\fff, U\fff)$\dss is a Banach space with the norm\qss \vspace*{3pt}
\[
\quad
\|\dff f\trf\|
\off =\off\sup_{\fff x\qff \in\qff S}\qff 
\|\qff f(x) \qff\|\qff.
\]

\vspace*{-9pt} 
A\qss \emph{mean}\qss on\dss $B\fff(\fff S\fff, U\fff)$\dss is
defined as a linear map\vspace*{3pt}
\[
\quad
m
\qff \colon\qff
B\fff(\fff S\fff, U\fff)
\qff \ttoo\qff
U
\]

\vspace*{-9pt}
such that\qss $\|\qff m \qff\|\qff \leq\qff 1$\qss and\qss
$m\dff(\dff \const_{\fff u} \fff)\off =\off u$\qss
for all\qss $u\qff \in\qff U$\dnsp,\oss
where $\const_{\fff u}$ is the constant function with the value $u$\nnsp,\oss
i.e.\qss $\const_{\fff u}\fff(s \fff)\off =\off u$\qss for all\qss $s\qff \in\qff S$\dnsp.\oss
For the trivial\dss $G$\dnsp-module\qss $U\qff =\qff \rrr$\qss 
this definition reduces to the definition of a mean on $B\dff(\dff S\fff)$\dnsp.\oss

Suppose that a group $G$ acts on a set $S$ from the right\halfff.\oss 
Then $G$ acts on the left on the space\dss $B\fff(\fff S\fff, U\fff)$\dss 
by the formula\qss 
$g\cdot f\dff (\fff s\fff )\qff =\qff f\dff (\fff s\cdot g\fff )$\dnsp,\oss 
where\qss $g\qff \in\qff G$\nnsp,\qss $f\qff \in\qff B\fff(\fff S\fff, U\fff)$\nnsp,\qss 
$s\qff \in\qff S$\nnsp.\oss 
The mean $m$ on\dss $B\fff(\fff S\fff, U\fff)$\dss is called\qss \emph{right invariant}\pss if\qss 
$m\fff(\fff g\cdot f\fff)\qff =\qff m\fff(\fff f\fff)$\qss 
for all\qss 
$g\qff \in\qff G$\nnsp,\qss 
$f\qff \in\qff B\fff(\fff S\fff, U\fff)$\dnsp.\oss 
As usual,\pss we  call the right invariant means simply\qss \emph{invariant means},\oss
and $G$ is always considered together with its action on itself by the right translations.\oss

Suppose now $G$ acts on a set $S$ on the left
and\dss this action is free and transitive.\oss
Then every\qss $s\qff \in\qff S$\qss defines a bijection\qss
$r_s\dff \colon\dff G\toto S$\qss by the rule\qss
$g\qff \longmapsto\qff g\cdot s$\nnsp,\oss where\qss $g\qff \in\qff G$\nnsp.\oss
If\qss $s\fff,\pff t\qff \in\qff S$\nnsp,\oss
then the bijections\dss $r_s$\dss and\dss $r_t$\dss differ by a right translation of $G$\nnsp.\oss 
Cf.\qss Section\qss \ref{amenable}.\oss
If $m$ is a mean on\dss $B\fff(\fff G\fff,\pff U\fff)$\dss and\qss $s\qff \in\qff S$\nnsp,\oss
then\qss
$f\qff \longmapsto\qff m\dff(\fff f\dff \circ\dff r_s\fff)$\qss
is a mean on\dss $B\fff(\fff S\fff,\pff U\fff)$\dnsp.\oss
By the same reason as in Section\qss \ref{amenable}\qss the mean
$f\qff \longmapsto\qff m\dff(\fff f\dff \circ\dff r_s\fff)$\qss
on\dss $B\fff(\fff S\fff,\pff U\fff)$\dss
does not depend on the choice of $s$\nnsp.\oss
This mean is said to be\qss 
\emph{induced by $m$ and the action of\qss $G$ on $S$\dnsp.}\oss
The\qss \emph{push-forwards}\qss of\dss invariant means by surjective homomorphisms
are defined in exactly the same way as in Section\qss \ref{amenable}.\oss

The following construction goes back at least to\dss B.E.\qss Johnson\qss \cite{j}.\oss
See\qss \cite{j},\oss Theorem\qss 2.5\qss and its proof\halfff.\oss
The key idea is to work with the\qss \emph{dual}\qss Banach spaces.\oss

\myuppar{A construction of\dss means of\dss vector-valued functions.}
Let\dss $A$\dss be a set and\dss let\qss 
$U\off =\off V^{\fff *}$\dss
be a Banach space dual to some other Banach space $V$\dnsp.\oss
A function\qss
$f\fff \in\qff B\fff(\dff A\fff,\pff U\fff)$\qss
can be considered as a function\qss
$f
\dff \colon\dff
A\dff \times\dff V\ttoo \rrr$\qss
linear by the second argument and such that\vspace*{3pt}
\begin{equation*}
\quad
|\qff f\dff(\fff a\fff,\pff v \fff) \qff|
\off \leq\off
c\dff \|\qff v \qff\|
\end{equation*}

\vspace*{-9pt}
for some\qss $c\qff \in\qff \rrr$\qss and all\qss $a\qff \in\qff A$\nnsp.\oss
Then\dss $\|\qff f \qff\|$\dss
is equal to the infimum of $c$ such that 
this inequality holds.\oss 
In particular\halfff,\pss
if\qss $a\qff \in\qff A$\qss and\qss $v\qff \in\qff V$\dnsp,\oss
then\vspace*{3pt}
\begin{equation}
\label{norm-dual}
\quad
|\qff f\dff(\fff a\fff,\pff v \fff) \qff|
\off \leq\off
\|\qff f \qff \|\cdot \|\qff v \qff\|
\end{equation}

\vspace*{-9pt}
Suppose that\qss  
$M\dff \colon\dff B\dff(\dff A\dff)\ttoo \rrr$\qss 
is a mean.\oss
Given a vector\qss $v\qff \in\qff V$\dnsp,\oss
let us define a function\qss $f_v\dff \colon\dff A\ttoo \rrr$\qss
by the rule\qss $f_v\fff(a \fff)\off =\off f\dff(\fff a\fff,\pff v\fff)$\dnsp.\oss
By\qss (\ref{norm-dual})\qss 
the function\dss
$f_v$\trs is bounded\dss and,\pss
more\-over\halfff,\pss
$\|\qff f_v \qff\|\qff \leq\qff \|\qff f \qff\|\cdot\|\qff v \qff\|$\nnsp.\oss
Let us define a map\qss
$m\fff(\dff f \dff)
\dff \colon\dff
V\ttoo \rrr$\qss by\vspace*{3pt}
\[
\quad
m\fff(\dff f \dff)\fff (v)
\off =\off
M\dff(\dff f_v\fff)\dff.
\]

\vspace*{-9pt}
Since $f\dff(a\fff,\pff \bullet \dff)$ is a linear functional,\pss
$f_v$ linearly depends on $v$\nnsp.\oss
Since $M$ is a linear functional,\oss this implies that\dss
$m\fff(\dff f \dff)$\dss is a linear functional on $V$\dnsp.\oss
It is bounded of the norm\qss $\leq\qff \|\qff f \qff\|$\qss because\vspace*{3pt}
\[
\quad
\|\qff m\fff(\dff f \dff)\fff (v\fff) \qff\|
\off \leq\off
\|\qff f_v \qff\|
\off \leq\off
\|\qff f \qff\|\cdot\|\qff v \qff\|\dff.
\]

\vspace*{-9pt}
It follows that\qss $m\fff(\dff f \dff)\qff \in\qff V^{\fff *}\off =\off U$\qss
and\qss
$\|\qff m\fff(\dff f \dff) \qff\|
\qff \leq\qff
\|\qff f \qff\|$\nnsp.\oss

\mypar{Lemma.}{vector-means}
\emph{The map\dss
$m\dff \colon\dff f\qff \longmapsto\qff m\dff(\dff f \dff)$\dss
is a mean\dss 
$B\fff(\dff A\fff,\pff U\fff)\ttoo U$\dnsp.\oss
If\oss $L\dff \colon\dff V\ttoo V$\oss is a bounded operator
and\oss $L^*\dff \colon\dff V^{\fff *}\qff =\qff U\ttoo U$\oss
is its adjoint operator\halfff,\oss
then}\qss\vspace*{3pt}
\[
\quad
m\fff(\qff L^*\fff \circ\fff f \dff)
\off =\off 
L^*(\dff  m\fff(\dff f \dff) \dff)\dff.
\]

\vspace{-9pt}
\emph{If\dss $A$\dss is a group and\dss $M$\dss is an invariant mean,\oss
then\dss $m$\dss is an invariant mean also.\oss}

\proof\qss
Since $f_v$\dss linearly depends on $f$ and $M$ is\dss linear\halfff,\pss
$m$\dss
is linear\halfff.\oss
By the inequality\qss
$\|\qff m\fff(\dff f \dff) \qff\|
\qff \leq\qff
\|\qff f \qff\|$\qss
its norm\qss $\leq\qff 1$\nnsp.\oss
Suppose that\oss 
$f\dff(\fff a \dff)
\off =\off 
u\qff \in\qff U\off =\off V^{\fff *}$\oss 
for all\qss $a\qff \in\qff A$\nnsp.\oss
Then\qss $f_v\dff(\fff a \dff)\off =\off u\dff(v\fff)$\qss
for all\qss $a\qff \in\qff A$\qss and\qss $v\qff \in\qff V$\dnsp,\oss
i.e.\qss all\dss $f_v$\dss are constant functions.\oss
It follows that\qss
$m\fff(\dff f \dff)\fff (v\fff)\off =\off u\dff(v\fff)$\qss
for all\qss $v\qff \in\qff V$\qss and\dss hence\qss 
$m\fff(\dff f \dff)\off =\off u$\nnsp.\oss
Therefore $m$ is indeed a mean.\oss
If\qss $a\qff \in\qff A$\qss and\qss $v\qff \in\qff V$\dnsp,\oss
then\qss
$L^*\fff \circ\fff f\dff(a\fff,\pff v\fff)
\off =\off
f\dff(a\fff,\pff L\dff(v\fff) \dff)$\dnsp.\oss
Therefore\qss
$(\qff L^*\fff \circ\fff f \dff)_v
\off =\off
f_{\qff L\dff(v\fff)}$\qss
and\vspace*{3pt}
\[
\quad
m\fff(\qff L^*\fff \circ\fff f \dff)\fff (v)
\off =\off
M\fff\left(\dff (\qff L^*\fff \circ\fff f \dff)_v \dff\right)
\off =\off
M\fff\left(\dff f_{\qff L\dff(v\fff)} \dff\right)
\]

\vspace*{-36pt}
\[
\quad
\phantom{m\fff(\qff L^*\fff \circ\fff f \dff)\fff (v)
\off =\off
M\fff\left(\dff (\qff L^*\fff \circ\fff f \dff)_v \dff\right)
\off }
=\off
m\fff(\dff f \dff)\fff (\qff L\dff(v\fff)\fff)
\off =\off
L^*(\dff  m\fff(\dff f \dff) \dff)\fff (v\fff)
\]

\vspace{-9pt}
It follows that\oss
$m\fff(\qff L^*\fff \circ\fff f \dff)
\off =\off 
L^*(\dff  m\fff(\dff f \dff) \dff)$\dnsp.\oss
The last statement of\dss the lemma is obvious.\oss  \eproof

\myuppar{Dual modules.}
Suppose that $A$ is a group and $V$ is a Banach space and a\qss \emph{right}\qss $A$\dnsp-module.\oss
Then\qss $U\off =\off V^{\fff *}$\qss has a canonical structure of a left
$A$\dnsp-module.\oss
Namely\halfff,\oss if\qss $u\dff \colon\dff V\ttoo \rrr$\qss is a bounded linear functional
and\qss $a\qff \in\qff A$\nnsp,\oss
then\dss $a\cdot u$\dss is the linear functional\qss
$v\qff \longmapsto\qff u\fff(\fff v\cdot a\dff)$\dnsp.\oss

For every\qss $a\qff \in\qff A$\qss
let\qss 
$L\dff(\fff a\fff)
\dff \colon\dff
V\ttoo V$\qss be the operator\qss
$v\qff \longmapsto\qff v\fff\cdot\fff a$\qss
of\dss the right multiplication by $a$ in $V$\dnsp.\oss
Then the adjoint operator\qss
$L\dff(\fff a\dff)^{*}
\dff \colon\dff
U\ttoo U$\qss
is the operator\qss
$u\qff \longmapsto\qff a\fff\cdot\fff u$\qss
of the left multiplication by $a$\dss in $U$\dnsp.\oss

\mypar{Lemma.}{mean-equivariance}
\emph{In the situation of\dss Lemma\qss \ref{vector-means},\oss
suppose that $U$ is a right $A$\dnsp-module.\oss
Then\qss
$m
\dff \colon\dff
B\fff(\dff A\fff,\pff U \dff)\ttoo U$\qss
is a morphism of\qss $A$\dnsp-modules.\oss}

\proof\qss
If\qss $a\qff \in\qff A$\qss
and\qss
$f\qff \in\qff B\fff(\fff A\fff,\pff U\fff)$\dnsp,\oss
then by definition\qss
$a\bullet f
\off =\off
L\dff(a\fff)^{*}\fff \circ\fff f\dff \circ\dff r_a$\nnsp,\oss
where\qss
$r_a
\dff \colon\dff
A\ttoo A$\qss
is the right shift\qss
$b\qff \longmapsto\qff b\halfff a$\qss
by $a$ in $A$\nnsp.\oss
By\dss Lemma\qss \ref{vector-means}\dss $m$ is an invariant mean and\dss hence\qss
$m\fff(\dff f\dff \circ\dff r_a \dff)
\off =\off
m\fff(\dff f \dff)$\dnsp.\oss
Moreover\halfff,\oss Lemma\qss \ref{vector-means}\qss implies that\vspace*{6pt}
\[
\quad
m\fff(\dff a\bullet f \dff)
\off =\off
m\fff(\dff L\dff(a\fff)^{*}\fff \circ\fff f\dff \circ\dff r_a \dff)
\]

\vspace*{-33pt}
\[
\quad
\phantom{m\fff(\dff a\bullet f \dff)
\off }
=\off
L\dff(a\fff)^{*}(\dff  m\fff(\dff f\dff \circ\dff r_a \dff) \dff)
\off =\off
L\dff(a\fff)^{*}(\dff  m\fff(\dff f \dff) \dff)
\off =\off
a\cdot  m\fff(\dff f \dff)\dff.
\]

\vspace*{-6pt}
It follows that $m$ is indeed a morphism of $A$\dnsp-modules.\oss  \eproof

\myuppar{Quotients by normal subgroups.}
Let $\Gamma$ be a discrete group,\pss
$A$ be a normal subgroup of\dss $\Gamma$\nnsp,\oss
and $U$ be a\dss $\Gamma$\nsp\dnsp-module.\oss
Let\qss $G\off =\off \Gamma/A$\qss and\qss
$\alpha\dff \colon\dff \Gamma\ttoo G$\qss
be the quotient map.\oss
Then $\alpha$ turns
$B\fff(\dff G\fff,\pff U \dff)$ into a $\Gamma$\nsp\dnsp-module,\oss
which can be described as follows.\oss
Let us consider $G$ as the group of\dss the right cosets.\oss
Then\dss $\Gamma$ acts on $G$ from the right by the rule\qss
$g\cdot h\off =\off g\halfff h$\nnsp,\oss
where\qss $h\qff \in\qff \Gamma$\qss and\qss $g\qff =\qff A\halfff \gamma$\qss 
for some\qss $\gamma\qff \in\qff \Gamma$\qss
is a right coset\halfff,\oss
and acts on $B\fff(\dff G\fff,\pff U \dff)$ from the left by the rule\vspace*{3pt}
\[
\quad
h\bullet f\dff(\dff A\halfff \gamma \dff)
\off =\off
h\cdot\left(\dff
f\dff(\dff A\halfff \gamma\halfff h \dff)
\dff\right)
\] 

\vspace*{-9pt}
The obvious map\oss
$\alpha^{\fff *}
\qff \colon\qff 
B\dff(\dff G\fff,\pff U \dff)
\qff \ttoo\qff 
B\dff(\dff \Gamma,\pff U \dff)$\oss
induced by\dss
$\alpha$\dss
is\dss a\dss $\Gamma$\dnsp-morphism.\oss

\myuppar{Quotients by normal amenable subgroups.}
Let $\Gamma$ be a group
and $V$ be a Banach space and a\qss \emph{right}\dss $\Gamma$\nsp\dnsp-module.\oss
Then\qss $U\off =\off V^{\fff *}$\qss is a left\dss $\Gamma$\nsp\dnsp-module.\oss
Suppose that $A$ is a normal amenable subgroup of\dss $\Gamma$\dnsp.\oss
Then $U$ is also an $A$\dnsp-module in a canonical way\halfff.\oss
Let\qss
$M
\dff \colon\dff
B\dff(\dff A \dff)\ttoo \rrr$\qss be an invariant mean 
and\dss let\qss
$m
\dff \colon\dff
B\dff(\dff A\fff,\pff U \dff)\ttoo \rrr$\qss
be the defined\dss by $M$ invariant mean on\dss $B\dff(\dff A\fff,\pff U \dff)$\dnsp.\oss
As before,\oss
let\qss $G\off =\off \Gamma/A$\qss and\qss
$\alpha\dff \colon\dff \Gamma\ttoo G$\qss
be the quotient map.\oss
The group $A$ acts freely and transitively from the left
on each right coset $A\halfff \gamma$\nnsp,\oss 
where\qss $\gamma\qff \in\qff \Gamma$\nnsp.\oss 
Given\qss $\gamma\qff \in\qff \Gamma$\dnsp,\oss
let the mean\qss\vspace*{3pt}
\[
\quad 
m_{\dff \gamma}
\dff \colon\dff
B\dff(\dff A\halfff \gamma\fff,\pff U \dff)
\ttoo
U
\]

\vspace{-9pt}
be defined by\qss
$m_{\dff \gamma}\dff(\dff f \dff)
\off =\off
m\dff(\dff f\dff \circ\dff r_{\fff \gamma} \dff)$\dnsp,\oss
where\qss
$r_{\fff \gamma}
\dff \colon\dff
A
\ttoo
A\halfff \gamma$\qss
is the map\qss
$a\qff \longmapsto\qff a\fff \gamma$\dnsp.\oss
The mean $m_{\dff \gamma}$ is nothing else but the mean on\dss
$B\dff(\dff A\halfff \gamma\fff,\pff U \dff)$\dss
induced by $m$ and the left action of $A$ on $A\halfff \gamma$\dnsp,\oss
which is free and transitive.\oss
In particular\halfff,\pss $m_{\dff \gamma}$ 
depends only on the coset $A\halfff \gamma$\dnsp,\oss
but not on the choice of a representative $\gamma$ of this coset\halfff.\oss
Cf.\qss the beginning of Section\qss \ref{amenable}.\oss
Let\qss\vspace*{3pt} 
\[
\quad
m_{\fff *}
\dff \colon\dff
B\dff(\dff \Gamma,\pff U \dff)
\ttoo
B\dff(\dff G\fff,\pff U \dff)
\]

\vspace{-9pt}
be the map defined by the formula\oss\vspace*{3pt} 
\[
\quad
m_{\fff *}\dff(\dff f \dff)(\dff A\halfff \gamma \dff)
\off =\off 
m_{\dff \gamma}\dff(\dff f \mid A\halfff \gamma \dff)\dff,
\]

\vspace{-9pt}
where $f \mid A\halfff \gamma$ 
is the restriction of $f$ to $A\halfff \gamma$\dnsp.\pss
It is well defined because
$m_{\dff \gamma}$ depends only on $A\halfff \gamma$\dnsp.

\mypar{Lemma.}{coset-means}
\emph{The map\dss $m_{\fff *}$\dss is\dss a\qss $\Gamma$\nsp\dnsp-morphism,\pss
$m_{\fff *}\dff \circ\dff \alpha^{\fff *}
\off =\off
\id$\nnsp,\oss
and\oss
$\|\qff m_{\fff *}\qff\|
\qff \leq\qff 1$\nnsp.\oss}

\proof\trs
Let us prove first that $m_{\fff *}$ commutes with the\sss $\Gamma$\nsp\dnsp-actions.\pss
Let\trs $h\trf \in\trf \Gamma$\trs
and\trs
$f\trf \in\qff B\fff(\fff \Gamma,\pff U\fff)$\dnsp.
Then by the definition of\dss the\dss $\Gamma$\nsp\dnsp-action\qss\vspace*{3pt}
\[
\quad
h\bullet f
\off =\off
L\dff(\halfff h\fff)^{*}\fff \circ\fff f\dff \circ\dff r_{\fff h}\dff,
\]

\vspace*{-9pt}
where\qss
$r_{\fff h}
\dff \colon\dff
\Gamma\ttoo \Gamma$\qss
is the right shift\qss
$\delta\qff \longmapsto\qff \delta\halfff h$\qss
by $h$ in $\Gamma$\dnsp.\oss

Let us consider an arbitrary\qss $\gamma\qff \in\qff \Gamma$\qss
and\dss prove that the values of\qss
$m_{\fff *}\fff(\dff h\bullet f \dff)$\qss
and\qss
$h\bullet m_{\fff *}\fff(\dff f \dff)$\qss
on the coset\qss $A\halfff \gamma$\qss are equal.\oss
The value of\qss
$m_{\fff *}\fff(\dff h\bullet f \dff)$\qss 
is\vspace*{5pt}
\[
\quad
m_{\fff *}\fff
(\dff h\bullet f \dff)
\fff
(\dff A\halfff \gamma \dff)
\off =\off
m_{\fff *}\fff
\left(\dff 
L\dff(\halfff h\fff)^{*}\fff \circ\fff f\dff \circ\dff r_{\fff h} 
\dff\right)
\fff
(\dff A\halfff \gamma \dff)
\]

\vspace*{-33pt}
\[
\quad
\phantom{m_{\fff *}\fff
(\dff h\bullet f \dff)
\fff
(\dff A\halfff \gamma \dff)
\off }
=\off
m_{\dff \gamma}\dff
\left(\dff
L\dff(\halfff h\fff)^{*}\fff \circ\fff f\dff \circ\dff r_{\fff h}
\qff \mid\qff
A\halfff \gamma
\qff\right)
\]

\vspace*{-33pt}
\[
\quad
\phantom{m_{\fff *}\fff
(\dff h\bullet f \dff)
\fff
(\dff A\halfff \gamma \dff)
\off }
=\off
m\dff
\left(\dff
L\dff(\halfff h\fff)^{*}\fff \circ\qff \varphi
\dff\right)\dff,
\]

\vspace*{-8pt}
where\qss $\varphi\dff \colon\dff A\ttoo U$\qss
is the map\qss
$a\qff \longmapsto\qff f\dff \circ\dff r_{\fff h}\dff(\dff a\halfff \gamma \dff)
\off =\off
f\dff(\dff a\halfff \gamma\fff h \dff)$\dnsp,\oss
and\dss therefore\vspace*{3pt}
\[
\quad
m_{\fff *}\fff
(\dff h\bullet f \dff)
\fff
(\dff A\halfff \gamma \dff)
\off =\off
L\dff(\halfff h\fff)^{*}\left(\dff  m\fff(\dff \varphi \dff) \dff\right)
\]

\vspace*{-9pt}
in view of\dss Lemma\qss \ref{vector-means}.\oss
On the other hand\vspace*{5pt}
\[
\quad
h\bullet m_{\fff *}\fff(\dff f \dff)\fff(\dff A\halfff \gamma \dff)
\off =\off
L\dff(\halfff h\fff)^{*}
\left(\dff m_{\fff *}\fff
(\dff f \dff)\fff(\dff A\halfff \gamma\halfff h \dff)
\dff\right)
\]

\vspace*{-33pt}
\[
\quad
\phantom{h\bullet m_{\fff *}\fff(\dff f \dff)\fff(\dff A\halfff \gamma \dff)
\off }
=\off
L\dff(\halfff h\fff)^{*}\left(\dff m_{\dff \gamma\halfff h}\fff
(\dff f 
\qff \mid\qff 
A\halfff \gamma\halfff h \dff)
\dff\right)
\]

\vspace*{-33pt}
\[
\quad
\phantom{h\bullet m_{\fff *}\fff(\dff f \dff)\fff(\dff A\halfff \gamma \dff)
\off }
=\off
L\dff(\halfff h\fff)^{*}
\left(\dff m
\fff(\dff \psi \dff) 
\dff\right)\dff,
\]

\vspace{-9pt}
where\qss $\psi\dff \colon\dff A\ttoo U$\qss
is the map\qss
$a\qff \longmapsto\qff f\dff \circ\dff r_{\gamma\halfff h}\dff(\dff a\dff)
\off =\off
f\dff(\dff a\halfff \gamma\halfff h \dff)$\dnsp,\oss
and\dss therefore\vspace*{3pt}
\[
\quad
h\bullet m_{\fff *}\fff(\dff f \dff)\fff(\dff A\halfff \gamma \dff)
\off =\off
L\dff(\halfff h\fff)^{*}
\left(\dff m
\fff(\dff \psi \dff) 
\dff\right)\dff.
\]

\vspace{-9pt}
Since,\oss obviously\halfff,\pss $\varphi\off =\off \psi$\nnsp,\oss it follows that\vspace*{3pt}
\[
\quad
m_{\fff *}\fff
(\dff h\bullet f \dff)
\fff
(\dff A\halfff \gamma \dff)
\off =\off
h\bullet m_{\fff *}\fff(\dff f \dff)\fff(\dff A\halfff \gamma \dff)
\]

\vspace*{-9pt}
for all\qss $\gamma\qff \in\qff \Gamma$\qss
and\dss hence\oss
$m_{\fff *}\fff
(\dff h\bullet f \dff)
\off =\off
h\bullet m_{\fff *}\fff(\dff f \dff)$\dnsp.\oss
This proves that $m_{\fff *}$ is a\dss $\Gamma$\nsp\dnsp-morphism.\oss
The last two claims of the lemma immediately 
hold\dss because all
$m_{\dff \gamma}$ are means.\oss  \eproof

\myuppar{Invariant means on\dss $B\fff(\dff A^n \dff)$\dnsp.}
We need an analogue of the coherent systems 
of invariant means from Section\qss \ref{amenable}.\oss
Let $A$ be an amenable group.\oss
Then $A^n$ are amenable for all $n$
by\qss \cite{gr},\oss Theorem\qss 1.2.6.\oss
Recall\dss that\vspace*{3pt}
\[
\quad
\pr_{\fff i\fff,\qff A}
\qff \colon\qff 
A^n
\qff \ttoo\qff  
A^{n\dff -\dff 1}
\]

\vspace*{-9pt}
is the projection onto the product of\dss all factors of\dss
$A^n$\dss
except the $i$\dnsp-th one.\oss
A sequence\qss\vspace*{3pt} 
\[
\quad
M_{\fff 0}\fff,\off M_{\fff 1}\fff,\off \ldots\fff,\off M_{\fff n}\fff,\off \ldots
\]

\vspace*{-9pt}
of\dss invariant means\qss
$M_{\fff n}\qff \in\qff B\fff(\dff A^n \dff)^*$\dnsp,\oss
where\oss
$n\off =\off 0\fff,\pff 1\fff,\pff 2\fff,\pff \ldots \off$\dnsp,\oss
is called\qss \emph{coherent}\pss if\vspace*{4pt}
\begin{equation}
\label{coherence-powers}
\quad
\left(\dff 
\pr_{\fff i\fff,\qff A}
\dff\right)_* M_{\fff n}
\off =\off
M_{\fff n\trf -\dff 1}
\end{equation}

\vspace*{-8pt} 
for every\qss $n\fff,\pff i$\qss as above.\oss
Equivalently,\oss  
$M_{\fff n}\fff(\dff f\dff \circ\qff \pr_{\fff i\fff,\qff A}\dff)
\off =\off
M_{\fff n\trf -\dff 1}\fff(\dff f\dff)
$\oss 
for every\qss $n\fff,\pff i$\qss as above.

\myuppar{\dnsp$\Sigma_{\dff n}$\nsp\dnsp-invariant means on\dss $B\fff(\dff A^n \dff)$\dnsp.}
The symmetric group $\Sigma_{\dff n}$ acts on $A^n$ by permuting the factors.\oss
This action induces actions of $\Sigma_{\dff n}$ on\dss $B\fff(\dff A^n \dff)^*$\dss
and on the set of invariant means on $B\fff(\dff A^n \dff)$\nnsp.\oss
We will say that an invariant mean\qss 
$M\qff \in\qff B\fff(\dff A^n \dff)^*$\dss
is a\dss 
\emph{$\Sigma_{\dff n}$\nsp\dnsp-invariant mean}\pss 
if\qss 
$M$ is fixed by the action of\trs 
$\Sigma_{\dff n}$\nnsp.\oss
If\qss $M\qff \in\qff B\fff(\dff A^n \dff)^*$\qss is an invariant mean,\oss 
then\vspace*{4pt}
\[
\quad
\frac{1}{\dff n\dff! \dff}\off
\sum\nolimits_{\qff \sigma\qff \in\qff \Sigma_{\dff n}}\dff 
M^{\dff \sigma}
\]

\vspace*{-8pt}
is\dss a\dss $\Sigma_{\dff n}$\nsp\dnsp-invariant mean,\oss 
where the action of\dss $\Sigma_{\dff n}$\dss is written as 
$(\dff M\fff,\pff \sigma\fff)\qff\longmapsto\qff M^{\dff \sigma}$\dnsp.\oss
Suppose that\dss $M_{\fff n}$\dss is\dss a\dss 
$\Sigma_{\fff n}$\nsp\dnsp-invariant mean on\qss 
$B\fff(\dff  A^n \dff)$\qss
and that\oss
$1\qff \leq\qff i\qff \leq\qff n$\nnsp.\oss
Let\vspace*{4pt} 
\[
\quad
M_{\fff n\dff -\dff 1}
\off =\qff\off 
\left(\dff \pr_{\fff i\fff,\qff A}\dff\right)_*\fff M_{\fff n}\dff.
\] 

\vspace*{-8pt}
Since\dss $M_{\fff n}$\dss is\dss 
$\Sigma_{\fff n}$\dnsp-invariant\halfff,\pss
$M_{\fff n\dff -\dff 1}$ is independent of the choice of\dss 
$i$\dss 
and\dss is\dss 
$\Sigma_{\fff n\dff -\dff 1}$\dnsp-invariant\halfff.\oss 
Since $M_{\fff n\dff -\dff 1}$ is independent of the choice of $i$\nnsp,\oss
the condition\qss (\ref{coherence-powers})\qss 
holds for every\dss 
$i$\nnsp.

\vspace{3pt}
\mypar{Theorem.}{coherent-exists-a}
\emph{There exists a coherent sequence of\dss invariant means\oss
$M_{\fff 0}\fff,\pff M_{\fff 1}\fff,\pff M_{\fff 2}\fff,\pff 
\ldots\off $\nnsp.}

\vspace{3pt}
\proof\qss 
The proof\dss is completely similar to the proof\dss of\dss 
Theorem\qss \ref{coherent-exists-infinite}.\oss  \eproof

\vspace{3pt}
\myuppar{Invariant means on\dss $B\fff(\dff A^n,\pff U \dff)$\dnsp.}
Let $U$ be a Banach space.\oss
A sequence\qss\vspace*{3pt} 
\[
\quad
m_{\fff 0}\fff,\off m_{\fff 1}\fff,\off \ldots\fff,\off m_{\fff n}\fff,\off \ldots
\]

\vspace*{-9pt}
of\dss invariant means\qss
$m_{\fff n}\dff \colon\dff B\fff(\dff A^n,\pff U \dff)\ttoo U$\dnsp,\oss
where\oss
$n\off =\off 0\fff,\pff 1\fff,\pff 2\fff,\pff \ldots \off$\dnsp,\oss
is called\qss \emph{coherent}\pss if\vspace*{4pt}
\begin{equation}
\label{coherence-powers-coeff}
\quad
\left(\dff 
\pr_{\fff i\fff,\qff A}
\dff\right)_* m_{\fff n}
\off =\off
m_{\fff n\trf -\dff 1}
\end{equation}

\vspace*{-8pt} 
for every\qss $n\fff,\pff i$\qss as above.\oss
Equivalently,\oss  
$m_{\fff n}\fff(\dff f\dff \circ\qff \pr_{\fff i\fff,\qff A}\dff)
\off =\off
m_{\fff n\trf -\dff 1}\fff(\dff f\dff)
$\oss 
for every\qss $n\fff,\pff i$\qss as above.

\vspace{3pt}
\mypar{Theorem.}{coherent-exists-a-u}
\emph{If\pss $U\off =\off V^{\fff *}$\qss
is the dual to a Banach space $V$\dnsp,\oss
then there exists an infinite coherent sequence\qss
$m_{\fff 0}\fff,\pff m_{\fff 1}\fff,\pff m_{\fff 2}\fff,\pff 
\ldots\off $\nnsp,\oss
where\dss $m_{\fff n}$ is an invariant mean\oss
$B\fff(\dff A^n,\pff U \dff)\ttoo U$\dnsp.\oss}

\proof\qss
The construction of an invariant mean\qss
$m$ on $B\fff(\dff A^n,\pff U \dff)$\qss
by an invariant mean\qss
$M$ on $B\fff(\dff A^n \dff)$\qss
commutes with push-forwards.\oss
Hence the theorem 
follows from 
Theorem\qss \ref{coherent-exists-a}.\oss  \eproof

\vspace{3pt}
\mypar{Theorem.}{chain-map-quotient}
\emph{Suppose that\dss $A$\dss is a normal amenable subgroup of\dss
a group\trs $\Gamma$\dss
and\dss that\trs $U$\dss is\dss a\dss left\qss $\Gamma$\nsp\dnsp-module 
dual to a right\qss $\Gamma$\nsp\dnsp-module.\oss
Let\oss $G\off =\off \Gamma/A$\oss
and\dss let\oss 
$\alpha\dff \colon\dff \Gamma\ttoo G$\oss 
be the quotient map.\oss 
Then there exists\dss a\qss $\Gamma$\nsp\dnsp-morphism}\vspace*{2pt} 
\[
\quad
\alpha_{\fff *}
\qff \colon\qff 
F\dff(\dff \Gamma^{\fff \bullet\dff +\dff 1},\pff U \dff)
\qff \ttoo\qff 
F\dff(\dff G^{\fff \bullet\dff +\dff 1},\pff U\dff)
\]

\vspace*{-10pt}
\emph{of\qss $\Gamma$\nsp\dnsp-resolutions
of\dss the\qss $\Gamma$\nsp\dnsp-module\dss $U$\dss
such that\oss 
$\alpha_{\fff *}\dff \circ\dff \alpha^{\fff *}
\off =\dff\off 
\id$\oss
and\oss
$\|\dff \alpha_{\fff *} \dff\|\off \leq\off 1$\nnsp.\oss}

\vspace{3pt}
\proof\qss
The proof follows the lines of\dss the proof\dss of\dss Theorem\qss \ref{chain-map}.\oss

Let\dss $\{\dff m_{\fff n} \dff \}$\dss be a 
coherent sequence of\dss invariant means\qss $B\fff(\dff A^n,\pff U \dff)\ttoo U$\dnsp.\oss
For every $n$ the group 
$A^n$ is a normal subgroup of\dss $\Gamma^{\fff n}$\dss and\dss $\alpha$\dss
induces an isomorphism\vspace*{4pt} 
\[
\quad
\alpha^{\fff n}
\qff \colon\qff
\Gamma^{\fff n}/A^n
\qff \ttoo\qff
G^{\fff n}.
\]

\vspace*{-8pt}
By applying\dss Lemma\qss \ref{coset-means}\qss to\qss $\Gamma^{\fff n}$\nsp\dnsp,\pss
$A^n$\nsp\dnsp,\pss $m_{\fff n}$\qss in the roles of\qss $\Gamma$\dnsp,\pss
$A$\nnsp,\pss $m$\qss respectively\halfff,\oss
we get maps\vspace*{4pt} 
\[
\quad
m_{\fff n\dff *}
\dff \colon\dff
B\dff(\dff \Gamma^{\fff n},\pff U \dff)
\ttoo
B\dff(\dff G^{\fff n},\pff U \dff)
\]

\vspace{-9pt}
such that\oss
$m_{\fff n\dff *}\dff \circ\pff \alpha^{\fff n\dff *}
\off =\off
\id$\nnsp,\oss
$\|\qff m_{\fff n\dff *}\qff\|
\qff \leq\qff 1$\nnsp,\oss
and each\dss $m_{\fff n\dff *}$\dss is a\dss $\Gamma^{\fff n}$\nsp\dnsp-morphism.\oss

The map\dss $m_{\fff n\dff *}$\dss can be also considered as a map\vspace*{3pt} 
\[
\quad
m_{\fff n\dff *}
\dff \colon\dff
F\dff(\dff \Gamma^{\fff n},\pff U \dff)
\ttoo
F\dff(\dff G^{\fff n},\pff U \dff)
\]

\vspace{-9pt}
with the same properties.\oss
The structure of\dss $\Gamma$\nsp\dnsp-modules on vector spaces\dss 
$F\dff(\dff \Gamma^{\fff n},\pff U \dff)\fff,\pff
F\dff(\dff G^{\fff n},\pff U \dff)$\dss
is induced from the structure of\dss
$\Gamma^{\fff n}$\nsp\dnsp-modules
by the diagonal\dss homomorphism\vspace*{3pt} 
\[
\quad
h
\off \longmapsto\off
(\dff h\fff,\pff h\fff,\pff \ldots\fff,\pff h \dff)
\qff \in\qff
\Gamma^{\fff n}\dff,
\]

\vspace{-9pt}
where\qss $h\qff \in\qff \Gamma$\dnsp.\oss
It\dss follows that the maps\dss $m_{\fff n\dff *}$\dss are morphisms of\dss $\Gamma$\nsp\dnsp-modules.\oss
Let us prove that the maps\dss $m_{\fff n\dff *}$\dss
commute with the differentials of the homogeneous standard resolutions.\oss

As the first step,\oss let us prove that\vspace*{6pt}
\begin{equation}
\label{inner-coherence}
\quad
m_{\fff n\dff +\dff 1\dff *}\fff 
\left(\dff
c\dff \circ\dff \pr_{\fff i\fff,\qff \Gamma}
\dff\right)
\off =\off
m_{\fff n\dff *}\fff 
(\dff c \dff)
\dff \circ\qff
\pr_{\fff i\fff,\qff G}\qff.
\end{equation}

\vspace*{-6pt}
for every\qss
$c\qff \in\qff F\dff(\dff \Gamma^{\fff n},\pff U \dff)$\qss
and\qss $0\qff \leq\qff i\qff \leq\qff n$\nnsp.\off\oss
Let us evaluate both sides of\qss (\ref{inner-coherence})\qss on an arbitrary\vspace*{4pt}
\[
\quad
g
\qff \in\qff
G^{\fff n\dff +\dff 1}
\off =\off
\Gamma^{\fff n\dff +\dff 1}/A^{n\dff +\dff 1}\dff.
\]

\vspace*{-8pt}
Let\qss
$\gamma
\qff \in\qff
\Gamma^{\fff n\dff +\dff 1}$
be some representative of\dss the coset $g$\nnsp,\oss
and\dss let\qss\vspace*{5pt}
\[
\quad
g\fff(\fff i\fff)
\off =\off
\pr_{\fff i\fff,\qff G}\dff (\dff g \dff)
\hspace*{1em}\mbox{ and }\hspace*{1em}
\gamma\dff(\fff i\fff)
\off =\off
\pr_{\fff i\fff,\qff \Gamma}\dff (\dff \gamma \dff)\dff.
\]

\vspace*{-7pt}
Then $\gamma\dff(\fff i\fff)$
is a representative of the coset $g\fff(\fff i\fff)$\dnsp.\oss
By the definition of\dss push-forwards
the left\dss hand side of\dss (\ref{inner-coherence})\dss 
is equal to $m_{\fff n\dff +\dff 1}\fff(\dff \varphi\dff)$\dnsp,\oss
where\qss
$\varphi\dff \colon\dff A^{n\dff +\dff 1}
\ttoo U$\qss
is the map\vspace{4pt}
\begin{equation}
\label{phi-function}
\quad
\varphi\dff(\dff {a}\dff)
\off =\off 
c\dff \circ\dff \pr_{\fff i\fff,\qff \Gamma}\qff
(\fff a\fff \gamma \dff)
\off =\off
c\dff(\qff \pr_{\fff i\fff,\qff A}\qff
(\dff a \dff)\dff \gamma\dff(\fff i\fff) \dff)\dff.
\end{equation}

\vspace{-8pt}
The right hand side
is equal to\qss $m_{\fff n}\fff(\dff \psi\dff)$\dnsp,\oss
where\qss
$\psi\dff \colon\dff A^{n}
\ttoo U$\qss
is the map\vspace*{4pt}
\begin{equation}
\label{psi-function}
\quad
\psi\dff(\dff b\dff)
\off =\off
c\dff(\dff 
b\dff \gamma\dff(\fff i\fff) \dff)\dff.
\end{equation}

\vspace{-8pt}
By comparing\qss (\ref{phi-function})\qss and\qss (\ref{psi-function})\qss
we see that\oss
$\varphi
\off =\off
\psi\qff \circ\qff
\pr_{\fff i\fff,\qff A}$\oss
and\dss hence\vspace*{5pt}
\[
\quad
m_{\fff n\dff +\dff 1}\fff(\dff \varphi\dff)
\off =\off
\Bigl(\dff
\left(\dff
\pr_{\fff i\fff,\qff A}
\dff\right)_* m_{\fff n}
\dff\Bigr)\dff (\dff \psi \dff)
\]

\vspace*{-7pt}
In view of the coherence condition\qss (\ref{coherence-powers-coeff})\qss
this implies that the values of both sides of\qss (\ref{inner-coherence})\qss
on $g$ are equal.\oss
Since $g$ is arbitrary\halfff,\oss
this proves\qss (\ref{inner-coherence}).\oss

Let us prove that\qss 
$m_{\fff n\dff +\dff 1\dff *}\dff \circ\dff d_{\dff n\dff -\dff 1}
\off =\off
d_{\dff n\dff -\dff 1}\dff \circ\dff m_{\fff n\dff *}$\nnsp.\oss
Let\qss $c\qff \in\qff F\dff(\dff \Gamma^{\fff n},\pff U \dff)$\qss 
and\qss
${g}
\qff \in\qff
G^{\fff n\dff +\dff 1}$\nnsp.\oss 
Then\vspace*{6pt} 
\[
\quad
m_{\fff n\dff +\dff 1\dff *}\fff(\dff d_{\dff n\dff -\dff 1}\dff(\fff c \fff) \dff)\dff
(\dff {g}\dff)
\off\off =\off\off
m_{\fff n\dff +\dff 1\dff *}\qff
\left(\qff 
\sum_{i\qff =\qff 0}^{n}\qff (-\qff 1)^{\dff i}\dff 
c\dff \circ\trf \pr_{\fff i\fff,\qff G}
\qff\right)
\dff
(\dff {g}\dff) 
\]

\vspace*{-24pt}
\[
\quad
\phantom{m_{\fff n\dff +\dff 1\dff *}\fff(\dff d_{\dff n\dff -\dff 1}\dff(\fff c \fff) \dff)\dff
(\dff {g}\dff)
\off\off }
=\off\off 
\sum_{i\qff =\qff 0}^{n}\qff (-\qff 1)^{\dff i}\dff
m_{\fff n\dff +\dff 1\dff *}\dff
\left(\dff 
c\dff \circ\trf \pr_{\fff i\fff,\qff G}
\dff\right)
\dff
(\dff {g}\dff) 
\]

\vspace*{-24pt}
\[
\quad
\phantom{m_{\fff n\dff +\dff 1\dff *}\fff(\dff d_{\dff n\dff -\dff 1}\dff(\fff c \fff) \dff)\dff
(\dff {g}\dff)
\off\off }
=\off\off 
\sum_{i\qff =\qff 0}^{n}\qff (-\qff 1)^{\dff i}\dff
m_{\fff n\dff *}\fff 
(\dff c \dff)\dff
(\qff {g}\fff(\fff i\fff) \qff) 
\]

\vspace*{-24pt}
\[
\quad
\phantom{m_{\fff n\dff +\dff 1\dff *}\fff(\dff d_{\dff n\dff -\dff 1}\dff(\fff c \fff) \dff)\dff
(\dff {g}\dff)
\off\off }
=\off\off 
\sum_{i\qff =\qff 0}^{n}\qff (-\qff 1)^{\dff i}\dff
m_{\fff n\dff *}\fff 
(\dff c \dff)
\dff \circ \qff 
\pr_{\fff i\fff,\qff G}\qff (\dff {g} \dff)  
\off\off 
=\off\off 
d_{\dff n}\dff
\left(m_{\fff n\dff *}\fff 
(\dff c \dff)
\right)\dff (\dff {g} \dff)\dff. 
\]

\vspace*{-6pt}
Therefore\qss 
$m_{\fff n\dff +\dff 1\dff *}\dff \circ\dff d_{\dff n\dff -\dff 1}
\off =\off
d_{\dff n\dff -\dff 1}\dff \circ\dff m_{\fff n\dff *}$\nnsp.\oss
Let\dss $\alpha_{\fff *}$\dss be equal to\dss $m_{\fff n\dff *}$\dss
on\dss $F\dff(\dff \Gamma^{\fff n},\pff U \dff)$\dnsp.\oss
By the above results,\pss $\alpha_{\fff *}$\dss
is a morphism of resolutions with the required properties.\oss  \eproof

\mypar{Theorem.}{normal-amenable}
\emph{In the situation of\qss Theorem\qss \ref{chain-map-quotient}\qss
there is a canonical isometric isomorphism}\oss\vspace*{3pt} 
\[
\quad
\widehat{H}^{\dff *}\fff(\dff G\fff,\pff U^{\fff A} \dff)
\qff \ttoo\qff
\widehat{H}^{\dff *}\fff(\dff \Gamma,\pff U \dff)\dff.
\]

\vspace*{-9pt}
\proof\qss
Let\oss 
$\alpha^{\fff *}
\dff \colon\dff
F\dff(\dff G^{\fff \bullet\dff +\dff 1},\pff U \dff)
\ttoo
F\dff(\dff \Gamma^{\fff \bullet\dff +\dff 1},\pff U\dff)$\oss be the morphism induced by $\alpha$
and\dss let $\alpha_{\fff *}$ be the morphism from\dss Theorem\qss \ref{chain-map-quotient}.\oss
Both of\dss them  
are\dss 
$\Gamma$\nsp\dnsp-morphisms and\dss hence induce maps\vspace*{1.5pt}
\begin{equation}
\label{alpha-alpha}
\quad
\begin{tikzcd}[column sep=normal, row sep=normal]\dis
\alpha^{\dff *}
\qff \colon\qff
F\dff(\dff G^{\fff \bullet\dff +\dff 1},\pff U \dff)^{\fff \Gamma}
\off 
\arrow[r, shift left=3pt]
&
\off 
F\dff(\dff \Gamma^{\fff \bullet\dff +\dff 1},\pff U\dff)^{\fff \Gamma}
\off\dff \colon\dff
\alpha_{\fff *}\dff.
\arrow[l, shift left=3pt]
\end{tikzcd}
\end{equation}

\vspace{-9pt}
Together with morphisms of resolutions,\oss
the latter maps have the norm\qss $\leq\qff 1$\nnsp.\oss
By the definition,\oss the cohomology of\dss the complex\qss
$F\dff(\dff \Gamma^{\fff \bullet\dff +\dff 1},\pff U\dff)^{\fff \Gamma}$\qss
are the bounded cohomology\qss
$\widehat{H}^{\dff *}\fff(\dff \Gamma,\pff U \dff)$\dnsp.

In order to identify the cohomology of\dss the complex\qss
$F\dff(\dff G^{\fff \bullet\dff +\dff 1},\pff U \dff)^{\fff \Gamma}$\dnsp,\oss
let us consider the complex\qss
$F\dff(\dff G^{\fff \bullet\dff +\dff 1},\pff U \dff)^{\fff A}$\qss
of\dss $A$\dnsp-invariants of\dss the complex\qss
$F\dff(\dff G^{\fff \bullet\dff +\dff 1},\pff U \dff)$\dnsp.\oss
Since the subgroup\qss $A\off \subset\off \Gamma$\qss acts
trivially on\qss $G\off =\off \Gamma/A$\nnsp,\oss
for every\qss $n\qff \geq\qff 0$\qss the subspace\qss
$F\dff(\dff G^{\fff n},\pff U \dff)^{\fff A}$\qss
is equal to\qss
$F\dff(\dff G^{\fff n},\pff U^{\fff A} \dff)$\dnsp.\oss
Obviously\halfff,\pss $U^{\fff A}$\dss is a\dss $G$\dnsp-module
and\dss hence\qss
$F\dff(\dff G^{\fff n},\pff U^{\fff A} \dff)$\qss
is also a\dss $G$\dnsp-module.\oss
Moreover\halfff,\vspace*{3pt}
\[
\quad
F\dff(\dff G^{\fff \bullet\dff +\dff 1},\pff U \dff)^{\fff \Gamma}
\off =\off
\left(\qff
F\dff(\dff G^{\fff \bullet\dff +\dff 1},\pff U \dff)^{\fff A}
\qff\right)^{\fff G}
\off =\off
F\dff(\dff G^{\fff \bullet\dff +\dff 1},\pff U^{\fff A} \dff)^{\fff G}\dff.
\]

\vspace{-9pt}
Therefore\qss
$F\dff(\dff G^{\fff \bullet\dff +\dff 1},\pff U \dff)^{\fff \Gamma}$\qss
is equal to the complex of\dss the\dss $G$\dnsp-invariants of\qss 
$F\dff(\dff G^{\fff \bullet\dff +\dff 1},\pff U^{\fff A} \dff)$\qss
of\dss the\dss $G$\dnsp-module\dss $U^{\fff A}$\dss
and\dss the morphisms\qss (\ref{alpha-alpha})\qss 
can be interpreted as morphisms\vspace*{1.5pt}
\begin{equation}
\label{alpha-alpha-m}
\quad
\begin{tikzcd}[column sep=normal, row sep=normal]\dis
\alpha^{\dff *}
\qff \colon\qff
F\dff(\dff G^{\fff \bullet\dff +\dff 1},\pff U^{\fff A} \dff)^{\fff G}
\off 
\arrow[r, shift left=3pt]
&
\off 
F\dff(\dff \Gamma^{\fff \bullet\dff +\dff 1},\pff U\dff)^{\fff \Gamma}
\off\dff \colon\dff
\alpha_{\fff *}\qff.
\arrow[l, shift left=3pt]
\end{tikzcd}
\end{equation}

\vspace{-9pt}
Theorem\qss \ref{chain-map-quotient}\qss implies that
the induced map in cohomology\oss\vspace*{3pt}
\[
\quad 
\alpha_{\fff *}\dff \circ\dff \alpha^{\fff *}
\qff \colon\qff 
\widehat{H}^{\dff *}\fff(\dff G\fff,\pff U^{\fff A} \dff)
\qff \ttoo\qff
\widehat{H}^{\dff *}\fff(\dff G\fff,\pff U^{\fff A} \dff)
\]

\vspace{-9pt}
is equal to the identity\halfff.\oss
In order to prove that
the induced map in cohomology\oss\vspace*{3pt}
\[
\quad 
\alpha^{\fff *}\dff \circ\dff \alpha_{\fff *}
\qff \colon\qff 
\widehat{H}^{\dff *}\fff(\dff \Gamma,\pff U \dff)
\qff \ttoo\qff
\widehat{H}^{\dff *}\fff(\dff \Gamma,\pff U \dff)
\]

\vspace{-9pt}
is equal to the identity\halfff,\oss 
it is sufficient to prove that the\dss $\Gamma$\nsp\dnsp-morphism\vspace*{3pt}
\[
\quad
\alpha^{\fff *}\dff \circ\dff \alpha_{\fff *}
\qff \colon\qff 
F\dff(\dff \Gamma^{\fff \bullet\dff +\dff 1},\pff U \dff)
\qff \ttoo\qff
F\dff(\dff \Gamma^{\fff \bullet\dff +\dff 1},\pff U \dff)
\] 

\vspace{-9pt}
is chain homotopic to the identity\halfff.\oss
But the\dss $\Gamma$\nsp\dnsp-resolution\dss 
$F\dff(\dff \Gamma^{\fff \bullet\dff +\dff 1},\pff U \dff)$\dss
of $U$
is strong and relatively injective and\dss $\alpha^{\fff *}\dff \circ\dff \alpha_{\fff *}$\dss
extends\dss $\id_{\dff U}$\nnsp.\oss
Therefore Lemma\qss \ref{homotopy-uniqueness}\qss 
implies that\dss $\alpha^{\fff *}\dff \circ\dff \alpha_{\fff *}$\dss
is chain homotopic to the identity morphism of\dss
$F\dff(\dff \Gamma^{\fff \bullet\dff +\dff 1},\pff U \dff)$\dnsp.\oss
It follows that the maps\qss (\ref{alpha-alpha-m})\qss
are mutually inverse.\oss 
Since both of them have the norm\qss $\leq\qff 1$\nnsp,\oss
they are isometric isomorphisms.\oss  \eproof

\myuppar{Relative injectivity of\dss $\Gamma$\nsp\dnsp-modules\dss
$B\fff(\dff G^{\fff n},\pff U \dff)$\dnsp.}
For the proof\dss of\dss Theorem\qss \ref{normal-amenable}\qss
there is no need to know if\dss $B\fff(\dff G^{\fff n},\pff U \dff)$\dss
is relatively injective as a\dss $\Gamma$\nsp\dnsp-module,\oss
in contrast with the\dss $\Gamma$\nsp\dnsp-modules\dss
$B\fff(\dff \Gamma^{\fff n},\pff U \dff)$\dnsp.\oss
The relative injectivity of\dss latter is 
needed for applying Lemma\qss \ref{homotopy-uniqueness}.\oss
Still,\pss $B\fff(\dff G^{\fff n},\pff U \dff)$\dss
is relatively injective as a\dss $\Gamma$\nsp\dnsp-module.\oss
This follows from Theorem\qss \ref{chain-map-quotient}\qss
and the next lemma.

\tikzcdset{column sep/my-size/.initial=6em}
\tikzcdset{row sep/my-size/.initial=9ex}

\mypar{Lemma.}{retracts}
\emph{Let\qss $V$\dss be a relatively injective\dss $\Gamma$\nsp\dnsp-module
and\qss $W$\dss be a retract of\qss $V$\dss in the sense that there exist\dss 
$\Gamma$\nsp\dnsp-morphisms\qss
$i\dff \colon\dff W\ttoo V$\qss
and\qss
$p\dff \colon\dff V\ttoo W$\qss
such that\qss
$p\dff \circ\dff i\off =\off \id_{\qff W}$\qss
and\qss
$\|\qff i \qff\|\fff,\pff \|\qff p \qff\|\qff \leq\qff 1$\nnsp.\oss
Then\qss $W$\dss is also relatively injective.\oss}

\proof\qss
Let us consider the diagram\vspace*{-3pt}
\begin{equation*}
\quad
\begin{tikzcd}[column sep=my-size, row sep=my-size]\dis
V_{1} \arrow[r, "\dis k"]
\arrow[d, "\dis \alpha"']
& 
V_{2} \arrow[ld, dashed, "\dis \beta"'] 
\arrow[d, dashed, "\dis \gamma"]
\\
W \arrow[r, shift left=3pt, "\dis i"]
& 
V\dff, \arrow[l, shift left=3pt, "\dis p"]
\end{tikzcd}
\end{equation*}

\vspace*{-12pt}
where\qss
$k\dff \colon\dff V_{1}\ttoo V_{2}$\qss is a given strongly injective $\Gamma$\nsp\dnsp-morphism,\pss
$\alpha$\dss is\dss a\dss given $\Gamma$\nsp\dnsp-morphism,\oss
and we need to find a $\Gamma$\nsp\dnsp-morphism $\beta$ such that\qss
$\beta\dff \circ\dff k\qff =\qff \alpha$\qss
and\qss
$\|\qff \beta \qff\|\qff \leq\qff \|\qff \alpha \qff\|$\nnsp.\oss
Since $V$ is  relatively injective,\oss
there exists a\dss $\Gamma$\nsp\dnsp-morphism $\gamma$  
such that\qss
$\gamma\dff \circ\dff k\qff =\qff i\dff \circ\dff \alpha$\qss 
and\vspace{2.5pt} 
\[
\quad
\|\qff \gamma \qff\|
\qff \leq\qff 
\|\qff i\dff \circ\dff \alpha \qff\|
\qff \leq\qff
\|\qff i \qff\| \dff \|\qff \alpha \qff\|
\qff \leq\qff
\|\qff  \alpha \qff\|\dff.
\]

\vspace{-9.5pt} 
Let\qss $\beta\qff =\qff p\qff \circ\qff \gamma$\nnsp.\oss
Then\qss 
$\beta\qff \circ\qff  k
\qff =\qff 
p\qff \circ\qff  i\qff \circ\qff \alpha
\qff =\qff
\alpha$\qss
and\oss  
$\|\qff \beta \qff\|
\qff \leq\qff
\|\qff p \qff\| \dff \|\qff \gamma \qff\|
\qff \leq\qff 
\|\qff \alpha \qff\|$\nnsp.\oss  \eproof

\myappend{The\qss differential\qss of\qss the\qss standard\qss resolution}{categories-classifying}

\vspace*{6pt}
\myuppar{Semi-simplicial sets.}
A\qss \emph{semi-simplicial\dss set}\pss consists of sets\oss
$S_{\fff 0}\fff,\pff S_{\fff 1}\fff,\pff S_{\fff 2}\fff,\pff \ldots \off$\oss
and the maps\vspace*{3pt}
\[
\quad
\partial_{\fff i}
\qff \colon\qff
S_{\fff n}\qff \ttoo\qff S_{\fff n\dff -\dff 1}
\hspace*{1.5em}\mbox{ and }\hspace*{1.5em}
\delta_{\fff i}\qff \colon\qff
S_{\fff n}\qff \ttoo\qff S_{\fff n\dff +\dff 1}
\]

\vspace*{-9pt}
defined for\qss $0\qff \leq\qff i\qff \leq\qff n$\qss
and called the\qss \emph{face}\qss and\qss \emph{degeneracy\dss operators}\pss respectively.\oss
They should satisfy some well known relations which we will not use explicitly.\oss
Neither will we use the degeneracy operators.\oss
The elements of\dss $S_{\fff n}$\dss are called\dss \emph{$n$\dnsp-simplices}.\oss

\myuppar{Nerve of a category.}
The\qss \emph{nerve}\dss $\mathcal{NC}$ of a category\dss $\mathcal{C}$\dss 
is a semi-simplicial set having as its $n$\dnsp-simplexes
the sequences of objects and morphisms of $\mathcal{C}$ of the form\vspace*{3pt}
\begin{equation*}
\quad
\begin{tikzcd}[column sep=large, row sep=normal]\dis
O_{\dff 0} 
& 
O_{\dff 1} \arrow[l, "\textstyle \off a_{\dff 0}"']
&   
O_{\dff 2} \arrow[l, "\textstyle \off a_{\dff 1}"']
&
\off \ldots \off \arrow[l, "\textstyle \off a_{\dff 2}"']
&
O_{\fff n\dff -\dff 1}
\arrow[l, "\textstyle \off a_{\dff n\dff -\dff 2}"']
&
O_{\fff n}\qff. 
\arrow[l, "\textstyle \off a_{\dff n\dff -\dff 1}"']
\end{tikzcd}
\end{equation*}

\vspace*{-9pt}
The $i$\dnsp-th face of such a simplex are obtained by removing the object $O_{\fff n\dff -\dff i}$\nnsp.\oss
The objects $O_{\dff 0}$ and $O_{\dff n}$  are removed together with the arrow
respectively ending or starting at the removed object\halfff.\oss
If an object $O_{\fff i}$ with\qss $0\qff <\qff i\qff <\qff n$\qss is removed,\oss 
then two morphisms\vspace*{3pt}
\begin{equation*}
\quad
\begin{tikzcd}[column sep=large, row sep=normal]\dis
O_{\fff i\dff -\dff 1}
& 
O_{\fff i} \arrow[l, "\textstyle \off a_{\dff i\dff -\dff 1}"']
&   
O_{\fff i\dff +\dff 1} \arrow[l, "\textstyle \off a_{\dff i}"']
\end{tikzcd}
\quad
\mbox{ are replaced by }
\quad
\begin{tikzcd}[column sep=mycolumnsize6, row sep=normal] 
O_{\fff i\dff -\dff 1}  & O_{\fff i\dff +\dff 1}\off. 
\arrow[l, "\textstyle \off a_{\dff i\dff -\dff 1}\fff \circ\fff a_{\dff i}"']
\end{tikzcd}
\end{equation*}

\vspace{-9pt}
The operator $\delta_{\fff i}$ acts by replacing the object $O_{\fff n\dff -\dff i}$ 
by the identity morphism\oss
$O_{\fff n\dff -\dff i}\qff \longleftarrow\qff O_{\fff n\dff -\dff i}$\nnsp.

\myuppar{A category associated to a group.}
A group $G$ gives rise to a category $\mathcal{G}$\sss having the elements of $G$ as objects
and exactly one morphism between any two objects.\oss
For every\qss $k\fff,\pff g\qff \in\qff G$\qss there is a morphism\qss
$k\fff g\qff \longleftarrow\qff g$\qss denoted also by $k$ by an abuse of notations.\oss
The composition\qss $g_{\fff 1}\dff \circ\dff g_{\dff 2}$\nnsp,\oss when defined,\oss
is the morphism\qss $g_{\fff 1}\fff g_{\dff 2}$\nnsp.\oss
The $n$\dnsp-simplices of $\mathcal{NG}$\sss
are diagrams of the form\vspace*{3pt}
\begin{equation*}
\quad
\begin{tikzcd}[column sep=large, row sep=normal]\dis
\bullet 
& 
\bullet \arrow[l, "\textstyle \off g_{\dff 0}"']
&   
\bullet \arrow[l, "\textstyle \off g_{\dff 1}"']
&
\off \ldots \off \arrow[l, "\textstyle \off g_{\dff 2}"']
&
\bullet
\arrow[l]
&
g_{\dff n\dff -\dff 1}\fff g_{\dff n}
\arrow[l, "\textstyle \off g_{\dff n\dff -\dff 2}"']
&
g_{\dff n}\qff. 
\arrow[l, "\textstyle \off g_{\dff n\dff -\dff 1}"']
\end{tikzcd}
\end{equation*}

\vspace*{-9pt}
The bullets $\bullet$ 
on this diagram stand for the objects\oss
$O_{\fff i}
\off =\off
g_{\dff i}\dff g_{\dff i\dff +\dff 1}\qff \ldots\qff g_{\dff n}$\nnsp.\oss
The $n$\dnsp-simplices of $\mathcal{NG}$ are in one-to-one correspondence
with the sequences\qss
$(\dff g_{\dff 0}\fff,\pff g_{\dff 1}\fff,\pff g_{\dff 2}\fff,\pff \ldots\fff,\pff g_{\dff n}\dff)
\qff \in\qff G^{\fff n\dff +\dff 1}$\dnsp,\oss
and we will identify them with such sequences.\oss 
Then the face operators take the form\vspace*{3pt}
\[
\quad
\partial_{\fff 0}\fff
(\dff g_{\dff 0}\fff,\pff g_{\dff 1}\fff,\pff g_{\dff 2}\fff,\pff \ldots\fff,\pff g_{\dff n}\dff)
\off =\off
(\dff g_{\dff 1}\fff,\pff g_{\dff 2}\fff,\pff 
\ldots\fff,\pff 
g_{\dff n\dff -\dff 2}\fff,\pff
g_{\dff n\dff -\dff 1}\fff g_{\dff n}\dff)
\dff,
\]

\vspace*{-36pt}
\[
\quad
\partial_{\fff i}\fff
(\dff g_{\dff 0}\fff,\pff g_{\dff 1}\fff,\pff g_{\dff 2}\fff,\pff \ldots\fff,\pff g_{\dff n}\dff)
\off =\off
(\dff g_{\dff 1}\fff,\pff 
\ldots\fff,\pff 
g_{\dff n\dff -\dff i\dff -\dff 1}\fff g_{\dff n\dff -\dff i}\fff,\pff 
\ldots\fff,\pff 
g_{\dff n}\dff)
\hspace*{1.5em}\mbox{ for }\hspace*{1.5em}
0\qff <\qff i\qff <\qff n\dff,
\]

\vspace*{-36pt}
\[
\quad
\partial_{\fff n}\fff
(\dff g_{\dff 0}\fff,\pff g_{\dff 1}\fff,\pff g_{\dff 2}\fff,\pff \ldots\fff,\pff g_{\dff n}\dff)
\off =\off
(\dff g_{\dff 1}\fff,\pff g_{\dff 2}\fff,\pff \ldots\fff,\pff g_{\dff n}\dff)\dff.
\]

\vspace{-9pt}
The group $G$ acts on $\mathcal{G}$ from the right\halfff.\oss
The action of\qss $h\qff \in\qff G$\qss
takes an object $g$ to the object $g\fff h$
and a morphism\qss
$k\fff g\qff \longleftarrow\qff g$\qss
to the morphism\qss
$k\fff g\fff h\qff \longleftarrow\qff g\fff h$\nnsp.\oss
i.e.\qss it\dss takes a morphism\oss\vspace*{-1pt}
\[
\quad
\begin{tikzcd}[column sep=large, row sep=normal] 
\bullet  & \bullet \arrow[l, "\textstyle \off k"']
\end{tikzcd}
\quad
\mbox{ to a morphism of the same form }
\quad
\begin{tikzcd}[column sep=large, row sep=normal] 
\bullet  & \bullet\off. \arrow[l, "\textstyle \off k"']
\end{tikzcd}
\]

\vspace{-7pt}
This action induces the right action of $G$ on the nerve $\mathcal{NG}$\dnsp.\oss
The action of\qss $h\qff \in\qff G$\qss on the $n$\dnsp-simplex\qss 
$(\dff g_{\dff 0}\fff,\pff g_{\dff 1}\fff,\pff g_{\dff 2}\fff,\pff \ldots\fff,\pff g_{\dff n}\dff)$\qss 
takes it to\qss 
$(\dff g_{\dff 0}\fff,\pff g_{\dff 1}\fff,\pff g_{\dff 2}\fff,\pff \ldots\fff,\pff g_{\dff n}\fff h\dff)$\nnsp.\oss

\myuppar{The nerve $\mathcal{NG}$ and the standard resolution.}
Let\dss $B\dff(\fff \mathcal{NG}_{n}\dff,\off \rrr \dff)$\dss
be the space of bounded functions on the  set $\mathcal{NG}_{n}$\dnsp.\oss
The right action of $G$ on $\mathcal{NG}$ induces a right action of $G$
on $\mathcal{NG}_n$ and\dss hence a left action on\dss
$B\dff(\fff \mathcal{NG}_{n}\dff,\off \rrr \dff)$\dnsp,\oss
turning\dss
$B\dff(\fff \mathcal{NG}_{n}\dff,\off \rrr \dff)$\dss
into a\dss $G$\dnsp-module.\oss
Let\vspace*{3pt}
\[
\quad
\partial^{\fff n}
\qff \colon\qff
B\dff(\fff \mathcal{NG}_{n}\dff,\off \rrr \dff)
\qff \ttoo\qff
B\dff(\fff \mathcal{NG}_{n\dff +\dff 1}\dff,\off \rrr \dff)
\]

\vspace*{-9pt}
be the restriction to the bounded cochains of the standard
coboundary map,\oss
i.e.\qss let\vspace*{3pt}
\[
\quad
\partial^{\fff n}\fff(\dff f \dff)\fff(\fff \sigma\fff)
\off\off =\off\off
\sum_{i\qff =\qff 0}^{n\dff +\dff 1}
\off (\dff -\qff 1)^{\fff i}
f\fff(\dff \partial_{\fff i}\fff \sigma \dff)\dff.
\]

\vspace*{-6pt}
for every\qss
$f\qff \in\qff B\dff(\fff \mathcal{NG}_{n}\dff,\off \rrr \dff)$\dnsp.\oss
The identification of $\mathcal{NG}_{n}$\dss with\dss $G^{\dff n\dff +\dff 1}$\dss
leads to an identification of\dss
$B\dff(\fff \mathcal{NG}_{n}\dff,\pff \rrr \dff)$\dss
and\dss
$B\dff(\fff G^{\dff n\dff +\dff 1},\off \rrr \dff)$\dss
as left\dss $G$\dnsp-modules.\oss
Comparing the formulas for\dss $\partial^{\fff n}$\dss
and\dss $d_{\dff n}$\dss shows that this identification turns
$\partial^{\fff n}$ into the differential\dss $d_{\dff n}$\nnsp.\oss
This interpretation of the standard resolution can be routinely extended
to the case of non-trivial coefficients.

\myuppar{A proof that the standard resolution\dss
$B\dff(\fff G^{\dff \bullet\dff +\dff 1},\off \rrr \dff)$\dss is a complex\halfff.}
It is sufficient to show that\vspace*{0pt}
\begin{equation*}
\quad
\begin{tikzcd}[column sep=large, row sep=normal]\dis
0 \arrow[r]
& 
\rrr \arrow[r, "\dis \partial_{\fff -\dff 1}"]
& 
B\dff(\fff \mathcal{NG}_{0}\dff,\pff \rrr \dff) \arrow[r, "\dis \partial_{\fff 0}"]
&   
B\dff(\fff \mathcal{NG}_{1}\dff,\pff \rrr \dff) \arrow[r, "\dis \partial_{\fff 1}"]
&
\off \ldots \off,
\end{tikzcd}
\end{equation*}

\vspace*{-9pt}
where\qss
$\partial_{\fff -\dff 1}\fff(a\fff)\fff(\fff \sigma\fff)
\off =\off
a$\qss
for every\qss $a\qff \in\qff \rrr$\nnsp,\pss $\sigma\qff \in\qff \mathcal{NG}_{0}$\nnsp,\oss
is a complex\halfff.\oss
Removing two objects 
from an\dss $(n\qff +\qff 2)$\dnsp-simplex\qss
$\sigma\qff \in\qff \mathcal{NG}_{n\dff +\dff 2}$\qss
leads to a simplex\qss
$\tau\qff \in\qff \mathcal{NG}_{n}$\qss
which
does not depend on the order in which these two objects were removed.\oss
Given\qss
$f\qff \in\qff B\dff(\fff \mathcal{NG}_{n}\dff,\off \rrr \dff)$\dnsp,\oss
the value\dss
$f\dff(\dff \tau \dff)$\dss
enters the tautological double sum expressing\vspace*{3pt}
\[
\quad
\partial^{\fff n\dff +\dff 1}\fff
\left(\qff
\partial^{\fff n}\fff(\dff f \dff)\fff(\fff \sigma\fff)
\qff\right)
\]

\vspace*{-9pt}
twice,\oss corresponding to the two orders of removing two objects.\oss
But removing the first object changes the number\qss
(the subscript)\qss 
of the second one by $1$ if\dss
the first is to the right of the second on the above diagram,\oss
and does not change the number if on the left\halfff.\oss
It follows that\dss
$f\dff(\dff \tau \dff)$\dss
enters the double sum twice,\oss but with different signs.\oss
Therefore all terms of this double sum cancel\halfff,\oss and\dss hence\qss
$\partial^{\fff n\dff +\dff 1}\dff \circ\dff \partial^{\fff n}
\off =\off
0$\qss
for all\qss $n\qff \geq\qff 0$\dnsp.\oss
The equality\qss
$\partial_{\fff 0}\dff \circ\qff \partial_{\fff -\dff 1}\off =\off 0$\qss
is immediate.\oss
This proof can be routinely generalized to the case of non-trivial coefficients.

\myappend{The\qss complex\qss $B\dff(\fff G^{\dff \bullet\dff +\dff 1}\dnsp,\pff U \dff)^{\fff G}$}{invariants}

\vspace*{6pt}
\myuppar{The invariant subspaces.}
The subspaces\dss
$B\dff(\fff G^{\dff \bullet\dff +\dff 1}\dnsp,\pff U \dff)^{\fff G}$\dss
admit an explicit description.\oss
Let\qss $f\qff \in\qff B\dff(\fff G^{\dff n\dff +\dff 1},\pff U \dff)$\dnsp.\oss
Then\qss $f$ is invariant with respect to the action of\sss $G$ 
if\dss and only\dss if\vspace*{4.5pt}
\[
\quad
h\cdot\left(\dff
f\dff(\fff g_{\fff 1}\fff,\pff  \ldots\fff,\pff g_{\fff n}\fff,\pff g_{\fff n\dff +\dff 1}\halfff h\fff)
\dff\right)
\off =\off
f\dff(\fff g_{\fff 1}\fff,\pff  \ldots\fff,\pff g_{\fff n}\fff,\pff g_{\fff n\dff +\dff 1}\fff)
\]

\vspace{-7.5pt}
for\dss all\oss
$g_{\fff 1}\fff,\pff \ldots\fff,\pff g_{\fff n}\fff,\pff g_{\fff n\dff +\dff 1}\fff,\pff h
\qff \in\qff G$\nnsp.\oss
By taking\qss $h\qff =\qff g_{\dff n\dff +\dff 1}^{\dff -\dff 1}$\nnsp,\oss 
we see that\vspace*{4.5pt}
\[
\quad
f\dff(\fff g_{\fff 1}\fff,\pff \ldots\fff,\pff g_{\fff n}\fff,\pff g_{\fff n\dff +\dff 1}\fff)
\off =\off
g_{\dff n\dff +\dff 1}^{\dff -\dff 1}\cdot
\left(\dff
f\dff(\fff g_{\fff 1}\fff,\pff \ldots\fff,\pff g_{\fff n}\fff,\pff 1 \fff)
\dff\right)
\]

\vspace{-6pt}
if\dss $f$\dss is\dss $G$\dnsp-invariant\halfff.\oss
It follows that the map\vspace*{4.5pt}
\[
\quad
i^{\fff n}
\qff \colon\qff 
B\dff(\fff G^{\dff n\dff +\dff 1},\pff U \dff)^{\fff G}
\qff \ttoo\qff
B\dff(\fff G^{\dff n}\fff,\pff U \dff)
\]

\vspace{-6pt}
defined by the formula\oss
$i^{\fff n}\fff(\dff f\dff)\fff(\fff g_{\fff 1}\fff,\pff \ldots\fff,\pff g_{\fff n}\fff)
\off =\off
f\dff(\fff g_{\fff 1}\fff,\pff \ldots\fff,\pff g_{\fff n}\fff,\pff 1 \fff)$\oss
is a canonical isomorphism.

\myuppar{The complex of invariant subspaces.}
By using isomorphisms $i^{\fff n}$
we see that the bounded cohomology\qss
$\widehat{H}^{\dff *}\fff(\dff G\fff,\pff U  \dff)$\qss
is equal to the cohomology of the complex\vspace*{3pt}
\begin{equation*}
\quad
\begin{tikzcd}[column sep=large, row sep=normal]\dis
0 \arrow[r]
& 
U \arrow[r, "\dis \delta_{\dff 0}\off"]
& 
B\dff(\fff G\fff,\pff U \dff) \arrow[r, "\dis \delta_{\dff 1}\off"]
&   
B\dff(\fff G^{\dff 2},\pff U \dff) \arrow[r, "\dis \delta_{\dff 2}\off"]
&
\off \ldots \off,
\end{tikzcd}
\end{equation*}

\vspace*{-3pt}
where\oss  
$\delta_{\dff 0}\dff(\fff v \fff)(\fff g\fff)\off =\off v$\oss 
for\dss all\qss
$v\qff \in\qff U$\dnsp,\qss $g\qff \in\qff G$\nnsp,\off\oss
and\dss\vspace*{4pt}
\[
\quad
\delta_{\dff n}\dff(\dff f \dff)\dff
(\dff g_{\fff 0}\fff,\pff g_{\fff 1}\fff,\pff \ldots\fff,\pff g_{\fff n} \dff)
\off =\off
(\dff -\qff 1\dff)^{n\dff +\dff 1}\dff
f\dff(\dff g_{\fff 1}\fff,\pff g_{\fff 2}\fff,\pff \ldots\fff,\pff g_{\fff n} \dff)
\]

\vspace*{-30pt}
\[
\quad
\phantom{\delta_{\dff n}\dff(\dff f \dff)\dff
(\dff g_{\fff 0}\fff,\pff g_{\fff 1}\fff,\pff \ldots\fff,\pff g_{\fff n} \dff)
\off = }
+\off\off
\sum_{i\qff =\qff 0}^{n\qff -\qff 1}\qff (-\qff 1)^{\dff n\dff -\dff i}\qff 
f\dff(\dff g_{\fff 0}\fff,\pff \ldots\fff,\pff
g_{\dff i}\fff g_{\fff i\dff +\dff 1}\fff,\pff \ldots\fff,\pff   g_{\dff n} \dff)
\]

\vspace*{-24pt}
\[
\quad
\phantom{\delta_{\dff n}\dff(\dff f \dff)\dff
(\dff g_{\fff 0}\fff,\pff g_{\fff 1}\fff,\pff \ldots\fff,\pff g_{\fff n} \dff)
\off = }
+\off\off
g_{\dff n}^{\dff -\dff 1}\cdot
\left(\dff
f\dff(\dff g_{\fff 0}\fff,\pff \ldots\fff,\pff g_{\fff n\dff -\dff 1} \dff)
\dff\right)\dff.
\]

\vspace*{-7pt}
for\qss $n\qff \geq\qff 1$\qss and\oss
$g_{\fff 0}\fff,\pff g_{\fff 1}\fff,\pff \ldots\fff,\pff g_{\fff n}
\qff \in\qff G$\nnsp.\off\oss
For the first differential $\delta_{\dff 1}$ this formula takes the form\vspace*{3pt}
\[
\quad
\delta_{\dff 1}\dff(\dff f \dff)\dff
(\dff g_{\fff 0}\fff,\pff g_{\fff 1} \dff)
\off =\off
f\dff(\dff g_{\fff 1} \dff)
\off -\off
f\dff(\dff g_{\dff 0}\fff g_{\fff 1} \dff)
\off +\off
g_{\dff 1}^{\dff -\dff 1}\cdot
\left(\dff
f\dff(\dff g_{\fff 0} \dff)
\dff\right)\dff.
\]

\vspace*{-9pt}
If\dss $U$ is a trivial\dss $G$\dnsp-module,\oss
then the formula for $\delta_{\dff 1}$ takes even simpler form\oss\vspace*{3pt}
\[
\quad
\delta_{\dff 1}\dff(\dff f \dff)\dff
(\dff g_{\fff 0}\fff,\pff g_{\fff 1} \dff)
\off =\off
f\dff(\dff g_{\fff 0} \dff)
\off +\off
f\dff(\dff g_{\fff 1} \dff)
\off -\off
f\dff(\dff g_{\dff 0}\fff g_{\fff 1} \dff)\dff.
\]

\myappend{The\qss second\qss bounded\qss cohomology\qss group}{second}

\vspace*{6pt}
\myuppar{Pseudo-limits.}
Let
$\nnn$ the set of all positive integers.\oss
The space\dss $B\dff(\dff \nnn \dff)$\dss
can be considered as the space of\dss bound\-ed sequences of real numbers.\oss
Let\dss $L\dff(\dff \nnn \dff)$\dss be the subspace of
sequences $f$ such that there exists  limit of\dss
$f\dff(n\fff)$\dss for\qss $n\qff \to\qff \infty$\nnsp.\oss 
A\qss \emph{pseudo-limit}\pss is defined as a bounded linear functional\qss
$l
\dff \colon\dff
B\dff(\dff \nnn \dff)
\ttoo
\rrr$\qss 
such that\vspace*{5pt}
\[
\quad
l\dff(\dff f \dff)
\off =\off
\lim_{n\dff \to\dff \infty}\qff f\dff(\fff n\fff)
\]

\vspace*{-7pt}
for all\qss $f\qff \in\qff L\dff(\dff \nnn \dff)$\qss
and\qss $\|\qff l \qff\|\qff \leq\qff 1$\nnsp.\oss
Let us prove that pseudo-limits exist\halfff.\oss

Let us define\dss $l\dff(\dff f \dff)$\dss for\qss
$f\qff \in\qff L\dff(\dff \nnn \dff)$\qss as the limit of $f$\dnsp.\oss
Then $l$ is a linear functional\qss
$L\dff(\dff \nnn \dff)\ttoo \rrr$\qss
such that\qss
$|\qff l\dff(\dff f \dff) \qff|\qff \leq\qff \|\qff f \qff\|$\qss
for all\qss $f\qff \in\qff L\dff(\dff \nnn \dff)$\dnsp.\oss
By applying Hahn--Banach theorem to $l$ 
and the norm\dss $\|\qff \bullet \qff\|$\dss on\dss $B\dff(\dff \nnn \dff)$\dss
we see that $l$ can be extended\dss to a linear functional on\dss
$B\dff(\dff \nnn \dff)$\dss
in such a way that\qss
$|\qff l\dff(\dff f \dff) \qff|\qff \leq\qff \|\qff f \qff\|$\qss
for all\qss $f\qff \in\qff B\dff(\dff \nnn \dff)$\dnsp.\oss
Obviously,\oss such an extension is a pseudo-limit\halfff.\oss

Let\qss $s\dff \colon\dff \nnn\ttoo \nnn$\qss be the map\qss
$n\qff \longmapsto\qff n\qff +\qff 1$\nnsp.\oss
One can prove that there exist pseudo-limits $l$
such that\qss
$l\dff(\dff f \dff)
\off =\off
l\dff(\dff f\dff \circ\dff s \dff)$\qss
for all\qss $f\qff \in\qff B\dff(\dff \nnn \dff)$\dnsp.\oss
Such a pseudo-limit is called a\qss \emph{Banach\dss limit\halfff}.\oss
See,\oss for example,\oss \cite{r},\oss 
Exercise\qss 4\qss to Chapter\qss 3.\oss
But for our purposes pseudo-limits are sufficient\halfff.\oss

\vspace*{3pt}
\myapar{Theorem\qss (see \cite{mm},\oss \cite{i2}).}{second-is-norm}
\emph{For every group\qss $G$\qss the semi-norm on\qss
$\widehat{H}^{\dff 2}\fff(\dff G \dff)$\qss
is a norm.\oss}

\vspace*{3pt}
\proof\qss
We will use the description of the complex\dss 
$B\dff(\fff G^{\dff \bullet\dff +\dff 1}\dnsp,\pff \rrr \dff)^{\fff G}$\dss
from Appendix\qss \ref{invariants}.\oss
It is sufficient to prove that\dss
$\image\dff \delta_{\dff 1}$\dss is closed.\oss
In order to do this it is sufficient\halfff,\oss
in turn,\oss 
to find a bounded left inverse to $\delta_{\dff 1}$\nnsp,\oss
i.e.\qss a bounded operator $P$ such that\qss
$P\dff \circ\dff \delta_{\dff 1}
\off =\off
\id$\nnsp.\oss

Indeed,\oss
suppose that\qss
$P\dff \circ\dff \delta_{\dff 1}
\off =\off
\id$\qss
and\dss let\qss
$Q
\off =\off
\id\qff -\pff \delta_{\dff 1}\dff \circ\dff P$\dnsp.\oss
Then\vspace*{5pt}
\[
\quad
Q\dff \circ\dff \delta_{\dff 1}
\off =\off
\bigl(\dff
\id\qff -\pff \delta_{\dff 1}\dff \circ\dff P
\dff\bigr)
\dff \circ\dff \delta_{\dff 1}
\off =\off
\delta_{\dff 1}\qff -\qff \delta_{\dff 1}\dff \circ\dff 
\bigl(\dff P\dff \circ\dff \delta_{\dff 1}
\dff\bigr)
\off =\off
\delta_{\dff 1}\qff -\qff \delta_{\dff 1}
\off =\off
0
\]

\vspace*{-7pt}
and\dss hence\qss
$\image\dff \delta_{\dff 1}
\qff \subset\qff
\kernel\dss Q$\dnsp.\oss
On the other hand,\oss
if\qss
$x\qff \in\qff \kernel\dff Q$\dnsp,\oss
then\qss
$x\qff -\qff \delta_{\dff 1}\dff \circ\dff P\dff(x\fff)
\off =\off
0$\qss
and\dss hence\qss
$x
\off =\off
\delta_{\dff 1}\fff(\dff P\dff(x\fff)\dff)$\qss
belongs to the image of\dss $\delta_{\dff 1}$\nnsp.\oss
Therefore\qss
$\kernel\dss Q
\qff \subset\qff
\image\dff \delta_{\dff 1}$\nnsp.\oss
By combining this with the already proved opposite inclusion,\oss we see that\qss
$\image\dff \delta_{\dff 1}
\off =\off
\kernel\dss Q$\dnsp.\oss
Since $Q$ is a bounded operator together with $\delta_{\dff 1}$ and $P$\dnsp,\oss
its kernel is closed\dss and\dss hence\dss $\image\dff \delta_{\dff 1}$
is closed.\oss

It remains to construct a bounded left inverse $P$\dnsp.\oss
Such an inverse should,\oss in particular\halfff,\oss
recover\qss $f\qff \in\qff B\dff(\fff G \dff)$\qss
if\qss
$\delta_{\dff 1}\dff(\dff f \dff)$\qss
is known.\oss
The kernel of\dss $\delta_{\dff 1}$\dss consists of the bounded homomorphisms\qss
$G\ttoo \rrr$\nnsp,\oss as it follows from the explicit formula for\dss $\delta_{\dff 1}$\nnsp.\oss
Since every bounded homomorphism\qss
$G\ttoo \rrr$\qss is obviously equal to zero,\pss
$\delta_{\dff 1}$\dss is injective.\oss
Therefore\qss $f\qff \in\qff B\dff(\fff G \dff)$\qss
is indeed determined\dss by\qss
$\delta_{\dff 1}\dff(\dff f \dff)$\dnsp.\oss
Let us try to find an explicit way to recover $f$\dss
from\dss $\delta_{\dff 1}\dff(\dff f \dff)$\dnsp.\oss

Let\qss $f\qff \in\qff B\dff(\fff G \dff)$\qss
and\trs let\qss
$A
\off =\off
\delta_{\dff 1}\dff(\dff f \dff)$\dnsp.\oss
Then\oss \vspace*{3pt}
\[
\quad
A\dff(\dff g\fff,\pff g^{\fff i\dff -\dff 1} \dff)
\off =\off
f\dff(\dff g \dff)
\off +\off
f\dff(\dff g^{\fff i\dff -\dff 1} \dff)
\off -\off
f\dff(\dff g^{\fff i} \dff)\
\oss \]

\vspace{-9pt}
for every\qss $g\qff \in\qff G$\qss and\qss $i\qff \geq\qff 0$\dnsp.\oss
Let us consider an arbitrary element\qss $g\qff \in\qff G$\nnsp.\oss
If\qss $m\qff \geq\qff 1$\nnsp,\oss then\vspace*{7pt}
\[
\quad
f\dff(\dff g^{\fff m} \dff)
\off =\off
f\dff(\dff g \dff)
\qff +\qff
f\dff(\dff g^{\fff m\dff -\dff 1} \dff)
\qff -\qff
A\dff(\dff g\fff,\pff g^{\fff m\dff -\dff 1} \dff)
\]

\vspace*{-30pt}
\[
\quad
\phantom{f\dff(\dff g^{\fff m} \dff)
\off }
=\off
f\dff(\dff g \dff)
\qff +\qff
f\dff(\dff g \dff)
\qff +\qff
f\dff(\dff g^{\fff m\dff -\dff 2} \dff)
\qff -\qff
A\dff(\dff g\fff,\pff g^{\fff m\dff -\dff 1} \dff)
\qff -\qff
A\dff(\dff g\fff,\pff g^{\fff m\dff -\dff 2} \dff)
\]

\vspace*{-39pt}
\[
\quad
\phantom{f\dff(\dff g^{\fff m} \dff)
\off }
\ldots
\]

\vspace*{-42pt}
\[
\quad
\phantom{f\dff(\dff g^{\fff m} \dff)
\off }
=\off
m\dff
f\dff(\dff g \dff)
\off -\off
\sum_{i\qff =\qff 1}^{m\qff -\qff 1}\qff
A\dff(\dff g\fff,\pff g^{\fff i} \dff)
\]

\vspace*{-12pt}
and\dss hence\vspace*{0pt}
\[
\quad
f\dff(\dff g \dff)
\off\off =\off\off
\frac{1}{\fff m\fff}\qff
f\dff(\dff g^{\fff m} \dff)
\off\off +\off\off
\frac{1}{\fff m\fff}\qff
\sum_{i\qff =\qff 1}^{m\qff -\qff 1}\qff
A\dff(\dff g\fff,\pff g^{\fff i} \dff)\qff.
\]

\vspace*{-6pt}
Since $f$\dss is bounded,\vspace*{0pt}
\[
\quad
\lim_{m\dff \to\dff \infty}\qff \frac{1}{\fff m\fff}\off
f\dff(\dff g^{\fff m} \dff)
\off\off =\off\off
0\dff.
\]

\vspace*{-6pt}
It follows that\vspace*{0pt}
\[
\quad
f\dff(\dff g \dff)
\off\off =\off\off
\lim_{m\dff \to\dff \infty}\qff
\frac{1}{\fff m\fff}\qff
\sum_{i\qff =\qff 1}^{m\qff -\qff 1}\qff
A\dff(\dff g\fff,\pff g^{\fff i} \dff)\dff.
\]

\vspace*{-6pt}
The last formula recovers $f$\dss by\qss
$A
\off =\off
\delta_{\dff 1}\dff(\dff f \dff)$\dnsp.\oss
In order to construct $P$\dnsp,\oss
it is sufficient to replace the limit in this formula by a pseudo-limit $l$\nnsp.\oss
Indeed,\oss if\qss $A\qff \in\qff B\dff(\fff G^{\dff 2} \fff)$\dnsp,\oss
then\vspace*{3pt}
\[
\left|\off
\frac{1}{\fff m\fff}\qff
\sum_{i\qff =\qff 1}^{m\qff -\qff 1}\qff
A\dff(\dff g\fff,\pff g^{\fff i} \dff)
\off\right|
\off\off \leq\off\off
\|\qff A \qff\|
\]

\vspace*{-9pt}
and\dss hence we can define\dss $P\dff(\dff A \dff)$\dss by the formula\vspace*{3pt}
\[
\quad
P\dff(\dff A\dff)\fff(\dff g \dff)
\off\off =\off\off
l\qff
\left(\off
\sum_{i\qff =\qff 1}^{m\qff -\qff 1}\qff
\frac{1}{\fff m\fff}\qff
A\dff(\dff g\fff,\pff g^{\fff i} \dff)
\off\right)\qff.
\]

\vspace*{-9pt}
Obviously,\pss $\|\qff P\dff(\dff A\dff) \qff\|\qff \leq\qff \|\qff A \qff\|$\qss
and\dss hence $P$ is a bounded operator\halfff.\oss
In fact\halfff,\pss $\|\qff P \qff\|\qff \leq\qff 1$\nnsp.\oss
Since $l$ is equal to the limit on $L\dff(\dff N\dff)$\dnsp,\oss
the above calculations show that $P$ is a left inverse to\dss $\delta_{\dff 1}$\nnsp.\oss  \eproof

\myappend{Functoriality\qss with\qss coefficients}{functoriality-coefficients}

\vspace*{6pt}
\myuppar{Change of groups.}
Let\qss 
$\alpha\dff \colon\dff \Gamma\toto G$\qss be a homomorphism.\oss 
Then any $G$\dnsp-module $V$ can be turned into a\dss $\Gamma$\nsp\dnsp-module
by defining the action of\dss $\Gamma$ by the rule\qss 
$(\dff \gamma\fff,\pff v \dff)\qff \longmapsto\qff \alpha\dff(\dff \gamma \dff)\cdot v$\nnsp.\oss
This $\Gamma$\nsp\dnsp-module is denoted by\dss $_\alpha V$\dss
and said to be the result of the\qss \emph{change of\dss groups}\qss 
by $\alpha$\dnsp.\oss
Obviously,\vspace*{3pt}
\[
\quad
V^{\dff G}\off \subset\off _\alpha V^{\dff \Gamma}\dff.
\]

\vspace{-9pt}
Suppose now that $U$ is a\dss $\Gamma$\nsp\dnsp-module and
$V$ is a $G$\dnsp-module.\oss
A linear map\qss
$u\dff \colon\dff V\ttoo U$\qss
is said to be an\qss \emph{$\alpha$\dnsp-morphism}\pss
if\trs it\trs is\dss a\dss $\Gamma$\nsp\dnsp-morphism of\dss $\Gamma$\nsp\dnsp-modules\qss
$_\alpha V\ttoo U$\dnsp.\oss

\myuppar{The induced maps.}
Let\qss $u\dff \colon\dff V\ttoo U$\qss be 
an $\alpha$\dnsp-morphism.\oss 
It\dss leads to the canonical maps\vspace*{4pt}
\[
\quad
(\dff \alpha\fff,\pff u \dff)^*
\qff \colon\qff
B\dff(\dff G^{\dff n\dff +\dff 1}\dnsp,\pff V \dff)
\qff \ttoo\qff
B\dff(\dff \Gamma^{\dff n\dff +\dff 1}\dnsp,\pff U \dff)
\]

\vspace*{-9pt}
defined by the formula\vspace*{3pt}
\[
\quad
(\dff \alpha\fff,\pff u \dff)^*\fff (\dff f \dff)
\fff
\left(\trf 
\gamma_{\fff 0}\fff,\pff \gamma_{\fff 1}\fff,\pff \ldots\fff,\pff \gamma_{n} 
\dff\right)
\off =\off
u\qff
\Bigl(\qff
f\dff
\left(\dff  \alpha\dff(\dff \gamma_{\fff 0}\dff)\fff,\pff  \alpha\dff(\dff \gamma_{\fff 1}\dff)\fff,\pff
\ldots\fff,\pff 
\alpha\dff(\dff \gamma_{n}\dff) 
\dff\right)
\qff\Bigr)\dff.
\]

\vspace*{-8pt}
A trivial check shows that the maps\dss $(\dff \alpha\fff,\pff u \dff)^*$\dss
are $\alpha$\dnsp-morphisms and commute with the differentials of\dss 
the standard resolutions.\oss
Therefore the maps\dss $(\dff \alpha\fff,\pff u \dff)^*$\dss 
define\dss a $\Gamma$\nsp\dnsp-morphism\vspace*{4pt}
\begin{equation*}
\quad
(\dff \alpha\fff,\pff u \dff)^*
\qff \colon\qff
_\alpha\dff B\dff(\dff G^{\dff \bullet\dff +\dff 1}\dnsp,\pff V \dff)
\qff \ttoo\qff
B\dff(\dff \Gamma^{\dff \bullet\dff +\dff 1}\dnsp,\pff U \dff)
\end{equation*}

\vspace*{-8.5pt}
of\dss $\Gamma$\nsp\dnsp-resolutions.\oss
Since\oss\vspace*{3.5pt}
\[
\quad
B\dff(\dff G^{\dff n}\dnsp,\pff V \dff)^{\fff G}
\off \subset\off
_\alpha\dff B\dff(\dff G^{\dff n}\dnsp,\pff V \dff)^{\fff \Gamma}\dff,
\]

\vspace*{-8pt}
for every\qss $n\qff \geq\qff 0$\dnsp,\oss
this $\Gamma$\nsp\dnsp-morphism leads to a morphism of complexes\vspace*{4pt}
\[
\quad
B\dff(\dff G^{\dff \bullet\dff +\dff 1}\dnsp,\pff V \dff)^{\fff G}
\qff \ttoo\qff
_\alpha\dff B\dff(\dff \Gamma^{\dff \bullet\dff +\dff 1}\dnsp,\pff V \dff)^{\fff \Gamma}
\]

\vspace*{-8pt}
and\dss then to homomorphisms\vspace*{4pt}
\[
\quad
(\dff \alpha\fff,\pff u \dff)^*
\qff \colon\qff
\widehat{H}^{\fff n}\fff(\dff G\fff,\qff V \dff)
\qff \ttoo\qff
\widehat{H}^{\fff n}\fff(\dff \Gamma\fff,\qff U \dff)
\]

\vspace*{-8pt}
of\dss the bounded cohomology spaces.\oss
If\qss $U\off =\off V\off =\off \rrr$\qss 
with the trivial action of\pss $\Gamma\fff,\off G$\qss respectively\dss
and\dss if\dss $u\off =\off \id_{\dff \rrr}$\nnsp,\oss
then the homomorphism\dss
$(\dff \alpha\fff,\pff u \dff)^*$\dss
is equal to 
$\alpha^{\fff *}$
from Section\qss \ref{spaces}.\oss

The norm of\dss 
the maps\dss
$(\dff \alpha\fff,\pff u \dff)^*$\dss 
is\dss obviously\qss $\leq\qff \|\qff u \qff\|$\nnsp.\oss

\myappend{Straight\qss and\qss Borel\qss straight\qss cochains}{borel}

\vspace*{6pt}
\myuppar{Borel straight cochains.}
Let $G$ be a discrete group 
and\dss let\qss 
$p\dff \colon\dff \mathcal{X}\ttoo X$\qss
be a locally trivial principal right\dss $G$\dnsp-bundle.\oss
A cochain\qss $f\qff \in\qff B^{\fff n}\fff(\dff \mathcal{X} \dff)$\qss
is called\qss \emph{straight}\pss if\dss $f\fff(\dff \sigma \dff)$\dss depends only
on the vertices of $\sigma$\nnsp.\oss
We will identify straight $n$\dnsp-cochains with functions\qss
$\mathcal{X}^{n\dff +\dff 1}\ttoo \rrr$\nnsp,\oss
and\dss the space of straight $n$\dnsp-cochains with the space\dss
$B\dff(\dff \mathcal{X}^{n\dff +\dff 1} \dff)$\dnsp.\oss
A straight cochain\dss $f\qff \in\qff B\dff(\dff \mathcal{X}^{n\dff +\dff 1} \dff)$\dss
is said to be a\dss \emph{Borel straight cochain}\pss
if\qss
$f\dff \colon\dff \mathcal{X}^{n\dff +\dff 1}\ttoo \rrr$\qss
is a Borel function,\oss
i.e.\qss if the preimage of every Borel subset of $\rrr$
is a Borel subset of\dss $\mathcal{X}^{n\dff +\dff 1}$\dnsp.\oss
The space of\dss Borel straight $n$\dnsp-cochains is denoted\dss by\dss
$\mathcal{B}\dff(\dff \mathcal{X}^{n\dff +\dff 1} \dff)$\dnsp.\oss
The image of\qss 
$d_{\dff -\dff 1}
\dff \colon\dff 
\rrr
\ttoo
B^{\fff 0}\fff(\dff \mathcal{X} \dff)$\qss
consists of constant functions and\dss hence
is contained in\dss
$\mathcal{B}\dff(\dff \mathcal{X} \dff)$\dnsp.\oss
Each space\dss
$\mathcal{B}\dff(\dff \mathcal{X}^{n\dff +\dff 1} \dff)$\dss
is in a natural way a\dss $G$\dnsp-module,\oss
and they form a subcomplex\dss
$\mathcal{B}\dff(\dff \mathcal{X}^{\fff \bullet\dff +\dff 1} \dff)$\dss 
of the complex\dss
$B^{\fff \bullet}\fff(\dff \mathcal{X} \dff)$\dnsp.\oss

\myapar{Lemma.}{borel-straight-ri}
\emph{The\qss $G$\dnsp-modules\qss
$\mathcal{B}\dff(\dff \mathcal{X}^{n\dff +\dff 1} \dff)$\dss
are relatively injective for all\dss $n$\nnsp.\oss}

\proof\qss
Let $F$ be a Borel\dss fundamental set for 
the action of $G$ on $\mathcal{X}$
and\dss let $V^{\fff n}$ be the space of\dss
bounded\dss Borel functions\qss
$F\dff \times\dff \mathcal{X}^{\fff n}\ttoo \rrr$\nnsp.\oss
Let\oss \vspace*{3pt}
\[
\quad
I^{n}
\qff \colon\qff
\mathcal{B}\dff(\dff \mathcal{X}^{n\dff +\dff 1} \dff)
\qff \ttoo\qff
B\dff (\dff G\fff,\pff V^{\fff n} \dff)\oss
\]

\vspace*{-9pt}
be the map given by the formula\oss\vspace*{3pt}
\[
\quad
I^{n}\fff(\dff f \dff)\fff(\dff g \dff)\dff
(\dff y_{\fff 0}\fff,\pff y_{\fff 1}\fff,\pff \ldots\fff,\pff y_{\fff n} \dff)
\off =\off
f\fff(\dff y_{\fff 0}\halfff g\fff,\pff y_{\fff 1}\halfff g\fff,\pff \ldots\fff,\pff y_{\fff n}\halfff g \dff)\fff,
\]

\vspace*{-9pt}
where\qss $g\qff \in\qff G$\nnsp,\qss
$y_{\fff 0}\qff \in\qff F$\nnsp,\oss
and\oss
$y_{\fff 1}\fff,\pff \ldots\fff,\pff y_{\fff n}\qff \in\qff \mathcal{X}$\dnsp.\oss
A routine check shows that $I^{n}$ is an isometric isomorphism of\dss $G$\dnsp-modules.\oss
Therefore\qss 
$\mathcal{B}\dff(\dff \mathcal{X}^{n\dff +\dff 1} \dff)$\dss 
is relatively injective
by Lemma\qss \ref{bgv}.\oss  \eproof

\myapar{Lemma.}{borel-straight-contracting}
\emph{The complex\qss
$\mathcal{B}\dff(\dff \mathcal{X}^{\fff \bullet\dff +\dff 1} \dff)$\qss
together with the map\qss
$d_{\dff -\dff 1}
\dff \colon\dff
\rrr
\ttoo
\mathcal{B}\dff(\dff \mathcal{X} \dff)$\qss
is a strong resolution of\dss the trivial\dss $G$\dnsp-module\dss $\rrr$\nnsp.\oss}

\proof\qss
Let\qss $b\qff \in\qff \mathcal{X}$\qss
and\dss let\qss\vspace*{3pt}
\[
\quad
K_{\dff n}
\qff \colon\qff
\mathcal{B}\dff(\dff \mathcal{X}^{n\dff +\dff 1} \dff)
\qff \ttoo\qff 
\mathcal{B}\dff(\dff \mathcal{X}^{n} \dff)
\hspace*{1em}\mbox{ and }\hspace*{1em}
\quad
K_{\dff 0}
\qff \colon\qff
\mathcal{B}\dff(\dff \mathcal{X} \dff)
\qff \ttoo\qff 
\rrr
\]

\vspace{-9pt}
be the maps defined\dss by the formulas\vspace*{3pt}
\[
\quad
K_{\dff n}\fff(\dff f \dff)\fff(\fff y_{\fff 1}\fff,\pff \ldots\fff,\pff y_{\fff n})
\off =\off
f\fff(\fff b\fff,\pff y_{\fff 1}\fff,\pff \ldots\fff,\pff y_{\fff n})
\hspace*{1em}\mbox{ and }\hspace*{1em}
K_{\dff 0}\fff(\dff f \dff)
\off =\off
f\fff(\dff b \fff)\dff.
\]

\vspace{-9pt}
Obviously,\pss $K_{\dff n}\fff(\dff f \dff)$
is a Borel function if\dss $f$\dss is,\oss
and\dss hence this definition is correct\halfff.\oss
A standard check shows that\dss $K_{\dff \bullet}$\dss
is a contracting homotopy.\oss
It remains to point out that\qss $\|\qff K_{\dff n} \qff\|\qff \leq\qff 1$\nnsp.\oss  \eproof

\tikzcdset{column sep/my-size-t/.initial=2em}
\tikzcdset{row sep/my-size-t/.initial=6em}

\myapar{Theorem.}{straight-cohomology}
\emph{If\trs the fundamental group\qss 
$\pi_{\dff 1}\dff(\dff \mathcal{X} \dff)$\qss
is amenable,\oss
then the map}\vspace*{6pt} 
\[
\quad
\mathcal{J}
\qff \colon\qff
\mathcal{H}^{\fff *}\left(\dff 
G\fff,\pff \mathcal{B}\dff(\dff \mathcal{X}^{\fff \bullet\dff +\dff 1} \dff) \dff\right)
\qff \ttoo\qff
\mathcal{H}^{\fff *}\left(\dff 
G\fff,\pff B^{\fff \bullet}\fff(\dff \mathcal{X} \dff) \dff\right)
\off\off =\off\off
\widehat{H}^{\dff *}\fff(\dff X \dff)
\]

\vspace{-6pt}
\emph{induced by the inclusion\qss
$\mathcal{B}\dff(\dff \mathcal{X}^{\fff \bullet\dff +\dff 1} \dff)
\ttoo
B^{\fff \bullet}\fff(\dff \mathcal{X} \dff)$\pss
is an isometric isomorphism.\oss}

\proof\qss
The resolutions\dss
$B^{\fff \bullet}\fff(\dff \mathcal{X} \dff)$\dss
and\dss
$\mathcal{B}\dff(\dff \mathcal{X}^{\fff \bullet\dff +\dff 1} \dff)$\dss
are strong resolutions by\dss Theorem\qss \ref{simply-connected-homotopy}\qss
and\dss Lemma\qss \ref{borel-straight-contracting}\qss respectively.\oss
Hence by\dss Theorem\qss \ref{comparing-to-standard}\qss 
there exist\dss $G$\dnsp-morphisms of resolutions\vspace*{4pt}
\[
\quad
u_{\dff \bullet}
\qff \colon\qff
B^{\fff \bullet}\fff(\dff \mathcal{X} \dff)
\qff \ttoo\qff 
B\dff(\fff G^{\dff \bullet\dff +\dff 1} \dff)
\hspace*{1em}\mbox{ and }\hspace*{1em}
v_{\dff \bullet}
\qff \colon\qff
\mathcal{B}\dff(\dff \mathcal{X}^{\fff \bullet\dff +\dff 1} \dff)
\qff \ttoo\qff 
B\dff(\fff G^{\dff \bullet\dff +\dff 1} \dff)
\]

\vspace{-8pt}
unique up to chain homotopy.\oss
By\dss Theorem\qss \ref{acyclic-bundles-isomorphism}\qss 
the map\vspace*{4pt}
\[
\quad
u\fff(\fff p \fff)_{*}
\qff \colon\qff
\mathcal{H}^{\fff *}\fff(\dff G\fff,\qff B^{\fff \bullet}\fff(\dff \mathcal{X} \dff) \dff)
\qff \ttoo\qff
\widehat{H}^{\dff *}\fff(\dff G \dff)
\]

\vspace{-8pt}
induced by $u_{\dff \bullet}$
is an isometric isomorphism.\oss
We claim that the map\vspace*{4pt}
\[
\quad
v\fff(\fff p \fff)_{*}
\qff \colon\qff
\mathcal{H}^{\fff *}\fff(\dff G\fff,\qff 
\mathcal{B}\dff(\dff \mathcal{X}^{\fff \bullet\dff +\dff 1} \dff) \dff)
\qff \ttoo\qff
\widehat{H}^{\dff *}\fff(\dff G \dff)
\]

\vspace{-8pt}
induced\dss by $v_{\dff \bullet}$
is also an isometric isomorphism.\oss
The proof is almost the same as the proof\dss of\dss
Theorem\qss \ref{acyclic-bundles-isomorphism}.\oss
The key part of\dss the latter is 
the construction of a morphism of resolutions $r_{\fff \bullet}$\nnsp.\oss
If\dss the fundamental set $F$ used in the construction of $r_{\fff \bullet}$
is a Borel set\halfff,\oss
then the image of $r_{\fff \bullet}$
consists of\dss Borel straight cochains.\oss
It follows that $r_{\fff \bullet}$ defines a morphism of resolutions\vspace*{4pt}
\[
\quad
s_{\fff \bullet}
\qff \colon\qff
B\dff(\fff G^{\dff \bullet\dff +\dff 1} \dff)
\ttoo 
\mathcal{B}\dff(\dff \mathcal{X}^{\fff \bullet\dff +\dff 1} \dff)
\]

\vspace{-8pt}
extending\dss $\id_{\dff \rrr}$\dss
and consisting of maps with the norm\qss $\leq\qff 1$\nnsp.\oss
The rest of the proof requires only replacing\qss
$B^{\fff \bullet}\fff(\dff \mathcal{X} \dff)$\qss
by\qss
$\mathcal{B}\dff(\dff \mathcal{X}^{\fff \bullet\dff +\dff 1} \dff)$\dnsp.\off\oss
Hence $v\fff(\fff p \fff)_{*}$ is an isometric isomorphism.\oss
The maps\dss $u\fff(\fff p \fff)_{*}$\nnsp,\pss 
$v\fff(\fff p \fff)_{*}$\nnsp,\pss
and\dss $\mathcal{J}$\qss form a triangle diagram\vspace*{0pt}
\begin{equation*}
\quad
\begin{tikzcd}[column sep=my-size-t, row sep=my-size-t]\dis
\mathcal{H}^{\fff *}\fff(\dff G\fff,\qff \mathcal{B}\dff(\dff \mathcal{X}^{\fff \bullet\dff +\dff 1} \dff) \dff)
\arrow[rr, "\dis \mathcal{J}"]
\arrow[dr, "\dis u\fff(\fff p \fff)_{*}"']
&
& 
\mathcal{H}^{\fff *}\fff(\dff G\fff,\qff B^{\fff \bullet}\fff(\dff \mathcal{X} \dff) \dff)
\arrow[ld, "\dis v\fff(\fff p \fff)_{*}"] 
\\ 
& 
\off\qff\widehat{H}^{\dff *}\fff(\dff G \dff)\dff.
&  
\end{tikzcd}
\end{equation*}

\vspace{-9pt}
This triangle is commutative because all maps are induced by morphisms of resolutions
and\dss the resolution\dss
$B\dff(\dff G^{\fff \bullet\dff +\dff 1} \dff)$\dss
is relatively injective.\oss
Cf.\dss the proofs of commutativity in Section\qss \ref{spaces}.\pss
Since\dss $u\fff(\fff p \fff)_{*}$\nnsp,\pss 
$v\fff(\fff p \fff)_{*}$\dss
are isometric isomorphisms,\pss
$\mathcal{J}$ is an isometric isomorphism also.\oss \eproof

\myappend{Double\qss complexes}{double-complexes}

\vspace*{6pt}
\myapar{Theorem.}{double-complex}
\emph{Let\qss
$(\dff K^{\fff p\fff,\dff q} \dff)\dff_{p,\trf q\qff \geq\qff 0}$\qss
be a double complex with the differentials}\vspace*{3pt}
\[
\quad
d\qff \colon\qff K^{\fff p\fff,\dff q}\ttoo K^{\fff p\fff,\dff q\dff +\dff 1}
\hspace*{1.5em}\mbox{ \emph{and} }\hspace*{1.5em}
\delta\qff \colon\qff K^{\fff p\fff,\dff q}\ttoo K^{\fff p\dff +\dff 1\fff,\dff q}
\]

\vspace*{-9pt}
\emph{and\dss let\qss
$(\dff T^{\fff n} \dff)\dff_{n\qff \geq\qff 0}$\qss
be its total complex.\oss
Let\qss $L^p$\qss be the kernel of the differential}\qss\vspace*{3pt}
\[
\quad
d\dff \colon\dff
K^{\fff p\fff,\dff 0}\ttoo K^{\fff p\fff,\dff 1}\dff.
\]

\vspace*{-9pt}
\emph{Then\qss $L^*$\qss together with the differential\qss $\delta$\qss
is a subcomplex of the total complex\qss $T^{\fff *}$\dnsp.\oss
If the complexes\qss
$(\dff K^{\fff p\fff,\dff *},\pff d \dff)$\qss 
are exact\halfff,\oss
then the homomorphism\qss
$H^{\fff *}\fff(\dff L^* \dff)\ttoo H^{\fff *}\fff(\dff T^{\fff *} \dff)$\qss
induced by the inclusion\qss
$L^*\ttoo T^{\fff *}$\qss
is an isomorphism.\oss}

\proof\qss
This is a special case of Theorem\qss 4.8.1\qss from the Chapter\qss I\qss of
the book\qss \cite{go}.\oss
Its standard proof is based on the properties of the spectral sequences 
associated with a double complex.\oss
But it can be also proved directly,\oss by using the diagram chase.\oss  \eproof

\myapar{Theorem.}{comparison}
\emph{Let\dss
$f^{\dff \bullet\fff,\dff \bullet}\dff \colon\dff
K^{\dff \bullet\fff,\dff \bullet}
\qff \ttoo\qff 
L^{\dff \bullet\fff,\dff \bullet}$\sss
be a morphism of\dss double complexes.\oss
If\qss for every\sss $p\qff \geq\qff 0$\sss
the morphism of\trs complexes\dss
$f^{\dff p\fff,\dff \bullet}\dff \colon\dff
K^{\dff p\fff,\dff \bullet}
\qff \ttoo\qff 
L^{\dff p\fff,\dff \bullet}$\sss
induces an\dss isomorphism\sss in\sss cohomology,\oss
then\dss
$T f^{\dff \bullet}\dff \colon\dff
T^{\dff \bullet}_{\dff K}
\pff \ttoo\qff 
T^{\dff \bullet}_{\dff L}$\qss
also induces an\dss isomorphism\sss in\sss cohomology.\oss}

\proof
The assumption means\sss that\sss $f$\sss induces an\sss isomorphism of\trs the\sss terms\sss
$E_{\dff 1}^{\dff \bullet\fff,\dff \bullet}$\sss of\trs the standard spectral\sss sequences
converging\dss to\sss the cohomology\sss of\dss $T^{\dff \bullet}_{\dff K}$\sss
and\dss $T^{\dff \bullet}_{\dff L}$\nsp.\oss
This\sss implies\sss that\dss
$T f^{\dff \bullet}\dff \colon\dff
T^{\dff \bullet}_{\dff K}
\pff \ttoo\qff 
T^{\dff \bullet}_{\dff L}$\qss
induces an\dss isomorphism\sss in\sss cohomology.\oss
This\sss result\sss also can\sss be proved directly\sss by\sss a diagram chase.\oss  \eproof

\myuppar{Remark.}
Elementary\dss proofs of\dss both\sss theorems are included\sss in\qss \cite{i3}.\oss

\begin{flushright}

August\qss {17},\oss {2017}.\oss
Revised\qss December\qss 14,\oss 2020.
 
https\halfff:/\!/\hspace*{-0.06em}nikolaivivanov.com

E-mail\halfff:\oss nikolai.v.ivanov{\fff}@{\dff}icloud.com

\end{flushright}

\end{document}